\documentclass[12pt]{article}
\usepackage{amssymb}

\def\theequation{\thesubsection.\arabic{equation}}
\def\theguess{\thesubsection.\arabic{guess}}
\newtheorem{theorem}{Theorem}

\def\thetheorem{\thesubsection.\arabic{theorem}}

\def\theprop{\thesubsection.\arabic{prop}}
\newtheorem{lemma}[theorem]{Lemma}
\def\thelemma{\thesubsection.\arabic{lemma}}
\newtheorem{cor}[theorem]{Corollary}
\def\thecor{\thesubsection.\arabic{cor}}
\newtheorem{exam}[theorem]{Example}
\def\theexam{\thesubsection.\arabic{exam}}
\newtheorem{remark}[theorem]{Remark}
\def\theremark{\thesubsection.\arabic{remark}}
\newcommand{\eqa}{\begin{eqnarray}}
\newcommand{\eeqa}{\end{eqnarray}}
\newcommand{\beq}{\begin{equation}}
\newcommand{\eeq}{\end{equation}}
\newcommand{\nn}{\nonumber}
\newcommand{\dbl}{\langle\hskip -0.09truecm \langle}
\newcommand{\dbr}{\rangle\hskip -0.09truecm \rangle}
\newcommand{\pal}{\partial}

\newcommand{\al}{\alpha}
\newcommand{\ga}{\gamma}
\newcommand{\de}{\delta}

\newcommand{\lm}{\lambda}

\newcommand{\ve}{\varepsilon}

\newcommand{\pf}{\noindent{\it Proof \ }}
\newcommand{\tr}{{\rm tr}}
\newcommand{\epf}{$\quad$\hfill
\raisebox{0.11truecm}{\fbox{}}\par\vskip0.4truecm}

\begin{document}

\title {Normal forms of hierarchies of integrable PDEs, Frobenius manifolds 
and Gromov - Witten 
invariants}
\author{ {Boris Dubrovin${}^*$ \ \ Youjin Zhang${}^{**}$}\\
{\small ${}^{*}$ \ SISSA, Via Beirut 2--4, 34014 Trieste, Italy, and Steklov
Math. Institute, Moscow}\\
{\small ${}^{**}$ \ Department of Mathematical Sciences,
Tsinghua University, Beijing 100084, P.R.China
}}
\maketitle
\begin{abstract}
We present a project of classification of a certain class of bihamiltonian 1+1
PDEs depending on a small parameter. Our aim is to embed the theory of Gromov -
Witten invariants of all genera into the theory of integrable systems. The
project is focused at describing normal forms of the PDEs and their local
bihamiltonian structures satisfying certain simple axioms. A Frobenius manifold
or its degeneration is associated to every bihamiltonian structure of our type.
The main result is a
universal loop equation on the jet space of a semisimple Frobenius manifold that
 can be
used for perturbative reconstruction of the integrable hierarchy. We show that
first few terms of the perturbative expansion correctly reproduce the universal
identities between intersection numbers of Gromov - Witten classes and their
descendents.
\end{abstract}

\medskip

\centerline{SISSA Preprint 65/2001/FM}

\vskip 0.5cm

\newpage

\tableofcontents

\def\theequation{\thesection.\arabic{equation}}
\def\theguess{\thesection.\arabic{guess}}
\def\thetheorem{\thesection.\arabic{theorem}}
\def\theprop{\thesection.\arabic{prop}}
\def\thelemma{\thesection.\arabic{lemma}}
\def\thecor{\thesection.\arabic{cor}}
\def\theexam{\thesection.\arabic{exam}}
\def\theremark{\thesection.\arabic{remark}}
\setcounter{equation}{0}
\setcounter{theorem}{0}

\newpage
\noindent
\section{Introduction}\par
In this paper we study the structure of integrable systems of evolutionary PDEs
with one spatial variable of the form
\beq\label{pde0}
u^i_t =\sum_{j=1}^n A^i_j(u) u^j_x +{\rm perturbation}, ~~i=1, \dots, n.
\eeq
The basic example of such a system is the celebrated Korteweg - de Vries
equation (in this example $n=1$)
\beq\label{kdv0}
u_t = u\, u_x +{\epsilon^2\over 12} u_{xxx}.
\eeq
Here $\epsilon$ is the parameter of the perturbation. Our aim is to develope an
approach to the problem of classification of such systems based on the
deep relationship between integrable systems and quantum field theory discovered
in the last decade \cite{BK, DS, GM, DVV,
FKN, Douglas1990, 
Witten1, DW, Witten2,
konts}. The amazing discovery of E. Witten and M. Kontsevich of the relationship
between KdV and the topology of the moduli spaces of stable algebraic curves
opened a new dimension in the theory of integrable systems. 
The surprising outcome of our classification project is that, at least at low
orders of the perturbative expansion, the topology of the moduli spaces
$\bar{\cal M}_{g,n}$ of stable 
algebraic curves is ``hidden'' in {\it every} integrable systems 
of our class.

Another motivation for our work was the wish to find a proper setup
for the general theory of Frobenius manifolds (cf. \cite{givental}).
As it has been suggested in \cite{cmp, selecta, taniguchi, icm} (see also
\cite{Witten2}) the right setting for the theory of semisimple Frobenius
manifolds is the theory of hierarchies of integrable PDEs.
 
Frobenius manifolds were introduced by one of the authors \cite{verdier}
as a coordinate-free form of the WDVV associativity equations \cite{Witten1,
DVV}. We refer the
reader to \cite{D3, D4, hitchin, manin} for the details of the theory
of Frobenius manifolds. In the mathematical literature Frobenius manifolds
are best known in quantum cohomologies, i.e., 
in the theory of the genus zero Gromov - Witten invariants of
smooth projective varieties or, more generally, of compact symplectic manifolds
\cite{Kon2, behrend, ruan-tian, mcduff-salamon}, although many ingredients of
the theory of Frobenius manifolds already appeared in the
singularity theory as a natural geometrical structure on the base of the
universal unfolding of an isolated hypersurface singularity \cite{saito, manin,
hertling}. Another source of Frobenius manifolds 
is the geometry of the orbit spaces of finite Coxeter groups
\cite{D2, coxeter}
and their generalizations \cite{weyl, bertola}. The notion of flat coordinates
on the orbit spaces discovered in \cite{sys, saito0} was important in these
constructions. The Frobenius manifolds of the singularity theory and of the
theory of reflection groups are always semisimple. The origin of semisimplicity
in quantum cohomology is still to be understood \cite{tian-xu, icm,
bayer-manin}. 

In certain cases mirror symmetry constructions, or the Arnold -
Brieskorn correspondence between ADE Weyl groups and simple hypersurface 
singularities
establish relationships between different classes of
examples of Frobenius manifolds. However, the general unifying principle
of the theory of Frobenius manifolds eventually covering also the theory of
Gromov - Witten invariants of higher genera is still missing.

Our suggestion is that, the right framework of the theory of Frobenius manifolds
of all genera is the theory of integrable PDEs along with all main ingredients
of this theory, i.e., bihamiltonian structures, tau-functions, Virasoro
symmetries, W-algebras etc. (an expert in the theory of integrable systems may
have his own opinion about how should this list of ``main ingredients''  be
continued). 

It has already been proven in \cite{D92} that all the genus zero
topological recursion relations for the so-called descendents can be correctly 
reproduced starting from
an arbitrary Frobenius manifold. Dispersionless integrable hierarchies
were crucial in this reconstruction theorem (see also \cite{Witten2}).
Bihamiltonian structure for these hierarchies was discovered in \cite{D2, D3}.
The next important step has been done in \cite{cmp}. It was shown that
also the genus one topological recursion relations \cite{Witten2} together
with the E. Getzler's defining relation for elliptic Gromov - Witten
invariants \cite{getzler1} can be reproduced starting from an arbitrary
semisimple Frobenius manifold. Using this result the Virasoro conjecture
of T. Eguchi {\it et al.} \cite{eguchi1, eguchi2, eguchi3} was proved
in \cite{selecta} up to genus one approximation
(see also \cite{liu-tian, liu-vir, getzler-vir} 
for an alternative 
approach to the
theory of Virasoro constraints in quantum cohomology).
 
In \cite{cmp} topology of \cite{Witten2, getzler1} was used to uniquely 
recover the first
order integrable perturbation of the dispersionless hierarchy of \cite{D92}.
In principle this approach can be extended also to higher genera (see
\cite{eguchi-getzler} where the genus two topological recursion relations
and other identities in the cohomology of $\bar {\cal M}_{g,n}$ for $g\leq 2$
\cite{belorus-panda, getzler2}
were used in order to compute the genus 2 free energy in topological sigma
models with two primaries and also \cite{liu} where a full system of equations
for the genus two Gromov - Witten potential has been obtained in the general
case).
We now want to change completely the setting. Instead of using topology for
constructing integrable hierarchies, as it was done in \cite{D92, cmp}
we want to develope an approach to the problem of classification of integrable
hierarchies eventually reproducing {\it all the universal identities between
Gromov - Witten invariants and their descendents of all genera} even for those
integrable hierarchies that {\it a priori} have nothing to do with topology.

The main result of this paper is a system of four axioms of the theory of
integrable hierarchies of the form (\ref{pde0}) that can be used as the
basis of the classification of these hierarchies.  We prove that, under
assumption of semisimplicity these axioms allow to uniquely reconstruct 
the whole
structure of the hierarchy starting from the dispersionless limit $\epsilon\to
0$. We are also able to correctly reproduce, starting from our axioms,
essentially all known universal identities for the Gromov - Witten cocycles and
their descendents
in $H^*(\bar {\cal M}_{g,n})$ written as differential constraints
for the tau-function of the hierarchy. In particular, we prove the
$3g-2$-conjecture of T. Eguchi and C.-S. Xiong \cite{eguchi3} and reproduce
the correct shape of the Virasoro constraints, derive the formulae for the genus
1 and genus 2 Gromov - Witten potential  expressing them in terms
of the genus 0 one (cf. \cite{givental, givental1}) etc. 
Let us emphasize again that,
all this has been done {\it for an arbitrary} semisimple Frobenius manifold.
It remains an open problem to figure out what could be the meaning of our
integrable hierarchies for other classes of Frobenius manifolds, e.g., in the
singularity theory.

To explain our axioms
let us first recall some features of KdV crucial for our classification project.

One of the starting point of the KdV theory was the discovery \cite{gkm68} of
an infinite family of commuting evolutionary PDEs commuting with (\ref{kdv0}),
\beq\label{hier-kdv}
{\pal u\over \pal t_j} = K_j( u, u_x, \dots, u^{(2j+1)}), ~~{\pal\over \pal t_i}
{\pal u\over \pal t_j}= {\pal\over \pal t_j}
{\pal u\over \pal t_i}
\eeq
with some polynomials $K_j$, $j\geq 0$. For $j=0$ one obtains the spatial
translations
$$
{\pal u\over \pal t_0} = {\pal u\over \pal x},
$$
for $j=1$ (\ref{hier-kdv}) coincides with (\ref{kdv0}), other equations of the
so-called {\it KdV hierarchy} (\ref{hier-kdv}) 
are more complicated. E.g., 
$$
{\pal u\over \pal t_2} = {1\over 2}{u^2 u'} +{\epsilon^2\over 12} (2 u' u'' + u\, u''') +
{\epsilon^4\over 240} u^{V}.
$$
They 
are obtained by a suitable
recursion procedure \cite{gkm68, ggkm74}. The latter has been represented
\cite{magri} in the {\it bihamiltonian form}: the equation of the hierarchy
are considered as flows on the space of functions $u(x)$ hamiltonian w.r.t.
two Poisson brackets $\{~,~\}_1$ and $\{~,~\}_2$
\beq\label{recurs}
{\pal u\over \pal t_j} = \{ u(x), H_j\}_1 = \left( j+{1\over 2}\right)^{-1}
\{ u(x), H_{j-1}\}_2,
\eeq
\beq\label{recurs1}
\{u(x), u(y)\}_1 = \delta'(x-y), ~~\{u(x), u(y)\}_2= 
u(x) \delta'(x-y) +{1\over 2} u'(x) \delta(x-y)
+{\epsilon^2\over 8}
\delta'''(x-y).
\eeq
The crucial property in the definition of a bihamiltonian structure 
is {\it compatibility} of the pair of the Poisson brackets:
$$
a_1\{~,~\}_1+a_2\{~,~\}_2
$$
must be a Poisson bracket for arbitrary constant coefficients $a_1$, $a_2$.
The Hamiltonians are local functionals
$$
H_j = \int h_j(u, u_x, \dots, u^{(2j+2)})dx
$$
to be determined recursively from (\ref{recurs}) starting from the Casimir
$$
H_{-1}=\int u\, dx.
$$
Explicitly,
$$
h_0 = {u^2\over 2} +\epsilon^2{u''\over 12}, ~~h_1={u^3\over 6} +
{\epsilon^2\over 24}({u'}^2 + 2 u\, u'') + \epsilon^4 {u^{IV}\over 240},
$$
$$
h_2 = {u^4\over 24} + {\epsilon^2\over 24} (u\, {u'}^2 +
u^2u'')+{\epsilon^4\over 480} (3{u''}^2+4u'u'''+2u\, u^{IV}) +{\epsilon^6\over
6720} u^{VI}.
$$
This observation is the starting point of our study: what we are classifying
is not just a single integrable PDE (\ref{pde0}) but a {\it hierarchy} of
integrable PDEs produced by a bihamiltonian recursion procedure. To say the same
thing in a shorter way: we want to classify bihamiltonian structures of integrable
PDEs of the form (\ref{pde0}). This distinguishes our approach from the
symmetry analysis technique 
(see \cite{mikh-shabat} and references therein) or from the Painlev\'e test
(see \cite{kruskal} and references therein) proved to be powerful
in classification of integrable PDEs of low orders. Our approach differs also
from the perturbative method of V.E. Zakharov {\it et al.} (see \cite{zakharov}
and references therein)
where nonlinear integrable perturbations of {\it linear} systems were studied.

We impose three additional constraints onto the bihamiltonian structure.
The first one is existence of a {\it tau-function} \cite{djkm83}. For the example of the KdV
hierarchy this means that, for an arbitrary solution $u=u(x+t_0, t_1, t_2,
\dots)$ of the hierarchy the densities of the Hamiltonians can be represented
in the form
\beq\label{recurs2}
h_j(u, u_x, \dots, u^{(2j+2)}) =\epsilon^2 {\pal^2 \log \tau(x+t_0, t_1, t_2, \dots)\over
\pal x \,\pal t_{j+1}}
\eeq
for some function $\tau(x+t_0, t_1, t_2, \dots)$.
In particular, for $j=-1$ one obtains the wellknown formula
$$
u=\epsilon^2{\pal^2\log\tau\over \pal x^2}.
$$
Existence of a {\it single} tau-function is a rather strong restriction onto
an integrable hierarchy (cf. \cite{miramontes} where such a tau-function was
constructed for generalized integrable hierarchies of KdV and affine Toda type). This restriction corresponds to the choice of a primitive
conjugacy class \cite{kac-peterson} in the Weyl group in the
setting of the theory of generalized Drinfeld - Sokolov hierarchies.

According to the conjecture by E.Witten \cite{Witten2} proved by M.Kontsevich
\cite{konts}, the logarithm of the tau-function of the particular solution of the KdV 
hierarchy specified by the initial data 
$$
u|_{t_j=0} =x,
$$
\eqa\label{tau-kw}
&&\log\tau=\frac1{\epsilon^2}\left(
\frac{t_0^3}{6} + \frac{t_0^3\,t_1}{6} + 
  \frac{t_0^3\,t_1^2}{6} + \frac{t_0^3\,t_1^3}{6} + 
  \frac{t_0^3\,t_1^4}{6} + \frac{t_0^4\,t_2}{24}+ 
  \frac{t_0^4\,t_1\,t_2}{8}\right.\nn\\
&&\quad\left. + 
  \frac{t_0^4\,t_1^2\,t_2}{4} + 
  \frac{t_0^5\,t_2^2}{40} + \frac{t_0^5\,t_3}{120} + 
  \frac{t_0^5\,t_1\,t_3}{30} + \frac{t_0^6\,t_4}{720}+\dots
\right)
\nn\\
&&\quad+\left(
\frac{t_1}{24} + \frac{t_1^2}{48} + \frac{t_1^3}{72} + 
  \frac{t_1^4}{96} + \frac{t_0\,t_2}{24} + 
  \frac{t_0\,t_1\,t_2}{12} +
 \frac{t_0\,t_1^2\,t_2}{8} + 
  \frac{t_0^2\,t_2^2}{24}\right.\nn\\
&&\quad \left. + \frac{t_0^2\,t_3}{48} + 
  \frac{t_0^2\,t_1\,t_3}{16} + 
\frac{t_0^3\,t_4}{144}+\dots\right)\nn\\
&&\quad+\epsilon^2\left(
\frac{7\,t_2^3}{1440} + \frac{7\,t_1\,t_2^3}{288} + 
  \frac{29\,t_2\,t_3}{5760} + 
  \frac{29\,t_1\,t_2\,t_3}{1440} + 
  \frac{29\,t_1^2\,t_2\,t_3}{576} + 
  \frac{5\,t_0\,t_2^2\,t_3}{144}\right.\nn\\ 
&&\quad+ 
  \frac{29\,t_0\,t_3^2}{5760} + 
  \frac{29\,t_0\,t_1\,t_3^2}{1152} + \frac{t_4}{1152} + 
  \frac{t_1\,t_4}{384} + \frac{t_1^2\,t_4}{192} + 
  \frac{t_1^3\,t_4}{96} + 
\frac{11\,t_0\,t_2\,t_4}{1440} \nn\\
&&\quad\left.+ 
  \frac{11\,t_0\,t_1\,t_2\,t_4}{288} + 
  \frac{17\,t_0^2\,t_3\,t_4}{1920}+\dots\right)+O(\epsilon^4)\nn
\eeqa
coincides with the generating function of the intersection numbers of 
the Mumford - Morita - Miller classes in $H^*(\bar {\cal M}_{g,n})$.
Here $\bar {\cal M}_{g,n}$ is the moduli space of stable algebraic curves of
genus $g$ with $n$ punctures, 
$$
\log \tau=\sum\epsilon^{2g-2} {\cal F}_g,
$$
$$
{\cal F}_g =\sum_n {1\over n!}t_{p_1}\dots t_{p_n}
\int _{\bar {\cal M}_{g,n}} c_1^{p_1}({\cal L}_1)\wedge \dots \wedge 
c_1^{p_n}({\cal L}_n)
$$
where ${\cal L}_i$ is the tautological line bundle over the moduli space
corresponding to the $i$-th puncture. In this setting different powers of 
the small dispersion
parameter $\epsilon$ in the KdV correspond to the contributions of different
genera $g$. In other words, the small dispersion expansion coincides with 
the genus
expansion. In the physical literature on topological field theory
the parameter $\epsilon$ is called
string coupling constant.

The first two assumptions, i.e., existence of a bihamiltonian structure
and of a tau-function of the integrable hierarchy imply that the {\it
dispersionless limit} $\epsilon\to 0$ is described by a {\it Frobenius manifold}
structure on the space of dependent variables of the hierarchy (see 
Section \ref{sec-3-5}
below) or by a {\it degenerate Frobenius manifold} structure. 
The dispersionless hierarchy itself is reconstructed by the Frobenius manifold
structure according to the construction of \cite{D92} (we call it Principal
Hierarchy in Section \ref{sec-3-6}
below). The Principal Hierarchy possesses all the universal properties 
observed in the theory of weak dispersive limits of integrable PDEs
\cite{dn83, dn89, pisa, cmp92, krichever, takasaki, bmrw}. (We consider here
only the formal geometric side of the theory of weakly dispersive integrable
PDEs. We refer the reader to the papers of P. Lax, D. Levermore, S. Venakides,
see \cite{lax-lever} and references therein,
for the analytic side of this theory. We hope, however, that our geometric 
analysis
could be useful also for the analytic theory.)

The next step is the main one: we are to reconstruct the full hierarchy
(\ref{pde0}) together with the bihamiltonian structure
starting from their dispersionless limit. (We do not consider in this paper the
hierarchies corresponding to degenerate Frobenius manifolds. We plan to do it
elsewhere.)
The assumption of semisimplicity is to
be added at this point. From the point of view of integrable systems
semisimplicity ensures complete integrability, i.e., completeness of the
family of commuting integrals (see Section \ref{sec-3-6-5} below). 
The last two axioms are
used to provide uniqueness of the reconstruction. 

The axiom 3 is the most disputable one. We call it {\it quasitriviality} of the
hierarchy. Before explaining this axiom we are to formulate in a more 
precise way
the classification problem. We study bihamiltonian PDEs depending on a formal
small parameter $\epsilon$ represented as a (formal) small dispersion expansion
\beq\label{pde00}
u^i_t =\sum_{j=1}^n A^i_j(u) u^j_x +\sum_{k>0} \epsilon^k K^i_{[k]}(u; u_x, \dots,
u^{(k+1)}), ~~i=1, \dots, n
\eeq
where $K^i_{[k]}(u; u_x, \dots,
u^{(k+1)})$ is a {\it polynomial  in the derivatives} weighted homogeneous 
of the degree $k+1$. It is understood that the $m$-th derivative ${u^i}^{(m)}$
has degree $m$. We classify these PDEs and their {\it local} 
bihamiltonian structures w.r.t. the {\it Miura group} of transformations of the form
\beq\label{miura0}
u^i\mapsto F^i_{[0]}(u) + \sum_{k>0} \epsilon^k F^i_{[k]}(u; u_x, \dots,
u^{(k)})
\eeq
where the coefficients $F^i_{[k]}(u; u_x, \dots,
u^{(k)})$ are {\it homogeneous polynomials in the derivatives} of the degree $k$ and
$$
\det \left( {\pal F^i_{[0]}(u)\over \pal u^j} \right) \neq 0.
$$
The problem of classification can be presented as the problem of description
of {\it normal forms} of integrable PDEs w.r.t. the transformations
(\ref{miura0}). The Miura group acts also on local translation invariant
Poisson brackets of systems
of the form (\ref{pde00}). We call them (0,n)-brackets (see the definition in
Section 2.4.3 below). It turns out that, at least over complex numbers all (0,n)
Poisson brackets are equivalent w.r.t. Miura group. This important technical
step of our theory is based on the differential-geometric theory, due to S.P.
Novikov and B. Dubrovin \cite{dn83} of the
so-called Poisson brackets of hydrodynamic type and also on the
triviality of the Poisson cohomology of these brackets proved by E. Getzler
\cite{getzler} and also by L. Degiovanni, F. Magri, V. Sciacca \cite{magri1}
\footnote{The problem of normal forms of Poisson brackets of
PDEs was studied also in \cite{medvedev}. However, a more general class of
admissible transformations was considered. This simplified the solution of the
classification problem.}. The main object of our study is {\it the problem of normal
forms of bihamiltonian structures of systems (\ref{pde00}) w.r.t. the Miura
group}. 

An integrable hierarchy (\ref{pde00}) is called trivial if it can be obtained,
together with the underlined bihamiltonian structure, from the dispersionless limit
$\epsilon=0$ by action of a transformation of the form (\ref{miura0}). It is
called {\it quasitrivial} if the same is true w.r.t. a transformation
\beq\label{quasi-miura0}
u^i\mapsto F^i_{[0]}(u) + \sum_{k>0} \epsilon^k F^i_{[k]}(u; u_x, \dots,
u^{(m_k)})
\eeq
where the coefficients are {\it rational functions} in the derivatives. The
quasitriviality property seems to be unobserved even in the theory of the KdV
equation (\ref{kdv0}). We prove it in Section \ref{sec-3-8}. 
The quasitriviality 
transforming the hierarchy of the dispersionless KdV
$$
u_t =u\, u_x
$$
to the full KdV together with the bihamiltonian
structure etc. reads
\beq\label{g-exp0}
u\mapsto u +{\epsilon^2\over 24} \left( \log u_x\right)_{xx}
+\epsilon^4 \left( {u^{IV}\over 1152\, {u'}^2} 
- {7\, u'' u'''\over 1920\, {u'}^3}
+{{u''}^3\over 360\, {u'}^4}\right)_{xx} + O(\epsilon^6).
\eeq
The reader easily recognizes in (\ref{g-exp0}) the genus expansion of the
topological gravity written in the form suggested by R. Dijkgraaf and E. Witten
\cite{DW} (see also \cite{itzykson-zuber, eguch0}). So, with respect to
the group of rational Miura transformations (\ref{quasi-miura0})
the normal forms of all our hierarchies are just the dispersionless
Principal Hierarchies; the reducing transformation (\ref{quasi-miura0})
is the candidate for the role of the genus expansion. 

Indeed, we prove that, choosing in a clever way the dependent variables of the
hierarchy the reducing transformation (\ref{quasi-miura0}) is expressed via
second derivatives of the $\epsilon$-expansion of the logarithm of the
tau-function of the hierarchy. Moreover, we prove, using that the
Poisson pencil depends polynomially in the
derivatives, that the latter must have the
form
\beq\label{logtau0}
\log \tau =\epsilon^{-2} {\cal F}_0 +\sum_{g\geq 1} \epsilon^{2g-2}
{\cal F}_g (u; u_x, \dots, u^{(3g-2)}).
\eeq
The terms with $g\geq 1$ of the expansion of the reducing transformation do not depend
on the choice of solution of the hierarchy.
In the setting of topological sigma-models (\ref{logtau0}) coincides with 
the so-called
$3g-2$-conjecture of \cite{eguchi3} (see also \cite{eguchi-getzler}).

The last axiom is used to uniquely fix the terms of the expansion
(\ref{logtau0}). It is based on study of symmetries of the integrable PDEs.
First we prove, in Section \ref{sec-3-10-1}
 that the Principal Hierarchy always admits 
a rich algebra of
infinitesimal symmetries isomorphic to the half of the Virasoro algebra. Due to
quasitriviality these symmetries can be lifted to Virasoro symmetries of the
full hierarchy. We require that the generators of the action of the half of 
the Virasoro
algebra by symmetries of the hierarchy {\it act linearly onto the tau-function}
of the hierarchy. For the KdV example the action of the
generators of the Virasoro algebra by symmetries of the hierarchy is given by
the following formulae \cite{DVV, FKN}
\beq\label{vir-sym-0}
\delta_m \tau = L_m \tau, ~~m\geq -1
\eeq
where
\eqa
&&
L_{-1} = \sum_{p\geq 1} t_p \pal_{p-1} +{1\over 2\epsilon^2} t_0^2
\nn\\
&&
L_0 =\sum_{p\geq 0} \left( p+{1\over 2}\right) t_p \pal_p +{1\over 16}
\nn\\
&&
L_1 = \sum_{p\geq 0} \left( p+{1\over 2}\right)\left( p+{3\over 2}\right)
 t_p \pal_{p+1} +{\epsilon^2\over 8} \pal_0^2
 \nn\\
 &&
 L_2= \sum_{p\geq 0} \left( p+{1\over 2}\right)\left( p+{3\over 2}\right)
 \left( p+{5\over 2}\right)
 t_p \pal_{p+2} +{3\,\epsilon^2\over 8} \pal_0 \pal_1
 \nn\\
&&
\dots\nn
\eeqa
with
$
\pal_k ={\pal\over \pal t_k}.
$
The Witten - Kontsevich tau-function (\ref{tau-kw}) is uniquely specified
\cite{KacSchwarz}
by the following system of {\it Virasoro constraints}
\eqa\label{kw-constraints}
&&
L_m \tau = \prod_{j=1}^{m+1} \left( j+{1\over 2}\right) \pal_{m+1}\tau, ~~m\geq 0
\\
&&
L_{-1} \tau =\pal_0 \tau.\nn
\eeqa

We call this last axiom {\it linearization of the Virasoro symmetries}. A
surprizingly looking consequence of this axioms says that, under certain
assumption of monotonicity the tau functions of {\it all} analytic solutions 
to our hierarchies
are annihilated by an appropriate linear combination of the Virasoro symmetries
and the flows of the hierarchy.

We also use the axiom of linearization of the Virasoro symmetries
to define a set of defining equations for our integrable hierarchies.
Using universality of the coefficients
of expansion of 
\beq\label{deltaf0}
\Delta {\cal F}(u; u_x, u_{xx}, \dots; \epsilon^2) = 
\sum_{g\geq 1} \epsilon^{2g-2}
{\cal F}_g (u; u_x, \dots, u^{(3g-2)})
\eeq
we derive a kind of {\it loop equation} (\ref{glavnoe}) 
for the function (\ref{deltaf0})
on the jet space (regarding loop equations and their applications
in matrix models and in topological gravity see \cite{AmbjornMakeenko, DVV,
FKN}). Universality of $\Delta{\cal F}$ as a function on the jet space is used
in order to develop a machinery of perturbative solution of loop equation. In
particular it allows to fix ambiguities inavoidable in the standard approach to
the loop equation
\cite{Makeenko} and to prove uniqueness of the reconstruction of the integrable
hierarchy starting from an arbitrary semisimple Frobenius manifold.

\def\theequation{\thesubsection.\arabic{equation}}
\def\theguess{\thesubsection.\arabic{guess}}
\def\thetheorem{\thesubsection.\arabic{theorem}}
\def\theprop{\thesubsection.\arabic{prop}}
\def\thelemma{\thesubsection.\arabic{lemma}}
\def\thecor{\thesubsection.\arabic{cor}}
\def\theexam{\thesubsection.\arabic{exam}}
\def\theremark{\thesubsection.\arabic{remark}}
\setcounter{equation}{0}
\setcounter{theorem}{0}

\section{Normal forms of Hamiltonian structures of evolutionary systems}
\label{sec-2}\par
\subsection{Brief summary of finite-dimensional Poisson geometry}\par
Let $P$ be a $N$-dimensional smooth manifold. A {\it Poisson bracket} on $P$ is
a structure of a Lie algebra on the ring of functions ${\cal F}:= {\cal
C}^\infty (P)$
$$
f, g \mapsto \{ f, g \},
$$
\beq
\{ g, f\} =-\{f,g\}, \\
\{ a f + b g, h\} = a \{f,h\} + b \{ g,h\}, ~~ a,\ b\in {\bf R}, f,\ g, \ h\in
{\cal F}
\eeq
\beq\label{jacobi}
\{\{f,g\},h\} + \{\{ h,f\}, g\} + \{\{ g,h\},f\} =0
\eeq
satisfying the Leibnitz rule
$$
\{ fg,h\} =f\{g,h\} + g\{ f,h\}
$$
for arbitrary three functions $f,\ g, \ h \in {\cal F}$. In a system of local
coordinates $x^1, \dots, x^N$ the Poisson bracket reads
\beq
\{ f,g\} =h^{ij}(x) \frac{\partial f}{\partial x^i} 
\frac{\partial g}{\partial x^j}
\eeq
(summation over repeated indices will be assumed) where the bivector
$h^{ij}(x)=-h^{ji}(x)$$\newline
 $$= \{ x^i, x^j\}$ satisfies the following system of
equations equivalent to the Jacobi identity (\ref{jacobi})
\beq\label{jacobi1}
\{\{x^i, x^j\}, x^k\} + \{\{ x^k, x^i\},x^j\} +\{ \{x^j, x^k\},x^i\} \\
\equiv \frac {\pal h^{ij}}{\pal x^s} h^{sk} + 
\frac {\pal h^{ki}}{\pal x^s} h^{sj}
+ \frac {\pal h^{jk}}{\pal x^s} h^{si} =0
\eeq
for any $i$, $j$, $k$. Such a bivector satisfying
(\ref{jacobi1}) is called {\it Poisson structure} on $P$.

Clearly any bivector constant in some coordinate system is a Poisson structure.
Vice versa \cite{lichn}, locally all solutions to (\ref{jacobi1}) of 
the constant
rank $2n = {\mathrm rk} (h^{ij})$ can be reduced, by a change of coordinates, 
to the 
following {\it normal form}
\beq\label{normal}
h=\left(\matrix{\bar h & 0 \cr 
0 & 0\cr}\right)
\eeq 
with a constant nondegenerate antisymmetric $2n \times 2n$ matrix 
$\bar h=\bar h^{ab}$. That means that locally there exist coordinates
$y^1, \dots , y^{2n},
c^1, \dots, c^k$, $2n+k=N$, s.t.
$$
\bar h^{ab} =\{ y^a, y^b\}
$$
and
\beq\label{casimir}
\{ f, c^j\} =0, ~~ j=1, \dots k
\eeq
for an arbitrary function $f$. 

For the case $2n =N$ the inverse matrix $\left( h_{ij}(x) \right) =\left(
h^{ij}(x) \right)^{-1}$ defines on $P$ a {\it symplectic structure}
$$
\Omega =\sum_{i<j} h_{ij}(x) dx^i \wedge dx^j, ~~ \Omega^n \neq 0.
$$
For $2n < N$ one obtains on $P$ a structure of {\it symplectic foliation}
$P=\cup_{c_0} P_{c_0}$, $c_0 = (c^1_0, \dots c^k_0)$, of the codimension $k=N-2n$
\beq\label{foliation}
P_{c_0} := \{ x \ |  \ c^1(x) = c^1_0, \dots , c^k(x) = c^k_0\}.
\eeq
The independent functions $c^1(x), \dots, c^k(x)$ defined in (\ref{casimir})
are called {\it Casimir functions}, or simply {\it Casimirs} of the Poissson
structure. Every leaf $P_{c_0}$ is a symplectic manifold, and the restriction
map ${\cal C}^\infty (P) \to {\cal C}^\infty (P_{c_0})$ is a homomorphism of
Lie algebras.

\begin{exam} Let ${\mathfrak  g}$ be $n$-dimensional Lie algebra. The {\it Lie -
Poisson} bracket on the dual space $P={\mathfrak  g}^*$ reads
\beq\label{lie-poiss}
\{ x^i, x^j\} =c^{ij}_k x^k.
\eeq
Here $c^{ij}_k$ are the structure constants of the Lie algebra. The Casimirs
of this bracket are functions on ${\mathfrak  g}^*$ invariant with respect to the
co-adjoint action of the associated Lie group $G$. The symplectic leaves
coincide with the orbits of the coadjoint action with the Berezin - Kirillov -
Kostant symplectic structure on them.
\end{exam}  

An arbitrary foliation $P=\cup_{\phi_0} P_{\phi_0}$ of a codimension $m$  
represented locally in the form
$$
P_{\phi_0} =\{ x \ | \ \phi^1(x) =\phi^1_0, \dots, \phi^m(x) = \phi^m_0\}
$$
will be  called {\it cosymplectic} if the $m\times m$ matrix $\{\phi^a,
\phi^b\}$ does not degenerate on the leaves. In this situation a new Poisson
structure $\{ \ , \ \}_D$ can be defined on $P$ s.t. the functions $\phi^a(x)$
are Casimirs of $\{ \ , \ \}_D$. This is the {\it Dirac bracket} given
explicitly by the formula
\beq\label{dirac}
\{ f,g\}_D =\{ f,g\} - \sum_{a, b} \{f, \phi^a\} \{ \phi^a, \phi^b\}^{-1}
\{\phi^b , g\}.
\eeq
It can be restricted in an obvious way to produce a Poisson structure on every
leaf. The restriction map
$$
\left( {\cal C}^\infty (P), \{ \ , \ \} \right) \to 
\left( {\cal C}^\infty (P_{\phi_0}), \{ \ , \ \}_D \right)
$$
is a homomorphism of Lie algebras.

A Poisson bracket defines an (anti)homomorphism 
$$
{\cal F} \to Vect(P)
$$
\beq\label{ham.vect.field}
H\mapsto X_H:= \{ \cdot , H\},
\eeq
$$
[X_{H_1}, X_{H_2}] =- X_{\{H_1, H_2\}}.
$$
$X_H$ is called {\it Hamiltonian vector field}. The corresponding dynamical
system
\beq\label{ham.system}
\dot x^i = h^{ij}(x) \frac{\pal H}{\pal x^j}
\eeq
is called {\it Hamiltonian system with the Hamiltonian $H(x)$}. It is a
{\it symmetry of the Poisson bracket}
\beq\label{symmetry}
Lie_{X_H} \{ \ , \ \} =0.
\eeq

The last one is the notion of Poisson cohomology of $\left( P, \{ \ , \ \}
\right)$ introduced by Lichnerowicz \cite{lichn}. We need to use the {\it
Schouten - Nijenhuis bracket}. Denote
$$
\Lambda^k =
H^0(P, \Lambda^k TP)
$$ 
the space of multivectors on $P$. The Schouten - Nijenhuis bracket
is a bilinear pairing $a, b \mapsto [a,b]$,
$$
\Lambda^k \times \Lambda^l \to \Lambda^{k+l-1}
$$ uniquely determined by the properties of supersymmetry
\beq\label{susy}
[b,a] =(-1)^{kl} [a,b], ~~ a\in \Lambda^k , ~ b\in \Lambda ^l
\eeq
the graded Leibnitz rule
\beq\label{schout1}
[c, a\wedge b] =[c,a]\wedge b + (-1)^{lk+k} a\wedge [c,b], ~~ a\in \Lambda^k , ~ c\in \Lambda ^l
\eeq
and the conditions $[f,g]=0$, $f, g \in \Lambda ^0 = {\cal F}$, 
$$
[v,f] =v^i \frac {\pal f}{\pal x^i}, ~~ v\in \Lambda^1 =Vect(P), ~~ f\in
\Lambda^0 = {\cal F},
$$
$[v_1, v_2]$ = commutator of vector fields for $v_1, v_2 \in \Lambda^1$. In
particular for a vector field $v$ and a multivector $a$
$$
[v,a] =Lie_v a.
$$
\begin{exam} For two bivectors $h=(h^{ij})$ and $f=(f^{ij})$ their Schouten - Nijenhuis
bracket is the following trivector
\beq\label{schout.explicit}
[h,f]^{ijk} =\frac{ \pal h^{ij} } {\pal x^s} f^{sk} +\frac{\pal f^{ij}}{\pal
x^s} h^{sk} + 
\frac{ \pal h^{ki} } {\pal x^s} f^{sj} +\frac{\pal f^{ki}}{\pal
x^s} h^{sj} +
\frac{ \pal h^{jk} } {\pal x^s} f^{si} +\frac{\pal f^{jk}}{\pal
x^s} h^{si}.
\eeq
Observe that the l.h.s. of the Jacobi identity (\ref{jacobi1})  
reads
$$
\{\{x^i, x^j\}, x^k\} + \{\{ x^k, x^i\},x^j\} +\{ \{x^j, x^k\},x^i\} =\frac 1
2 [h,h]^{ijk}.
$$
\end{exam}

The Schouten - Nijenhuis bracket satisfies the graded Jacobi identity \cite{nijen}
\beq\label{schout2}
(-1)^{km} [[a,b],c] +(-1)^{lm}[[c,a],b] +(-1)^{kl} [[b,c],a]=0, ~~ a\in
\Lambda^k, \ b\in \Lambda^l, \ c \in \Lambda^m.
\eeq
It follows that, for a Poisson bivector $h$ the map
\beq\label{pdif}
\pal : \Lambda^k \to \Lambda^{k+1}, \ \pal a =[h,a]
\eeq
is a differential, $\pal^2 =0$. The cohomology 
of the complex $\left(
\Lambda^*, \pal \right)$ is called {\it Poisson cohomology} of $\left( P, \{
\ , \ \}\right)$. We will denote it
$$
H^*(P,\{ \ , \} ) = \oplus_{k\geq 0} H^k (P,\{ \ , \} ).
$$  
In particular,
$H^0 (P,\{ \ , \} )$ coincides with the ring of Casimirs of the Poisson
bracket, $H^1 (P,\{ \ , \} )$ is the quotient of the Lie algebra of
infinitesimal symmetries 
$$
v\in Vect(P), \ Lie_v \{ \ , \ \} =0
$$
over the subalgebra of Hamiltonian vector fields, $H^2 (P,\{ \ , \} )$ is the
quotient of the space of infinitesimal deformations of the Poisson bracket by
those obtained by infinitesimal changes of coordinates (i.e., by those of the
form $Lie_v \{ \ , \ \}$ for a vector field $v$).

On a symplectic manifold $(P,\{ \ , \} )$ Poisson cohomology coincides with the
de Rham one. The isomorphism is established by ``lowering the indices'':
for a cocycle $a=(a^{i_1\dots i_k}) \in \Lambda^k$ the $k$-form
$$
\sum_{i_1< \dots <i_k} \omega_{i_1 \dots i_k} dx^{i_1}\wedge \dots \wedge
dx^{i_k}, \quad
\omega_{i_1\dots i_k} = h_{i_1 j_1} \dots h_{i_k j_k} a^{j_1 \dots j_k}
$$
is closed. In particular, for $P$ = ball the Poisson cohomology is trivial. In
the general case ${\mathrm rk} (h^{ij}) < { \dim} ~ P$ the Poisson cohomology does not
vanish even locally (see \cite{lichn}).  We will prove now a simple criterion
of triviality of 1- and 2-cocycles. 

\begin{lemma}\label{first-lemma}
Let $h=(h^{ij}(x))$ be a Poisson structure of a constant rank
$2n<N$ on a sufficiently small ball $U$. 1). A one-cocycle $v=(v^i(x))\in
H^1(U,h)$ is trivial {\rm iff} the vector field $v$ is tangent to the leaves of
the symplectic foliation (\ref{foliation}).  2). A 2-cocycle $f=(f^{ij}(x))\in
H^2(U,h)$  is trivial {\rm iff} 
\beq\label{2-trivial}
f(dc', dc'')=0
\eeq
for arbitrary two Casimirs of $h$.
\end{lemma}

Of course, the statement of the lemma can be easily derived from
the results of \cite{lichn}. Nevertheless, we give a proof since we will
use similar arguments also in the infinite dimensional situation.

\pf 1). For a coboundary $v=\pal f$
and  for any Casimir $c$ of $h$ we have
$$
\pal_v c = \{ c, f\} =0.
$$
This means that $v$ is tangent to the symplectic leaves (\ref{foliation}). To
prove the converse statement let us choose the canonical coordinates 
$x=(y^1,
\dots, y^{2n}, c^1, \dots, c^k)$ reducing the bracket to the constant form
(\ref{normal}).
Here $c^1$, \dots, $c^k$ are independent Casimirs (\ref{casimir}). 
In these coordinates $v=(v^1, \dots, v^{2n}, 0, \dots, 0)$. The 1-form
$\omega=(\omega_1, \dots, \omega_{2n}, 0, \dots, 0)$ given by
$$
\omega_i =\sum_{j=1}^{2n} \bar h_{ij} v^j
$$
has the property
$$
d\omega|_{P_{c_0}} =0.
$$
Therefore a function $g$ locally exists s.t. 
$$
dg =\sum_{i=1}^{2n} \omega_idy^i + \sum_{a=1}^k \phi_a dc^a
$$
for some functions $\phi_1$, \dots, $\phi_k$. This function is the Hamiltonian
for the vector field $v$.

2). We will again use the canonical coordinates for $h$ as in the proof of the
first part.
For an exact 2-cocycle $f=\pal v$  and arbitrary two functions $c'$, $c''$
$$
f(dc',dc'') =-\{ c', v^i\} \pal_i c'' - \pal_i c' \{ v^i, c''\} .
$$
This is equal to zero if $c'$ and $c''$ are Casimirs of the bracket $h$.

To prove the converse statement we first consider, for every $a=1$, \dots, $k$,
the vector field $w$ (depending on $a$)
\beq\label{correction}
w^i=f^{ia}, ~~i=1, \dots, N.
\eeq
From (\ref{2-trivial}) it follows that $w$ is tangent to the symplectic leaves
of $h$.
The cocycle condition
\beq\label{2-computation}
0=[h,f]^{aij} =\pal_k f^{ai} h^{kj} + \pal_k f^{ja} h^{ki} =(\pal w)^{ij}
\eeq
implies $\pal w=0$. According to the first part of the lemma, there
exists a function $q^a(x)$ s.t. $w=\pal q^a$:
\beq\label{2-1}
f^{ia} =h^{ik} \pal_kq^a, ~~a=1, \dots, k.
\eeq
Let us now change the cocycle by a coboundary
$$
f\mapsto f+\pal z
$$
where the vector field $z$ is given by
\beq
z=\sum_{a=1}^k q^a \frac{\pal}{\pal c^a}.
\eeq
After such a change, due to (\ref{2-1}), we obtain
$$
f^{ia}=f^{ai}=0, ~~i=1, \dots, N.
$$
The rest of the proof repeats the arguments of the first part. The 2-form
$$
\omega_{ij} = \sum_{i,\, j=1}^{2n} \bar h_{ik}\bar h_{lj} f^{kl}
$$
is closed along the symplectic leaves. Hence there exists a 1-form
$\phi=(\phi_i)$ s.t.
$$
\omega =d\phi+\tilde\omega
$$
where every monomial in $\tilde\omega$ contains at least one $dc^a$ for some
$a$. Therefore
$$
f=\pal u
$$
for the vector field 
$$
u^i=\sum_{k=1}^{2n}\bar h^{ik} \phi_k, ~i=1, \dots, 2n, ~~~u^i=0, ~
i> 2n.
$$
The lemma is proved.\epf
\setcounter{equation}{0}
\setcounter{theorem}{0}

\vskip 0.5cm
\subsection{Formal loop spaces}\par

Let $M$ be a $n$-dimensional smooth manifold. Our aim  is to describe an
appropriate class of Poisson brackets on the loop space
$$
{\cal L}(M) = \{ S^1 \to M\}.
$$
In our definitions we will treat ${\cal L}(M)$ formally in the spirit of formal
variational calculus of \cite{dickey, dt}. We define the formal loop
space ${\cal L}(M)$ in terms of ring of functions on it. We also describe
calculus of differential forms and vector fields on the formal loop space. In
the next section we will also deal with multivectors on the formal loop space.

Let $U\subset M$ be a chart on $M$ with the coordinates $u^1, \dots, u^n$.
Denote ${\cal A} = {\cal A}(U)$ the space of polynomials in the independent
variables $u^{i,s}$, $i=1, \dots, n$, $s=1, 2, \dots$
\beq\label{dif.pol}
f(x; u; u_x, u_{xx}, \dots):=\sum_{m\geq 0} f_{i_1 s_1; \dots ;i_m s_m} (x;u) u^{i_1, s_1} \dots u^{i_m,
s_m}
\eeq
with the coefficients $f_{i_1 s_1; \dots ; i_ms_m}(x;u)$ being smooth functions
on $S^1\times M$. Such an expression will be called {\it differential
polynomial}.
We will often use an alternative notation for the independent
variables 
$$
u^i_x = u^{i,1}, ~~ u^i_{xx} = u^{i,2} , \dots
$$
Observe that polynomiality w.r.t. $u=(u^1, \dots, u^n)$ is not assumed.

The operator $\pal_x$ is defined as follows
\beq\label{delx}
\pal_x f =\frac{\pal f}{\pal x} +\frac {\pal f}{\pal u^i} u^{i,1} + \dots+
\frac{\pal f}{\pal u^{i,s} }u^{i,s+1}+\dots. 
\eeq
We will often use the notation
\beq\label{superscript}
f^{(k)}:= \pal_x^k f.
\eeq
The following identities will be
useful
\beq
\frac{\pal}{\pal u^i} \pal_x =\pal_x \frac{\pal}{\pal u^i} 
\eeq
\beq\label{comm-der}
\frac{\pal}{\pal u^{i,s}} \pal_x =\pal_x \frac{\pal}{\pal u^{i,s}}
+\frac{\pal}{\pal u^{i,s-1}}.
\label{commute}
\eeq
We define the space
$$
{\cal A}_{0,0} ={\cal A}/ {\bf R}, ~~{\cal A}_{0,1} = {\cal A}_{0,0}\, dx ,
$$
the operator
\beq\label{diff}
d: {\cal A}_{0,0} \to {\cal A}_{0,1}, ~~df := \pal_x f\, dx
\eeq
and the quotient
\beq\label{functionals}
\Lambda_0 ={\cal A}_{0,1}/d {\cal A}_{0,0}.
\eeq
The elements of the space $\Lambda_0$ will be written as integrals
over the circle $S^1$
\beq\label{functional}
\bar f :=\int f(x; u; u_x, u_{xx}, \dots) dx \in \Lambda_0
\eeq
We will
 use below the following simple statement.
\begin{lemma}\label{feqzero}
 If $\int fg\, dx=0$ for an arbitrary $g\in {\cal A}$ then $f\in
{\cal A}$ is equal to zero.
\end{lemma}
 
 The expressions (\ref{functional}) 
 are also called {\it local functionals} with the
{\it density} $f$. The space of local functionals  is the main building block 
of the ``space of functions''
on the formal loop space. The full ring ${\cal F} = {\cal F}({\cal L}(U))$ 
of functions on the formal loop space
by definition coincides with the tensor algebra of $\Lambda_0$
\beq\label{ring}
{\cal F} = {\bf R} \oplus \Lambda_0 \oplus \Lambda_0 \hat\otimes
\Lambda_0 \oplus \Lambda_0\hat\otimes \Lambda_0 \hat\otimes \Lambda_0
\oplus \dots
\eeq
Elements of $\Lambda_0^{\hat\otimes \, k}$ will be written as multiple integrals
of 
differential
polynomials of $k$ copies of the variables that we denote $u^i(x_1)$, \dots,
$u^i(x_k)$, $u_x^i(x_1)$, \dots, $u_x^i(x_k)$ etc.
\beq\label{k-functional}
\int f(x_1, \dots, x_k; u(x_1), \dots, u(x_k); u_x(x_1), \dots, u_x(x_k),
\dots) dx_1 \dots dx_k \in 
\Lambda_0^{\hat\otimes \, k}.
\eeq

The short exact sequence
$$
0\to {\cal A}_{0,0} \stackrel{d}\to {\cal A}_{0,1} \stackrel{\pi}\to
\Lambda_0 \to 0
$$
($\pi$ is the projection) is included in the {\it variational bicomplex}
$$
\begin{array}{ccccccccc}
 & & \Big\uparrow\vcenter{%
 \rlap{$\scriptstyle{\delta}$}} & & \Big\uparrow\vcenter{%
 \rlap{$\scriptstyle{\delta}$}} & & \Big\uparrow\vcenter{%
 \rlap{$\scriptstyle{\delta}$}} & & \\
0 & \to & {\cal A}_{2,0} & \stackrel{d}\rightarrow & {\cal A}_{2,1} &
\stackrel{\pi}\rightarrow & \Lambda_2 & \to & 0 \\
 & & \Big\uparrow\vcenter{%
 \rlap{$\scriptstyle{\delta}$}} & & \Big\uparrow\vcenter{%
 \rlap{$\scriptstyle{\delta}$}} & & \Big\uparrow\vcenter{%
 \rlap{$\scriptstyle{\delta}$}} & & \\
 0 & \to & {\cal A}_{1,0} & \stackrel{d}\rightarrow & {\cal A}_{1,1} &
\stackrel{\pi}\rightarrow & \Lambda_1 & \to & 0 \\
 & & \Big\uparrow\vcenter{%
 \rlap{$\scriptstyle{\delta}$}} & & \Big\uparrow\vcenter{%
 \rlap{$\scriptstyle{\delta}$}} & & \Big\uparrow\vcenter{%
 \rlap{$\scriptstyle{\delta}$}} & & \\
 0 & \to & {\cal A}_{0,0} & \stackrel{d}\rightarrow & {\cal A}_{0,1} &
\stackrel{\pi}\rightarrow & \Lambda_0 & \to & 0 \\
 & & \Big\uparrow\vcenter{%
 \rlap{$\scriptstyle{\delta}$}} & & \Big\uparrow\vcenter{%
 \rlap{$\scriptstyle{\delta}$}} & & \Big\uparrow\vcenter{%
 \rlap{$\scriptstyle{\delta}$}} & & \\
 & & 0 & & 0 & & 0 & &
 \end{array}
 $$
Here ${\cal A}_{k,l}$ are elements of the total degree $k+l$ in the Grassman
algebra
with the generators $\delta u^{i,s}$, $i=1, \dots, n$,
$s=0, 1, 2, \dots$ (observe the difference in the range of the second index
of $u^{i,s}$ and $\delta u^{i,s}$) and $dx$ with the coefficients in ${\cal A}$
having the degree $l$ in $dx$. We will
often identify
$$
\delta u^{i,0}=\delta u^i.
$$
For example, every $k$-form $\omega\in{\cal A}_{k,0}$ is a finite sum
\beq\label{omega}
\omega= \frac{1}{k!} \omega_{i_1s_1; \dots; i_ks_k} \delta u^{i_1, s_1}\wedge
\dots \wedge \delta u^{i_k s_k}
\eeq
where the coefficients $\omega_{i_1s_1; \dots; i_ks_k}\in {\cal A}$ are assumed
to be antisymmetric w.r.t. permutations of pairs $i_p, s_p \leftrightarrow i_q,
s_q$.

The exterior differential in the Grassman algebra is decomposed into a sum
$d+\delta$. 
The horizontal differential
$$
d: {\cal A}_{k,0} \to {\cal A}_{k,1}
$$
is defined by
\beq\label{d}
d\omega =dx\wedge\pal_x\omega
\eeq
where the derivation $\pal_x$,
$$
\pal_x (\omega_1 \wedge \omega_2) = \pal_x \omega_1 \wedge \omega_2 + \omega_1
\wedge \pal_x \omega_2
$$
is given by (\ref{delx}) on the coefficients of the differential form and
by
$$
\pal_x \delta u^{i,s} = \delta u^{i,s+1}.
$$
The elements of the quotient
$$
\Lambda_k = {\cal A}_{k,1}/d {\cal A}_{k,0}
$$
will be called {\it (local) $k$-forms} on the loop space. $k$-forms will also be
written by integrals
$$
\int dx\wedge\omega\in \Lambda_k, ~~ \omega \in {\cal A}_{k,0}.
$$
\begin{exam} Any one-form has a unique representative
\beq\label{1-form}
\phi = \int dx\wedge\phi_i \delta u^i 
\eeq
(use integration by parts).
\end{exam}

More generally, for every $k$-form $\omega$ written as in (\ref{omega}) 
$$
dx\wedge \omega \sim dx\wedge \tilde\omega ~({\mathrm mod} ~d({\cal A}_{k,0}))
$$
where
\beq\label{k-reduced}
\tilde\omega = \frac{!}{(k-1)!} \tilde\omega_{i_1; i_2s_2; \dots ; i_ks_k}
\delta u^{i_1}\wedge \delta u^{i_2, s_2} \wedge \dots \wedge \delta u^{i_k,
s_k}
\eeq
\eqa 
&&\tilde\omega_{i_1; i_2s_2; \dots ; i_ks_k} \nn\\
&&
= \frac{1}{k} \sum_{\begin{array}{c}0\leq r_l\leq
s_l\\ 2\leq l\leq k\end{array}}
\sum_{s\geq r_2+\dots +r_k} (-1)^s \left(\begin{array}{ccc}
 & s & \\ r_2 & \dots & r_k\end{array}\right) 
\omega_{i_1s; i_2,s_2-r_2; \dots;
 i_k,s_k-r_k}^{(s-r_2-\dots-r_k)}\nn
 \eeqa
 here
 \beq\label{multinomial}
 \left(\begin{array}{ccc} & s & \\ r_2 & \dots & r_k \end{array}\right)
 =\frac{s!}{r_2! \dots r_k! (s-r_2-\dots -r_k)!}
 \eeq
 stands for the multinomial coefficients. The coefficients 
 $\tilde\omega_{i_1; i_2s_2; \dots ;
 i_ks_k}$ will be called {\it reduced components} of $\omega$. They are
 antisymmetric w.r.t. pairs $i_2, s_2$, \dots, $i_k, s_k$ but with the
 permutation of $i_1$ and $i_2$ they behave as
 \eqa
&&\tilde\omega_{i_2; i_1s_2; \dots ; i_ks_k} 
\nn\\
&&
= \sum_{\begin{array}{c} 0\leq t_l\leq s_l\\ 3\le l\le k\end{array}}
 \sum_{t\geq s_2+t_3+\dots + t_k}
  (-1)^{t+1} \left(\begin{array}{ccc} & t  & \\
 s_2\, t_3 & \dots & t_k\end{array}\right) 
  \tilde\omega_{i_1;i_2 s_2; i_3,s_3-t_3;\dots;
 i_k,s_k-t_k}^{(t-s_2-t_3-\dots -t_k)}\nn\\
\label{k-antisymmetry}
\eeqa

We now define vertical arrows of the bicomplex. For a monomial 
$$
\omega = f\, \delta u^{i_1, s_1} \wedge \dots \wedge \delta u^{i_k, s_k}
$$
put
\beq
\delta\omega =\sum_{s\geq 0} \frac{\pal f}{\pal u^{j,t}} \delta u^{j,t} \wedge 
\delta u^{i_1,
s_1} \wedge \dots \wedge \delta u^{i_k, s_k},
\eeq
where we denote
$$
\frac{\pal}{\pal u^{j,0}} := \frac{\pal}{\pal u^j}.
$$
This gives vertical differential $\delta: {\cal A}_{k,0} \to {\cal A}_{k+1,
0}$,
$$
\delta^2 =0.
$$
The map $\delta$ on ${\cal A}_{k,1}$ is defined by
essentially same formula, $\delta dx =0$. Anticommutativity
$$
\delta d = - d\delta
$$
justifies action of $\delta$ on the quotient $\Lambda_k$.

\begin{exam} On $\Lambda_0$ the differential $\delta$ acts as follows
\beq\label{euler-lagrange}
\delta \int f\, dx =\int dx \wedge\left( \sum_s (-1)^s \pal_x^s \frac {\pal f} {\pal
u^{i,s}} \right)\delta u^i 
\eeq
(the Euler - Lagrange differential). We will use the notation
\beq\label{var.der}
\frac{\delta \bar f}{\delta u^i(x)} := \sum_s (-1)^s \pal_x^s \frac {\pal f} {\pal
u^{i,s}} 
\eeq
for the components of the 1-form, $\bar f = \int f\, dx$.
\end{exam}

\begin{theorem}[\cite{dt}] For $M$ = {\rm ball} both arrows and columns of the
variational bicomplex are exact.
\end{theorem}

\begin{exam}\label{zero-fun} 
A necessary and sufficient condition for 
$$
\frac{\delta\bar f}{\delta u^i(x)}=0, ~~i=1, \dots, n.
$$
is the existence of a differential polynomial $g=g(x;u;u_x;\dots)$ such that
$f=\pal_x g$.
\end{exam}

Let us now consider the space $\Lambda^1$ of vector fields on the formal loop space. These will be
formal infinite sums
\beq\label{vector}
\xi = \xi^0 \frac{\pal}{\pal x} + 
\sum_{k\geq 0} \xi^{i,k} 
\frac{\pal}{\pal u^{i,k}}, ~~ \xi^{i,k} \in {\cal A}
\eeq
where we denote
$$
\frac{\pal }{\pal u^{i,0}} := \frac{\pal}{\pal u^i}.
$$
The derivative of a functional $\bar f=\int f(x; u; u_x, \dots) dx \in \Lambda_0$
along $\xi$ reads
$$
\xi \,\bar f := \int \left( \xi^0 \frac{\pal f}{\pal x} 
+ \sum \xi^{i,k} 
\frac{\pal f}{\pal u^{i,k}}\right) dx.
$$
The Lie bracket of two vector fields is defined by
\eqa
&&[\xi,\eta] =(\xi^0 \eta^0_x - \eta^0 \xi^0_x
+\xi^{j,t}\frac{\pal\eta^0}{\pal u^{j,t}}-
\eta^{j,t}\frac{\pal\xi^0}{\pal u^{j,t}}) \frac {\pal}{\pal x}
\nn\\
&&\quad
+ \sum _{s\geq 0}\left( \xi^0 \frac{\pal \eta^{i,s}}{\pal x} 
- \eta^0 \frac{\pal\xi^{i,s}}{\pal x}
+ \xi^{j,t}\frac {\pal \eta^{i,s}} {\pal u^{j,t}}
-\eta^{j,t} \frac{\pal \xi^{i,s}}{\pal u^{j,t}}\right) \frac{\pal}{\pal
u^{i,s}}\label{commutator}
\eeqa

{\it Evolutionary vector fields} $a$ are defined by the conditions 
of vanishing of the $\pal/\pal x$-component and the commutativity
$$
[\pal_x, a]=0
$$
They are parameterized by $n$-tuples $a^1$, \dots, $a^n$ of elements of
${\cal A}$ as follows 
\beq\label{evolution}
a= \sum _{s\geq 0} \pal_x^s a^i \frac {\pal}{\pal u^{i,s}}.
\eeq
The corresponding system of evolutionary PDEs reads
\beq\label{pde}
u^i_t =a^i(x; u; u_x, u_{xx}, \dots), ~~ i=1, \dots, n.
\eeq
In particular, an evolutionary vector field $a$ is called {\it translation
invariant} if the coefficients $a^i$ do not depend explicitly on $x$,
$$
\frac{\pal a^i}{\pal x}=0, ~~ i=1, \dots, n.
$$

The contraction $i_\xi \omega$ of a $k$-form $\omega\in {\cal A}_{k,0}$ given
by (\ref{omega})
and a vector field $\xi$
is a $(k-1)$-form defined by  
\beq\label{contraction}
i_\xi\omega =\frac{1}{(k-1)!}\xi^{j,t} \omega_{jt; i_1 s_1; \dots ;i_{k-1} s_{k-1}}
\delta u^{i_1, s_1}\wedge \dots \wedge \delta u^{i_{k-1}, s_{k-1}}.
\eeq
As usual
$$
i_\xi i_\eta =-i_\eta i_\xi
$$
for two vector fields $\xi$, $\eta$.
For a form $\omega\in {\cal A}_{k,1}$ the contraction $i_\xi \omega\in {\cal
A}_{k-1,1}$ is defined by essentially same formula provided the vector field
$\xi$ contains no $\pal / \pal x$-term. It is an easy exercise to check, using
(\ref{commute}), that for an evolutionary vector field $a$
\beq\label{dick}
d i_a + i_a d =0.
\eeq
It readily follows that contraction with evolutionary vector fields is
well-defined on the quotient $i_a: \Lambda_k \to \Lambda_{k-1}$. A more strong
statement holds true
\begin{lemma} Let $\omega \in {\cal A}_{k,1}$. It belongs to $d({\cal
A}_{k,0})$ {\rm iff} $i_a \omega \in d( {\cal A}_{k-1,0})$ for an arbitrary
evolutionary vector field $a$.
\end{lemma}

\pf We use induction in $k$. For $k=1$ we can choose a representative
of the class of $\omega\wedge dx\in \Lambda_1$ with the 1-form $\omega$
given by (\ref{1-form}). The contraction reads
$$
i_a(dx\wedge\omega)=- \int dx \, a^i\omega_i.
$$
Using Lemma \ref{feqzero} we obtain $\omega_i=0$ for all $i$. 

Let us assume validity of the lemma for any $(k-1)$-form. We will prove
that the condition $i_a \omega\wedge dx =0 \in \Lambda_{k-1}$ implies
vanishing of all the reduced components (\ref{k-reduced}). By induction the
above condition is equivalent to
$$
i_{b_2} \dots i_{b_k}i_a \omega\wedge dx\in d({\cal A}_{0,0})
$$
for arbitrary evolutionary vector fields $b_2$, \dots, $b_k$. Integrating by
parts we rewrite the last line in the form
\beq\label{i1}
\int a^i \phi_i \, dx=0 
\eeq
where
$$
\phi_i =k\,\tilde \omega_{i; i_2 s_2; \dots ;i_k s_k} \pal_x^{s_2}
b_2^{i_2} \dots \pal_x^{s_k} b_k^{i_k}.
$$
From (\ref{i1}) it follows that $\phi_i=0$ for all $i$. Since $b_2^i$, \dots,
$b_k^i$ are arbitrary differential polynomials, this implies 
$\tilde \omega_{i; i_2 s_2; \dots ;i_k s_k}=0$. That means that the form
$\omega\wedge dx$ is equivalent to zero, modulo $d({\cal A}_{k,0})$. The lemma
is proved.\epf

\begin{cor} A form $\omega\in {\cal A}_{k,1}$ belongs to $d({\cal A}_{k,0})$
{\rm iff} 
$$
i_{a_1} \dots i_{a_k} \omega \in d({\cal A})
$$
for arbitrary evolutionary vector fields $a_1$, \dots, $a_k$.
\end{cor}

\begin{exam} For the one-form $\omega =\delta \int f\, dx$ the contraction $i_a
\omega$ reads
$$
i_a\omega = \int a^i \frac {\delta \bar f}{\delta u^i(x)} dx \in \Lambda_0.
$$
This is the time derivative of the functional $\bar f =\int f\, dx$ w.r.t. the
evolutionary system (\ref{pde}).
\end{exam}

\begin{exam} For a one-form
$\omega =\omega_i \delta u^i\wedge dx\in \Lambda_1$
the condition of closedness 
$$
\delta \omega =0\in \Lambda_2
$$
reads  
\beq\label{one-main}
\frac{\pal \omega_i} {\pal u^{j,s}} =\sum_{t\geq s} (-1)^t 
\left(\begin{array}{c}  t \\ s\end{array}\right)
\pal_x^{t-s} \frac{\pal \omega_j}{\pal u^{i,t}}
\eeq
for any $i,\, j=1, \dots, n$, $s=0,\, 1, \, \dots$. This is the classical
Volterra's criterion \cite{volterra} for the system of ODEs
$$
\omega_1(x; u; u_x, u_{xx}, \dots)=0, \dots, \omega_n(x; u; u_x, u_{xx}, \dots
)=0
$$
to be locally representable in the Euler - Lagrange form 
$$
\omega_i =\frac{\delta \bar f}{\delta u^i(x)}, ~~i=1, \dots, n
$$
(use exactness of the
variational bicomplex).
\end{exam}

\begin{exam} For a 2-form
$$
\omega=\frac{1}{2}\omega_{is; jt} \delta u^{i,s}\wedge \delta u^{j,t},
$$
the contraction with two evolutionary vector fields $a$ and $b$ can be
represented in the form
\beq\label{2-form.1}
i_a i_b \omega =-2 \int a^i \tilde \omega _{i; js} \pal_x^s b^j \, dx
\eeq
where
$$
\tilde\omega_{i;js} =\frac{1}{2}\sum_{r=0}^s \sum_{t\geq s-r} (-1)^t 
\left(\begin{array}{c}t \\ 
s-r\end{array}\right)
\pal_x^{t-s+r} \omega _{is; jr}
$$
are the reduced components (\ref{k-reduced}).
The tilde will be omitted in the subsequent formulae.
According to this we will often represent 2-forms in the
reduced form
\beq\label{2-reduced}
dx\wedge\omega =  \omega_{i;js} dx\wedge\delta u^i \wedge \delta u^{j,s}.
\eeq
The reduced coefficients must satisfy the antisymmetry conditions
(\ref{k-antisymmetry}). They are spelled out as follows
\beq\label{2-antisymmetry.forms}
\omega_{i;js} =\sum_{t\geq s} (-1)^{t+1} \left(\begin{array}{c}t \\
s\end{array}\right) \partial_x^{t-s}
\omega_{j;it}
\eeq
(integrate by parts in (\ref{2-form.1}) and use arbitraryness of $a$ and $b$).
\end{exam}

\begin{exam} A 2-form $dx\wedge\omega =\delta dx\wedge \phi$ for 
$$
\phi=\phi_i \delta u^i
$$
has the reduced representative (\ref{2-reduced}) with
\beq\label{2-form.exact}
\omega_{i;js} = \frac12\,\left(\frac{\pal \phi_i}{\pal u^{j,s}} +\sum_{t\geq s} (-1)^{t+1}
\left(\matrix{ t \cr s\cr}\right) \pal_x^{t-s} \frac{\pal \phi_j}{\pal
u^{i,t}}\right).
\eeq
\end{exam}

\begin{exam} A 2-form (\ref{2-reduced}) is closed, $\delta \omega =0$, {\rm iff}
\beq\label{2-form.closed}
\left( \sum_{m=s}^{t+s}\sum_{r=0}^{m-s} +\sum_{m\geq t+s+1} \sum_{r=0}^t\right)
(-1)^m \left(\begin{array}{c} m \\ r \, s\end{array}\right) \pal_x^{m-r-s} \frac{\pal \omega_{j; k, t-r}}{\pal
u^{i,m}}
+ \frac{\pal \omega_{i;j,s}}{\pal u^{k,t}} -\frac{\pal \omega _{i;k,t}}{\pal
u^{j,s}} =0
\eeq
for any $i$, $j$, $k=1, \dots n$, $s=0, \, 1,\, 2, \dots$.
\end{exam}
\pf
By definition 
$$
\delta (\omega)=\sum\frac{\pal \omega_{i;js}}{\pal u^{k,l}} 
\delta u^i\wedge \delta u^{j,s}\wedge
\delta u^{k,l}\wedge dx.
$$
So $\delta\omega=0$ means that for any three evolutionary vector fields
$$
a=\sum (a^i)^{(s)}\frac{\pal}{\pal u^{i,s}},\quad
b=\sum (b^i)^{(s)}\frac{\pal}{\pal u^{i,s}},\quad
c=\sum (c^i)^{(s)}\frac{\pal}{\pal u^{i,s}}
$$
the contraction $i_ai_bi_c \delta (dx\wedge \omega) \in d({\cal A}_{0,0})$, i,e,
\eqa
&&\int \frac{\pal \omega_{i;j,s}}{\pal u^{k,l}}\left[a^i\,(b^k)^{(l)}\,(c^j)^{(s)}-
a^i\,(b^j)^{(s)}\,(c^k)^{(l)}-(a^k)^{(l)}\,b^i\,(c^j)^{(s)}+
(a^k)^{(l)}\,(b^j)^{(s)}\,c^i\right.\nn\\
&&\left.
+
(a^j)^{(s)}\,b^i\,(c^k)^{(l)}-
(a^j)^{(s)}\,(b^k)^{(l)}\,c^i\right]dx=0
\eeqa
Using integration by parts we get
\eqa
&&\int \frac{\pal \omega_{i;j,s}}{\pal u^{k,l}}\left[a^i\,(b^k)^{(l)}\,(c^j)^{(s)}-
a^i\,(b^j)^{(s)}\,(c^k)^{(l)}\right]\nn\\
&&+
\sum (-1)^{m+1}
\left( \begin{array}{c}  m  \\ l \, r
\end{array}\right)
\pal_x^{m-l-r}\left(\frac{\pal\omega_{i;j,s}}{\pal u^{k,m}}
\right)
\left(a^k\,(b^i)^{(l)}\,(c^j)^{(s+r)}-a^k\,(b^j)^{(s+l)}\,(c^i)^{(r)}\right)
\nn\\
&&+
\sum (-1)^{m}
\left(\begin{array}{c}  m  \\ s \, r \end{array}\right)
\pal_x^{m-s-r}\left(\frac{\pal\omega_{i;j,m}}{\pal u^{k,l}}
\right)
\left(a^j\,(b^i)^{(s)}\,(c^k)^{(l+r)}-a^j\,(b^k)^{(l+s)}\,(c^i)^{(r)}\right)dx
\nn\\
&&=0
\nn
\eeqa
The above identity is equivalent to
\eqa
&&\frac{\pal \omega_{i;j,s}}{\pal u^{k,l}}-\frac{\pal \omega_{i;k,l}}{\pal u^{j,s}}\nn\\
&&+
\left(\sum_{m=l}^{l+s}\sum_{r=0}^{m-l}+\sum_{m\ge l+s+1}\sum_{r=0}^s\right)
(-1)^{m+1}
\left( \matrix{
m\cr l \, r}\right)
\pal_x^{m-l-r}\left(\frac{\pal\omega_{k;j,s-r}}
{\pal u^{i,m}}\right)\nn\\
&&+
\left(\sum_{m=s}^{l+s}\sum_{r=0}^{m-s}+\sum_{m\ge l+s+1}\sum_{r=0}^l\right)
(-1)^{m}
\left(\begin{array}{c} m  \\ s \, r\end{array}\right)
\pal_x^{m-s-r}\left(\frac{\pal\omega_{j;k,l-r}}
{\pal u^{i,m}}\right)\nn\\
&&+
\left(\sum_{m=l}^{l+s}\sum_{r=0}^{m-l}+\sum_{m\ge l+s+1}\sum_{r=0}^s\right)
(-1)^{m}
\left(\begin{array}{c}  m  \\ l \, r\end{array}\right)
\pal_x^{m-l-r}\left(\frac{\pal\omega_{k;i,m}}
{\pal u^{j,s-r}}\right)
\nn\\
&&+
\left(\sum_{m=s}^{l+s}\sum_{r=0}^{m-s}+\sum_{m\ge l+s+1}\sum_{r=0}^l\right)
(-1)^{m+1}
\left(\begin{array}{c}  m  \\ r \, s\end{array}\right)
\pal_x^{m-s-r}\left(\frac{\pal\omega_{j;i,m}}
{\pal u^{k,l-r}}\right)=0\nn
\eeqa
Now by using the antisymmetry condition (\ref{2-antisymmetry.forms})
we see that the third term in the 
above sum is equal to the forth term, and by using the identity
(\ref{comm-der}) and the antisymmetry condition
(\ref{2-antisymmetry.forms}) again
we see that the last two terms equal to the second and first term 
respectively. Thus we
arrive at the proof of (\ref{2-form.closed}).\epf

\begin{remark} 
The equations (\ref{2-form.closed}) were derived by Dorfman in the theory of
the so-called
local symplectic structures \cite{dorf1}, see also the book
\cite{dorf2}.
\end{remark}

\begin{cor}\label{exact-2form} Any solution to (\ref{2-form.closed}) satisfying
(\ref{2-antisymmetry.forms}) can be locally represented in the form
(\ref{2-form.exact}).
\end{cor}
This follows from exactness of the variational bicomplex.

We will briefly outline necessary points of the global picture of functionals,
differential forms and vector fields on the formal loop space ${\cal L}(M)$
for a general smooth manifold $M$ (i.e., not only for a ball). The
corresponding objects must be defined for any chart $U\subset M$ as it was
explained above. On the intersections $U\cap V$ they must satisfy certain
consistency conditions. For example, functionals in the charts
$U$, $V$ with the coordinates $u^1$, \dots, $u^n$ and $v^1$, \dots, $v^n$ 
are defined by densities
$f_U(x; u; u_x, \dots)$ and $f_V(x; v; v_x, \dots)$ s.t.
$$
f_V(x; v(u); \frac{\pal v}{\pal u}u_x, \dots ) dx = f_U (x; u; u_x, \dots )dx 
~ ({\mathrm mod} ~{\mathrm Im}\, d).
$$
Such objects comprise the space $\Lambda_0(M)$. As above, we obtain the 
ring of 
functions on the formal loop space taking the tensor algebra of $\Lambda_0(M)$.
One-forms in the charts $U$, $V$ are described by their reduced components
$\omega_i^U$ and $\omega_a^V$ s.t., on $U\cap V$
transform as components of a covector
$$
\omega_a^V(x; v(u); \frac{\pal v}{\pal u}u_x, \dots) =\omega_i^U(x; u; u_x,
\dots) \frac{\pal u^i}{\pal v^a}
$$
etc. The components of evolutionary vector fields transform like vectors
\beq\label{1-transform}
a^k_V(x; v(u); \frac{\pal v}{\pal u} u_x, \dots) = \frac{\pal v^k}{\pal u^i}
a^i_U (x; u; u_x, \dots).
\eeq
The contraction $i_a dx\wedge \omega$ of a 1-form with an evolutionary vector
field is well-defined as an element of $\Lambda_0(M)$. 

The global theory of the variational bicomplex was developed in \cite{takens},
\cite{vinograd}.

\setcounter{equation}{0}
\setcounter{theorem}{0}

\vskip 0.5cm
\subsection{Local multivectors and local Poisson brackets}\label{sec-2-3}\par 

We first define more general, i.e. non-local $k$-vectors as elements of
$(\Lambda^1)^{\wedge k}$. They will be written as infinite sums of expressions
of the form
\eqa
&&\alpha =\frac1{k!}\,
\alpha^{i_1s_1; \dots i_ks_k} (x_1, \dots, x_k; u(x_1), \dots , u(x_k);
u_x(x_1), \dots, u_x(x_k), \dots) 
\nn\\
&&\quad\times\frac{\pal}{\pal u^{i_1, s_1}(x_1)}
\wedge \dots \wedge\frac{\pal}{\pal u^{i_k, s_k}(x_k)}
\label{k-vector}
\eeqa
(in this subsection we will consider only multivectors not containing
$\pal/\pal x$). The coefficients must satisfy the antisymmetry condition
w.r.t. simultaneous permutations
$$
i_p, s_p, x_p \leftrightarrow i_q, s_q, x_q.
$$
The exterior algebra structure on multivectors is introduced in a usual way:
the product of a $k$-vector  $\alpha$ by a $l$-vector $\beta$ is a
$(k+l)$-vector
\eqa
&&(\alpha\wedge\beta)^{i_1 s_1; \dots ;i_k s_k;i_{k+1}s_{k+1}; \dots ;
i_{k+l}s_{k+l}}
(x_1,  \dots, x_{k+l}; u(x_1), \dots, u(x_{k+l}); \dots)
\nn\\
&&=\frac{1}{k!l!} \sum_{\sigma\in S_{k+l}}(-1)^{{\mbox {sgn}}\sigma}\,
\alpha^{i_{\sigma(1)}s_{\sigma(1)};\dots; i_{\sigma(k)}
s_{\sigma(k)}}(x_{\sigma(1)}, \dots, x_{\sigma(k)}; u(x_{\sigma(1)}), \dots,
u(x_{\sigma(k)}); \dots) 
\nn\\
&&\times\beta^{i_{\sigma(k+1)}s_{\sigma(k+1)}; \dots; i_{\sigma(k+l)}
s_{\sigma(k+l)}}
(x_{\sigma(k+1)}, \dots, x_{\sigma(k+l)}; u(x_{\sigma(k+1)}), \dots,
u(x_{\sigma(k+l)}); \dots )\nn\\
\eeqa

\begin{exam} Lie derivative of a $k$-vector $\alpha$ (\ref{k-vector})
along a vector field (\ref{vector}) reads
\eqa
&&Lie_\xi \alpha^{i_1s_1; \dots ;i_ks_k}=
\nn\\&&\quad
\sum_{p=1}^k \left[\xi^0 (x_p; u(x_p); \dots) \frac{\pal}{\pal x_p} 
\alpha^{i_1s_1; \dots ;i_ks_k} + \xi^{j_p,t_p}(x_p; \dots) \frac{\pal}{\pal
u^{j_p, t_p}(x_p)} \alpha^{i_1s_1; \dots ;i_ks_k}\right]
\nn\\
&&\quad
-\sum_{p=1}^k \frac{\pal \xi^{i_p,s_p}(x_p; \dots)}{\pal u^{j_p, t_p}(x_p)}
\alpha^{i_1s_1; \dots; i_{p-1}s_{p-1};j_p t_p; \dots ;i_ks_k}.
\label{k-lie}
\eeqa
here we assume that $\xi^0$ does not depend on $u^{j,t}$.
\end{exam}

\noindent
{\bf Definition}. A $k$-vector $\alpha$ is called {\it translation invariant}
if
$$
Lie _{\pal_x}\alpha=0
$$
and
$$
\left( \frac{\pal}{\pal x_1} + \dots +\frac{\pal}{\pal x_k}\right) \alpha=0.
$$

\begin{lemma} Every translation invariant $k$-vector $\alpha$ has coefficients
of the form
\eqa
&&\alpha^{i_1s_1; \dots ;i_ks_k}(x_1, \dots, x_k; u(x_1), \dots, u(x_k); \dots)
\nn\\
&&\quad
=\pal_{x_1}^{s_1} \dots \pal_{x_k} ^{s_k} A^{i_1\dots i_k}(x_1, \dots, x_k;
u(x_1), \dots u(x_k); \dots)\label{k-tr.inv.}
\eeqa
where the differential polynomials $A^{i_1\dots i_k}(x_1, \dots, x_k;u(x_1),
\dots, u(x_k); \dots)$ are antisymmetric w.r.t. simultaneous permutations
$$
i_p, x_p \leftrightarrow i_q, x_q
$$
and also they satisfy
$$
A^{i_1\dots i_k}(x_1+t, \dots, x_k+t; u(x_1), \dots, u(x_k); \dots)=
A^{i_1\dots i_k}(x_1, \dots, x_k; u(x_1), \dots, u(x_k); \dots)
$$
for any $t$.
\end{lemma}

The functions $A^{i_1\dots i_k}(x_1, \dots; x_k; u(x_1), \dots, u(x_k); \dots)$
will be called {\it components} of the translation invariant $k$-vector
$\alpha$.

Translation invariant multivectors form a graded Lie subalgebra of the full
graded Lie algebra of multivectors
closed w.r.t.
Schouten - Nijenhuis bracket. 

\begin{exam} The Lie derivative of a bivector $\alpha$ with the components
\newline
$A^{ij}(x-y; u(x), u(y); \dots)$ along a translation invariant vector field $a$
with the components $a^i(u; u_x, \dots)$ has the components
\eqa
&&
Lie_a \alpha^{ij} =\pal_x^t a^k(u(x); \dots) 
\frac{\pal A^{ij}}{\pal u^{k,t}(x)}
+\pal_y^ta^k (u(y); \dots) \frac{\pal A^{ij}}{\pal u^{k,t}(y)}
\nn\\
&&\quad
-\frac{\pal a^i(u(x); \dots)}{\pal u^{k,t}(x)} \pal_x^t A^{kj} -\frac{\pal
a^j(u(y); \dots)}{\pal u^{k,t}(y)} \pal_y^t A^{ik}.
\label{lie-der}
\eeqa
\end{exam}

\begin{exam} Let $\alpha$, $\beta$ be two translation invariant bivectors
with the components $A^{ij}(x-y; u(x), u(y); u_x(x),u_x(y), \dots)$ and 
$B^{ij}(x-y; u(x),u(y);
u_x(x),u_x(y); \dots)$ that we redenote resp. $A^{ij}_{x,y}$ and $B^{ij}_{x,y}$
for brevity. The Schouten - Nijenhuis bracket $[\alpha,\beta]$ is a translation invariant
trivector with the components
\eqa
&&[\alpha,\beta]^{ijk}_{x,y,z} =
\frac{\pal A^{ij}_{x,y}}{\pal u^{l,s}(x)} \pal_x^s B^{lk}_{x,z}
+ \frac{\pal B^{ij}_{x,y}}{\pal u^{l,s}(x)} \pal_x^s A^{lk}_{x,z}
+\frac{\pal A^{ij}_{x,y}}{\pal u^{l,s}(y)} \pal_y^s B^{lk}_{y,z}
+ \frac{\pal B^{ij}_{x,y}}{\pal u^{l,s}(y)} \pal_y^s A^{lk}_{y,z}
\nn\\
&&\quad
+\frac{\pal A^{ki}_{z,x}}{\pal u^{l,s}(z)} \pal_z^s B^{lj}_{z,y}
+ \frac{\pal B^{ki}_{z,x}}{\pal u^{l,s}(z)} \pal_z^s A^{lj}_{z,y}
+\frac{\pal A^{ki}_{z,x}}{\pal u^{l,s}(x)} \pal_x^s B^{lj}_{x,y}
+ \frac{\pal B^{ki}_{z,x}}{\pal u^{l,s}(x)} \pal_x^s A^{lj}_{x,y}
\nn\\
&&\quad
+\frac{\pal A^{jk}_{y,z}}{\pal u^{l,s}(y)} \pal_y^s B^{li}_{y,x}
+ \frac{\pal B^{jk}_{y,z}}{\pal u^{l,s}(y)} \pal_y^s A^{li}_{y,x}
+\frac{\pal A^{jk}_{y,z}}{\pal u^{l,s}(z)} \pal_z^s B^{li}_{z,x}
+ \frac{\pal B^{jk}_{y,z}}{\pal u^{l,s}(z)} \pal_z^s A^{li}_{z,x}.
\label{schouten}
\eeqa
\end{exam}

For a translation invariant $k$-vector $\alpha$ and k 1-forms $\omega^1$,
\dots, $\omega^k$,
$$
\omega^j= \omega^j_{is}\delta u^{i,s}\wedge dx \in {\cal A}_{1,1}, ~~j=1,
\dots, k
$$
the contraction
\eqa
&&<\alpha, \omega^1\wedge \dots \wedge \omega^k>
\nn\\
&&\quad
:=\frac{1}{k!}\int
\sum_{\sigma \in S_k} (-1)^{{\mbox{sgn}} \, \sigma}
\omega_{i_1s_1}^{\sigma(1)} (x_1; u(x_1); \dots) \dots
\omega_{i_ks_k}^{\sigma(k)}(x_k; u(x_k); \dots)
\nn\\
&&\quad
\alpha^{i_1s_1; \dots; i_ks_k}(x_1, \dots, x_k; u(x_1), \dots, u(x_k); \dots)
dx_1 \dots dx_k
\in \Lambda_0^{\hat\otimes k}
\label{k-bracket.1}
\eeqa
is well defined on $\Lambda_1^{\otimes k}$.

\begin{exam} The value of a translation invariant $k$-vector $\alpha$ with the
components $A^{i_1 \dots i_k}$ on the 1-forms $\delta \bar f^1$, \dots, $\delta
\bar f^k$ equals
\eqa
&&<\alpha, \delta\bar f^1\wedge \dots \wedge \delta \bar f^k>
\nn\\
&&\quad
=\int \frac{\delta \bar f^1}{\delta u^{i_1}(x_1)} \dots \frac{\delta \bar
f^k}{\delta u^{i_k}(x_k)} A^{i_1 \dots i_k} (x_1, \dots , x_k; u(x_1), \dots,
u(x_k);\dots) dx_1 \dots dx_k\nn\\
&&\quad  \in \Lambda_0^{\hat\otimes k}.
\label{k-bracket.2}
\eeqa
\end{exam}

The transformation law of components of translation invariant multivectors
w.r.t. changes of coordinates on the intersection of two coordinate charts
$(U, u^1, \dots, u^n)$ and $(V, v^1, \dots, v^n)$ is analogous to the
transformation law of components of multivectors on a finite-dimensional
manifolds:
\eqa
&&
A_V^{a_1 \dots a_k} (x_1, \dots x_k; v(u(x_1)), \dots, v(u(x_k));
\frac{\pal v}{\pal u} u_x(x_1), \dots, \frac{\pal v}{\pal u}u_x(x_k), \dots)
\nn\\
&&
= \frac{\pal v^{a_1}}{\pal u^{i_1}}(x_1)\dots 
\frac{\pal v^{a_k}}{\pal u^{i_k}}(x_k)
A_U^{i_1 \dots i_k} (x_1, \dots, x_k; u(x_1), \dots, u(x_k); u_x(x_1),\dots,
u_x(x_k), \dots).\nn\\
\label{transf.multivector}
\eeqa

We now proceed to the main definition of {\it local multivectors}. They are
translation invariant multivectors $\alpha$  such that their dependence on
$x_1$, \dots, $x_k$ is given by a finite order distribution with the support
on the diagonal $x_1=x_2=\dots=x_k$
\beq\label{k-local}
A^{i_1 \dots i_k}= \sum_{p_2,p_3,\dots , p_k \geq 0} 
B^{i_1 \dots i_k}_{p_2,\dots,p_k}(u(x_1);
u_x(x_1), \dots) \delta^{(p_2)}(x_1-x_2) \delta^{(p_3)}(x_1-x_3) \dots
\delta^{(p_k)} (x_1-x_k).
\eeq
The coefficients $B^{i_1 \dots i_k}{p_2,\dots,p_k}(u(x_1);
u_x(x_1), \dots) $ are differential polynomials in ${\cal A}$ not depending
explicitly on $x$.
All the sums at the moment are assumed to be finite. In the next section we
will relax this condition. Delta functions and their derivatives and products are defined by
the formulae
\beq\label{delta}
\int f(y) \delta(x-y) dy=f(x), ~~\int f(y) \delta^{(p)}(x-y) dy = f^{(p)}(x)
\eeq
\eqa
&&\int f(x_1, \dots, x_k) \delta^{(p_2)}(x_1-x_2) \delta^{(p_3)}(x_1-x_3) 
\dots
\delta^{(p_k)} (x_1-x_k) dx_2 \dots dx_k \nn\\
&&\quad =\pal_{x_2}^{p_2}\dots \pal_{x_k}^{p_k}
f(x_1, \dots, x_k)|_{x_1=x_2=\dots =x_k}.\nn
\eeqa

\begin{lemma} The value (\ref{k-bracket.1}) of a local $k$-vector $\alpha$ on k
1-forms $\omega^1$, \dots, $\omega^k$ is given by
\eqa
&&
<\alpha, \omega^1\wedge\dots\wedge\omega^k>
\nn\\
&&\quad
=\int B^{i_1\dots i_k}_{p_2 \dots p_k} (u; u_x, u_{xx}, \dots)
\omega^1_{i_1}(x; u; u_x, \dots)\pal_x^{p_2} \omega_{i_2}^2(x; u; u_x, \dots)
\nn\\
&&\quad\quad\dots \pal_x^{p_k} \omega_{i_k}^k(x; u; u_x, \dots)\, dx
\in \Lambda_0.\label{value-kform}
\eeqa
It gives a well-defined polylinear map
$$
\alpha:\Lambda_1^{\otimes k} \to \Lambda_0.
$$
\end{lemma}

In calculations with local multivectors various simple identities for
delta-functions will be useful. All of them are simple consequences of the
definition (\ref{delta}). First,
\beq\label{binom}
f(y) \delta^{(p)}(x-y) =\sum_{q=0}^p \left(\matrix{ p\cr q\cr}\right)
f^{(q)}(x)\delta ^{(p-q)}(x-y).
\eeq
Next,
\eqa
&&
\delta(x_1-x_2) \dots \delta(x_1-x_k) = \delta (x_2-x_1)\delta (x_2-x_3) \dots
\delta (x_2-x_k) = \dots \nn\\
&&= \delta(x_{k}-x_1) \dots \delta(x_{k}-x_{k-1}).\label{diagonal}
\eeqa
Differentiating (\ref{diagonal}) w.r.t. $x_1$, \dots, $x_k$ we will obtain
relations between products of derivatives of delta-functions.

We leave as a simple exercise for the reader to prove that the space of local
multivectors  that we denote
$$
\Lambda^*_{loc} =\oplus \Lambda^k_{loc}
$$
is closed w.r.t. the Schouten - Nijenhuis bracket. Warning: this space is not
closed w.r.t. the exterior product! Because of this we were to introduce
a wider algebra of multivectors to introduce the definition of the Schouten -
Nijenhuis bracket according to the rules one uses in the finite dimensional case.

\begin{exam} The component of a local bivector $\varpi$  has the form
\beq\label{bivector}
\varpi^{ij}=\sum_{s\geq 0} A^{ij}_s(u(x); u_x(x), \dots) \delta^{(s)}(x-y).
\eeq
The value of the bivector on two 1-forms $\phi=\phi_i \delta u^i\wedge dx$
and $\psi = \psi_i \delta u^i \wedge dx$ equals
\beq\label{2-local}
\int \phi_i A^{ij}_s \pal_x^s \psi_j dx.
\eeq
The conditions of antisymmetry of the bivector reads
\beq\label{2-antisymmetry}
A^{ji}_s = \sum_{t\geq s} (-1)^{t+1} \left(\matrix{t \cr s\cr}\right)
\pal_x^{t-s}A^{ij}_t.
\eeq
\end{exam}

\pf Let us explain how to prove the antisymmetry condition
(\ref{2-antisymmetry}). We must have
$$
\sum_s A^{ji}_s(u(y); \dots) \delta^{(s)}(y-x) = -
\sum A^{ij}_s(u(x); \dots) \delta^{(s)}(x-y).
$$
Using $\delta^{(s)}(y-x) =(-1)^s \delta^{(s)}(x-y)$ and (\ref{binom}) we obtain
(\ref{2-antisymmetry}).\epf

\begin{remark} One can represent the bivector as
\beq\label{operator}
A^{ij}(u(x); u_x(x), \dots; \frac{d}{dx}) \delta(x-y).
\eeq
Here the differential operators $A^{ij}$ are
$$
A^{ij}(x; u(x); u_x(x), \dots; \frac{d}{dx}) =\sum_s A^{ij}_s \frac{d^s}{dx^s}.
$$
For local multivectors of higher rank the language of differential operators was used
by Olver \cite{olver}.
\end{remark}

\begin{exam}\label{schouten-bivector-sym} 
The Schouten - Nijenhuis bracket of the bivector
$$
\varpi=h^{ij}\delta(x-y),
$$
where $h^{ij}$ is a constant antisymmetric matrix,
with $\alpha$ of the form (\ref{bivector}) 
reads
\eqa
&&
\hskip -1.0 truecm 
[\varpi,\alpha]^{ijk}_{x,y,z} = \left[\frac{\pal A^{ij}_t}{\pal u^{l,s}} h^{lk}
+\sum (-1)^{q+r+s}\left(\matrix{q+r+s\cr q ~~ r\cr}\right) \left(\frac{\pal
A^{ki}_{q+r+s}}{\pal u^{l,t-q}}\right)^{(r)} h^{lj}\right.
\nn\\
&&\hskip -1.0truecm
\left.
+\sum (-1)^{q+r+t} \left(\matrix{q+r+t\cr q ~~ r\cr}\right)
\left( \frac{\pal A^{jk}_{s-q}}{\pal
u^{l,q+r+t}}\right)^{(r)} h^{li}\right] \delta^{(v)}(x-y)\delta^{(s)}(x-z).
\label{ultralocal.cocycle}
\eeqa
\end{exam}

\pf Substituting into (\ref{schouten}) we obtain
\eqa
&&\hskip -1.0truecm
[\varpi,\alpha]^{ijk}_{x,y,z}=\frac{\pal A^{ij}_t(x)}{\pal u^{l,s}(x)}h^{lk}
\delta^{(v)}(x-y)\delta^{(s)}(x-z)
\nn\\
&&\hskip -1.0truecm
+\frac{\pal A^{ki}_t(z)}{\pal u^{l,s}(z)}h^{lj}
\delta^{(v)}(z-x)\delta^{(s)}(z-y) +\frac{\pal A^{jk}_t(y)}{\pal u^{l,s}(y)}
h^{li}\delta^{(v)}(y-z) \delta^{(s)}(y-x).
\eeqa 
Here $A^{ij}_t(x)$, $A^{jk}_t(y)$, $A^{ki}_t(z)$ stand for $A^{ij}_t(u(x);
\dots)$, $A^{jk}_t(u(y); \dots)$, $A^{ki}_t(u(z);\dots)$ resp.  
Use the identities (\ref{diagonal}) 
$$
\delta^{(v)}(z-x)\delta^{(s)}(z-y) = (-\pal_x)^t (-\pal_y)^s [ \delta(z-x)
\delta(z-y)] 
$$
$$=(-\pal_x)^t (-\pal_y)^s [\delta(x-y)\delta(x-z)] =(-1)^t
\sum_{q=0}^t \left(\matrix{t\cr q\cr}\right)
\delta^{(s+q)}(x-y)\delta^{(t-q)}(x-z),
$$
$$
\delta^{(v)}(y-z) \delta^{(s)}(y-x) = (-\pal_x)^s
(-\pal_z)^t[\delta(y-z)\delta(y-x)] 
$$
$$
= (-\pal_x)^s(-\pal_z)^t
[\delta(x-y)\delta(x-z)] =
(-1)^s \sum_{q=0}^s \left(\matrix{s\cr q\cr}\right)
\delta^{(s-q)}(x-y)\delta^{(t+q)}(x-z)
$$
and also (\ref{binom}) to arrive at (\ref{ultralocal.cocycle}).\epf

\noindent 
{\bf Definition.} A {\it local Poisson structure} on the formal loop space is a local
bivector $\varpi\in\Lambda^2_{loc}$ (\ref{bivector}) satisfying $[\varpi,\varpi]=0$.

Adopting the notations common in the physical literature we will represent the
Poisson structure in the form
\beq\label{pb}
\{ u^i(x), u^j(y)\} 
=\sum_s A^{ij}_s(u(x); u_x(x), u_{xx}(x),
\dots)\delta^{(s)}(x-y).
\eeq
The Poisson bracket of two local functionals
$\bar f= \int f(x; u; u_x, \dots) dx$ and\newline $\bar g = \int g(x; u; u_x, \dots)
dx$ can be written in the following equivalent forms (see above the general
theory of multivectors)
\eqa
&&\{ \bar f, \bar g\}=<\varpi,\delta{\bar f}\wedge \delta{\bar g}> 
= \int\int dx dy \frac{\delta \bar f}{\delta u^i(x)}
\{ u^i(x), u^j(y)\} \frac{\delta \bar g}{\delta u^j(y)}
\nn\\
&&\quad
=\sum_s\int dx \frac{\delta \bar f}{\delta u^i(x)} A^{ij}_s(u; u_x, u_{xx}, \dots)
\left( \frac{\delta \bar g}{\delta u^j(x)}\right)^{(s)}\in \Lambda_0.
\eeqa
Therefore it is again a local functional. 

The crucial property of local Poisson brackets is that, the Hamiltonian
systems
\beq\label{ham.pde}
u^i_t =-i_{\delta{\bar H}}\varpi
=\{ u^i(x), \bar H\} = A^{ij}_s(u; u_x, u_{xx}, \dots)\pal_x^s 
\frac{\delta \bar H}{\delta u^j(x)}
\eeq
with local translation invariant Hamiltonians
$$
\bar H= \int H(u; u_x, \dots) dx
$$
are translation invariant evolutionary PDEs (\ref{pde}).

Living in the infinite dimensional loop space we will not impose conditions
on the rank of the Poisson bracket. However, in the main examples the corank of
the bivector will be finite.
\begin{exam} For a constant antisymmetric matrix $h^{ij}$
the bivector
\beq\label{ultralocal}
\{ u^i(x), u^j(y)\} =h^{ij}\delta (x-y)
\eeq
is a local Poisson structure. It is called {\it ultralocal Poisson bracket}.
This is a symplectic structure on the loop space {\it iff} ${\det}
h^{ij}\neq 0$. The Hamiltonian evolutionary PDEs read
$$
u^i_t = h^{ij} \frac{\delta \bar H}{\delta u^j(x)}.
$$
Reducing the nodegenerate matrix  $h^{ij}$ to the canonical form
we arrive at the Hamiltonian formulation of 1+1 dimensional variational
problems 
$$
q^i_t = \frac{\delta \bar H}{\delta p_i(x)},\quad {p^i}_t = -
\frac{\delta \bar H}{\delta q^i(x)}.
$$
\end{exam}

\begin{exam} For a constant symmetric matrix $\eta^{ij}$ the bivector
\beq\label{delta'}
\{u^i(x), u^j(y)\} =\eta^{ij}\delta'(x-y)
\eeq
is a local Poisson structure. Under the assumption 
${\det} \,\eta^{ij}\neq
0$ this Poisson bracket has $n$ independent Casimirs
\beq\label{delta'.casimir}
\bar u^1 = \int u^1\, dx, \dots, \bar u^n = \int u^n \, dx.
\eeq
The annihilator of (\ref{delta'}) is generated by the above Casimirs.

The Hamiltonian evolutionary PDEs read
\beq\label{sht}
u^i_t =\eta^{ij} \pal_x \frac{\delta \bar H}{\delta u^j(x)}.
\eeq
\end{exam}

We finish this section with spelling out the transformation law of coefficients
of local bivectors imposed by the general formula (\ref{transf.multivector}).
If $A^{ij}_s(u; u_x, \dots)$ and $A^{ab}_t(v; v_x, \dots)$ are the coefficients
of a local bivector in two coordinate charts $(U, u^1, \dots, u^n)$ and
$(V, v^1, \dots, v^n)$ resp. then, on $U\cap V$ one has
\beq\label{transf.pb}
A^{ab}_t(v; \frac{\pal v}{\pal u} u_x, \dots) =
\sum_{s\geq t} \left(\matrix{ s \cr t\cr}\right) \frac{\pal v^a}{\pal u^i}
\left( \frac{\pal v^b}{\pal u^j}\right)^{(s-t)} A^{ij}_s(u; u_x, \dots).
\eeq.

\begin{exam}
Applying the transformation
\beq
u=\frac{1}{4} v^2\nn
\eeq
to the bivector
\beq\label{magri.dless}
\{ u(x), u(y)\} =u(x) \delta'(x-y) +\frac{1}{2} u'(x) \delta(x-y)
\eeq
we obtain a constant Poisson bracket of the form (\ref{delta'})
\beq\label{gzf}
\{ v(x), v(y)\} = \delta'(x-y).
\eeq
Hence (\ref{magri.dless}) is itself a Poisson structure. This is the Lie -
Poisson bracket on the space dual to the Lie algebra of vector fields on the
circle \cite{dn89}.
\end{exam}
 
\setcounter{equation}{0}
\setcounter{theorem}{0}

\vskip 0.5cm
\subsection{Problem of
classification of local Poisson brackets}\label{sec-2-4}\par

The last example of the previous section is the simplest issue of the problem
of reduction of local Poisson brackets to the simplest (possibly, to the
constant one) form. In this example the reduction to the constant form was
achieved by a change of coordinates in the target space $M$ 
($M$ was one-dimensional). We give now another well-known example: to transform
the bivector (the {\it Magri bracket} for the KdV equation)
\beq\label{magri}
\{ u(x), u(y)\} =u(x) \delta'(x-y) +\frac{1}{2} u'(x) \delta (x-y) -
\delta'''(x-y)
\eeq
to the constant form (\ref{gzf}) one is to use the celebrated {\it Miura
transformation}
\beq\label{miura}
u=\frac{1}{4} v^2 + v'.
\eeq
Our strategy will be to classify local Poisson brackets on the loop space
${\cal L}(M)$ (the target space $M$ will be a ball in this section)
with respect to the action of the group of
Miura-type
transformations.

The problem of reduction of certain classes of Poisson brackets to a canonical
form by coordinate transformations was first investigated in \cite{dn83}
for the Poisson brackets of hydrodynamic type and in \cite{dn84}
for the so-called differential geometric Poisson brackets (see also \cite{dn89}
and the references therein). Some results regarding reduction of the local
Poisson brackets to the canonical form by using Miura - B\"acklund
transformations were obtained in\cite{astvin}, \cite{geldorf}, \cite{olver1},
\cite{polyak} (in the latter non translation invariant  Poisson brackets were
studied).

We want to classify local Poisson brackets w.r.t. general Miura type
transformations of the form
\beq\label{miura-type}
u^i\to \tilde u^i = F^i(u; u_x, u_{xx}, \dots).
\eeq
The problem is that these transformations do not form a group. The main trouble
is with inverting such a transformation. E.g., to invert the 
Miura transformation   one is to solve Riccati equation (\ref{miura}) 
w.r.t. $v$. To resolve this problem we will extend the class of Miura-type
transformations. Simultaneously we will be to also extend the class of local
functionals, vector fields, and Poisson brackets.
 
\subsubsection{Extended formal loop space}\par
 
Let us introduce gradation on the ring ${\cal A}$ of differential polynomials
putting
\beq\label{grading} {\deg} u^{i,k}= k, ~~k\geq 1, ~~ {\deg}
f(x;u) =0.
\eeq
We extend the gradation onto the spaces ${\cal A}_{k,l}$ of differential forms
by
$$
{\deg} dx = -1, ~~{\deg} \delta u^{i,s}=s, ~s\geq 0.
$$
The differentials $d$ and $\delta$ preserve  the gradation. Introduce a formal
indeterminate $\epsilon$ of the degree
$$
{\deg} \epsilon = -1.
$$
Let us define a subcomplex
$$
\hat{\cal A}_{k,l}\subset {\cal A}_{k,l}\otimes {\mathbb C}[[\epsilon],\epsilon^{-1}], ~~
\hat\Lambda_k\subset \Lambda_k\otimes {\mathbb C}[[\epsilon],\epsilon^{-1}]
$$
collecting all the elements of the total degree $k-1$. In particular, the space
of local functionals $\hat\Lambda_0$ consists of integrals of the form
\eqa
&&
\bar f = 
\int f(u; u_x, u_{xx}, \dots; \epsilon)dx, 
\nn\\
&&
 f(u; u_x, u_{xx},
\dots)=
\sum_{k=0}^\infty \epsilon^k f_k(u; u_x, \dots, u^{(k)}), 
~~f_k\in {\cal A}, ~~{\deg}
f_k =k.
\label{e-local}
\eeqa
We will still call such a series differential polynomials when it will not cause
confusions.
Taking the tensor algebra of $\hat\Lambda_0$ we obtain the ring of functionals
on the {\it extended } formal loop space that we will denote $\hat{\cal L }
(M)$. The vertical differential $\delta$ of the bicomplex must be renormalized
$$
\delta\mapsto \hat\delta={1\over \epsilon} \delta.
$$
As it follows from Theorem 1.2.1 the bicomplex 
$(\hat {\cal A}_{k,l}, d, \hat \delta)$ is exact

The gradation on the vector and multivector fields is defined by
$$
{\deg} \frac{\pal}{\pal x} =1, ~~{\deg} \frac{\pal}{\pal
u^{i,s}}=-s, ~s\geq 0.
$$
Observe that $\pal_x$ increases degrees by one:
$$
 {\deg} \pal_x f = {\deg} f + 1.
$$
The space $\hat\Lambda^1$ of vector fields on $\hat{\cal L }(M)$ is
obtained by collecting  
all the elements in $\Lambda^1 \otimes {\bf
C}[[\epsilon],\epsilon^{-1}]$
of the total degree 1.  
In particular,
the translation invariant evolutionary vector fields are
\beq\label{e-vector}
a^i =\sum_{k=0}^\infty \epsilon^{k-1} 
a^i_k(u; u_x, \dots, u^{(k)}), ~~a^i_k \in
{\cal A}, ~~
{\deg}
a^i_k = k.
\eeq
The corresponding evolutionary system of PDEs reads
\eqa
&&
u^i_t =\epsilon^{-1}\,a^i_0(u) + a^i_1(u; u_x) + \epsilon\, a^i_2(u; u_x,
u_{xx}) + O(\epsilon^2)
\nn\\
&&
a^i_1(u;u_x) = v^i_j(u) u^j_x,
\nn\\
&&
a^i_2(u; u_x, u_{xx}) = b^i_j(u) u^j_{xx} +\frac{1}{2} c^i_{jk}(u) u^j_x u^k_x
\label{e-pde}
\eeqa
etc. 
 
Proceeding in a similar way we introduce the subspace
$$
\hat \Lambda^k \subset \Lambda^k \otimes {\mathbb C}[[\epsilon],\epsilon^{-1}]
$$
of $k$-vectors of the total degree k.

\begin{lemma} The Schouten - Nijenhuis bracket gives a well-defined map
$$
\epsilon\,[\ ,\ ]:\
\hat \Lambda^k \times \hat \Lambda ^l \to \hat \Lambda^{k+l-1}.
$$
\end{lemma}

There is an important subtlety with the grading of the local multivectors.
Indeed, a local $k$-vector is a map
$$
\Lambda_1^{\otimes k} \to \Lambda_0
$$
not
$$
\Lambda^{\otimes k}_1 \to \Lambda_0^{\hat\otimes k}
$$
(cf. the formulae (\ref{k-bracket.1}), (\ref{k-local}) and 
(\ref{value-kform})). Such a map does not
respect the grading.
 So we must assign
a nonzero degree to delta-function
\beq\label{deg-delta}
{\deg} \delta(x-y) = 1, ~~ {\deg} \delta^{(s)}(x-y) =s+1
\eeq
Therefore a general local bivector in $\hat\Lambda_{loc}$ will be represented by an infinite sum
\eqa\label{e-pb}
&&\{ u^i(x), u^j(y)\} =\sum_{k=-1}^\infty \epsilon^k \{ u^i(x), u^j(y)\}^{[k]}
\\
&&
\{ u^i(x), u^j(y)\}^{[k]} =\sum_{s=0}^{k+1} A^{ij}_{k,s}(u; u_x, \dots,
u^{(s)}) \delta^{(k-s+1)}(x-y),
\nn\\
&&
A^{ij}_{k,s}\in {\cal A}, ~~
{\deg} A^{ij}_{k,s} = s, ~~s=0, 1, \dots, k+1.\nn
\eeqa

More explicitly, the first three terms in the expansion (\ref{e-pb})
read
\eqa\label{hydro-1}
&&\{ u^i(x), u^j(y)\} ^{[-1]} = h^{ij}(u(x)) \delta(x-y)
\\&&
\{ u^i(x), u^j(y)\}^{[0]} = g^{ij}(u(x))\delta'(x-y) +\Gamma^{ij}_k(u(x)) u^k_x
\delta(x-y)\label{hydro0}
\\&&
\{ u^i(x), u^j(y)\}^{[1]} = a^{ij}(u(x)) \delta''(x-y) +b^{ij}_k(u(x)) u^k_x
\delta'(x-y) 
\nn\\
&&\quad
+[c^{ij}_k(u(x))u^k_{xx} +\frac{1}{2} d^{ij}_{kl}(u(x)) u^k_x
u^l_x]\delta(x-y)
\label{hydro1}
\eeqa
where $h^{ij}$, $g^{ij}(u)$, $\Gamma^{ij}(u)$, $a^{ij}(u)$, $b^{ij}_k(u)$, $c^{ij}_k(u)$,
$d^{ij}_{kl}(u)$ are some functions on the manifold $M$.

\begin{remark}
Our rules of introducing the gradation can be memorized
using the following simple trick. Do a rescaling of the independent variable
$x$,
\beq\label{rescale.1}
x\mapsto\epsilon x.
\eeq
The $x$-derivatives $u^{i,k} = d^k u^i/dx^k$ will change
\beq\label{rescale.2}
u^{i,k} \mapsto \epsilon^k u^{i,k}.
\eeq
We also have
\beq\label{rescale.dx}
dx\mapsto \epsilon^{-1} dx.
\eeq
The delta-function, according to the definition (\ref{delta}) must be rescaled
as
\beq\label{rescale.3}
\delta(x) \mapsto \epsilon\,\delta(x).
\eeq
In other words, ``delta-function'' is not a function but a density.
Simultaneously with the rescaling we will also redefine the integrals
$$
\int ~ . ~ dx_1\dots dx_k \mapsto \epsilon^k \int ~.~dx_1 \dots dx_k.
$$
After such a rescaling we expand all the formulae of the previous two sections
in a power series in $\epsilon$ to arrive at our grading conventions.  
\end{remark}

\begin{exam} Rescaling (\ref{rescale.1}) the KdV equation $u_t=uu_x+u_{xxx}$
one obtains
$$
u_t=\epsilon(uu_x + \epsilon^2 u_{xxx}).
$$
One usually introduces {\it slow time variable} $t\mapsto \epsilon t$ to recast
the last equation into the form
\beq\label{e-kdv}
u_t =uu_x + \epsilon^2 u_{xxx}.
\eeq
This is the {\it small dispersion expansion} of the KdV equation. The smooth solutions
of (\ref{e-kdv}) describe solutions to KdV slow varying in space and time.
\end{exam}

\begin{remark}
Besides the rescaling procedure, one can arrive at the above
series
(\ref{e-local}), (\ref{e-vector}), (\ref{e-pb}) considering the continous
limits of differential-difference systems. E.g., for the well-known example of
Toda lattice
\beq\label{toda-df}
\dot u_n =v_n - v_{n-1}, ~~
\dot v_n = e^{u_{n+1}} - e^{u_n}, ~~ n\in {\bf Z} 
\eeq
the continuous limit $u_n=u(\epsilon n)=u(x)$, 
$v_n=v(\epsilon n)=v(x)$ gives
an evolutionary system of the form (\ref{e-pde}) with an infinite series in the
r.h.s.
\eqa
&&u_t =(v(x) - v(x-\epsilon)) = \epsilon[ v' -\frac{1}{2} \epsilon v''
+O(\epsilon^2)], \nn\\
&&v_t= e^{u(x+\epsilon)} - e^{u(x)} =\epsilon[(e^u)'
+\frac{1}{2}\epsilon (e^u)'' + O(\epsilon^2)].\nn
\eeqa
Replacing in the Hamiltonian structure
$$
\{ u_m, u_n\} = \{ v_m,v_n\}=0, ~~ \{ u_m, v_n\} =\delta_{mn}-\delta_{m,n+1}
$$
the Kronecker symbols $\delta_{mn}$  by $\epsilon^{-1} \delta(x-y)$, $\delta_{m, n+1}$
by $\epsilon^{-1} \delta(x-y-\epsilon)$,
we obtain a Poisson bracket of the (\ref{e-pb}) form
$$
\{ u(x), u(y)\} = \{ v(x), v(y)\} =0, 
$$
$$\{  u(x), v(y)\} =\frac{1}{\epsilon}
[\delta(x-y) -\delta(x-y-\epsilon)] = \delta'(x-y) -\frac{\epsilon}{2}
\delta''(x-y) +O(\epsilon^2).
$$
\end{remark}

\subsubsection{Miura group}

The last definition will be that to extend the class of 
Miura-type transformations (\ref{miura-type}).
Let us consider the transformations
\eqa
&&
u^i \mapsto \tilde u^i =\sum_{k=0}^\infty \epsilon^k F^i_k(u; u_x, \dots,
u^{(k)}), ~~ i=1, \dots, n \label{miura-group}
\nn\\
&&
F^i_k\in {\cal A}, ~~{\deg} \, F^i_k =k,
\\
&&
{\det}\left( \frac{\pal F^i_0(u)}{\pal u^j}\right) \neq 0.
\nn
\eeqa

\begin{lemma}\label{lemma-miura-group}
 The transformations of the form (\ref{miura-group}) form a group.
The Lie algebra of the group is isomorphic to 
the subalgebra $\hat \Lambda_{ev}^1$ of all translation invariant evolutionary 
vector fileds in $\hat \Lambda^1$ with the Lie
bracket operation.
\end{lemma}

\begin{exam} Let us invert the clasical Miura transformation
$$
u=\frac{1}{4}v^2 + \epsilon v'
$$
using successive approximations. Rewriting the equation in the form
$$
v= 2\sqrt{u-\epsilon v'} = 2\sqrt{u} -\epsilon\frac{v'}{\sqrt{u}} +O(\epsilon^2)
=2\sqrt{u} -\epsilon\frac{u'}{u} + O(\epsilon^2)
$$
we obtain first two terms of the solution $v=F(u; u', \dots; \epsilon)$.
\end{exam}

\begin{remark}
This way of solving the Riccati equation is essentially equivalent
to the classical WKB method of solving the related linear second order ODE
$$
\epsilon^2 y'' =\frac{1}{4} u\, y, ~~ v=4\epsilon\frac{y'}{y}.
$$
Substituting the above series solution to Riccati we obtain the WKB 
asymptotic solution
to the second order ODE with the small parameter $\epsilon \to 0$
$$
y= u^{-1/4} \exp{\frac{1}{2\epsilon}\int \sqrt{u}dx} \left( 1+O(\epsilon)\right).
$$
\end{remark}

\noindent
{\bf Definition}. The group ${\mathcal G}$ of all the transformations of the
form (\ref{miura-group}) is called {\it Miura group}.

The Miura group ${\mathcal G}$ 
looks to be a natural candidate for the role of the group of ``local
diffeomorphisms'' of the extended formal loop space $\hat{\cal L }(M)$
(recall that at the moment $M$ is a ball). ${\mathcal G}$ 
contains the group of diffeomorphisms $Diff(M)$
of the
manifold $M$ as a subgroup. It coincides with the semidirect product of
$Diff(M)$ and the pro-unipotent subgroup ${\mathcal G}_0$ of Miura-type
transformations close to identity,
\beq
u^i\mapsto u^i +\epsilon\,A^i_j(u) u^j_x +\epsilon^2 \left( B^i_j(u) u^j_{xx}
+{1\over 2} C^i_{jk}(u) u^j_x u^k_x\right) +\dots
\eeq
The produc in the group ${\mathcal G}_0$ reads
\eqa
&&
A^i_j = {A_1}^i_j +{A_2}^i_j 
\nn\\
&&
B^i_j ={B_1}^i_j + {B_2}^i_j +{A_2}^i_k {A_1}^k_j
\nn\\
&&
C^i_{jk} ={C_1}^i_{jk} +{C_2}^i_{jk} +{1\over 2} \left[ \pal_s {A_2}^i_j
{A_1}^s_k+\pal_s{A_2}^i_k{A_1}^s_j +{A_2}^i_s \pal_k{A_1}^s_j +{A_2}^i_s \pal_j{A_1}^s_k\right]
\nn\\
&&
\dots
\nn\\
\eeqa

The following simple statements
immediately follow from Lemma \ref{lemma-miura-group}.

\begin{lemma} The class of local functionals (\ref{e-local}), evolutionary PDEs
(\ref{e-pde}), and local translation invariant multivectors (see the formula
(\ref{e-pb}) for the bivectors) on the extended formal loop space $\hat{\cal
L}(M)$ is invariant w.r.t. the action of the Miura group. 
\end{lemma}

\begin{lemma}\label{l-1-4-9} 
An arbitrary vector field $a$ of the form (\ref{e-vector}) 
with $a_0\neq0$ can be reduced, by a transformation of the Miura
group, to a constan form
$$
a_0={\rm const}, ~~a_i=0 ~{\rm for}~ i>0.
$$
\end{lemma}

This is an infinite dimensional analogue of the theorem of ``rectifying
of a vector field''.

\pf By using the theorem of ``rectifying of a vector field'' on a 
finite dimensional manifold, we can reduce the vector field 
$(a_0^1,\dots,a^n_0)$ to a constant one by performing a 
Miura transformation of the form 
(\ref{miura-group}) with $F_k=0,\ k\ge 1$. We prove the lemma by induction. 
Let us assume  the vector field $a$ to be of the form
$$
a^i=a_0^i+\ve^k\,a^i_k(u,u_x,\dots,u^{(k)})+{\cal O}(\ve^{k+1})
$$
with $\deg a_k=k$. 
Since $a_0\ne 0$, we can find differential polynomials $F^i_k(u,
\dots,u^{(k)})$ of degree $k$ such that
$$
a_0^l\,\frac{\pal F^i_k}{\pal u^l}+a^i_k=0.
$$
Then the Miura transformation
$$
{\bar u}^i=u^i+\ve^k\,F^i_k(u,\dots,u^{(k)})
$$
reduces the vector field $a$ to the form
$$
a^i=a_0^i+{\cal O}(\ve^{k+1}).
$$
\epf

\subsubsection{$(p,q)$-brackets on the extended formal loop space}\par

Let us write explicitly down the transformation law of the coefficients of a local
Poisson bracket w.r.t. transformations from the Miura group (cf. \cite{olver}).
Let $A^{kl}$ be the differential operator of the Poisson bracket given in 
(\ref{operator}). In the new ``coordinates'' $\tilde u^i$ of the form
(\ref{miura-type}) the Poisson bracket will be given by the operator
\beq\label{olver-transform}
\tilde A^{ij} ={L^*}^i_k A^{kl} L^j_l
\eeq
where the matrix-valued operator $L^i_k$ and the adjoint one ${L^*}^i_k$ are
given by
$$
L^i_k = \sum_s (-\pal_x)^s \circ \frac {\pal\tilde u^i}{\pal u^{k,s}} ,
~~
{L^*}^i_k =\sum_s  \frac {\pal\tilde u^i}{\pal u^{k,s}} \pal_x^s.
$$
{\bf Main Problem.} To describe the orbits of the action of the Miura group
${\mathcal G}$ on ${\hat\Lambda}^2$.
\medskip

To our opinion this problem is the natural setup of the problem of
classification of Poisson structures of one-dimensional evolutionary PDEs
with respect to Miura-type transformations (they are called also Darboux,
or Bianchi, or B\"acklund transformations. Our transformations do not
involve a change of the independent variable $x$ since we consider only
translation invariant PDEs).

Our conjecture is that, for a reasonable class of Poisson brackets to be defined
below, 
 the orbits are labelled by certain finite-dimensional
geometrical structures on the underlined manifold $M$. Below we will illustrate
this claim describing two orbits being, in a certain sense, generic. The full
problem remains open.

We first  explain how a local Poisson bracket from $\hat\Lambda^2_{loc}$  
induces certain finite-\newline
dimensional geometrical structures
on $M$.

\begin{lemma} The subgroup $Diff(M)\subset{\mathcal G}$ acts independently on
every term $\{ ~,~\}^{[k]}$ of the expansion (\ref{e-pb}), $k\geq -1$. In
particular, the leading term $A^{ij}_{k,0}(u)$ is a $(2,0)$-tensor 
field on $M$,
symmetric/antisymmetric for even/odd $k$.
\end{lemma}

\pf This follows from the transformation law (\ref{transf.pb}). The
symmetry/antisymmetry of the coefficients follows from the general antisymmetry
condition (\ref{2-antisymmetry}). The lemma is proved.\epf

\begin{lemma} The subspaces ${\rm span}\left( \{~ , \}^{[-1]}, \{~,
~\}^{[0]}, \dots, \{~, ~\}^{[k]}\right)$ for every $k$ remain invariant
w.r.t. to the action of the Miura group ${\mathcal G}$.
\end{lemma}

This follows from the explicit formula (\ref{olver-transform}).

\begin{lemma} The first non-zero term in the expansion (\ref{e-pb})
is itself a local Poisson bracket.
\end{lemma}

This is obvious.

\begin{cor} The coefficient $h^{ij}(u)$ in (\ref{hydro-1}) is a Poisson
structure on $M$. This Poisson structure is invariant w.r.t. the action
of ${\mathcal G}$ on ${\hat\Lambda}_{loc}^2$.
\end{cor}

We obtain a map
\beq\label{classify-1}
{\hat\Lambda}^2_{loc} /{\mathcal G} \to ~{\rm Poisson~structures~on~}M.
\eeq
Let us assume that the Poisson structure $h^{ij}(u)$ on $M$ has constant rank
$p=2p_1$. Denote $q:=n-p$ the corank of $h^{ij}(u)$.

{\bf Definition.} (\ref{e-pb}) is called {\it $(p,q)$-bracket} if the
coefficient $g^{ij}(u)$ in (\ref{hydro0}) does not degenerate on $Ker\,
h^{ij}(u)\subset T^*_uM$ on an open dense subset in $M\ni u$.
\smallskip

\begin{exam} The ultralocal Poisson bracket (\ref{ultralocal}) with a
non-degenerate matrix $h^{ij}$ is a $(n,0)$-bracket.
\end{exam}

\begin{exam} The Poisson bracket (\ref{delta'}) with a non-degenerate matrix
$\eta^{ij}$ is a $(0,n)$-bracket.
\end{exam}

\begin{exam} Let $c^{ij}_k$ be the structure constants of a semisimple
$n$-dimensional Lie algebra ${\mathfrak  g}$. The Killing form $\eta^{ij}$ on 
${\mathfrak
g}$ defines a central extension $\hat {\mathfrak  g}$ of ${\mathfrak  g}$ called Kac - Moody
Lie algebra \cite{kac}. The Lie - Poisson bracket (\ref{lie-poiss}) on the dual
space $\hat {\mathfrak  g}^*$ 
\beq\label{current}
{1\over \epsilon} \{ u^i(x), u^j(y)\} = \eta^{ij}\delta'(x-y) +{1\over \epsilon}
c^{ij}_k u^k\delta(x-y)
\eeq
can be considered as a Poisson bracket of the form (\ref{e-pb}) on the loop
space $\hat {\cal L}({\mathfrak  g}^*)$. Here $\epsilon$ is the central charge. This
is a $(p,q)$-bracket with
$$
q={\rm rk}\, {\mathfrak  g}, ~~ p={\rm dim}\, {\mathfrak  g} - {\rm rk}\, 
{\mathfrak  g}.
$$
\end{exam}

Let a $(p,q)$-Poisson bracket on $\hat{\cal L}(M)$ of the form (\ref{e-pb}) -
(\ref{hydro0}) be given. We will now construct a {\it flat metric} on the base
of the symplectic foliation of $M$ defined by the finite-dimensional Poisson
bracket $h^{ij}(u)$. Let us assume that $M$ is a small ball such that the
symplectic foliation defines a fibration
\beq\label{fibr}
M\to N, ~~{\rm dim}\, N=q.
\eeq
Functions on $N$ are Casimirs of the finite-dimensional Poisson bracket
$h^{ij}(u)$. Define first a symmetric bilinear form $(~,~)^*$ on $T^*N$ putting
\beq\label{fibr1}
(df_1, df_2)^* := {\pal f_1\over \pal u^i}{\pal f_2\over \pal u^j}g^{ij}(u)
\eeq
for any two Casimirs of $h^{ij}(u)$. By the assumption this bilinear form does
not degenerate. Define the non-degenerate symmetric tensor $(~,~)$ on $TN$ by
\beq\label{fibr2}
(~,~):= \left[ (~,~)^*\right]^{-1}.
\eeq

\begin{theorem} The metric (\ref{fibr2}) on $TN$ is well-defined and flat.
\end{theorem}

\pf The simplest way to prove that the metric (\ref{fibr2}) is constant along
symplectic leaves and also prove vanishing of the curvature of the metric is
the following one. Choose local coordinates $u=(w^a, v^\alpha)$ on $M$, $a=1,
\dots, p=2 p_1$, $\alpha=1, \dots, q$, $p+q=n$, such that $w^1, \dots, w^p$ are
canonical coordinates on the symplectic leaves and $v^1, \dots, v^q$ is a system
of independent Casimirs of $h^{ij}(u)$. The $v$'s can be considered as
coordinates on $N$. In the local coordinates the Poisson bracket (\ref{e-pb})
reads
\eqa\label{fibr5}
&&
\{ w^a(x), w^b(y)\} = {1\over \epsilon} h^{ab} \delta(x-y) +A^{ab}_{00}(u)
\delta'(x-y) +A^{ab}_{01}(u, u_x) \delta(x-y) +O(\epsilon)
\nn\\
&&
\{v^\alpha(x), w^a(y)\} = A^{\alpha a}_{00}(u) \delta'(x-y) +A^{\alpha
a}_{01}(u,u_x) \delta(x-y) +O(\epsilon)
\nn\\
&&
\{ v^\alpha(x), v^\beta(y)\} = g^{\alpha\beta}(u) \delta'(x-y) +\left(
\Gamma^{\alpha\beta}_a(u) w^a_x +\Gamma^{\alpha\beta}_\gamma(u)
v^\gamma_x\right) \delta(x-y) +O(\epsilon)
\nn\\
\eeqa
Here $h^{ab}$ is a constant antisymmetric nondegenerate matrix, the $q\times q$
matrix $g^{\alpha\beta}(u)$ coincides with the Gram matrix of the bilinear form
(\ref{fibr1}) in the basis $dv^1$, \dots, $dv^q$. Let us consider the following
foliation on the loop space $\hat{\cal L}(M)$
\beq\label{fibf3}
w^a(x)\equiv w^a_0, ~~a=1, \dots, p
\eeq
for arbitrary given numbers $w_0^1$, \dots, $w_0^p$. Due to nondegeneracy of
$h^{ab}$ this foliation is cosymplectic. The corresponding Dirac bracket
$\{~,~\}_D$ can be considered as a Poisson bracket on $\hat{\cal L}(N)$ since,
by definition
$$
\{v^\alpha(x), w^a(y)\}_D = \{w^a(x), w^b(y)\}_D =0.
$$
Let us show that the Dirac bracket is a $(0,q)$-bracket on $\hat{\cal L}(N)$
with the same leading term
$$
\{ v^\alpha(x), v^\beta(y)\}_D = \{ v^\alpha(x), v^\beta(y)\}+O(\epsilon).
$$
Introduce the differential operator 
$$
\Pi_{ab} = h_{ab} -\epsilon\, h_{a\, a'} \left( A^{a'b'}_{00}(u) {d\over dx} +
A^{a'b'}_{01}(u, u_x)\right) \, h_{b'\, b} +O(\epsilon^2)
$$
inverse to the operator $\epsilon\, A^{ab}$. Then the Dirac bracket has the form
\beq\label{fibf4}
\{v^\alpha(x), v^\beta(y)\}_D = \{v^\alpha(x), v^\beta(y)\}
-\epsilon\, \sum_{a, b=1}^p A^{\alpha a} \Pi_{ab} A^{b\beta} \delta(x-y).
\eeq
We obtain a $(0,q)$ Poisson bracket on $\hat {\cal L}(N)$ eventually depending
on the parameters $w_0^1$, \dots, $w_0^p$. The leading term
$$
\{v^\alpha(x), v^\beta(y)\}_D^{[0]} =g^{\alpha\beta}(v,w_0) \delta'(x-y)
+\Gamma^{\alpha\beta}_\gamma(v,w_0) v^\gamma_x(x-y)
$$
is itself a Poisson bracket (the so-called Poisson bracket of hydrodynamic
type). According to the theory of such brackets \cite{dn83} one can choose local
coordinates $v^\alpha$ on $N$ in such a way that
$$
g^{\alpha\beta}={\rm const}, ~~\Gamma^{\alpha\beta}_\gamma=0.
$$
This proves the theorem.
\epf

An alternative way to prove of the Theorem is to write explicitly down the terms of
the order $\epsilon^{-1}$ in the Jacobi identity
$$
\{ \{ v^\alpha(x), v^\beta(y)\}, w^a(z)\} +({\rm cyclic}) =0
$$
in order to prove that the leading term in $\{ v^\alpha(x), v^\beta(y)\}$
does not depend on $w$, $w_x$. Then from the leading term in the Jacobi identity
$$
\{ \{ v^\alpha(x), v^\beta(y)\}, v^\gamma(z)\} +({\rm cyclic})=0
$$
it follows, as in \cite{dn83}, vanishing of the curvature of the metric
$g^{\alpha\beta}(v)$.
\smallskip

We believe that the problem of classification of $(p,q)$-brackets (and also 
classification of pencils of $(p,q)$-brackets) is very important in the
Hamiltonian theory of integrable PDEs. In this paper we will mainly consider
$(0,n)$-brackets leaving the general case for a subsequent publication.

To illustrate our technique we will begin with a more simple example of
$(n,0)$-brackets.

\subsubsection{Classification of $(n,0)$-brackets}\par
Our first result is

\begin{theorem}\label{sym-case}
If $M$ is a ball then all  $(n,0)$ Poisson brackets in
$\hat\Lambda^2_{loc}$ 
are equivalent w.r.t. the action (\ref{olver-transform}) of the Miura group
${\mathcal G}$.
\end{theorem} 

\pf First we choose the Darboux coordinates for the symplectic structure
on $M$. The Poisson bracket in question will read
\beq\label{step0}
\{ u^i(x), u^j(y)\} = \frac{1}{\epsilon} h^{ij} \delta(x-y) +\sum_{k=0}^\infty
\epsilon^k \{u^i(x),u^j(y)\}^{[k]}.
\eeq
Next we will try to kill all the terms of the expansion (\ref{step0})
by transformations of the form (\ref{miura-group}) with $F^i_0(u)=u^i$, $i=1, \dots,
n$. To this end an appropriate version of Poisson cohomology will be 
useful. We define the Poisson cohomology $H^*(\hat{\cal L}(M), \varpi)$ for a Poisson
structure $\varpi\in \hat\Lambda^2_{loc}$ as the cohomology of the complex
\beq\label{complex}
0\to\hat\Lambda^0_{loc} \stackrel{\pal}\to \hat \Lambda^1_{loc} 
\stackrel{\pal}\to 
\hat \Lambda^2_{loc} \stackrel{\pal}\to \dots
\eeq
with the differential $\pal\beta:= [\varpi,\beta]$. In the present proof 
$\varpi=
h^{ij}\delta(x-y)$. The cohomology is naturally decomposed into the
direct sum
\beq\label{decomposition}
H^k=\oplus_{m\geq -1} H^{k,m}
\eeq
with respect to monomials in $\epsilon$,
where $H^{k,m}$ consists of the cocycles proportional to $\epsilon^m$.
Denote
\beq\label{decomposition1}
\tilde H^k := \oplus_{m\geq 0} H^{k,m}.
\eeq

The following obvious statement holds true.

\begin{lemma} The first non-zero term in the expansion (\ref{step0}) is a
2-cocycle in the Poisson cohomology $\tilde H^2$ of the ultralocal Poisson bracket 
(\ref{ultralocal}). 
\end{lemma}
So, we will be able to kill this first nonzero term of the expansion if we
prove that, for any 2-cocycle $\in \hat\Lambda^2_{loc}$ of the ultralocal bracket  
$\varpi:= \epsilon\{~, ~\}^{[-1]}$ , there exists a vector field $a$ of the form
(\ref{e-vector}) such that the Lie derivative $Lie_a\varpi$ gives the cocycle
and $a|_{\epsilon=0}=0$. This will follow from the following general statement
about triviality of cohomology of the ultralocal Poisson bracket.

\begin{lemma} For $M$ = ${\rm ball}$ and the ultralocal Poisson bracket $\varpi$
(\ref{ultralocal}) with det $h^{ij}\neq 0$ the Poisson cohomology 
$\tilde H^1(\hat{\cal L}(M), \varpi)$, $\tilde H^2(\hat{\cal L}(M), \varpi)$ vanish.
\end{lemma}
 
The first proof of the lemma (and of the lemma \ref{H1H2-der}
 below) was obtained
by E. Getzler \cite{getzler} (also triviality of the higher cohomology has been
proved). Independently,  L. Degiovanni, F. Magri, V. Sciacca obtained another proof
\cite{magri1}.
We have decided to present here  our own proofs of triviality of cohomologies
that closely follow the finite-dimensional case. Our proof will also be useful
in the study of bihamiltonian structures below.

Let us prove first  triviality of $H^1$. Let an evolutionary vector field
$a$ with the components $a^1$, \dots, $a^n$  be a cocycle. Denote
$$
\omega_i =h_{ij}a^j
$$
where the constant matrix $h_{ij}$ is inverse to $h^{ij}$. The condition 
$\pal a=Lie_a \varpi=0$ reads
$$
\frac{\pal \omega_i} {\pal u^{j,s}} =\sum_{t\geq s} (-1)^t 
\left(\begin{array}{c}  t \\ s\end{array}\right)
\pal_x^{t-s} \frac{\pal \omega_j}{\pal u^{i,t}}.
$$
Using (\ref{one-main})  
we conclude that there exists a local functional
$ \bar f = \int f \, dx$ such that
$$
\omega_i =  \frac {\delta \bar f}{\delta u^i(x)}.
$$
Therefore the vector field is a Hamiltonian one,
$$
a^i = h^{ij}\frac {\delta \bar f}{\delta u^j(x)}.
$$

Let us now proceed to the proof of triviality of $H^2$. 
The idea is very simple:
the bivector
$$
 \varpi + \varepsilon  \alpha =h^{ij}\delta(x-y)  + \varepsilon \sum_s
A^{ij}_s\delta^{(s)} (x-y) 
$$
satisfies the Jacobi identity
$$
[\varpi+\varepsilon \alpha, \varpi+\varepsilon\alpha]=0 ({ \mathrm  
mod}\, \varepsilon^2)
$$
{\it iff}\ the inverse matrix is a closed differential form
$$
\frac{1}{2} h_{ij}dx\wedge\delta u^i\wedge\delta u^j
+\frac12\, \ve\, \omega_{i; js}\, dx\wedge\delta u^i\wedge\delta u^{j,s}
\ ({\mathrm mod}\,\ve^2)
$$
where
\beq\label{om}
\omega_{i;js}:= h_{ip}h_{jq} A^{pq}_s.
\eeq
Denote
$$
\omega =\frac{1}{2} \omega_{i;js} \delta u^i\wedge\delta u^{j,s}.
$$
From the condition of closedness $\delta(dx\wedge\omega)=0\in \Lambda_3$ we
derive, due to Corollary \ref{exact-2form},
 existence of a one-form $dx\wedge\phi$,
$\phi = \phi_i \delta u^i$ 
such that $\delta(dx\wedge\phi) = dx\wedge \omega$.
The vector field $a$ with the components
$$
a^i = h^{ij}\phi_j
$$
gives a solution to the equation
$$
[\varpi, a] = \alpha.
$$

To be on the safe side we will now show, by straightforward calculations,
that, indeed, the above geometrical arguments work. First, from the
antisymmetry condition (\ref{2-antisymmetry}) for the bivector $\alpha$
it readily follows the antisymmetry condition (\ref{2-antisymmetry.forms})
for the 2-form $\omega$ with the reduced components (\ref{om}). Next, we are to
verify that from the cocycle condition $[\varpi, \alpha]=0$ 
where the Schouten - Nijenhuis bracket $[\varpi, \alpha]$ is written
in (\ref{ultralocal.cocycle}), it follows closedness (\ref{2-form.closed})
of $dx\wedge\omega\in \Lambda_2$. First we will rewrite the formula for the 
bracket in a slightly modified form. Differentiating the antisymmetry condition
$$
A^{ik}_s =-\sum(-1)^m \left(\matrix{m\cr s\cr}\right) \pal_x^{m-s}A^{ki}_m
$$
w.r.t. $u^{l,t}$ and using the commutators (\ref{commute}) we obtain
$$
\frac{\pal A^{ik}_s}{\pal u^{l,t}} =- \sum (-1)^{q+r+s} \left(\matrix{q+r+s\cr
q ~~~r\cr}\right) \left( \frac{\pal A^{ki}_{q+r+s}}{\pal u^{l,t-q}}\right)^{r}.
$$
So the coefficients of the Schouten - Nijenhuis bracket 
(\ref{ultralocal.cocycle}) can be rewritten as follows
\eqa
&&
[\varpi, \alpha]^{ijk}_{x,y,z}=
\left[\frac{\pal A^{ij}_t}{\pal u^{l,s}} h^{lk}
 -\frac{\pal A^{ik}_s}{\pal u^{l,t}}\right.
\nn\\
&&\hskip -1.0 truecm\left.
+\sum (-1)^{q+r+t} \left(\matrix{q+r+t\cr q ~~ r\cr}\right)
\left( \frac{\pal A^{jk}_{s-q}}{\pal
u^{l,q+r+t}}\right)^{(r)} h^{li}\right] \delta^{(t)}(x-y)\delta^{(s)}(x-z).
\eeqa
The coeffcient of  $ \delta^{(t)}(x-y)\delta^{(s)}(x-z)$ must vanish for every
$t$ and $s$. 
Multiplying this coefficient by $h_{ia}h_{jb}h_{kc}$ we arrive at the condition
of closedness (\ref{2-form.closed})
of the 2-form (\ref{om}).
Using Corollary \ref{exact-2form} we establish existence of
differential polynomials $\phi_1$, \dots, $\phi_n$ representing the 2-form
as in (\ref{2-form.exact}). The translation invariant vector field
$$
a^i =h^{ij}\phi_j
$$
will satisfy $\pal a =\alpha$. 
This proves the lemma, and also the theorem.\epf

\subsubsection{Classification of $(0,n)$-brackets}\label{sec-2-4-5}\par

Let us now proceed to considering the Poisson structures in $\hat\Lambda^2$
with identically vanishing leading term $\{~,~\}^{[-1]}$. 
As above, the first nonzero term (\ref{hydro0}) is itself
a Poisson bracket. The leading coefficient $g^{ij}(u)$ of it determines a symmetric
tensor field on $M$ invariant w.r.t. the action of the Miura group. This gives a
map
$$
\hat\Lambda^2/{\mathcal G} \to {\rm symmetric ~tensors ~on~}M.
$$

\begin{theorem} Let $M$ be a ball. Then the only invariant of a $(0,n)$
 Poisson bracket
in $\hat\Lambda^2$ with respect to the action of the Miura group is the
signature of the quadratic form $g^{ij}(u)$.
\end{theorem}
Proof. The symmetric nondegenerate tensor $g^{ij}(u)$
defines a {\it pseudoriemannian metric} 
$$
g_{ij}(u) du^i du^j, ~~\left(g_{ij}\right) = \left( g^{ij}\right)^{-1}
$$
on the manifold $M$. 
From the general theory of \cite{dn83} of the Poisson brackets
of the form (\ref{hydro0})  it follows that the Riemann curvature of the metric
vanishes, and that the coefficient $\Gamma^{ij}_k(u)$ in (\ref{hydro0})
is related to the Christoffel coefficients $\Gamma_{ij}^k(u)$ of the Levi-Civita connection
for the metric by
$$
\Gamma^{ij}_k = -g^{is}\Gamma_{sk}^j.
$$
Using standard arguments of differential geometry we deduce that, locally
coordinates $v^1(u)$, \dots, $v^n(u)$ exists such that, in the new coordinates
the metric becomes constant
$$
\frac{\pal v^k}{\pal u^i}\frac{\pal v^l} {\pal u^j} g^{ij}(u) =\eta^{kl} = {\rm
const}.
$$
The Christoffel coefficients in these coordinates vanish. Of course,
all constant symmetric matrices $\eta^{kl}$ of a given signature
are equivalent w.r.t. 
linear changes of coordinates.

We have reduced the proof of the theorem to reducing to the normal form
(\ref{delta'}) the 
Poisson bracket 
\beq\label{step1}
 \{u^i(x), u^j(y)\} = \eta^{ij}\delta'(x-y) + 
\sum_{k=1}^\infty \epsilon^k\{ u^i(x),
u^j(y)\}^{[k]}
\eeq
by the transformations of the form (\ref{miura-group}) 
with $F^i_0=$id. As in the
proof of Theorem \ref{sym-case}
the latter problem is reduced to proving triviality of the
second Poisson cohomology of the Poisson bracket $\alpha$ of the form 
(\ref{delta'}). Triviality of
the cohomology is somewhat surprising from the point of view of
finite-dimensional Poisson geometry. Indeed, as we have seen above 
this Poisson bracket degenerates. So we will be to also prove that
all the cocycles are tangent to the leaves of the symplectic foliation
(see the end of section \ref{sec-2-3}).

\begin{lemma}\label{H1H2-der}
 For $M$ = ball all the cocycles in $H^1(\hat {\cal L}(M))$
and $H^2(\hat {\cal L}(M))$ vanishing at $\epsilon=0$ are trivial.
\end{lemma}

Denote, like in (\ref{decomposition1}), 
\beq\label{decomposition2}
\tilde H^k =\oplus_{m>0} H^{k,m}.
\eeq
We are to prove that $\tilde H^1=\tilde H^2=0$.

Let us begin with proving triviality of $\tilde H^1$. Using (\ref{lie-der})
we obtain
\eqa
&&
\pal_a\varpi =Lie_a\varpi^{ij} =-\pal_x \sum_{t\geq 0}(-1)^t \eta^{ip} \left(
\frac{\pal a^j}{\pal u^{p,t}}\right)^{(v)} \delta(x-y)
\nn\\
&&\hskip -1.0truecm
- \sum_{r\geq 0}\left[ \frac{\pal a^i}{\pal u^{p,r}} \eta^{pj} +\sum_{t\geq r}
(-1)^t \left(\matrix{t+1\cr r+1\cr}\right) \eta^{ip} \left(\frac{\pal a^j}{\pal
u^{p,t}}\right)^{(t-r)}\right]\delta^{(r+1)}(x-y).
\label{later}
\eeqa
Since
$$
\frac{\delta\bar a^j}{\delta u^p(x)} =\sum_{t\geq 0}(-1)^t \left(\frac{\pal
a^j}{\pal u^{p,t}}\right)^{(v)}
$$
is a differential polynomial in $\Lambda_0\otimes {\mathbb C}[[\epsilon]]$ of the
degree 0 vanishing at $\epsilon=0$, from vanishing of the coefficient in front
of $\delta(x-y)$ we derive that
$$
\frac{\delta\bar a^j}{\delta u^p(x)}=0, ~~j, \, p=1, \dots, n.
$$
Using Example \ref{zero-fun}
we derive existence of differential polynomials $b^j$ s.t.
$$
a^j=\pal_xb^j, ~~j=1, \dots, n.
$$
This is the crucial point in the proof: we have shown that the vector field $a$
is tangent to the level surface of the Casimirs (\ref{delta'.casimir}). The
remaining part of the proof is rather straightforward. Using (\ref{commute})
and also the Pascal triangle identity
$$
\left(\matrix{m\cr n\cr}\right) +\left(\matrix{m\cr n-1\cr}\right) 
=\left(\matrix{m+1\cr n\cr}\right) 
$$
we rewrite the coefficient of $\delta^{(r+1)}(x-y)$ in the form
$$
\pal_x\left[ \frac{\pal \omega_k}{\pal u^{l,r}}-\sum_{t\geq r}(-1)^t
\left(\matrix{t+1\cr r\cr}\right) \left(\frac{\pal \omega_l}{\pal
u^{k,t}}\right)^{(t-r)}\right]
$$
$$
+\frac{\pal \omega_k}{\pal u^{l,r-1}} +(-1)^r \frac{\pal \omega_l}{\pal
u^{k,r-1}}=0.
$$
Here
$$
\omega_k=\eta_{li}b^i, ~~(\eta_{ij})=(\eta^{ij})^{-1},
$$
the last two terms are not present for $r=0$.
As above, for $r=0$ we derive that
$$
\frac{\pal \omega_k}{\pal u^l} =\sum_{t\geq 0}(-1)^t
\left(\frac{\pal\omega_l}{\pal u^{k,t}}\right)^{(v)}.
$$
Proceeding by induction in $r$ we prove that the 1-form 
$\int dx\wedge \omega_i\,\delta u^i$ is closed. Using the Volterra
criterion we derive existence of a differential polynomial $f$ s.t. $\omega =
\delta\int f\,dx$. Hence
$$
a^i=\eta^{ij} \pal_x \frac{\delta\bar f}{\delta u^j(x)}.
$$
We proved triviality of $\tilde H^1$.

Let us proceed to prove the triviality of $\tilde H^2$.
 The condition $\pal \alpha=0$ for $\alpha$ of the form
(\ref{bivector}) can be computed similarly to 
Example \ref{schouten-bivector-sym}. We obtain a system of
equations
\eqa
&&
\frac{\pal A^{ij}_t}{\pal u^{l,s-1}}\eta^{lk}
+\sum (-1)^{q+r+s}\left(\matrix{q+r+s\cr q ~~ r\cr}\right) 
\left(\frac{\pal A^{ki}_{q+r+s}}{\pal u^{l,t-q-1}}\right)^{(r)}\eta^{lj}
\nn\\
&&\quad
+\sum (-1)^{q+r+t} \left(\matrix{q+r+t\cr q ~~ r\cr}\right)
\left( \frac{\pal A^{jk}_{s-q}}{\pal
u^{l,q+r+t-1}}\right)^{(r)}\eta^{li}=0 \nn\\
&&\quad ~~{\rm for~any}~i,j,k,s,t
\label{eq_s_t}
\eeqa
(it is understood that the terms with $s-1$, $t-q-1$ or $t+q+r-1$ negative
do not appear in the sum). Recall that the crucial point in the proof
of triviality of the 2-cocycle is to establish validity of 
(\ref{2-trivial}) for the Casimirs (\ref{delta'.casimir}) of $\varpi$.
Explicitly, we need to show that
\beq\label{eq.main}
\alpha(\delta {\bar u}^i, \delta {\bar u}^j)
=\int A^{ij}_0dx=0 ~~{\rm for~any}~i,j.
\eeq
We first use (\ref{eq_s_t}) for $s=t=0$ to prove that
$$
\pal_x \sum_r (-1)^r \left( \frac{\pal A^{jk}_0}{\pal
u^{l,r}}\right)^{(r)}\eta^{li}=0.
$$
Hence
$$
A^{jk}_0=\pal_x B^{jk}
$$
for some differential polynomial $B^{jk}$. This implies (\ref{eq.main}).
The rest of the proof is identical to the proof of Lemma \ref{first-lemma}. 
We first construct
the vector field $z$ (see the proof of the lemma \ref{first-lemma}). 
To this end we use the
equation (\ref{eq_s_t}) for $s=0$, $t>0$:
$$
\sum_{q,r} (-1)^{q+r}\left(\matrix{ q+r\cr r}\right) \left(
\frac{\pal A^{ki}_{q+r}}{\pal u^{l,t-q-1}}\right)^{(r)} \eta^{lj}
+\sum_r (-1)^{t+r} \left(\matrix{ t+r\cr r\cr} \right)
\left(\frac{\pal A^{jk}_0}{\pal u^{l,t+r-1}}\right)^{(r)}\eta^{li}=0.
$$
Differentiating the antisymmetry condition 
$$
A^{ik}_0 = \sum (-1)^{r+1} \left(A^{ki}_r\right)^{(r)}
$$
w.r.t. $u^{l,t-1}$ we identify the first term of the previous equation with 
$$
-\frac{\pal A^{ik}_0}{\pal u^{l,t-1}} \eta^{lj}.
$$
The resulting equation coincides with the condition $\pal a^k =0$
of closedness of the 1-cocycle 
$$
(a^k)^i  = A^{ik}_0
$$
for every $k=1, \dots, n$
(see (\ref{later}) for the explicit form of this condition). Using the first
part of Lemma we arrive at existence of $n$ differential polynomials
$q^1$, \dots, $q^n$ s.t.
\beq\label{z-field}
A^{ik}_0 =\eta^{is}\pal_x\frac{\delta\bar q^k}{\delta u^s(x)}.
\eeq

The last step, as in the proof of Lemma \ref{first-lemma}, 
is to change the cocycle $\alpha$
to a cohomological one to obtain a closed 2-cocycle
$$
\alpha \mapsto \alpha + \pal z=:\alpha'
$$
for 
$$
z = q^i\frac{\pal}{\pal u^i}.
$$
The new 2-cocycle $\alpha'$ will have the same form as above with $A^{ij}_0=0$.
Denote
$$
g_{i;js} :=\eta_{ip}\eta_{jq} A^{ij}_s, ~~s\geq 1.
$$
We will now show existence of differential polynomials $\omega_{i;j0}$,
$\omega_{i;j1}$,\dots
s.t.
\eqa
&&
g_{i;j1}=\pal_x\omega_{i;j0},
\nn\\
&&
g_{i;js}=\pal_x\omega_{i; j,s-1}+\omega_{i;j,s-2} ~~{\rm for}~s\geq 2.
\label{int}
\eeqa
From (\ref{eq_s_t}) for $s=1$, $t=0$ we obtain
$$
\pal_x\sum_r (-1)^r \left(\frac{\pal A^{jk}_1}{\pal u^{l,r}}\right)^{(r)}=0.
 $$
As we already did many times, from the last equation it follows
that 
$$
\sum_r (-1)^r \left(\frac{\pal A^{jk}_1}{\pal u^{l,r}}\right)^{(r)}=0.
$$
This shows existence of $\omega_{i;j0}$. Using (\ref{eq_s_t}) for $s=1$
and $t>0$ we inductively prove existence of the differential polynomials
$\omega_{i;j,t-1}$. Actually, we can obtain 
\beq\label{yz}
\omega_{i;jl} =\sum_{s\geq l+2} \pal_x^{s-l-2} g_{i;js}.
\eeq
From this it readily follows that the coefficients $\omega_{i;js}$ satisfy
the antisymmetry conditions (\ref{2-antisymmetry.forms}). 
Thus they determine a 2-form
$\omega$.

Let us prove that the 2-form $\omega$ is closed. Denote 
$$
J_{ijk;st}:=
\left( \sum_{m=s}^{t+s}\sum_{r=0}^{m-s} +\sum_{m\geq t+s+1} \sum_{r=0}^t\right)
(-1)^m \left(\begin{array}{c} m \\ r \, s\end{array}\right) \pal_x^{m-r-s} \frac{\pal \omega_{j; k, t-r}}{\pal
u^{i,m}}
$$
$$
+ \frac{\pal \omega_{i;j,s}}{\pal u^{k,t}} -\frac{\pal \omega _{i;k,t}}{\pal
u^{j,s}}
$$
the l.h.s. of the equation (\ref{2-form.closed}) of closedness of a 2-form.
Let us show that the coefficient of $\delta^{(v)}(x-y)\delta^{(s)}(x-z)$
in (\ref{eq_s_t}) is equal to
\beq\label{coef}
\pal_x J_{ijk;t-1,s-1}+J_{ijk;t-1,s-2} +J_{ijk;t-2,l-1}.
\eeq
To this end we replace the second sum in (\ref{eq_s_t}) by
$$
-\frac{\pal A^{ik}_s}{\pal u^{l,t-1}}\,\eta^{lj}.
$$
Lowering the indices by means of $\eta_{ij}$ and using (\ref{int})
we obtain (\ref{coef}). From vanishing of (\ref{coef}) we inductively deduce
that $J_{ijk;st}=0$ for all $i,j,k=1, \dots, n$ and all $s,t\geq 0$
(observe that the coefficients $J_{ijk;t0}=J_{ijk;0s}=0$ due to our assumption
$A^{ij}_0=0$. This proves that the 2-form $\omega$ is closed. So $\omega
=\delta \int dx\wedge \phi$ for some 1-form $\phi =\phi_i\delta u^i$. Introducing
the vector field 
$$
a^i = \eta^{ik}\phi_k
$$
we finally obtain, for the original cocycle $\alpha$,
$$
\alpha = \pal (a-z).
$$
Theorem is proved.

\newpage
\setcounter{equation}{0}
\setcounter{theorem}{0}
\section{Bihamiltonian geometry of loop spaces}\par
\subsection{Bihamiltonian structures and hierarchies of commuting flows}\par

\subsubsection{Poisson pencils and bihamiltonian recursion procedure: summary
of the finite-dimensional case}\par 

{\bf Definition.} A {\it bihamiltonan structure} on the manifold $P$ is a
2-dimensional linear subspace in the space of Poisson structures on $P$. 

Choosing two points $\{ ~, ~\}_1$ and $\{ ~, ~ \}_2$ of the subspace we
obtain that the linear combination
\beq\label{biham1}
a_1 \{ ~, ~\}_1 + a_2 \{ ~, ~\}_2
\eeq
with arbitrary constant coefficients $a_1$, $a_2$ is again a Poisson bracket.
This reformulation is usually referred to as {\it the compatibility condition}
of the two Poisson brackets. It is spelled out as vanishing of the Schouten -
Nijenhuis bracket
\beq\label{biham2}
\left[ \{ ~, ~\}_1, \{ ~, ~\}_2\right] =0.
\eeq

An importance of bihamiltonian structures for recursive constructions of
integrable systems was discovered by F.Magri \cite{magri} in the analysis of the
so-called Lenard scheme of constructing the KdV integrals. The basic idea of these
constructions is given by the following simple

\begin{lemma} \label{magri-lemma} Let $H_0$, $H_1$, \dots, be a sequence of
functions on $P$ satisfying the recursion relation
\beq\label{recursion-a}
\{ ~.~, H_{p+1}\}_1=\{~.~,H_{p}\}_2, ~~~p=0,\, 1, \, \dots
\eeq
Then
$$
\{ H_p, H_q\}_1=\{ H_p, H_q\}_2=0, ~~p, \, q=0, \, 1, \, \dots
$$
\end{lemma}

For convenience of the reader we reproduce the proof of the lemma. Let $p<q$
and $q-p=2m$ for some $m>0$. Using the recursion and antisymmetry of the
brackets we obtain
$$
\{ H_p, H_q\}_1=\{ H_p, H_{q-1}\}_2 =-\{ H_{q-1}, H_p\}_2=
-\{ H_{q-1}, H_{p+1}\}_1= \{ H_{p+1}, H_{q-1}\}_1.
$$
Iterating we arrive at
$$
\{ H_p, H_q\}_1= \dots = \{ H_{p+m}, H_{q-m}\}_1=0
$$
since $p+m=q-m$. Doing similarly in the case $q-p=2m+1$ we obtain
$$
\{ H_p, H_q\}_1 = \dots = \{ H_n, H_{n+1}\}_1=\{ H_n, H_{n}\}_2=0
$$
where $n=p+m=q-m-1$. The commutativity $\{ H_p, H_q\}_2=0$ easily follows from
the recursion. The Lemma is proved.

We have not used yet the compatibility condition of the two brackets. It turns
out to be crucial in contructing the Hamiltonians
satisfying the recursion relation (\ref{recursion-a}). 
There are two essentially
different realizations of the recursive procedure.

The first one applies to the case when the bihamiltonian structure is { \it 
symplectic} , i.e. $N=2n$ and the Poisson structures of the affine line
(\ref{biham1}) do not degenerate for generic $a_1$, $a_2$. Without loss of
generality one may assume nondegeneracy of $\{ ~, ~\}_1$. The {\it recursion
operator} 
$$
{\cal R}: TP \to TP
$$
is defined by
\beq\label{roperator} 
{\cal R}:= \{ ~,~\}_2 \cdot \{ ~, ~\}_1^{-1}.
\eeq
The main recursion relation (\ref{recursion-a}) can be rewritten in the form
\beq\label{recursion2}
dH_{p+1} ={\cal R}^* dH_p, ~~p=0, \, 1, \dots
\eeq
where
$$
{\cal R}^*: T^*P \to T^*P
$$
is the adjoint operator.

\begin{theorem} \cite{magri, mckean} The Hamiltonians
$$
H_p : ={1\over p+1 }\tr \, {\cal R}^{p+1}, ~~p\geq 0
$$
satisfy the recursion (\ref{recursion2}).
\end{theorem}

Clearly there are at most $n$ independent of these commuting functions. 
We say that the
bihamiltonian symplectic structure is {\it generic} if exactly $n$ of these
functions are independent. Let us denote $\lambda_i=\lambda_i(x)$ the
eigenvalues of the recursion operator.
Since the characteristic polynomial of ${\cal R}$
is a perfect square
$$
\det \left( {\cal R} -\lambda\right) =\prod_{i=1}^n (\lambda-\lambda_i)^2.
$$
only $n$ of these eigenvalues can be distinct, say,
$\lambda_1=\lambda_1(x)$, \dots, $\lambda_n=\lambda_n(x)$. For generic
bihamiltonian symplectic structure these
are independent functions on $P\ni x$.

\begin{theorem}\cite{magri, magri-8, mckean} 
Let $\{ ~, ~\}_{1,2}$ be a generic symplectic
bihamiltonian structure. Then

1) All the commuting Hamiltonians
$$
H_p ={1\over p+1} \tr\, {\cal R}^{p+1} = {1\over p+1} \sum_{i=1}^n
\lambda_i^{p+1}(x), ~~p=0, 1, \dots , n-1
$$
generate completely integrable systems on $P$.

2) The eigenvalues $\lambda_i(x)$ can be included in a coordinate system
$\lambda_1, \mu_1, \dots, \lambda_n, \mu_n$ in order to reduce the two Poisson
structures to a block diagonal form where the $i$-th block in $\{~,~\}_1$
and in $\{~,~\}_2$ reads, respectively
$$
\left( \matrix{ 0 & 1\cr -1 & 0\cr}\right), ~~~
\left( \matrix{ 0 & \lambda_i\cr -\lambda_i & 0\cr}\right), ~~i=1, \dots, n.
$$
\end{theorem}

The last formula gives the normal form of a generic symplectic bihamiltonian
structure. Therefore all such structures are equivalent w.r.t. the group of
local diffeomorphisms.

Let us now consider the degenerate situation. We assume that the Poisson
structure (\ref{biham1}) has constant rank for generic $a_1$ and $a_2$. Without
loss of generality we may assume that
\beq\label{corank}
k={{\rm corank}}\{ ~, \}_1 = {{\rm corank}}(\{ ~, \}_1+ \epsilon \{ ~, \}_2)
\eeq
for an arbitrary sufficiently small $\epsilon$. 

Let us first prove the following useful property of bihamiltonian structures of
the constant rank.

\begin{lemma}\label{2-1-4} Let the bihamiltonian structure satisfy (\ref{corank}). 
Then the Casimirs of $\{ ~, \}_1$ commute w.r.t. $\{ ~, \}_2$.
\end{lemma}

\pf Let $2m$ be the rank of $\{ ~, \}_1$. We first reduce the matrix
of this bracket to the canonical constant block diagonal form. Denote
$(h^{ab})$ the matrix of the second Poisson bracket in these coordinates.
Let us now choose two integers $i$, $j$ such that $2m<i<j\leq N=2m+k$ and form 
a $(2m+1)\times (2m+1)$ minor of the matrix $\{ ~, \}_1 +\epsilon \{ ~, \}_2$ 
by adding $i$-th column and $j$-th row to the principal $2m\times 2m$
minor standing in the
first $2m$ columns and
first $2m$ rows.
The condition (\ref{corank}) is equivalent to vanishing of the determinants of
all these minors. It is easy to see that the determinant in question is equal
to $-\epsilon \,h^{ij}+O(\epsilon^2)$. Therefore $h^{ij}=0$
for all pairs $(i,j)$ greater than $2m$. The lemma is proved.
\epf

\begin{cor}\label{2-1-5} 
For a compatible pair of Poisson brackets of the constant
${{\rm rank}} \left(\{~,~\}_2\right.$ 
$\left. -\lambda \{~,~\}_1 \right)=
{{\rm rank}}\{~,~\}_1$, $\lambda\to\infty$,
$$
\{~,~\}_2\in H^2(P, \{~,~\}_1)
$$
is a trivial cocycle.
\end{cor}

\pf What $\{~,~\}_2$ is a cocycle w.r.t. the Poisson cohomology of $\left( P,
\{~,~\}_1\right)$ follows from (\ref{biham2}). To prove triviality use
commutativity of the Casimirs of the first Poisson bracket
and also Lemma \ref{first-lemma}. \epf

A bihamiltonian structure with a marked line $\lambda\, \{~,~\}_1$ is called
{\it Poisson pencil}. Choosing another Poisson bracket $\{~,~\}_2$
of the pencil one can represent the brackets of the bihamiltonian structure  
in the
form
\beq\label{lambda-pencil}
\{~,~\}_\lambda := \{~,~\}_2-\lambda \{~,~\}_1.
\eeq
The representation (\ref{lambda-pencil}) is well-defined up to an affine change
of the parameter $\lambda$
$$
\lambda\mapsto a\, \lambda + b.
$$

In the case of Poisson pencils of constant rank the corank 
of $\{~,~\}_\lambda$ equals $k$ for $\lambda\to\infty$. The recursive
construction of the commuting flows in this case is given by the following
simple statement (cf. \cite{magri, gz1, panas}).

\begin{theorem}\label{2-1-6} Under the assumption (\ref{corank}) the coefficients of the
Taylor expansion 
\beq\label{hier-int}
c^\alpha(x,\lambda) =c_{-1}^\alpha(x) +{c_0^\alpha(x)\over \lambda} +
{c_1^\alpha(x)\over \lambda^2}+\dots,  ~~\lambda\to\infty
\eeq
of the Casimirs $c^\alpha(x,\lambda)$,
$\alpha=1, \dots, k$ of the Poisson bracket $\{~,~\}_\lambda$ 
commute with respect to both the Poisson brackets
$$
\{ c_p^\alpha, c_q^\beta\}_{1,2}=0, ~~\alpha, \beta=1, \dots, k, ~~p, q\geq -1.
$$
\end{theorem}

\pf Spelling out the definition of the
Casimirs
$$
\{ ~.~, c^\alpha\}_\lambda=0
$$
for the coefficients of the expansion (\ref{hier-int}) we must have first that
\beq\label{hier-casi}
\{ ~.~, c_{-1}^\alpha\}_1=0.
\eeq
That is, the leading coefficients of the Taylor expansions are Casimirs of
$\{~,~\}_1$. For the subsequent coefficients we get the recursive relations
\beq\label{hier_recur} 
\{~.~, c_{p+1}^\alpha\}_1 = \{~.~, c_{p}^\alpha\}_2, ~~~p=-1, \, 0, \, 1, \dots
\eeq
From (\ref{hier_recur}) and Theorem 1 it follows that
$$
\{ c^\alpha_p, c^\alpha_q\}_{1,2}=0, ~~p, q\geq -1.
$$
The commutativity $\{ c^\alpha_p, c^\beta_q\}_{1,2}=0$ for $\alpha\neq\beta$
easily follows from the same recursion trick and from commutativity of the
Casimirs 
\beq\label{casi-comm}
\{ c^\alpha_{-1}, c^\beta_{-1}\}_2=0
\eeq
proved in Lemma \ref{2-1-4}. The theorem is proved. \epf

\begin{exam}\label{e-2-1-7} 
According to Corollary \ref{2-1-5} there exists a vector field $Z$ such
that 
$$
Lie_Z \{~,~\}_1 =\{~,~\}_2.
$$
We say, following \cite{casati} that the bihamiltonian structure is 
{\rm exact} if the vector field
$Z$ can be chosen in such a way that
\beq\label{exact}
\left( Lie_Z\right)^2 \{~,~\}_1 =0.
\eeq
For an exact bihamiltonian structure the generating functions (\ref{hier-int})
of the commuting Hamiltonians $c^\alpha_p(x)$ have the form
\beq\label{exact-int}
c^\alpha(x;\lambda) =\exp\left(-Z/\lambda\right) c^\alpha_{-1}(x)
=c^\alpha_{-1}(x) -{1\over \lambda} \pal_Z c^\alpha_{-1}(x) +{1\over \lambda^2}
\pal_Z^2 c^\alpha_{-1}(x)  \dots
\eeq
for every $\alpha=1, \, \dots, k$. 
\end{exam}
This formula can be easily proved by choosing
a system of local coordinates $x^1, \dots, x^N$ on the phase space $P$ such that
the vector field $Z$ corresponds to the shift along $x^1$. In these coordinates
the tensor of the first Poisson bracket depends linearly on $x^1$ and the 
second Poisson bracket is $x^1$-independent. The Poisson pencil $\{~,~\}_\lambda$
is obtained from $\{~,~\}_1$ by the shift $x^1 \mapsto x^1 -1/\lambda$ and by 
multiplication by $-\lambda$.

Conversely, if, in a given coordinate system, $\{~,~\}_1$ depends linearly
on one of the coordinates and $\{~,~\}_2$ does not depend on this coordinate
then the bihamiltonian structure is exact. In particular this trick can be
applied to the standard linear Lie - Poisson structures on the dual spaces
to Lie algebras (cf \cite{manakov}). In this case it was called in
\cite{mischenko-fomenko} the method of argument translation.

All our bihamiltonian structures on the loop spaces to be studied below
will be exact. However, at the moment we do not see their Lie algebraic
origin.

The construction of Theorem \ref{2-1-6} 
for a bihamiltonian structure of the constant
corank $k$
produces $k$ chains of pairwise commuting bihamiltonian flows
\beq\label{hier-flows}
{dx\over dt^{\alpha,p} } =\{ x, c^\alpha_p\}_1 = \{ x, c^\alpha_{p-1}\}_2,
~~\alpha=1, \dots, k, ~~~p=0, 1, 2, \dots
\eeq
The chains are labeled by the Casimirs $c^\alpha_{-1}$ of the first Poisson
bracket. The level $p$ in each chain corresponds to the number of iterations of
the recursive procedure (we will keep using this expression although the
recursion operator is not defined in the degenerate case). All the family of
commuting flows organized by the above recursion procedure is called {\it
the hierarchy} determined by the bihamiltonian structure. 

The hierarchy structure of the constructed family of commuting flows depends 
nontrivially on the choice of $\{~,~\}_1$ in the
Poisson pencil (\ref{biham1}). On the contrary, a different choice of the
second Poisson bracket in the pencil produces a triangular linear transformation
of the commuting Hamiltonians, i.e., to the Hamiltonians of the level $p$ it
will be added a linear combination of the Hamiltonians of the lower levels.

In the finite dimensional case we are discussing now all the chains of the
hierarchy will be finite. In other words, the generating functions
(\ref{hier-int}) of the commuting flows will become polynomials after
multiplication by a suitable power of $\lambda$ (the degrees of an appropriate 
system of these polynomials correspond to the {\it type} of the bihamiltonian
structure \cite{panas}).
A simple necessary
and sufficient condition of complete integrability of the flows of the hierarchy
was found by A.Brailov and A.Bolsinov (see in \cite{bolsinov}). The problem
of normal forms of degenerate bihamiltonian structures has been studied by
I.M.Gelfand and I.Zakharevich \cite{gz1}, \cite{gz2} for the case of the corank 1
and by I.Zakharevich \cite{zakhar} and A.Panasyuk \cite{panas} for higher coranks.

\subsubsection{Construction of bihamiltonian hierarchies on the extended 
loop \\ spaces}\label{sec-3-1-2}\par

In the remaining part of the paper we will study bihamiltonian structures
on the extended loop spaces $\hat{\cal L}(M)$ assuming $M$ to be a
$n$-dimensional ball. Moreover we restrict ourselves at considering $(0,n)$
bihamiltonian structures, i.e., of the form 
\beq\label{mainclass1}
\{ u^i(x), u^j(y)\}_{1,\, 2} =\sum_{k=0}^\infty \epsilon^k \{ u^i(x),
u^j(y)\}_{1,\, 2}^{[k]}, ~~i, j=1, \dots, n
\eeq
with
\beq\label{mainclass2}
\{ u^i(x), u^j(y)\}_{1,\, 2}^{[k]}=\sum_{s=0}^{k+1} {A^{ij}_{k,s}}_{1,2}(u; u_x,
\dots, u^{(s)})\delta^{(k-s+1)}(x-y)
\eeq
$$
{A^{ij}_{k,s}}_{1,2}(u; u_x,
\dots, u^{(s)})\in {\cal A}, ~~\deg {A^{ij}_{k,s}}_{1,2}=s, 
$$
both satisfying the condition
\beq\label{mainclass3}
\det {A^{ij}_{0,0}}_1(u)\neq 0, ~~ \det {A^{ij}_{0,0}}_2(u)\neq 0.
\eeq

\begin{exam}
 For $n=1$ take 
\beq\label{kdv-ham2}
\{ u(x), u(y)\}_2 = u(x) \delta'(x-y) +{1\over 2} u' \delta (x-y)
-{1\over 4}\epsilon^2\delta'''(x-y) 
\eeq
This is the Lie - Poisson bracket on the dual space to the Virasoro algebra
(see Example \ref{lie-poiss} above). 
It depends linearly on $u(x)$. Taking the Lie
derivative of this bracket along the vector field $\pal / \pal u$ we obtain
another Poisson bracket
\beq\label{kdv-ham1} 
\{ u(x), u(y)\}_1 =\delta'(x-y).
\eeq
This is an exact bihamiltonian structure in the sense of Example \ref{e-2-1-7}
(the roles of $\{~,~\}_1$ and $\{~,~\}_2$ have been interchanged, $Z=\pal /\pal
u$). The 
Casimir $c([u];\lambda)$ of the Poisson pencil $\{~,~\}_2 -\lambda\{~,~\}_1$
is determined from the following third order equation
\beq\label{3d-order}
-{\epsilon^2\over 4} y''' + u \, y' + {1\over 2} u' y  = \lambda\, y'
\eeq
where
$$
y:={\delta c([u]; \lambda)\over \delta u(x)}.
$$
Starting from the Casimir
$$
I_{-1}[u]=\int u \, dx
$$
of (\ref{kdv-ham1}) and applying the recursive procedure of the previous 
section
we produce an infinite sequence of commuting Hamiltonians of the
KdV hierarchy
\beq\label{2-2-6-b}
{\pal u\over \pal{t^k}} = \pal_x {\delta I_k[u]\over \delta u(x)} =\left[ 
-{1\over 4} \epsilon^2 \pal_x^3 + u \pal_x +{1\over 2} u_x\right] {\delta
I_{k-1}\over \delta u(x)}, 
~~k=0, 1, \dots
\eeq
where the first few equations of the hierarchy and their Hamiltonians have the
form 
\eqa\label{kdv-ham}
&&I_0[u]=\int {u^2\over 2}dx,\quad 
I_1={1\over 4}\int \left( \epsilon^2{{u'}^2 \over 2} +u^3\right)\, dx,\nn\\
&&
I_2[u] = {1\over 16}\int \left( \epsilon^4{{u''}^2 \over 2} 
+ 5\epsilon^2 u\, {u'}^2 + {5\over 2}
u^4\right)\, 
dx \, , \dots,
\eeqa 
\eqa\label{kdv-eq}
&&{\pal u\over \pal{t^0}} =u_x, \quad{\pal u\over\pal{t^1}}=
{1\over 4}( 6 u\, u' - \epsilon^2 u'''), \nn\\
&&
{\pal u\over\pal {t^2}} = {1\over 16}[
30 u^2 u'-10 \epsilon^2 ( 2 u' u''
+ u \, u''') + \epsilon^4 u^V]\, , \dots
\eeqa
\end{exam}
The above algorithm of constructing the KdV hierarchy is due to Lenard (see in
\cite{ggkm74}). The bihamiltonian interpretation of it is due to Magri
\cite{magri}. The generating function of the Hamiltonian densities 
\beq\label{resolvent}
{1\over 2 \, \chi_R} = {1\over 2\, \sqrt{\lambda}} + \sum_{j=0}^\infty 
{R_j\over \lambda^{2j+1\over 2}}
\eeq 
\beq
H_k = -4 \int R_{k-2}dx
\eeq
that is, the density of the Casimir of the Poisson pencil $\{~,~\}_2
-\lambda\{~,~\}_1$
coincides, according to I.M.Gelfand and L.A.Dickey \cite{gelfand-dickey}
with the diagonal  of the resolvent of the Lax operator
$$
L=-\epsilon^2 \pal_x^2 + u(x).
$$
The series in (\ref{resolvent}) is to be understood in the sense of the
asymptotic expansion at $\lambda \to\infty$. The coefficients of the asymptotic
expansion do not depend on the choice of the boundary conditions.

Let us return back to the general case (\ref{mainclass1})--(\ref{mainclass3}).
We will apply now the results of the first part of the paper 
to prove existence of the commuting hierarchy for an arbitrary bihamiltonian
structure of the above form. Indeed, according to these results for every $\lambda$
the Poisson bracket
\beq\label{pencil}
\{u^i(x),u^j(y)\}_\lambda:= \{u^i(x),u^j(y)\}_2 - \lambda \{u^i(x),u^j(y)\}_1
\eeq
can be reduced to the constant form (\ref{delta'}). We will now prove that the
coefficients of the Taylor expansion of the reducing transformation are
densities of the commuting Hamiltonians. 

\begin{theorem}\label{t-2-2-2}
 There exists a transformation of the form
\beq\label{reducing} 
u^i \mapsto \tilde u^i =F^i(u; u_x, u_{xx}, \dots; \epsilon;\lambda)= 
\sum_{p=-1}^\infty {f^i_p(u; u_x, u_{xx}, \dots ;\epsilon) \over \lambda^{p+1}}
\eeq
reducing the Poisson bracket (\ref{pencil}) to the normal form
\beq\label{reduced}
\{ \tilde u^i(x), \tilde u^j (y)\}_\lambda = - \lambda \, \eta^{ij} \delta'
(x-y).
\eeq
The coefficients of the expansion
$$
f^i_p (u; u_x, u_{xx}, \dots ; \epsilon) \in {\cal A}\otimes
{\mathbb C}[[\epsilon]]
$$
are densities of pairwise commuting integrals
\beq\label{coeff}
\bar f^i_p := \int f^i_p (u; u_x, u_{xx}, \dots ; \epsilon)\, dx
\eeq
$$
\{\bar f^i_p, \bar f^j_q\}_{1, 2} =0, ~~i, j=1, \dots, n, ~~~ p, q= -1, 0, 1, 2,
\dots
$$
\end{theorem}

\pf Let us prove first that the reducing transformation for the Poisson pencil
(\ref{pencil})
can be chosen in the form of the series (\ref{reducing}). Let 
$$
u^i\mapsto \hat
u^i=f^i_{-1}(u; u_x,
\dots; \epsilon)
$$ 
be the Miura-type transformation reducing the first Poisson
bracket to the normal form $\eta^{ij}\delta'(x-y)$ with a constant nodegenerate
symmetric matrix $\eta^{ij}$. The compatibility condition of the Poisson pencil says
that the second Poisson bracket is a 2-cocycle of the first one. Due to
triviality of the 2-cohomology there exists an infinitesimal transformation
$$
\hat u^i \mapsto \hat u^i + {\hat f^i_0(u; u_x, \dots ;\epsilon)\over \lambda}
$$
reducing the pencil to the form
$$
-\lambda \left[ \eta^{ij} \delta' (x-y) +O\left( {1\over
\lambda^2}\right)\right].
$$
Iterating this procedure we will kill all the terms in the $1/\lambda$ expansion
of the pencil (\ref{pencil}). The superposition of the Miura-type
transformations gives the series (\ref{reducing}).

To prove commutativity of the functionals (\ref{coeff}) it suffices to observe
that 
$$
\int \tilde u^i dx = \sum_{p=-1}^\infty {\bar f^i_p \over \lambda^{p+1}}, ~~i=1,
\dots, n
$$
are the Casimirs of the Poisson pencil (\ref{pencil}). This immediately follows from the
reduced form of it (\ref{reduced}). The theorem is proved. \epf


\begin{cor} Every bihamiltonian structure of the form (\ref{mainclass1}) -
(\ref{mainclass3}) generates a commuting bihamiltonian hierarchy of evolutionary
PDEs 
\beq\label{mainclass4}
{\partial u^i\over\partial t^{j,p}} = \{ u^i(x), \bar f^j_p\}_1= 
\{ u^i(x), \bar f^j_{p-1}\}_2, ~~j=1, \dots, n, ~~p=0, 1, 2, \dots
\eeq
\end{cor}

\begin{exam}\label{e-2-2-4}
 Let us apply the above procedure to the bihamiltonian structure\newline
(\ref{kdv-ham1}), (\ref{kdv-ham2}) of the KdV hierarchy. We already know 
that (\ref{kdv-ham2}) is reduced to the constant form by means of the Miura
transformation $u\mapsto \chi$,
$$
i\epsilon\chi' -\chi^2 =u.
$$
The whole Poisson pencil (\ref{pencil}) built of (\ref{kdv-ham1}) and (\ref{kdv-ham2})
is obtained from (\ref{kdv-ham2}) by a shift $u\mapsto u-\lambda$
(i.e., this is an exact bihamiltonian structure). Therefore the
reducing transformation for the pencil has the form
$$
u\mapsto \tilde u = -2 \sqrt{\lambda} [\chi - \sqrt{\lambda}],
$$
\beq\label{2-2-16}
\{ \tilde u(x), \tilde u(y)\}_\lambda = - \lambda \, \delta'(x-y).
\eeq
Here $\chi$ is the unique solution to the Riccati equation 
$$
i \epsilon\chi' -\chi^2 = u-\lambda
$$
of the form
$$
\chi = k + \sum_{m=1}^\infty {\chi_m \over k^m}, ~~k=\sqrt{\lambda}.
$$
The coefficients $\chi_m$ are polynomials in $u$, $\epsilon u'$, 
$\epsilon^2 u''$ etc.,
$$
\chi_1 =-u/2, ~~\chi_2=-i \epsilon u'/4, ~~\chi_3={1\over 8}(\epsilon^2 u''-u^2), \dots
$$
The even coefficients all are total derivatives; the odd ones give the
Hamiltonians of the KdV hierarchy
$$
{\pal u\over \pal t^k} = \pal_x {\delta I_k\over \delta u(x)}
$$
where
$$
I_k =-4\int \chi_{2k+3} dx, ~~k=0, 1, \dots
$$
(we redenote $f^1_k=-4 \chi_{2k+3}$, $t^k=t^{1,k}$).
\end{exam}

This algorithm of constructing the conservation laws of the KdV equation was
found by Gardner, Kruskal and Miura in \cite{gkm68}. The equivalence of the
family of KdV integrals produced by the Lenard scheme to those produced by the
Gardner, Kruskal and Miura algorithm was established by G.Wilson \cite{wilson}.
 To our best knowledge the 
Hamiltonian nature of
this algorithm was not discussed in the literature.

\begin{remark} It can be shown that the even part
$$
\sum {\chi_{2\,l}\over \lambda^l}
$$
of the reducing transformation can be gauged out by a Miura-type transformation
preserving the canonical form (\ref{2-2-16}).
\end{remark}

\setcounter{equation}{0}
\setcounter{theorem}{0}
\subsection{The leading order of $(0,n)$ bihamiltonian structures}
\label{sec-3-2}\par

The main problem we address in this paper is the classification of Poisson
pencils of the form (\ref{mainclass1})-(\ref{mainclass3}) on the extended
loop space w.r.t. the action of Miura group. Actually we will add below further
restrictions on the class of bihamiltonian structures to be classified. 

As we already did above in solving the problem of normal forms of a single
Poisson bracket on the extended loop space, we can try to classify Poisson
pencils by successive approximations, i.e., first reducing to a normal form
the leading term of the series (\ref{mainclass1}), then studying 2-cocycles on
this leading term etc.

It is clear that the leading term
\beq\label{lead1}
\{u^i(x), u^j(y)\}_\lambda ^{[0]} =\{u^i(x), u^j(y)\}_2^{[0]} -
\lambda \{u^i(x), u^j(y)\}_1^{[0]}
\eeq
is itself a Poisson pencil. Let us redenote the coefficients of the
pencil as follows
\eqa\label{pbht}
&& 
\{u^i(x), u^j(y)\}_1^{[0]}= g^{ij}_1(u(x)) \delta' (x-y) + \Gamma_{1\, k}^{ij}(u)
u^k_x \delta (x-y), \nn\\
&& 
\{u^i(x), u^j(y)\}_2^{[0]} =g^{ij}_2(u(x)) \delta' (x-y) 
+ \Gamma_{2\, k}^{ij}(u)
u^k_x \delta (x-y).
\eeqa
According to our assumption (\ref{mainclass3}) the leading coefficients
$g^{ij}_1(u)$ and $g^{ij}_2(u)$ define two flat geometries on $M$. Under what
conditions these two flat geometries correspond to a bihamiltonian structure
(\ref{pbht}) on $M$? It is clear that the two flat geometries must be members of
a {\it flat pencil}, i.e., the linear combination $\Gamma_{2\, k}^{ij}
-\lambda \Gamma^{ij}_{1 \, k}$ of the Christoffel coefficients of the two
metrics gives the Christoffel coefficients of the linear combination of the
metrics $g^{ij}_2-\lambda g^{ij}_1$ and the latter is a flat (contravariant)
metric for an arbitrary $\lambda$.  The notion of the flat pencil was introduced
by one of the authors in \cite{coxeter} (see also \cite{D3}, \cite{taniguchi}).
It was shown in \cite{coxeter} that flat pencils are parametrized by solutions
to a certain system of nonlinear PDEs that resembles the equations of
associativity. Let us briefly recall this parametrization.

{\bf Definition.} A {\it quasi-Frobenius structure} on a manifold $M$ is 
a pair $(f, <\ ,\ >)$ where $f=f_\alpha(v)dv^\alpha$ is a one-form and 
$<\ ,\ >$
is a symmetric bilinear nondegenerate form on the tangent planes. 
This form will be called metric on $M$. The following conditions must hold true.

1) the metric $<\ ,\ >$ is flat

2) in the flat coordinates $v^1$, \dots, $v^n$ for the flat metric,
$$
<\pal_\alpha,\pal_\beta> = \eta_{\alpha\beta}={{\rm const}}, ~~\pal_\alpha := {\pal\over\pal
v^\alpha}
$$
the multiplication table 
\beq
\pal_\alpha\cdot\pal_\beta = \pal_\alpha\pal_\gamma f_\beta(v)\eta^{\gamma\mu}
\pal_\mu
\eeq
(here $(\eta^{\alpha\beta}) =
(\eta_{\alpha\beta})^{-1}$)
defines on the tangent spaces $T_vM$ a structure of algebra satisfying
the following conditions

1)
\beq\label{2-3-4}
(a\cdot b)\cdot c=(a\cdot c)\cdot b.
\eeq

2) The bilinear form $<\ ,\ >$ is  {\it invariant} on the algebra $T_vM$, i.e.
\beq\label{2-3-5}
<a\cdot b, c> = <a, c\cdot b>.
\eeq
Here $a$, $b$, $c$ are arbitrary three vectors in $T_vM$.

3) The bilinear form on the cotangent bundle defined by
\beq
( dv^\alpha, dv^\beta) = 
\eta^{\alpha\lambda}\eta^{\beta\mu}(\pal_\lambda f_\mu +
\pal_\mu f_\lambda)
\eeq
is an invariant form  for the dual algebra
$$
dv^\alpha\cdot dv^\beta = \eta^{\alpha\lambda} \eta^{\beta\mu} \pal_\lambda
\pal_\gamma f_\mu(v) \, dv^\gamma
$$
on $T^*_vM$.
\smallskip 

Warning: the name ``quasi-Frobenius manifold'' was also used in different
senses in \cite{verdier, manin}!

\begin{theorem} For an arbitrary quasi-Frobenius manifold the two metrics
$<\ ,\ >$ and $(~,~)$ form a flat pencil. 
Conversely, any flat pencil can
locally be obtained in such a way.
\end{theorem}

\begin{remark} The metric $(~,~)$ could be degenerate. However, the linear
combination $(~,~)-\lambda<\ ,\ >$ does not degenerate for generic $\lambda$.
\end{remark}

\begin{exam} In the particular case of quadratic functions in $v$
$$
f_\beta(v) = {1\over 2} c_{\beta\lambda\mu} v^\lambda v^\mu
$$
the structure constants
$$
c^\gamma_{\alpha\beta} =\eta^{\gamma\nu}
c_{\nu\alpha\beta}
$$
define on the linear space $M=T_vM$ a structure of
{\rm Novikov algebra}, i.e., an algebra satisfying (\ref{2-3-4}) 
and also {\em
left-symmetric}, i.e.,
\beq
a\cdot (b\cdot c) - (a\cdot b)\cdot c = b\cdot (a\cdot c) - (b\cdot a)\cdot c
\eeq
(the algebras satisfying the last identity are called Vinberg algebras or also
pre-Lie algebras). The inner product $<\ ,\ >$ satisfying (\ref{2-3-5})
is called an invariant bilinear form on the Novikov algebra.
The above second Poisson bracket in this case is a linear Poisson bracket of
hydrodynamic type. The theory of such Poisson brackets was first studied by
A.Balinsky and S.P.Novikov \cite{bal-nov}. 
\end{exam}
 
Proof of the theorem see in \cite{coxeter}, \cite{D3}.

The problem of local classification of quasi-Frobenius manifolds can be reduced
to the theory of the following system of nonlinear PDEs
\beq
\pal_\alpha\pal_\lambda f_\beta \eta^{\lambda\mu} \pal_\mu \pal_\delta f_\gamma
= 
\pal_\alpha\pal_\lambda f_\gamma \eta^{\lambda\mu} \pal_\mu \pal_\delta f_\beta
\eeq
\beq
(\pal_\alpha f_\lambda+ \pal_\lambda f_\alpha) \eta^{\lambda\mu} \pal_\mu
\pal_\beta f_\gamma 
= (\pal_\beta f_\lambda+ \pal_\lambda f_\beta) \eta^{\lambda\mu} \pal_\mu
\pal_\alpha f_\gamma.
\eeq
Very recently O. Mokhov \cite{mokhov} and E.Ferapontov \cite{ferapontov}
proved integrability of this system.

As it was essentially shown in \cite{coxeter} (see the precise formulation in
\cite{taniguchi}) that, under certain quasihomogeneity assumption the
quasi-Frobenius structure on $M$ reduces to a Frobenius one. In the next 
section
we impose a somewhat different additional constraint onto the 
bihamiltonian structure
(\ref{pbht}) that will also give a correspondence between the normal forms 
of the bihamiltonian structures
(\ref{mainclass1}) on the extended loop spaces ${\cal L}(M)$ and Frobenius     
structures on $M$.

\setcounter{equation}{0}
\setcounter{theorem}{0}

\subsection{Tau-structures, tau-covers and normal coordinates of bihamiltonian hierarchies
of PDEs}\par

The densities of the commuting hamiltonians of the hierarchy (\ref{mainclass4})
constructed in the
Section \ref{sec-3-1-2} 
admit certain freedom in their definition. Indeed, instead of
the densities $f^i_p$ one can take 
\beq\label{change1}
h_{i,p} =\sum_{q=-1}^p \sum_{j=1}^n a_{ij, pq} f^j_q +{{\rm total}} ~{{\rm
derivative}}
\eeq
with arbitrary constant coefficients $a_{ij,pq}$ satisfying the following
nondegeneracy condition
\beq\label{change2}
\det\left(a_{ij,pp}\right)\neq 0, ~~p=-1, 0, 1, \dots.
\eeq
The flows of the new hierarchy
\beq\label{change3}
{\partial\over\pal \tilde t^{i,p}} =\{ ~.~, \bar h_{i,p}\}_1, ~~\bar h_{i,p} :=\int
h_{i,p}dx, ~~ i=1, \dots, n, ~~p=0, 1, 2, \dots
\eeq
still commute.

{\bf Definition.} 1). A {\it tau-structure} for a Poisson pencil
(\ref{pencil}) is a choice of the densities of the commuting hamiltonians in the
class of equivalence w.r.t. the transformations (\ref{change1}) s.t.
the 1-form
\beq\label{closed}
\omega =\sum_{p=0}^\infty \sum_{i=1}^n h_{i,p-1} d\tilde t^{i,p}
\eeq
is closed, i.e.,
\beq\label{sc}
{\pal h_{i,p-1}(u; u_x, \dots; \epsilon)\over \pal t^{j,q}}=
{\pal h_{j,q-1}(u; u_x, \dots; \epsilon)\over \pal t^{i,p}}.
\eeq
A Poisson pencil admitting a tau-structure is called {\it
tau-symmetric}.

2). We say that the tau-symmetric Poisson pencil is {\it compatible with
spatial translations} if
\beq\label{translation}
\{ ~.~, {\bar {h}}_{1,0}\}_1 
= a_{11, 00}\{ ~.~, {\bar {f}}^1_0\}_1 = {\pal\over \pal x}.
\eeq 

Clearly tau-symmetry is a geometric property of a Poisson pencil. 
Indeed, if a
Poisson pencil $\{ ~,~\}_\lambda^{\tilde{}}$ is obtained from $\{ ~,~\}_\lambda$
by means of a Miura-type transformation
$$
u^i\mapsto \tilde u^i =F^i(u; u_x, \dots; \epsilon)
$$
then the densities $h_{i,p}(u; u_x, \dots;\epsilon)$ satisfying (\ref{change1}),
(\ref{sc}) considered as functions
of the new coordinates $\tilde u^i$ will give a tau-structure for the pencil
$\{~,~\}_\lambda^{\tilde{}}$ with the same coefficients $a_{ij,pq}$.

The functions
\beq\label{normal1}
h_{i,-1}=h_{i,-1}^{[0]}(u) +\epsilon\, h_{i,-1}^{[1]}(u; u_x) + \dots, ~~i=1,
\dots, n
\eeq
will be of particular importance for working with tau-symmetric bihamiltonan
hierarchies. 

\begin{lemma} The functions (\ref{normal1}) define a Miura-type transformation
\beq\label{normal2}
u^i\mapsto \tilde u_i:= h_{i,-1}, ~~i=1, \dots, n .
\eeq
The functionals $\int\tilde u_1 dx$, \dots, $\int\tilde u_n dx$
are Casimirs of $\{~,~\}_1$.
\end{lemma}

\pf From the definition we have
$$
h_{i,-1}(u; u_x, \dots;\epsilon) = \sum_j a_{ij, -1 \, -1} f^j_{-1}(u; u_x,
\dots; \epsilon)+{{\rm total}} ~ {{\rm derivative}}
$$
where the functions $f^j_p(u;u_x,u_{xx},\dots;\epsilon)
\in {\cal A}_{0,0}(\epsilon)$ were constructed in the proof of 
Theorem \ref{t-2-2-2}.
In particular
$$
h_{i,-1}^{[0]}(u) =  \sum_j a_{ij, -1 \, -1} {f^j_{-1}}^{[0]}(u)
$$
where
$$
f^j_{-1}(u; u_x, \dots) = {f^j_{-1}}^{[0]}(u)+ O(\epsilon).
$$
Since $u^j\mapsto f^j_{-1}(u; u_x, \dots)$ is the transformation reducing
$\{ ~,~\}_1$ to the normal form, the nondegeneracy
$$
\det \left( {\pal {f^j_{-1}}^{[0]}(u)\over \pal u^k}\right) \neq 0
$$
holds true. This implies the nondegeneracy 
$$
\det \left( {\pal {h_{i,-1}}^{[0]}(u)\over \pal u^k}\right) \neq 0
$$
also for the transformation (\ref{normal2}). It remains to observe that
the functionals 
$$
\int f^j_{-1}(u; u_x, \dots; \epsilon)\, dx
$$
are Casimirs of the first Posson bracket. Taking an invertible linear
combination
of them and adding a total derivative will still give a system of Casimirs. 
\epf

{\bf Definition.} The dependent variables $\tilde u_1$ \dots, $\tilde u_n$
are called {\it normal coordinates} on ${\cal L}(M)$ w.r.t. the given
tau-structure.

The first Poisson bracket in the normal coordinates has the form
\beq\label{normal3}
\{\tilde u_i(x), \tilde u_j(y\}_1 = \eta_{ij}\delta'(x-y) + O(\epsilon)
\eeq
with a constant invertible matrix $(\eta_{ij})$. The variables 
$$
\tilde u^i :=\eta^{ij}\tilde u_j, ~~ \left( \eta^{ij}\right) :=\left(
\eta_{ij}\right)^{-1}
$$
will also be called the normal coordinates. The tilde over the normal
coordinates will often be omitted. The equations of the hierarchy in the normal
coordinates are written in the following form
\beq\label{2-4-9}
{\pal u_i\over \pal t^{j,q}} = \pal_x \Omega_{i,0; j,q}(u; u_x, \dots
;\epsilon)
\eeq
where the matrix $\Omega_{i,p;j,q}(u; u_x, \dots; \epsilon)$ is defined
in the formula (\ref{stress}) below.

We will settle below the problem of
ambiguity in the choice of tau-structure of a given Poisson pencil under
certain additional assumption about the pencil.
In particular we will describe the freedom in the choice of normal coordinates.

Let us denote 
$$
{\cal K}={\cal K}(M; \{~,~\}_1,  \{~,~\}_2):=
\cap_\lambda {\rm Ker}\, \left(  \{~,~\}_2 -\lambda \,  \{~,~\}_1\right)\in 
\Lambda_0(M)
$$
the subspace spanned by the commuting Hamiltonians
$\bar h_{i,p}$ of the hierarchy. This can be lifted to a subspace
$$
\hat {\cal K}=\in {\cal A}_{0,1}(M)
$$
w.r.t. to the factorization map ${\cal A}_{0,1}(M)\stackrel{\pi}\rightarrow 
\Lambda_0(M)$
(all the notations as in Section \ref{sec-2}). 
Namely, the Hamiltonian $\bar h_{i,p}$ lifts to the density $h_{i,p}$ satisfying (\ref{change3}).
We will now define an important symmetric {\it product map}
\beq\label{productmap}
*:\hat{\cal K} \times \hat{\cal K} \to {\cal A}_{0,0}(M)
\eeq
as follows.
Due to commutativity of the flows there exists an infinite matrix of the
densities of the fluxes of the conserved quantities 
\eqa
&&\Omega_{i,p; j,q}(u; u_x, u_{xx}, \dots;\epsilon)
=\sum_{k=0}^\infty \epsilon^k 
\Omega_{i,p; j,q}^{[k]}(u; u_x, \dots, u^{(k)}) \in {\cal A}_{0,0}, \\
&&\quad 1\leq i, j\leq n, ~~p, q =0, 1, 2, \dots\nn
\eeqa
such that
\beq\label{stress}
{\pal h_{i,p-1}(u; u_x, \dots; \epsilon)\over \pal t^{j,q}}
=\pal_x \Omega_{i,p; j,q}(u; u_x, u_{xx}, \dots;\epsilon)
\eeq
(we will omit tilde over the time variables always assuming the hierarchy to be
written in the normal coordinates). 
Observe the shift $p\mapsto p-1$ in the level of the conserved quantity
$\bar h_{i, p-1}$. The matrix $\Omega_{i,p;j,q}$ is determined up to adding
of constants. 

\begin{lemma} The collection of the densities $h_{i,p}(u;u_x, \dots;\epsilon)$
of the commuting \newline
Hamiltonians $\bar h_{i,p}$ of the form (\ref{change1}),
(\ref{change2}) is a tau-structure {{\rm iff}} the matrix $\Omega_{i,p; j,q}
= \Omega_{i,p; j,q}(u;
u_x, \dots;\epsilon)$ is symmetric
\beq\label{sym1}
\Omega_{i,p; j,q}=\Omega_{j,q;i,p}.
\eeq
\end{lemma} 

\pf The condition of closedness of the 1-form (\ref{closed}) 
written in (\ref{sc}) 
is equivalent to the symmetry (\ref{sym1}). \epf

We define the product map (\ref{productmap}) on the basis of the space
${\cal K}$ putting
\beq\label{productmap1}
(h_{i,p}, h_{j,q}) \mapsto h_{i,p} * h_{j,q} := \Omega_{i,p+1; j,q+1
}.
\eeq

All the equations of the constructed hierarchy have the form (\ref{pde00}). 
We define 
solutions to the hierarchy as vector functions
$$
u^i=u^i({\bf t};\epsilon)=\sum_{k=0}^\infty \epsilon^k u^i_{[k]}({\bf t}),
 ~~i=1, \dots, n, ~~{\bf t}= \left(
t^{i,p}\right)_{1\leq i \leq n, ~ p=0, 1, \dots}
$$
satisfying all the equations of the hierarchy as formal series in $\epsilon$.
Compatibility of the equations of the hierachy gives a possibility to solve them
simultaneously. Actually there is a problem of the definition of an appropriate
class of functions of infinite number of variables. This will be done below
where we will also construct the solution to the Cauchy problem for the
hierarchy with the initial data in a suitable class of functions of $x$.

\begin{cor} For every solution 
$$
u^i=u^i({\bf t};\epsilon), ~~i=1, \dots, n
$$
of a tau-symmetric hierarchy there exists a function 
$\tau=\tau({\bf t};\epsilon)$ such
that
$$
\Omega_{i,p; j,q}(u({\bf t};\epsilon); u_x({\bf t};\epsilon), \dots;\epsilon)
=\epsilon^2 {\pal^2\log \tau\over \pal t^{i,p}\pal t^{j,q}}.
$$
The function $\tau({\bf t}; \epsilon)$ is determined uniquely up to a
transformation of the form
$$
\tau({\bf t};\epsilon)\mapsto e^{{1\over 2 \epsilon^2}\sum 
Q_{i,p; j,q}(\epsilon) t^{i,p}t^{j,q}} \tau
({\bf t};\epsilon)
$$
with a constant symmetric matrix $Q_{i,p; j,q}(\epsilon)$.
\end{cor}

The function $\tau$ will be called {\it tau-function} of a given tau-structure
of the bihamiltonian hierarchy. 

For a tau-structure compatible
with spatial translations
we have
$$
{\pal\over \pal t^{1,0}} = {\pal\over\pal x}.
$$
Therefore the solution to the hierarchy,
when written in the normal coordinates reads
\beq\label{normal4}
u_i({\bf t};\epsilon) = \epsilon^2 {\pal^2 \log \tau\over \pal x \pal t^{i,0}}, ~~i=1,
\dots, n.
\eeq

In the study of symmetries of integrable PDEs it will be useful
the following extension of the PDEs  to systems with an infinite family of
dependent variables $f$, $f_{i,p}$, $i=1, \dots, n$, $p=0$, $1$, \dots,
$u_1$, \dots, $u_n$:
\eqa
&&
\epsilon\, {\partial f\over \partial t^{j,q}} = f_{j,q}
\nn\\
&&
\epsilon\, {\partial f_{i,p}\over \pal t^{j,q}} = \Omega_{i,p; j,q}(u; u_x,
\dots; \epsilon)
\label{2-4-16}\\
&&
{\pal u_i\over \pal t^{j,q}} = \pal_x \Omega_{i,0; j,q}(u; u_x, \dots
;\epsilon).
\nn
\eeqa

{\bf Definition.} The system (\ref{2-4-16}) of PDEs is called {\it tau-cover}
of the hierarchy (\ref{2-4-9}).

It is clear that the hierarchy (\ref{2-4-16}) is still commutative.
\medskip

The notion of tau-function was discovered by Date, Jimbo, Kashiwara and Miwa
\cite{djkm83} in their study of KP hierarchy. In the next section we will
construct the tau-structure for the KdV hierarchy.

\setcounter{equation}{0}
\setcounter{theorem}{0}

\subsection{Example: tau-structure of KdV and heat kernel expansion}\par

Let us define the differential polynomials 
$$
h_k=h_k(u, \epsilon\, u_x, \epsilon^2u_{xx}, \dots, \epsilon^{2k+2}
u^{(2k+2)})
$$
as the Seeley coefficients \cite{seeley} of the Lax operator 
$L=-\epsilon^2 \pal_x^2 + u$, i.e., as the
coefficients of the asymptotic expansion at $z\to 0$ of the diagonal of 
the heat kernel
\beq\label{2-5-1}
< x\, | \, e^{-z\,( -\epsilon^2\pal_x^2 + u(x))}\, |\, x> \sim {1\over \sqrt{4\, \pi\, z}} \sum_{k=0}^\infty
(-z)^k h_{k-2}, ~~h_{-2}=1.
\eeq
Note that the coefficients of the asymptotic expansion (\ref{2-5-1}) 
do not depend
on the choice of the boundary conditions for the Lax operator
\cite{seeley}.
They are related to the gradients of the old Hamiltonians by
\beq\label{h-kdv}
h_k:=\prod_{i=1}^{k+1}\left( i+{1\over 2}\right)^{-1}{\delta I_{k+1}\over \delta u(x)}, ~~k=-1, 0, 1, \dots
\eeq
We want to prove that these differential polynomials define a tau-structure
of the KdV hierarchy, i.e. the following symmetry holds true
\beq\label{2-5-3}
\{ h_{k-1}, \bar h_l\}_1 = \{h_{l-1}, \bar h_k\}_1.
\eeq
Observe that, redefining the flows of the hierarchy
as follows
\beq\label{tau-kdv}
{\pal u\over \pal \tilde t^k}:= \{ u(x), \bar h_k\}_1 = \pal_x {\delta \bar h_k\over \delta u(x)}, ~~
\bar h_k = \int h_k dx, ~~k=0, 1,
\dots
\eeq
we obtain the same hierarchy up to a change of normalization of the flows
\beq\label{normalize}
{\pal u\over \pal \tilde t^k}=\prod_{i=1}^{k}\left( i+{1\over 2}\right)^{-1}
{\pal u\over \pal t^k}.
\eeq
The symmetry (\ref{2-5-3}) means that the new densities are coefficients of 
a closed 1-form
\beq\label{kdv-1-form}
\omega =h_{-1} d\tilde t^0+ h_0 d\tilde t^1 + h_1 d\tilde t^2 + \dots.
\eeq
Explicitly the first few Hamiltonians and the corresponding flows of the
hierarchy are
$$
h_{-1}=u, ~~h_0 = {u^2\over 2} -\epsilon^2{u''\over 6}, ~~h_1={u^3\over 6} - 
{\epsilon^2\over 12}({u'}^2 + 2 u\, u'') + \epsilon^4 {u^{IV}\over 60},
$$
$$
h_2 = {u^4\over 24} - {\epsilon^2\over 12} (u\, {u'}^2 +
u^2u'')+{\epsilon^4\over 120} (3{u''}^2+4u'u'''+2u\, u^{IV}) - {\epsilon^6\over
840} u^{VI}
$$
$$
u_{\tilde t^0}=u', ~~u_{\tilde t^1}= u\, u' -{\epsilon^2\over 6} u''',
$$
$$
u_{\tilde t^2}={u^2 u'\over 2} -{\epsilon^2\over 6} (2 u' u'' + u\, u''') +
{\epsilon^4\over 60} u^{V}.
$$

To prove that the Seeley coefficients define a tau-structure of the KdV
hierarchy let us recall some identities for the generating
functional
$$
p=p(\lambda) = k + \sum_{m=1}^\infty {\bar \chi_m\over k^m}, ~~k=\sqrt{\lambda}
$$
(see, e.g., \cite{D75}). Let us introduce the ``real part''
$$
\chi_R := k + \sum_{j=0}^\infty {\chi_{2j+1}\over k^{2j+1}}.
$$
Then we have
$$
{\delta p\over \delta u(x)} = -{1\over 2 \chi_R}
$$
and
$$
{dp\over d\lambda} =\int {1\over 2 \chi_R}dx.
$$
From the last two formulae it immediately follows that
$$
\int {\delta \bar \chi_{2j+3}\over \delta u(x)} dx= \left( j+{1\over 2}\right)
\bar \chi_{2 j+1}.
$$
This proves that the new flows (\ref{tau-kdv}) 
are proportional to the equations of the KdV
hierarchy with the coefficient of proportionality given in (\ref{normalize}). 
From the same
formula it follows that the new densities of the Hamiltonians satisfy the
recursion similar to (\ref{2-2-6-b}) but with a different normalization
$$
\left( -{1\over 4} \epsilon^2 \pal_x^3 + u\pal_x +{1\over 2} u_x\right)
{\delta \bar h_{j-1}\over \delta u(x)} = \left(j+{1\over 2}\right) \pal_x
{\delta \bar h_{j}\over \delta u(x)}.
$$
All this means that the new Hamiltonians  
can be obtained from
the old ones by multiplying a nonzero constant and by adding of a total
derivative.

The generating function of the densities $h_k$ coincides with
$$
w(x,\lambda):= {\sqrt{\lambda}\over \chi_R(x,\lambda)}=1+{h_{-1}\over 2\lambda}
+\sum_{k=2}^\infty {(2k-1)!! h_{k-2}\over (2\lambda)^k}.
$$

We are now to prove the symmetry
\beq\label{symkdv}
{\pal h_{i-1}\over \pal \tilde t^j} = {\pal h_{j-1}\over \pal \tilde t^i} 
\eeq
for each pair $(i,j)$. We will use the following formula \cite{D75}
\beq\label{d75}
\left\{ \chi_R(x,\lambda), p(\mu)\right\}_1 ={1\over 8\,(\mu-\lambda)}
\left({\chi_R(x,\lambda)\over \chi_R(x,\mu)}\right) '
\eeq
(here $'=d/dx$, $\mu$ is an indeterminate). The both sides of the formula 
are understood as formal series
in inverse powers of $\sqrt{\mu}$. From this it easily follows that
$$
-\frac14\,\sum_{i,j=0}^\infty {(2i+1)!!\,(2j+1)!!\over (2\,\mu)^{i+1}\,
(2\,\lambda)^{j+1}}\,
{\pal h_{j-1}\over \pal \tilde t^i} ={1\over \mu-\lambda}
[w(x,\mu)w'(x,\lambda)-w'(x,\mu)w(x,\lambda)].
$$
The symmetry of the r.h.s. proves  (\ref{symkdv}).

Using the formula (28) from \cite{gelfand-dickey} we can represent the r.h.s. as the total derivative of the
function
\eqa 
&&{2\over \mu-\lambda}
[w(x,\mu)w'(x,\lambda)-w'(x,\mu)w(x,\lambda)] \nn\\
&&= \pal_x {1\over (\mu-\lambda)^2} \left[ w''(\mu) w(\lambda) +
w(\mu) w''(\lambda) - w'(\mu) w'(\lambda)\right.\nn\\
&&\quad \left. + 2 (\lambda+\mu - 2\,u) \, w(\mu) w(\lambda) -
2(\lambda+\mu)\right].
\eeqa

It would be interesting to derive the symmetry (\ref{symkdv}) directly from the properties of
the  heat kernel (cf. the recent paper \cite{mwz} where an analogous symmetry
of the derivatives of the Green function of the Dirichlet boundary value problem
on the plane has been derived from the Hadamard variational formula for the
Green function).

\begin{remark} In \cite{magri-et-al} it was suggested an approach to the theory
of tau-functions of the KP hierarchy and of its reductions. This approach is
based on the theory of exact Poisson pencils (see Example \ref{e-2-1-7} above).
Let 
$$
\{~,~\}_\lambda = \{~,~\}_2 -\lambda\, \{~,~\}_1
$$
be an exact Poisson pencil
$$
Lie_Z \{~,~\}_2=\{~,~\}_1, ~~ Lie_Z \{~,~\}_1=0
$$
for some vector field $Z$. Let us assume that the vector field $Z$ satisfies 
the conditions of the lemma \ref{l-1-4-9}. Doing if necessary a Miura-type
transformation we may assume that
$$
Z={\pal\over\pal u^1}.
$$
Let
$$
u^i\mapsto \tilde u^i = F^i(u^1, \dots, u^n; u^1_x, \dots, u^n_x,
\dots;\epsilon)
$$
be the reducing transformation for the Poisson bracket $\{~,~\}_2$. Then the
$\lambda$-dependent reducing transformation for the pencil described in
Theorem \ref{t-2-2-2} must have the form
\beq\label{2-5-10}
u^i\mapsto \tilde u^i = F^i(u^1-\lambda, \dots, u^n; u^1_x, \dots, u^n_x,
\dots;\epsilon)
\eeq
assuming that the r.h.s. is analytic at $\lambda=\infty$. If this is the case
then we can derive, following \cite{magri-et-al} that
\beq\label{2-5-11}
{\delta \bar f^i_p\over \delta u^1(x)} = p\, f^i_{p-1}+{\rm total ~
derivative}.
\eeq
Indeed, from
$$
{\pal\over \pal \lambda} \int  F^i(u^1-\lambda, 
\dots, u^n; u^1_x, \dots, u^n_x,
\dots;\epsilon)\, dx =-\int {\delta \bar F^i\over \delta u^1(x)} dx
$$
it follows that
$$
{\pal\over \pal \lambda} F^i(u^1-\lambda, \dots, u^n; u^1_x, \dots, u^n_x,
\dots;\epsilon) =-{\delta \over \delta u^1(x)}\int F^i\,  dx + 
{\rm total ~
derivative}.
$$
From (\ref{2-5-11}) it follows that the coefficients 
\beq
m^i_{p-1}:= {\delta\bar f^i_p\over\delta u^1(x)}
\eeq
of the expansion of $\delta F^i/\delta
u^1(x)$ are densities of conserved quantities of the bihamiltonian hierarchy.
They were used in \cite{magri-et-al} in order to define the tau-function. 

In the particular case of KdV one is to use, following the above prescription,
the differential polynomials
$$
m_{k-1} := {\delta I_k\over \delta u(x)},
$$
$$
\int m_{k-1}dx ={2k+1\over 2} I_{k-1}
$$
in order to define the tau-function of the KdV hierarchy
$$
m_{k-1} = \epsilon^2 {\pal^2\log\tau\over \pal x\, \pal t_k}.
$$
Instead of the symmetry (\ref{2-5-3}) one obtains
\beq
{2k+3\over 2} \{m_{k-1}, \bar m_l\}_1 = {2l+3\over 2} \{m_{l-1}, \bar m_k\}_1.
\eeq

The main problem with extending this approach to more general class of exact
Poisson pencils is that to describe the analytic properties for {\it large} 
$\lambda$ of the reducing transformation defined for {\it small} $\lambda$
by the formula (\ref{2-5-10}). 
Our axioms imply that this transformation always has a
regular singularity at $\lambda=\infty$ described by the formulae 
(\ref{period1})
-(\ref{period21}) below.
\end{remark}

\setcounter{equation}{0}
\setcounter{theorem}{0}

\subsection{From tau-structures to Frobenius manifolds}\label{sec-3-5}\par

Here we will construct the main invariant of an arbitrary tau-symmetric 
Poisson pencil of $(0,n)$ brackets on $\hat{\cal L(M)}$ with respect
to the action of Miura group. We will
prove that every such a bihamiltonian structure satisfying certain genericity
assumption  determines
a Frobenius structure on $M$. 

Actually our invariant will depend only on the leading order (\ref{pbht}) of the
bihamiltonian structure. We will prove
the following general result: classes of equivalence of 
tau-symmetric bihamiltonian structures of the form (\ref{pbht}) on ${\cal L}(M)$
satisfying certain genericity assumption are in one-to-one correspondence with 
Frobenius structures on $M$. 

All the calculations will be done in the system of flat coordinates $v^1, \dots,
v^n$ for the metric $g_1^{ij}$. Denote $\eta^{\alpha\beta}$ the (constant)
components of the metric $g_1^{ij}$ and $g^{\alpha\beta}(v)$ and
$\Gamma^{\alpha\beta}_\gamma(v)$ the components of the metric $g_2^{ij}$ and the
Levi-Civita connection for this metric. Recall (see Section \ref{sec-3-2}
 above) that,
in these coordinates $g^{\alpha\beta}(v)$ and
$\Gamma^{\alpha\beta}_\gamma(v)$ are expressed in terms of certain functions
$f_1(v)$, \dots, $f_n(v)$ as follows
\eqa
&&\Gamma^{\alpha\beta}_\gamma (v) = \pal^\alpha\pal_\gamma f^\beta(v),\nn\\
&&
g^{\alpha\beta}(v) =\pal^\alpha f^\beta(v) + \pal^\beta f^\alpha(v).
\label{2-6-1}
\eeqa
Here all raising and lowering of the Greek indices is to be done with the help
of the matrix $(\eta^{\alpha\beta})$ and the inverse one
$(\eta_{\alpha\beta})=(\eta^{\alpha\beta})^{-1}$. E.g.,
$$
\pal^\alpha := \eta^{\alpha\lambda}\pal_\lambda.
$$

We first observe that the densities of the commuting flows of the
hierarchy determined by the bihamiltonian structure (\ref{pbht}) are just
functions on $M$. They are determined by the following procedure.

\begin{lemma}\label{l-2-6-1}
 The bihamiltonian hierarchy determined by (\ref{pbht}) in the flat
coordinates $(v^1, \dots, v^n)$ for $\{~,~\}_1^{[0]}$ has the form 
$$
{\pal v^\beta\over \pal t^{\alpha,p}} = \pal_x \left( \pal^\beta
\phi_{\alpha,p}(v)\right)
$$
where the functions $\phi_{\alpha,p}(v)$ on $M$ are the coefficients of the
series solution 
$$
\phi_\alpha(v;\lambda) = 
\sum_{p=-1}^\infty {\phi_{\alpha,p}(v)\over \lambda^{p+1}}
$$
to the linear system 
\beq\label{2-6-2}
\left( g^{\beta\epsilon}(v)-\lambda\, \eta^{\beta\epsilon}\right)
\pal_\epsilon\pal_\gamma \phi + 
\Gamma^{\beta\epsilon}_\gamma(v) \pal_\epsilon \phi=0,
~~\beta, \, \gamma = 1, \dots, n
\eeq
such that
$$
\phi_{\alpha, -1}(v) = v_\alpha:= \eta_{\alpha\epsilon}v^\epsilon.
$$
\end{lemma}

\pf The system (\ref{2-6-2}) is just the spelling of the definition
$$
\{ ~.~, \bar \phi(v;\lambda)\}_\lambda^{[0]}=0
$$
of the Casimir of the pencil. \epf

Observe that the system (\ref{2-6-2}) 
coincides with the equations for the flat coordinates
for the flat pencil $g^{ij}_2-\lambda\, g^{ij}_1$ written in the flat coordinates
for $g^{ij}_1$.

\begin{cor} The flows of the lower level of the hierarchy have the form
$$
{\pal v^\beta\over \pal t^{\alpha, 0}} = \pal_x (\pal^\beta f_\alpha(v)).
$$
\end{cor}

\pf For the coefficient $\phi_{\alpha,0}(v)$ of $\lambda^0$ 
we obtain
$$
\pal^\beta\pal_\gamma \phi_{\alpha, 0} = \Gamma^{\beta\epsilon}_\gamma
\eta_{\epsilon\alpha} = \pal^\beta\pal_\gamma f_\alpha(v)
$$
where the functions $f_\alpha(v)$ were defined in (\ref{2-6-1}). So
$$
\phi_{\alpha,0}(v) = f_\alpha(v) + {{\rm linear}}.
$$
Adding of a function linear in $v$ does not affect the flat pencil neither
the Hamiltonian equations. \epf

By the assumption the densities $h_{\alpha,p}^{[0]}(v)$ are related to the
densities $\phi_{\alpha,p}(v)$ constructed in Lemma \ref{l-2-6-1} 
by an invertible 
triangular transformation with constant coefficients. In particular we may
assume that
$$
h_{\alpha,-1}^{[0]}(v) = v_\alpha
$$ 
and
\beq\label{h0}
h_{\alpha,0}^{[0]}(v) = a_\alpha^\beta f_{\beta}(v) + {{\rm linear}} ~
{{\rm function}}.
\eeq
Actually the linear function in the last formula can be absorbed by a
redefinition of $f_\beta(v)$. 

\begin{theorem}\label{t-2-6-3} 
Let (\ref{pbht}) be a tau-symmetric bihamiltonian structure.
Introduce the matrices 
$$
a =(a^\beta_\alpha), ~~b:=a^{-1}=(b^\beta_\alpha),
~~ b^{\alpha\beta} := \eta^{\alpha\lambda}b_\lambda^\beta.
$$
Then there exists a function $F(v)$ s.t.
\eqa
&&\Gamma_{\gamma}^{\alpha\beta}(v) = b^{\beta\epsilon} \pal^\alpha\pal_\gamma
\pal_\epsilon F(v),\label{2-6-3}\\
&&
g^{\alpha\beta} (v) = b^{\alpha\epsilon} \pal_\epsilon \pal^\beta F(v) +
b^{\beta\epsilon} \pal_\epsilon \pal^\alpha F(v) + g_0^{\alpha\beta}.
\label{2-6-4}
\eeqa
The function $F(v)$ satisfies associativity equations
\eqa\label{ass.eq}
&&\pal_\alpha\pal_\beta \pal_\lambda F(v) \eta^{\lambda\mu}
\pal_\mu\pal_\gamma\pal_\delta F(v) = \pal_\delta\pal_\beta \pal_\lambda F(v) 
\eta^{\lambda\mu}
\pal_\mu\pal_\gamma\pal_\alpha F(v),\\
&&\qquad \alpha, \beta, \gamma, \delta =1, \dots, n.
\eeqa
The function $F(v)$ is determined uniquely up to adding of at most quadratic
polynomial.
\end{theorem}

\pf The first of the symmetry conditions (\ref{sc}) reads
$$
{\pal v_\alpha\over \pal t^{\beta,0}} ={\pal v_\beta\over \pal t^{\alpha,0}}.
$$
Due to (\ref{h0}) it implies
$$
a_\alpha^\epsilon \pal_\beta f_\epsilon(v) = 
a_\beta^\epsilon \pal_\alpha f_\epsilon(v).
$$
Therefore there exists a function $F(v)$ s.t.
$$
a_\alpha^\epsilon f_\epsilon(v) =\pal_\alpha F(v).
$$
Hence
$$
f_\alpha(v) = b_\alpha^\lambda \pal_\lambda F(v).
$$
The freedom in adding of linear functions to $f_\alpha(v)$ corresponds to the
freedom in adding of a quadratic polynomial to $F(v)$.

It remains to prove that the algebra with the structure constants
$$
c_{\alpha\beta}^\gamma(v) : = \pal_\alpha\pal_\beta\pal^\gamma F(v)
$$
is associative for any $t\in M$. Indeed, from the quasi-Frobenius property
$$
\Gamma^{\alpha\beta}_\epsilon\Gamma^{\epsilon\gamma}_\delta
=\Gamma^{\alpha\gamma}_\epsilon\Gamma^{\epsilon\beta}_\delta
$$
for the Christoffel coefficients
$$
\Gamma^{\alpha\beta}_\gamma (v) = c_{\epsilon\gamma}^\alpha(v) b^{\beta\epsilon}
$$
we obtain
$$
b^{\beta\lambda} b^{\gamma\mu}\left[
c_{\epsilon\lambda}^\alpha(v) c_{\mu\delta}^\epsilon(v) -
c^\alpha_{\epsilon\mu}(v) c_{\lambda\delta}^\epsilon(v)\right] =0.
$$
This proves associativity. \epf

Due to the Theorem the formula
\beq\label{2-6-5}
\pal_\alpha\cdot \pal_\beta := c_{\alpha\beta}^\gamma(v) \pal_\gamma
\eeq
where
\beq\label{2-6-6}
c_{\alpha\beta}^\gamma(v) :=
\eta^{\gamma\epsilon}\pal_\epsilon\pal_\alpha\pal_\beta F(v)
\eeq
defines on the tangent planes $T_vM$ a structure of commutative associative
algebra.

We will now put into the game the condition of compatibility with spatial
translations.

\begin{lemma} If the tau-structure is compatible with spatial translations
then the function $F(v)$ constructed in Theorem \ref{t-2-6-3}
satisfies
$$
\pal_1\pal_\alpha\pal_\beta F(v)=\eta_{\alpha\beta}.
$$
\end{lemma}

\pf By definition we must have
$$
a_1^\epsilon \pal_\beta f_\epsilon(v) = v_\beta + {{\rm const}},
$$
i.e.,
\beq\label{2-6-7}
\pal_1\pal_\beta F(v) = v_\beta + v_\beta^0
\eeq
for a constant shift $v_\beta^0$.
Adding if necesary a quadratic polynomial to $F(v)$ we end the proof. \epf

According to the lemma the vector field
\beq\label{2-6-8}
e:= \pal /\pal v^1
\eeq
is the unity of the Frobenius algebra on $T_vM$ at every point $v\in M$. 
Actually we have a somewhat stronger condition 
\beq\label{constc}
{\pal\over \pal \tilde t^{1,0}} 
= c {\pal\over \pal t^{1,0}} ={\pal\over \pal x}, ~~
a_1^\alpha = c\,\delta_1^\alpha
\eeq
for some nonzero constant $c$. 

We are now to construct the Euler vector field on $M$. 
Let us introduce the linear function
$\varphi(v)= v_1 := \eta_{1\epsilon}v^\epsilon$. The gradient of this function
w.r.t. to the second metric we will denote $E$. 

\begin{lemma} The components $E^\alpha(v)$
of the vector field $E=E^\alpha(v)\pal_\alpha$
are the following linear functions on $M$
\beq\label{2-6-8b}
E^\alpha(v) = c \, v^\alpha + \bar b^{\alpha}_\beta
(v^\beta+v^\beta_0)+\eta_{1\epsilon} g_0^{\epsilon\alpha}
\eeq
where $\bar b^\alpha_\beta$ is the matrix of the operator adjoint to
$b^\alpha_\beta$,
$$
\bar b^\alpha_\beta := \eta^{\alpha\lambda}\eta_{\beta\mu} b^\mu_\lambda.
$$
\end{lemma}

\pf By definition
$$
E^\alpha=g^{\alpha\epsilon}(v)\pal_\epsilon \varphi(v) =\eta_{1\epsilon}
g^{\alpha\epsilon}(v)= \pal^\alpha f_1(v) + \pal_1 f^\alpha(v) +
\eta_{1\epsilon}g_0^{\epsilon\alpha}.
$$
Using $b_1^\alpha=c\, \delta_1^\alpha$ we obtain 
$$
\pal^\alpha f_1 = c \, v^\alpha.
$$
Representing
$$
f^\alpha(v) = \bar b^\alpha_\beta \pal^\beta F(v) 
$$
and using (\ref{2-6-7}) we end the proof. \epf

\begin{lemma}\label{l-2-6-6} 
Denote $(~,~)$ the bilinear form on $T_v^*M$
corresponding to the second 
metric. Then
\beq\label{intersectionform}
(\omega_1, \omega_2) =i_E \omega_1\cdot \omega_2.
\eeq
Here $\omega_1$, $\omega_2\in T^*_vM$ are two 1-forms on $M$. There product
is the 1-form defined via the product of tangent vectors on $M$
with the help of the isomorphism
$$
\eta: T_vM\to T_v^* M.
$$
\end{lemma}

\pf From the condition of symmetry of the Levi-Civita connection
$$
g^{\alpha\epsilon} \Gamma_\epsilon^{\beta\gamma}=
g^{\beta\epsilon} \Gamma_\epsilon^{\alpha\gamma}
$$
we derive that
$$
g^{\alpha\epsilon} c_\epsilon^{\beta\gamma}=
g^{\beta\epsilon} c_\epsilon^{\alpha\gamma}.
$$
Multiplying the last equation by $\pal _\beta\varphi(v)$ and taking the sum
w.r.t. $\beta$ we obtain
$$
E^\epsilon c_\epsilon^{\alpha\gamma} = g^{\alpha\gamma}.
$$
\epf

{\bf Definition.} We say that a solution $F(v)$ to the associativity equations
(\ref{ass.eq}) 
is {\it  rigid} if every constant invariant symmetric bilinear form
on the algebras $T_vM$ is proportional to $<\ , >$. A bihamiltonian structure
(\ref{mainclass1}) - (\ref{mainclass3})  possesing of a tau-structure
is called rigid if the corresponding solution to the associativity equations is.

Every cubic solution to the associativity equations is not rigid: taking an
arbitrary constant linear form $l$ we define a constant invariant symmetric
bilinear form by
$$
(a,b)_l:= l(a\cdot b), ~~a, ~b\in T_v M.
$$ 
Observe that the structure constants $c_{\alpha\beta}^\gamma(v)$ of such
a Frobenius manifold do not depend on $v$. Conversely, nonrigid solutions to
the associativity equations can be characterized by the following statement.

\begin{lemma}\label{l-2-6-7}
A solution $F(v)$ to the associativity equations is not rigid
if and only if there exists a constant (in the flat coordinates $v^\alpha$)
vector field $w$ non proportional to $e$ such that
\beq\label{delv}
\partial_w F(v) = {\rm at ~ most ~ quadratic} ~ {\rm polynomial ~ in} ~ v.
\eeq
\end{lemma}

\pf. Let $(~,~)$ be a constant invariant symmetric bilinear form on $T_vM$
non proportional to $<\ ,\ >$. We introduce a constant linear form
$$
l(a):= (a,e)
$$
and a constant vector $w=(w^\alpha)$ dual to $l$, i.e.
$l(a) =<w,a>$,
$$
w^\alpha=\eta^{\alpha\beta} l_\beta, ~~l_\beta:=l(\partial_\beta).
$$
Since
$$
(\partial_\alpha,\partial_\beta)= c_{\alpha\beta}^\gamma(v) l_\gamma,
$$
the sum
$$
w^\epsilon c_{\alpha\epsilon}^\gamma(v)
$$
does not depend on $v$. Differentiating the last equation along $v^\beta$ we
obtain
$$
\partial_w c_{\alpha\beta}^\gamma(v)=0.
$$
That means validity of (\ref{delv}).
Inverting the arguments, we also prove the converse statement. \epf

There is a simple way to generate new solutions to the equations of
associativity starting from a non rigid solution.

\begin{lemma} Let $F(v)$ be a solution to the equations of associativity and 
$w$ be the constant vector field satisfying (\ref{delv}). There
locally exists a function $F^w(v)$ and a constant antisymmetric linear operator
$\rho=(\rho^\alpha_\beta)$ such that
\beq\label{fv}
\nabla F^w(v)=w\cdot \nabla F(v) +{1\over 2} \rho(v) 
\eeq
where
$$
(\rho(v))^\alpha =\rho^\alpha_\beta\, v^\beta.
$$
The function $F^w(v)$ satisfies associativity equations (\ref{ass.eq}) with the
same $\eta^{\alpha\beta}$.
\end{lemma}

\pf Let us introduce the following 1-form
$$
\omega_\gamma:=q_{\gamma\nu}\nabla^\nu F(v) = w^\epsilon
c_{\gamma\epsilon}^\nu(v) \partial_\nu F(v).
$$
Here $q_{\gamma\nu}:=w^\epsilon c_{\gamma\epsilon\nu}(v)$ by assumption is a
constant symmetric matrix.  
Differentiating the
above equation twice along $v^\alpha$ and $v^\beta$ and using the associativity
we obtain
$$
\partial_\alpha\partial_\beta \omega_\gamma = w^\nu c_{\nu\alpha}^\epsilon(v)
c_{\epsilon\beta\gamma}(v).
$$
Due to the symmetry of the r.h.s. in $\beta$ and $\gamma$, there exists a
constant skew-symmetric matrix $\rho_{\beta\gamma}$ such that
$$
\partial_\beta \omega_\gamma - \partial_\gamma \omega_\beta =
\rho_{\beta\gamma}.
$$
The one-form
$$
\omega_\gamma-{1\over 2} \rho_{\gamma\nu}v^\nu
$$
is closed. Therefore it is locally equal to the differential of a function
we denote $F^w(v)$. 

To prove the last statement of Lemma we observe that
$$
\partial_\alpha\partial_\beta\partial_\lambda F^w(v) \eta^{\lambda\mu}
\partial_\mu\partial_\gamma \partial_\delta F^w(v) = 
\partial_\alpha\partial_\beta\partial_\lambda F(v) \eta^{\lambda\mu}_w
\partial_\mu\partial_\gamma \partial_\delta F(v)
$$
where
$$
\eta^{\lambda\mu}_w = q^\lambda_{\lambda'} q^\mu_{\mu'} \eta^{\lambda'\mu'}=
(w\cdot w)^\nu c_\nu^{\lambda\mu}(v)
$$
is a constant invariant symmetric bilinear form on $T_vM$. This proves validity
of the associativity equations for the new function $F^w(v)$.
\epf

\begin{theorem} For a  solution $F(v)$ to the associativity equations 
corresponding to a tau-symmetric bihamiltonian structure compatible with
spatial translations the constant vector $w=(w^\alpha)$
\beq\label{vectorw}
w^\alpha = 2c\, \delta^\alpha_1 + \bar b^\alpha_1
\eeq
satisfies (\ref{delv}). The function $F(v)$ satisfies
\beq\label{delE}
\partial_E F(v) = 2F^w(v) +{\rm quadratic}
\eeq
The Lie bracket of the vector fields $E$, $e$ equals
\beq\label{liebracket}
[e,E]= w-c\, e.
\eeq 
\end{theorem}

\pf From the equation (\ref{2-6-1}) we obtain
$$
E^\epsilon\pal_\epsilon\pal_\alpha\pal_\beta F = b_\alpha^\epsilon \pal_\epsilon
\pal_\beta F + b_\beta^\epsilon \pal_\epsilon\pal_\alpha F + {g_0}_{\alpha\beta}.
$$
We can rewrite it as
\beq\label{2ndder}
\pal_\alpha\pal_\beta [E^\epsilon \pal _\epsilon F] = Q_\alpha^\epsilon 
\pal_\epsilon\pal_\beta
F + Q_\beta^\epsilon \pal_\epsilon\pal_\alpha F + {g_0}_{\alpha\beta}
\eeq
where
$$
Q_\alpha^\beta = c\, \delta_\alpha^\beta + b_\alpha^\beta + \bar b_\alpha^\beta,
$$
the constant $c$ was defined in (\ref{constc}).
The compatibility conditions
$$
\pal_\gamma [ Q_\alpha^\epsilon \pal_\epsilon\pal_\beta
F + Q_\beta^\epsilon \pal_\epsilon\pal_\alpha F] = \pal_\beta [ 
Q_\alpha^\epsilon \pal_\epsilon\pal_\gamma
F + Q_\gamma^\epsilon \pal_\epsilon\pal_\alpha F]
$$
implies
$$
Q_\beta^\epsilon c_{\epsilon\alpha\gamma}(v) = Q_\gamma^\epsilon 
c_{\epsilon\alpha\beta}(v).
$$
This means invariance of the constant symmetric bilinear form
$$
<Q\, x, y> = c <x,y> + <b\, x, y> + <x, b\, y>.
$$
Therefore the matrix $Q_\beta^\alpha$ can be represented 
in the form
$$
Q^\alpha_\beta= w^\epsilon c_{\epsilon\beta}^\alpha(v), ~~w^\alpha =Q^\alpha_1.
$$
This gives the formula (\ref{vectorw}) since $b^\alpha_1=c\,\delta^\alpha_1$. 
Using Lemma \ref{l-2-6-7} we arrive at the proof
of (\ref{delv}).

To prove (\ref{delE}) we differentiate (\ref{2ndder}) along $v^\gamma$ and use
associativity to obtain
$$
\pal_\alpha\pal_\beta\pal_\gamma (\pal_E F(v)) = 2 Q_\gamma^\epsilon
c_{\epsilon\alpha\beta}(v) = \pal_\alpha\pal_\beta [2 (w\cdot \nabla
F(v))_\gamma]= 2 \partial_\alpha\pal_\beta \pal_\gamma F^w(v).
$$ 
The Theorem is
proved.\epf

Recall \cite{D3} that a solution $F(v)$ to the associativity equations 
(\ref{ass.eq})
possesing a constant unity vector field $e$ defines on $M\ni v$ a structure of
Frobenius manifold if a linear vector field $E$ exists on $M$ s.t.
\beq\label{fm}
\pal_E F(v) =(3-d) F(v) +{\rm quadratic}, ~~~[e,E]=e.
\eeq
Here $d$ is a constant called {\it the charge} of the Frobenius manifold.
The linear vector field is called {\it Euler vector field} of the Frobenius
manifold.
If a linear vector field $E$ satisfies 
\beq\label{degfm}
\pal_E F(v) =k\, F(v) +{\rm quadratic}, ~~~[e,E]=0
\eeq 
for some constant $k$ then we will call $M$ {\it degenerate Frobenius manifold}.
The theory of degenerate Frobenius manifolds is simpler than the theory of
Frobenius manifolds due to presence of a commutative group of algebra
automorphisms generated by $e$, $E$. In the particular semisimple case a
degenerate Frobenius manifold can be described in terms of Prym theta-functions
of plane algebraic curves of the degree $n=\dim M$ with an involution (see
Appendix below).

\begin{lemma}\label{l-2-6-10} 
If the constant vector field (\ref{vectorw}) is proportional to
$e$,
$$
w=(c+\kappa)e
$$
for some constant $\kappa$ then, for $\kappa\neq 0 $ the function $\kappa\,
F(v)$, defines on $M$ a Frobenius structure with the Euler vector field
$\kappa^{-1}E$, the unity $e$ and the charge
\beq\label{charge}
d=1-{2\,c\over \kappa}.
\eeq
If $\kappa=0$ then $(F(v), E)$ define on $M$ a structure of degenerate Frobenius
manifold.
\end{lemma}

\pf For $w=(c+\kappa)e$
$$
F^w(v) =(c+\kappa) F(v) +{\rm quadratic}.
$$
So (\ref{delE}) gives
$$
\pal_EF(v) =2(c+\kappa) F(v) +{\rm quadratic}, ~~[e,E]=\kappa\,e.
$$
For $\kappa=0$ the last equation coincides with (\ref{degfm}) with $k=2c$.
For $\kappa\neq 0$ the renormalization $E\mapsto \kappa^{-1} E$, $F\mapsto
\kappa\, F$ gives (\ref{fm}). Observe that, due to the renormalization of
 $F(v)$
the formula (\ref{intersectionform}) for the second metric remains unchanged.
\epf

\begin{cor} Every rigid tau-symmetric bihamiltonian structure on ${\cal L}(M)$ of the form 
(\ref{pbht}) compatible with
spatial translations in some coordinates can be reduced to one of the following
two normal forms.

1)
\eqa
&&\{ v^\alpha(x), v^\beta(y)\}_1 =\eta^{\alpha\beta}\delta'(x-y),\nn\\
&&
\{ v^\alpha(x), v^\beta(y)\}_2 = g^{\alpha\beta}
(v(x))\delta'(x-y)+\Gamma^{\alpha\beta}_\gamma(v(x)) v^\gamma_x \delta(x-y)
\label{normalfm}
\eeqa
where 
$$
g^{\alpha\beta} (v) =E^\epsilon(v) c_\epsilon^{\alpha\beta}(v)
$$
is the intersection form of a Frobenius structure on $M$
\eqa
&&\Gamma^{\alpha\beta}_\gamma(v) 
=c^{\alpha\epsilon}_\gamma (v) \left({1\over 2}
-\mu\right)^\beta_\epsilon\nn\\
&&
\mu:= {2-d\over 2} -\nabla E.\nn
\eeqa

2)
\eqa
&&\{ v^\alpha(x), v^\beta(y)\}_1=\eta^{\alpha\beta}\delta'(x-y),\nn\\
&&\{ v^\alpha(x), v^\beta(y)\}_2 =  g^{\alpha\beta}
(v(x))\delta'(x-y)+\Gamma^{\alpha\beta}_\gamma(v(x)) v^\gamma_x \delta(x-y)
\label{normaldegfm}
\eeqa
where
$$
g^{\alpha\beta} (v) =E^\epsilon(v) c_\epsilon^{\alpha\beta}(v)
$$
is the intersection form of a degenerate Frobenius structure on $M$
$$
\Gamma^{\alpha\beta}_\gamma(v) =- c^{\alpha\epsilon}_\gamma (v) 
\mu^\beta_\epsilon
$$
$$
\mu:={k\over 2}-\nabla E.
$$
\end{cor}

We will now show that, under the assumption of semisimplicity the bihamiltonian
structure can be reduced to the direct sum of the structures (\ref{normalfm}),
(\ref{normaldegfm}) even without assumption of rigidity.

{\bf Definition.} The bihamiltonian structure (\ref{mainclass1}) -
(\ref{mainclass3}) is called {\it semisimple} if the characteristic equation
$$
\det (g_2^{ij}(u) -\lambda\, g_1^{ij}(u))=0
$$
has pairwise distinct roots for any $u\in M$.

We will see below that semisimplicity guarantees integrability of the
bihamiltonian hierarchy (\ref{mainclass4}).

\begin{lemma} For a semisimple tau-symmetric bihamiltonian structure compatible
with spatial translations the Frobenius algebra constructed in 
Theorem \ref{t-2-6-3} is
semisimple.
\end{lemma}

\pf According to Lemma \ref{l-2-6-6} the linear operator
$$
(g_2^{ij})\cdot (g_1^{ij})^{-1}
$$
coincides with the operator of multiplication by the vector field $E$. Therefore
the Frobenius algebra on $T_vM$ contains at least one element with simple
spectrum. Hence all the elements of the algebra are semisimple.
\epf

\begin{lemma} Every semisimple solution to the associativity equations
possesing of a constant unity can be decomposed into a direct sum of
rigid solutions.
\end{lemma}

\pf Let us consider the subspace $W$ of all constant vector fields $w$ 
satisfying (\ref{delv}). It contains the unity vector field $e$. Observe
that the operator of multiplication of vectors of $T_vM$
by any vector $w\in W$ has constant matrix in the basis of flat coordinates
on $M$. Therefore $W$ is a finite dimensional subalgebra in the Frobenius
algebra of vector fields on $M$. 

Due to semisimplicity the commuting operators of multiplication by the vectors
$w$ from $W$
can be simultaneously reduced to the diagonal form. Let $\lambda_1(w)$, \dots,
$\lambda_m(w)$ be the pairwise distinct eigenvalues of these operators.
We consider these eigenvalues as elements of the dual space $W^*$. We define
$$
TM= TM_1\oplus \dots \oplus TM_m
$$
as the decomposition of the tangent bundle into the direct sum of the
corresponding
eigensubspaces. It is easy to see that this is an orthogonal decomposition
w.r.t. the bilinear form $<\ ,\ >$ and that the products of vectors from $TM_i$
and $TM_j$ are all zeroes for $i\neq j$. In this way we obtain the submanifolds
$M_1$, \dots, $M_m$. The solution $F(v)$ decomposes into the sum
$$
F(v) = F_1(v_1)+\dots +F_m(v_m), ~~v_s\in M_s
$$
up to adding of at most quadratic  polynomial.

Denote $e_s$, $E_s$, $w_s$ the projections of the vectors $e$, $E$, $ w$ resp.
onto the $s$-th factor. The vector $e_s$ is the unity of the Frobenius algebra
on $M_s$. The
eigenvalues of the operator of multiplication by $ w_s-\lambda_s(w) e_s$
restricted onto $TM_s$ are zeroes. Due to semisimplicity this implies that
$$
w_s = \lambda_s(w) e_s.
$$
Therefore $M_s$ is rigid.
\epf

\begin{cor} Let $F(v)$, $v\in M$ be a semisimple solution to the associativity
equations possesing a constant unity vector field $e$ and a linear vector field $E$
satisfying (\ref{delE}) where the constant vector field $w$ satisfies
(\ref{delv}).  Then $M$ is isomorphic to the direct product 
\beq\label{decompose}
M= M_1 \times M_2 \times \dots \times M_m
\eeq
of Frobenius manifolds or degenerate Frobenius manifolds. The vector field 
$w$
must have the form
$$
w=c\,e + \sum_{s=1}^m \kappa_s e_s,
$$
where $e=e_1+\dots + e_m$ is the decomposition of the unity vector fields
w.r.t. the product structure (\ref{decompose}), and $\kappa_1$, \dots,
$\kappa_m$ are some constants. The factors with $\kappa_s=0$ correspond to
degenerate Frobenius manifolds and the factors with $\kappa_s\neq 0$
correspond to the Frobenius manifold $M_s$ with 
$$
d_s=1-{2\,c\over \kappa_s}
$$
and with the Euler vector field
$$
E_s := \kappa_s^{-1} {\rm pr}_s E.
$$
\end{cor}

\pf Decomposing the constant vector field $w$ given in (\ref{vectorw}) w.r.t.
the decomposition (\ref{decompose}) we obtain, due to rigidity of the factors
$$
w=c+\sum_{i=1}^m \kappa_s e_s.
$$
where $\kappa_1$, \dots, $\kappa_m$ are some constants. Decomposing also
(\ref{delE}), (\ref{liebracket}) we obtain 
$$
\pal_{E_s} F_s= 2(c+\kappa_s) F_s+{\rm quadratic}, ~~[e_s, E_s]=\kappa_s e_s.
$$
The factors with $\kappa_s=0$ give degenerate Frobenius manifolds; those
with $\kappa_s\neq 0$ after a renormalization like in Lemma \ref{l-2-6-10} 
are Frobenius
manifolds.
\epf

It remains to only investigate the freedom in the choice 
of the tau-structure. We will do it for the rigid 
Poisson pencil (\ref{pbht}).

\begin{lemma} For a rigid tau-symmetric Poisson pencil of the form 
(\ref{pbht}) the tau-structure is determined uniquely up to simultaneous
linear transformations
\beq\label{2-6-26}
h_{\alpha,p}(v)\mapsto A_\alpha^\beta h_{\beta,p}, ~~p=-1, 0, 1, \dots
\eeq
and transformations of the form
\beq\label{2-6-27}
h_{\alpha,p}(v)\mapsto \rho^{p+1} h_{\alpha,p}(v)+\sum_{q=0}^p
\rho^{p-q}B_{\alpha,q}^\gamma h_{\gamma, p-q-1}, ~~p=0, 1, \dots 
\eeq
with some nonzero $\rho$ and a collection of $n\times n$ constant matrices
$B_{\alpha,p}^\beta$, $p=0, 1, \dots$.
\end{lemma}

\pf Indeed, let $\tilde h_{\alpha,p}(v):=\tilde h_{\alpha,p}^{[0]}$ 
be another tau-structure.
Let us first assume that
$$
\tilde h_{\alpha,-1}=h_{\alpha,-1}=v_\alpha.
$$
Let
$$
\tilde h_{\alpha,0}=A_\alpha^\beta h_{\beta,0}+ B_\alpha^\beta v_\beta.
$$
From
$$
{\pal v_\alpha\over \pal t^{\beta,0}}= c_{\alpha\beta\gamma}(v) v^\gamma_x
$$
using the symmetry
$$
{\pal v_\alpha\over \pal \tilde t^{\beta,0}} = 
{\pal v_\beta\over \pal \tilde t^{\alpha,0}}
$$
we derive
$$
A_\alpha^\lambda c_{\lambda\beta\gamma}(v) = A_\beta^\lambda
c_{\lambda\alpha\gamma}(v).
$$
For $\gamma=1$ the last equation implies that the constant bilinear form $<A\, x, y>$
on $T_vM$ is symmetric. From the whole equation it follows that this form is
invariant. Due to rigidity it must be proportional to $<\ ,\ >$. Hence
$A_\alpha^\beta = \rho\, \delta _\alpha^\beta$ for some nonzero constant $\rho$.

Let us redenote $B_{\alpha,0}^\beta:= B_\alpha^\beta$. From
the next
symmetry condition
$$
{\pal v_\alpha\over \pal \tilde t^{\beta,1}}=
{\pal \tilde h_{\beta,0}\over \pal
\tilde t^{\alpha,0}}
$$
we derive that
$$
{\pal\over \pal \tilde t^{\alpha,1}}= \rho^2 
{\pal\over \pal  t^{\alpha,1}}+ \rho\, B_{\alpha,0}^\gamma 
{\pal\over \pal t^{\gamma,0}}.
$$
Hence
$$
\tilde h_{\alpha,1}(v) =\rho^2 h_{\alpha,1}(v) + \rho\, B_{\alpha,0}^\gamma
h_{\gamma,0}(v) + B_{\alpha,1}^\gamma v_\gamma
$$
where $B_{\alpha,1}^\gamma$ is a new constant matrix. The proof of the lemma
can be finished using induction.
\epf

Observe that, conversely, an arbitrary transformation of the form (\ref{2-6-26}),
(\ref{2-6-27}) map
a tau-structure to another one.
 
We postpone for a subsequent publication the study of integrable hierarchies
corresponding to degenerate Frobenius manifolds. In the remaining part of this
paper we will assume that the $(0,n)$ Poisson pencil (\ref{pbht}) defines on $M$
a structure of a semisimple Frobenius manifold. Actually, we will see that not
an arbitrary Frobenius manifold can be obtained starting from a tau-symmetric 
Poisson pencil (\ref{pbht}). The restriction is that,
 the spectrum of the Frobenius manifold must contain no half-integers. However,
 an arbitrary semisimple Frobenius manifold generates an integrable hierarchy
according to the construction of the next section. If the spectrum contains no
half-integers then the hierarchy admits a tau-symmetric bihamiltonian structure.
Integrable hierarchies corresponding to the Frobenius manifolds with
half-integers in the spectrum can be considered as the closure of our
construction. For example, these hierarchies are conjecturally in the theory
of Gromov - Witten invariants of smooth projective varieties of odd complex
dimension (see below).

\setcounter{equation}{0}
\setcounter{theorem}{0}

\subsection{From Frobenius manifold to the Principal Hierarchy}
\label{sec-3-6}

Let $M$ be a $n$-dimensional Frobenius manifold.
In this section we construct the integrable hierarchy on ${\cal L}(M)$ of the
first order quasilinear systems
\beq
{\pal v^i\over \pal t} = A^i_j(v) {\pal v^j\over \pal x}, ~~i=1, \dots, n.
\eeq
Under certain assumptions about the eigenvalues of the gradient $\nabla E$
of the Euler vector field $E$ we will show that the hierarchies are generated
by the bihamiltonian structure of the form (\ref{normalfm}). 
We will also construct
the conservation laws for the hierarchy and, for an appropriate class of its 
solutions we compute the tau-function. Finally, for the case of semisimple
Frobenius manifold we will prove completeness of our system of conservation
laws.

\subsubsection{Commuting bihamiltonian flows on the loop space of a Frobenius manifold}

Let us concentrate first at the hamiltonian systems w.r.t. the Poisson bracket
$\{~,~\}_1$. Recall that, in the flat coordinates $v^1$, \dots, $v^n$ 
for the metric $<\ ,\ >$ on $M$ the Poisson bracket has the form
$$
\{ v^\alpha(x), v^\beta(y)\}_1 = \eta^{\alpha\beta} \delta'(x-y).
$$
The flow with the hamiltonian $\bar f =\int f(v)\,dx$ reads
$$
v_t = \{ v(x), \bar f \}_1 = \pal_x \nabla f.
$$
This is a first order quasilinear PDE with the matrix of coefficients
$$
A^\alpha_\beta (v) =\nabla^\alpha\nabla_\beta f(v).
$$

Denote $A(M)$ the space of smooth functions $f(v)$ on $M$ satisfying
\beq\label{star}
\nabla_a\nabla_b\,f =\nabla_{a\cdot b} Lie_e f
\eeq
for arbitrary two vector fields $a$, $b$ on $M$. Here $e$ is the unity vector
field on $M$.

\begin{theorem}\label{t-2-7-1} 
1) For arbitrary two functions $f$, $g\in A(M)$ the hamiltonian 
flows
\beq\label{2fields}
v_t = \{ v(x), \bar f\}_1, ~~v_s =\{ v(x), \bar g\}_1
\eeq
commute.
2) For any $f\in A(M)$ there exists $g\in A(M)$ such that
$$
\{~.~, \bar f\}_1 = \{~.~, \bar g\}_2.
$$
\end{theorem}

Here the second Poisson bracket has the form (\ref{normalfm}).

\pf  The commutator $(v_t)_s - (v_s)_t$ of the vector fields (\ref{2fields})
reads
$$
(v^\alpha_t)_s - (v^\alpha_s)_t=\pal_x \left[ \pal_\beta \nabla^\alpha f \pal_x
\nabla^\beta g - \pal_\beta \nabla^\alpha g \pal_x
\nabla^\beta f\right].
$$
Using (\ref{star}) we rewrite the expression in the brackets as follows
\eqa
&&\pal_\beta \nabla^\alpha f \pal_x
\nabla^\beta g - \pal_\beta \nabla^\alpha g \pal_x
\nabla^\beta f\nn\\
&&= \left[c^\alpha_{\beta\lambda}c^\beta_{\gamma\mu} (Lie_e \nabla^\lambda
f Lie_e \nabla^\mu g - Lie_e \nabla^\lambda
g Lie_e \nabla^\mu f\right] v^\gamma_x =0\nn
\eeqa
due to associativity.

To prove the bihamiltonian property of the flows with the densities of
hamiltonians in $A(M)$ it suffices, due to Lemma \ref{H1H2-der}, 
to prove that these flows
are symmetries of the second Poisson bracket. We leave this calculation as an
exercise for the reader.
\epf

\begin{remark} To find the hamiltonian $\bar g$ for the hamiltonian flow
$$
v_t=\{ v(x), \bar f\}_1= \{ v(x), \bar g\}_2,
$$ 
$f(v)\in A(M)$, w.r.t. the second Poisson bracket (\ref{normalfm})
one is to find a solution $g(v)$ 
to (\ref{star})
satisfying also
\beq\label{2nd-ham}
Lie_E \nabla g + {3-d\over 2} \nabla g = \nabla Lie_e f.
\eeq
\end{remark}

Our nearest goal is to construct a basis in a suitable subspace of $A(M)$.
Such a basis is just the hierarchy we promised to construct.

\subsubsection{Spectrum of Frobenius manifold, Levelt basis of deformed flat coordinates, and 
Hamiltonians of the Principal Hierarchy}\label{sec-3-6-2}\par

Denote $\tilde \nabla$ the deformed flat connection on $M\times {\mathbb C}^*$
\eqa\label{def-con}
&&
\tilde \nabla_u\, v = \nabla_u v + z\, u\cdot v, ~~u, v\in TM, ~~z\in{\bf
C}^*\nn\\
&&
\tilde \nabla _{d\over dz} v = \pal_z v +E\cdot v -{1\over z} {\cal V} v
\nn
\eeqa
where
\beq\label{2-7-4b}
{\cal V} := {2-d\over 2} -\nabla E
\eeq
is an antisymmetric operator w.r.t. $<\ ,\ >$. The unity vector field $e$ is an
eigenvector of this operator with the eigenvalue
$$
{\cal V} e = -{d\over 2} e.
$$

\begin{lemma} Let $f=f(v,z)$ be a {\rm horizontal function} for the connection
(\ref{def-con}), i.e., 
\beq\label{hor}
\tilde\nabla\, df=0.
\eeq
Then 
\beq\label{hor1}
f(v,z)\in A(M) ~~{\rm for ~any}~z.
\eeq
\end{lemma}

\pf Spelling the first half of the horizontality equation $\tilde \nabla_\alpha
df=0$ one obtains
$$
\pal_\alpha\pal_\beta f = z \, c_{\alpha\beta}^\gamma(v) \pal_\gamma f.
$$
In particular,
$$
\pal_1\pal_\gamma f = z \pal_\gamma f.
$$
These two equations imply (\ref{star}). \epf

To construct a basis in $A(M)$ we will use the coefficients of expansion at
$z=0$ of the {\it deformed flat coordinates}. By definition these are $n$ 
independent horizontal
functions $\tilde v_1(v;z)$, \dots, $\tilde v_n(v;z)$, that is,
$$
\det \left( {\partial \tilde v_\alpha(v;z)\over \partial v^\beta}\right) \neq 0
$$
 and the functions satisfy
\beq\label{horizont}
\tilde \nabla d\tilde v_\alpha=0, ~~\alpha = 1, \dots, n.
\eeq
We now describe a particular system of solutions to (\ref{horizont}).

First, using the Levi-Civita connection for the flat metric $<\ ,\ >$ 
we obtain a
natural trivialization of the tangent bundle
$$
TM\simeq M\times V
$$
where $V$ is a $n$-dimensional complex space with a symmetric nondegenerate
bilinear form that we denote by the same symbol $<\ ,\ >$. For an arbitrary linear
operator $A:V\to V$ denote $A^*: V\to V$ the adjoint operator
$$
<A^* x, y> =<x, Ay> ~~~{\rm for} ~ {\rm any} ~~x, y\in V.
$$
The horizontal $(1,1)$-tensor ${\cal V}$ becomes a linear operator on $V$ satisfying
$$
{\cal V}^*=-{\cal V}.
$$
The horizontal unity vector field $e$ gives a distinguished eigenvector of ${\cal V}$
in $V$ that we also denote $e$. After such a trivialization the deformed flat
connection becomes a flat connection on the trivial bundle $M\times {\mathbb C}^*
\times V$. The equation (\ref{horizont}) for the gradients
$$
Y(v;z):=\nabla\tilde v(v;z)
$$
of deformed flat coordinates
takes the form of a system
\beq\label{lin-system1}
\pal_\alpha Y =z C_\alpha(v) Y, ~~\alpha=1, \dots, n
\eeq
\beq\label{lin-system2}
\pal_z Y =\left( {\cal U}(v) +{{\cal V}\over z}\right)Y.
\eeq
Here  $C_\alpha(v)$ and ${\cal U}(v)$ are  the linear operators in $V$ of multiplication by
$\pal/\pal v^\alpha$ and $E$ resp. These operators are symmetric,
$$
C_\alpha^* = C_\alpha, ~~{\cal U}^* = {\cal U}.
$$
To fix a system of the deformed flat coordinates we are to choose a basis in the space
of solutions to the system (\ref{lin-system1}), (\ref{lin-system2}). 
Such a basis corresponds 
to a choice of a representative in the equivalence class of {\it normal forms} of the system 
(\ref{lin-system2}) near $z=0$ (see details in \cite{D4}). 
The parameters of such a normal form are called {\it monodromy at the origin} 
of the Frobenius manifold. Let us first recall the description of the 
parameters.

{\bf Definition.} The {\it spectrum} of a Frobenius manifold is a quadruple
$\left( V,<~,~>, \hat\mu, R\right)$ where
$V$ is a $n$-dimensional linear space over ${\mathbb C}$ equipped 
with a symmetric non-degenerate bilinear form $<\ ,\ >$, 
semisimple antisymmetric 
linear operator $\hat \mu: V\to V$, $<\hat\mu\, a, b>+<a, \hat\mu\,b>=0$
and a nilpotent linear operator $R: V\to V$ satisfying the following properties.
First, 
\beq\label{R2}
R^* =-e^{\pi\,i\,\hat\mu} R e^{-\pi\,i\,\hat\mu}
\eeq
Observe the following consequence of (\ref{R2})
\beq\label{levelt51}
R e^{2\pi\,i\,\hat\mu} = e^{2\pi\,i\,\hat\mu} R.
\eeq 
The operator $R$ must also be {\it $\hat\mu$-nilpotent}, i.e., it must preserve
a natural flag in $V$ associated to the operator $\hat\mu$. The flag is
constructed as follows.

For a given $\rho\in 
Spec \,e^{2\pi\, i\,\hat\mu}$ denote $\mu_{\rm max}=\mu_{\rm max}(\rho)\in 
Spec \,\hat\mu$ the eigenvalue with the
maximal real part satisfying $e^{2\pi\,i\,\mu_{\rm max}}=\rho$. 
For every nonnegative integer $m$ define the subspace 
\beq\label{flag3}
V^m := 
\oplus_{0\leq k \leq m} \oplus_{\rho\in Spec \,e^{2\pi\,i\,\hat\mu} }
{\rm Ker} \left[\hat \mu-(\mu_{\rm
max}(\rho)-k)\cdot 1\right] \subset V.
\eeq
Clearly
\beq\label{flag2}
0=V^{-1}\subset V^0\subset V^1\subset
\dots V
\eeq
For sufficiently large $m$ one has $V^m = V$. 
Let $0=k_1\leq k_2< k_3<\dots< k_l$
be all the integers such that
$$
V^{k_i-1}\neq V^{k_i}.
$$
Denoting
\beq\label{flag4}
F_i := V^{k_i}, i=1, \dots, l
\eeq
we obtain a flag
\beq\label{flag5}
0=F_0\subset F_1 \subset \dots \subset F_l=V.
\eeq
{\bf Definition.} The flag (\ref{flag5}) is called {\it Levelt flag} associated
with ($V$, $<\ ,\ >$, $\hat\mu$). The operator $R: V\to V$ is called
{\it $\hat\mu$-nilpotent} if the Levelt flag is invariant
\beq \label{flag6}
R(F_j)\subset F_j, ~~j=0, 1, \dots, l.
\eeq

By the construction the operator $R$ satisfies
\beq\label{R3}
z^{\hat\mu} R z^{-\hat\mu} = R_0 + R_1 z + R_2 z^2+\dots
\eeq
where the coefficients of the matrix valued polynomial are nilpotent operators
$R_0$, $R_1$, \dots such that
\beq\label{R4}
R=R_0+ R_1 + \dots
\eeq
and
\beq\label{flag7}
R_k(V^j) \subset V^{j-k}.
\eeq

Observe the following useful identity
\beq\label{R5}
z^{\hat\mu} R_k z^{-\hat\mu} = z^k R_k, ~~k=0, 1, \dots 
\eeq
\beq\label{R6}
[\hat\mu, R_k] = k\, R_k, ~~k=0, 1, \dots .
\eeq
The spelling of the equation (\ref{R2}) for the coefficients $R_k$ reads
\beq\label{R7}
R_k^* =(-1)^{k+1} R_k, ~~k=0, 1, \dots .
\eeq
Any polynomial of the matrices $R_k$ can be uniquely decomposed as
follows
\beq\label{R8}
P(R_0, R_1, \dots)= [P(R_0, R_1, \dots)]_0+[P(R_0, R_1, \dots)]_1+\dots
\eeq
\beq\label{R8-b}
z^{\hat\mu} [ P(R_0, R_1, \dots)]_k z^{-\hat\mu} = z^k [P(R_0, R_1,
\dots)]_k.
\eeq 
The last ingredient of the monodromy at the origin is an eigenvector $e\in V$
of $\hat\mu$ satisfying $R_0\,e=0$. It will be needed later on.

We will now explain how to associate a 5-tuple $(V, <\ ,\ >, \hat\mu, R, e)$ 
to a Frobenius manifold. The linear space $V$ with a symmetric nondegenerate bilinear form $<\ ,\ >$ 
and a vector $e\in V$ have already been constructed above.
Denote $\hat\mu:V\to V$ the semisimple part of the operator ${\cal V}$, i.e., 
\beq\label{mu}
\hat\mu:= \oplus_{\mu\in Spec\, {\cal V}} \mu\, P_\mu
\eeq
where $P_\mu:V\to V_\mu$ is the projector of $V$ onto the root subspace of 
${\cal V}$
$$
V=\oplus_{\mu\in Spec \, {\cal V}} V_\mu, ~~~V_\mu := {\rm Ker}\, ({\cal V}-\mu\cdot 1)^n,
$$
$P_\mu (V_{\mu'})=0$ for $\mu\neq \mu'$, $P_\mu |_{V_\mu}={\rm id}_{V_\mu}$.
Clearly the operator $\hat\mu$ is antisymmetric, $\hat\mu^*=-\hat\mu$. Denote $R_0$ the nilpotent part of ${\cal V}$
$$
{\cal V}=\hat\mu+R_0.
$$ 
Other operators $R_1$, $R_2$, \dots are not determined by ${\cal V}$ only. They appear only
in
presence of resonances, i.e., pairs of eigenvalues $\mu$, $\mu'$ of ${\cal V}$ 
such that $\mu-\mu' \in {\bf Z}_{>0}$.

Let us choose a basis $e_1, \dots, e_n$ in $V$ such that $e_1=e$. The matrices
of the linear operators $\hat\mu$ and $R$ we will denote by the same symbols.

\begin{theorem} For a sufficiently small ball $B\in M$ there exists a
fundamental matrix of solutions to the system (\ref{lin-system1}), 
(\ref{lin-system2}) of the form
\beq\label{fundamental}
Y(v;z) = \Theta(v;z)z^{\hat\mu}z^R
\eeq
such that
the matrix valued function $\Theta(v;z):V\to V$ is analytic on $B\times {\mathbb C}$ satisfying
\beq\label{normalize-theta}
\Theta(v;0)\equiv 1
\eeq
\beq\label{orthogonal}
\Theta^*(v;-z)\Theta(v;z)\equiv 1.
\eeq
\end{theorem}

This was proved in \cite{D4} for the case of diagonalizable ${\cal V}$. The general
case can be settled in a similar way. Note that a branch of logarithm $\log z$
is to be fixed in order to define the matrices $z^{\hat\mu}:=e^{\hat\mu\log z}$
and $z^R: = e^{R\, \log z}$. The latter matrix is polynomial in $\log z$. The
fundamental matrix (\ref{fundamental}) is analytic on the universal covering
$B\times \widetilde{\mathbb C}^*$.

\begin{remark} Forgetting about the first part (\ref{lin-system1}) of the linear
system and also about the orthogonality (\ref{orthogonal}) we arrive, for a
fixed $v$, at a distinguished fundamental matrix for the system of linear differential
equations (\ref{lin-system2}) with rational coefficients. It essentially coincides with the
fundamental matrix constructed by F.R.Gantmakher \cite{gantmakher} 
and A.H.M.Levelt
\cite{levelt}. The decomposition (\ref{R4}) corresponds to the Levelt's flag
in the space of solutions ${\bf V}$ to (\ref{lin-system2}) determined by the following
non-archimedian valuation function
\beq\label{levelt}
\nu(Y) := {\rm maximal~integer} ~m~ {\rm such ~ that}~\lim_{z\to 0} z^{-r} Y(z)=0
~{\rm for ~any ~ real ~}r<m
\eeq
for a non-zero solution $Y=Y(z)$ (it is  understood that, during the limit in
(\ref{levelt}) $z$ goes to zero within an arbitrary  fixed sector of the
universal covering of ${\mathbb C}^*$), and $\nu(0)=\infty$ (see details in
\cite{levelt}). If 
$$
\infty=\nu_0> \nu_1>\dots >\nu_l
$$
are all the values of the valuation function then
the Levelt's flag 
$$
0={\bf F}_0\subset {\bf F}_1\subset \dots \subset {\bf F}_l={\bf V}
$$
in the space of solutions is defined by
\beq\label{flag25}
{\bf F}_k :=\left\{ y\in {\bf V} \, | \, \nu (y) \geq \nu_k\right\} .
\eeq 
The flag is invariant w.r.t. the monodromy around $z=0$ transformation
given in the basis (\ref{fundamental}) by the matrix
\beq\label{monodromy}
M_0 =e^{2\pi\,i\,\hat\mu} e^{2\pi\,i\,R}.
\eeq
The fundamental matrix (\ref{fundamental}) maps the flag (\ref{flag5}) to the
flag (\ref{flag25}).
\end{remark}

Let us now describe, following \cite{D4}, the ambiguity in the choice
of the fundamental matrix (\ref{fundamental}). Denote $P(V,<\ ,\ >, \hat \mu,
e)\subset Aut\, V$ the group of linear transformations $\Delta: V\to V$
satisfying
\beq\label{Delta1}
z^{\hat\mu} \Delta z^{-\hat\mu} = {\rm polynomial ~in}~z
\eeq
\beq\label{Delta2} 
\Delta^* e^{\pi\, i\, \hat \mu} \Delta =  e^{\pi\, i\, \hat \mu}
\eeq
\beq\label{Delta3}
\Delta\, e = e.
\eeq
The group acts on the fundamental matrices of the form (\ref{fundamental}) by
the formulae
\beq\label{Delta4}
R\mapsto \Delta^{-1} R\Delta
\eeq
\beq\label{Delta5}
\Theta(v;z)\mapsto \Theta(v;z) (\Delta_0 + z \Delta_1 + \dots)
\eeq
where
\beq\label{Delta6}z^{\hat\mu} 
\Delta z^{-\hat\mu} = \Delta_0+ z \Delta_1 +\dots
\eeq

\begin{theorem} Two fundamental matrices of the form (\ref{fundamental})
correspond to the same Frobenius manifold {\rm iff} they are related by the
transformation (\ref{Delta4}) - (\ref{Delta5}).
\end{theorem}

{\bf Definition.} Two 5-tuples $(V_i, <\ ,\ >_i, \hat \mu_i, e_i, R_i)$, $i=1,2$,
are called {\it equivalent} if there exists an isomorphism $\phi: V_1\to V_2$ and
$\Delta\in P(V_1, <\ ,\ >_1, \hat\mu_1)$ such that
$$
\phi^* <\ ,\ >_2 = <\ ,\ >_1
$$
$$
\hat \mu_2 \phi = \phi \hat \mu_1
$$
$$
\phi(e_1) = e_2
$$
$$
\phi^{-1} R_2 \phi = \Delta^{-1} R_1 \Delta.
$$

{\bf Definition.} Class of equivalence of the parameters in (\ref{fundamental})
is called {\it monodromy at $z=0$} of the Frobenius manifold.

Columns of the fundamental matrix (\ref{fundamental}) are gradients of a system
of deformed flat coordinates. Due to constancy of $\hat \mu$, $R$ also columns
of $\Theta(v;z)=\left( \Theta^\alpha_\beta(v;z)\right)$ are gradients of some 
analytic functions on $B\times {\mathbb C}$
that we denote
$\theta_1(v;z)$, \dots, $\theta(v;z)$
$$
\Theta^\alpha_\beta(v;z) =\nabla^\alpha\theta_\beta(v;z).
$$
We obtain a system of deformed flat coordinates of the form
\beq\label{vtilde}
(\tilde v_1(v;z), \dots, \tilde v_n(v;z)) = (\theta_1(v;z), \dots,
\theta_n(v;z)) z^{\hat\mu} z^R.
\eeq

{\bf Definition.} We will call (\ref{vtilde}) {\it Levelt system of deformed
flat coordinates} on $M$ at $z=0$.

Denote $\theta_{\alpha,p}$ the coefficients of the Taylor expansions of the
analytic part of a Levelt
system of deformed flat cordinates
\beq\label{theta-alpha-p}
\theta_\alpha(v;z) =\sum_{p=0}^\infty \theta_{\alpha,p}(v) z^p, ~~\alpha=1,
\dots, n.
\eeq
The coefficients $\theta_{\alpha,p}$ are determined from the recursion procedure
\beq\label{recursion1}
\pal_\lambda\pal_\mu \theta_{\alpha,p}(v) =c_{\lambda\mu}^\nu(v) \pal_\nu
\theta_{\alpha,p-1}(v), ~~p>0,
\eeq
\beq\label{recursion-theta}
\theta_{\alpha,0}= v_\alpha\equiv \eta_{\alpha\epsilon} v^\epsilon
\eeq
with the additional constraints for the matrices $\Theta_p(v):
=(\nabla^\alpha\theta_{\beta,p}(v))$ given by 
\beq\label{recursion3}
(p+1) \Theta_{p+1}(v) + [\Theta_{p+1}(v), {\cal V}] 
= {\cal U}(v) \Theta_p(v) -\sum_{k\geq 1}
\Theta_{p-k+1}(v) R_k, ~~p=0, 1, \dots .
\eeq
In particular,
\beq\label{theta-alpha-1}
\theta_{\alpha,1}(v) ={\pal F(v)\over \pal v^\alpha}
\eeq
\beq\label{theta-1-2}
\theta_{1,2}(v) = {\pal F(v)\over \pal v^\epsilon} v^\epsilon - 2 F(v).
\eeq
In these formulae $F(v)$ is the potential of the Frobenius manifold.

Denote $A_0(M)\subset A(M)$ the subspace of all solutions to (\ref{star})
polynomial in the first coordinate $v^1$. It is a dense subspace when
restricting to functions on a compact in $M$.

\begin{lemma}\label{l-2-7-4} The coefficients $\theta_{\alpha,p}(v)$ form a basis in $A_0(M)$,
i.e., every solution $f(v)$ to (\ref{star}) polynomial in $v^1$ can be uniquely represented
as a finite linear combination 
$$
f(v)=\sum_{\alpha,p} c^{\alpha,p} \theta_{\alpha,p}(v)
$$
with constant coefficients $c^{\alpha,p}$.
\end{lemma}

\pf Let us first prove that $\theta_{\alpha,p}(v)$ are polynomials in $v^1$.
Indeed, from (\ref{recursion1}) it follows that
$$
\pal_1 \theta_{\alpha,p}(v) = \theta_{\alpha,p-1}(v), ~~p>0, 
$$
$$
\pal_1 \theta_{\alpha,0}(v) =\eta_{\alpha\,1}.
$$
Polynomiality follows from these equations. Moreover, we can compute the leading
terms of the polynomials:
\beq\label{leading}
\theta_{\alpha,k+1}(v) =\eta_{\alpha \, 1}{(v^1)^{k+1}\over (k+1)!} + 
\sum_{\ga\ne 1}\eta_{\al\gamma}v^\gamma\,
{(v^1)^k\over k!} + {\rm terms ~of ~lower~ degrees ~in ~} v^1.
\eeq
From the above equations and from (\ref{recursion1}) it also follows that
$\theta_{\alpha,p}(v)$ satisfies (\ref{star}). So, $\theta_{\alpha,p}(v)\in
A_0(M)$ for $\alpha=1, \dots, n$, $p=0, 1, \dots $. 

Let us now prove that these functions form a basis in $A_0(M)$. We use induction
w.r.t. the degree of $f(v)\in A_0(M)$ as a polynomial in $v^1$. For the
polynomials of the degree 0 the equation (\ref{star}) gives
$$
\pal_\alpha\pal_\beta f(v)=0.
$$
So $f(v)$ is a linear function of the flat coordinates
$$
f(v) = \sum c^{\alpha,0} \theta_{\alpha,0}(v) + {\rm const}.
$$
Assuming the Lemma already proved for the polynomials of the degree $\deg_{v^1}
f(v) \leq k-1$ consider
$$
f(v) = f_0(\bar v) + f_1(\bar v) v^1 + \dots +f_k(\bar v) {(v^1)^k\over k!}
$$
where
$$
\bar v = (v^2, \dots v^n).
$$
Then from (\ref{star}) we deduce that $f_k(\bar v)$ is a linear function
$$
f_k(\bar v) = a^\epsilon v_\epsilon + b
$$
where $a^1$, \dots, $a^n$, $b$ are some constant coefficients. Due to
independence of $f_k$ on $v^1$ the coefficients $a^\epsilon$ satisfy
$$
a^\epsilon\eta_{\epsilon\, 1} =0.
$$
Using (\ref{leading}) we show that the polynomial 
$$
f'(v) := f(v) - a^\epsilon \theta_{\epsilon, k+1}(v) -b\, \eta^{1\,
\alpha} \theta_{\alpha,k}(v) \in A_0(M)
$$
has degree in $v^1$ less than $k$. 

It remains to prove linear independence of the functions $\theta_{\alpha,p}(v)$.
Assume that a nontrivial linear combination
$$
\sum_{p=0}^m\sum_{\alpha=1}^n c^{\alpha,p} \theta_{\alpha,p}(v)=0
$$
and not all the coefficients $c^{\alpha,m}$ vanish. Applying the operator
$\pal_1^{m}$ we obtain
$$
\sum_{\alpha=1}^n c^{\alpha,m} v_\alpha ={\rm const}.
$$
This contradicts independency of the flat coordinates.
\epf

We arrive at the main construction of this section, i.e., at an infinite
family of commuting flows
\beq\label{F-hierarchy}
{\pal v\over \pal t^{\alpha,p}} =\{ v(x), H_{\alpha,p}\}_1
=\pal_x \nabla \theta_{\alpha,p+1}(v)=  \nabla\theta_{\alpha,p}(v)\cdot v_x,
~~~H_{\alpha,p} := \bar \theta_{\alpha,p+1}, ~~p=0, 1, \dots .
\eeq
In particular from (\ref{theta-alpha-1}) it follows that
\beq\label{t-1-0}
{\pal v\over \pal t^{1,0}} = {\pal v\over \pal x}
\eeq 
\beq\label{t-alpha-0}
{\pal v^\alpha\over \pal t^{\beta,0}} = c_{\beta\gamma}^\alpha(v) v_x^\gamma.
\eeq
From (\ref{theta-1-2}) we also obtain that
\beq\label{t-1-1}
{\pal v\over \pal t^{1,1}} = v\cdot v_x
\eeq
In the last formula we identify the vector $v\in M$ of the flat coordinates with
the tangent vector $v\in TM$ having the same components.

{\bf Definition.} The hierarchy (\ref{F-hierarchy}) of the first order quasilinear evolutionary PDEs
on ${\cal L}(M)$ is called {\it the Principal Hierarchy} corresponding to the Frobenius manifold $M$.

The product map (\ref{productmap}), (\ref{productmap1}) is given by the following multiplication table
\beq\label{productmap2}
\theta_{\al,p}*\theta_{\beta, q} =\Omega_{\al,p; \beta,q}(v)
\eeq
where the generating function of the coefficients $\Omega_{\al,p; \beta,q}(v)$
is given by
\eqa
&&
\sum \Omega_{\al,p; \beta,q}(v) z^p w^q 
= {<\nabla \theta_{\alpha}(v; z), \nabla \theta_\beta(v;w)> 
-\eta_{\alpha\beta}\over z+w}\nn\\
&&
=\sum_{k=1}^\infty (-1)^k {(w+z)^{k-1}\over k!} < \nabla \theta_\al (v; z),
\pal_z^k \nabla \theta_\beta (v; -z)>.\nn\\
\label{productmap3}
\eeqa

\begin{exam}\label{e-kdv-2-7}
For one-dimensional Frobenius manifold $F(v) ={1\over 6} v^3$.
Here
\beq\label{theta-kdv}
\Theta(v,z)= e^ {z\, v}, ~~ \theta_{1,p} = {v^{p+1}\over (p+1)!}\,,\quad  p=0, 1, \dots .
\eeq
Redenoting the times $t^{1,p}=: t^p$ we obtain the hierarchy
\beq\label{burgers}
{\pal v \over \pal t^p} = {v^p\over p!} v_x
\eeq
In particular
$$
{\pal v\over \pal t^1} = v\, v_x.
$$
This equation (after the change of the sign of the time $t^1\mapsto - t^1$)
is sometimes called dispersionless KdV or nonviscous Burgers equation.
It also coincides with the Riemann simple wave equation.
The hierarchy (\ref{burgers}) is the simplest example of the hierarchies
in our considerations. We suggest to call it {\it Riemann hierarchy}.

The product map of two polynomials (or power series) in $v$ 
is given by the following formula
$$
f(v)* g(v) = \int dv\, f'(v) g'(v).
$$
\end{exam}

\begin{exam}\label{e-I2-k-2-7} For the  two-dimensional Frobenius manifolds
the only parameter is the charge $d$. It is convenient to introduce parameter
$\kappa$ s.t.
$$
d=1-{2\over \kappa}.
$$
For generic $\kappa\neq -1, 0, 1$
\beq\label{I2-k}
F=\frac12\,(v^1)^2\,v^2+{(v^2)^{\kappa+1}\over \kappa^2-1}.
\eeq
The deformed flat coordinates can
be expressed via modified Bessel functions.
The normalized system of deformed
flat coordinates (\ref{vtilde}) reads
\eqa
&&{\tilde v}_1=z^{-\frac12}\sqrt{v^2}\,e^{z\, v^1}\,  
\Gamma(1+\kappa^{-1})\,(\kappa-\kappa^{-1})^{-{1\over 2\,\kappa}}
I_{1\over\kappa}[{2\,z\,\sqrt{\kappa^2-1}\over \sqrt{\kappa}} 
{(v^2)}^{\kappa\over2}]\nn\\
&&=e^{z\, v^1} \left[\sum_{m\geq 0} {\Gamma(1+\kappa^{-1})\over
\Gamma(m+1+\kappa^{-1})} (\kappa-\kappa^{-1})^m (v^2)^{m\, \kappa+1} {z^{2\,
m}\over m!} \right] \,z^{-{1\over 2}+{1\over \kappa}}, \nn\\
&&{\tilde v}_2= z^{-\frac12}\sqrt{v^2}\,e^{z\, v^1} \Gamma(1-\kappa^{-1})
(\kappa-\kappa^{-1})^{{1\over 2\, \kappa}}
I_{-{1\over \kappa}}[{2\,z\,\sqrt{\kappa^2-1}\over \sqrt{\kappa}} 
(v^2)^{\kappa\over2}]
-{z^{-\frac1{\kappa}-\frac12}}
\nn\\
&&={z^{-1}}
\left[ e^{z\, v^1} \sum_{m\geq 0} 
{\Gamma(1-\kappa^{-1})\over \Gamma(m+1-\kappa^{-1})}(\kappa-\kappa^{-1})^{m}
(v^2)^{\kappa\, m} {z^{2\, m}\over m!}
 -1\right] \, z^{{1\over 2} -{1\over \kappa}}.\nn
\eeqa
The matrices $\hat\mu$ and $R$ in (\ref{vtilde}) are given by
$$
\hat\mu=\left(\matrix{ -\frac12+\frac1{\kappa}&0 \cr 0&\frac12-\frac1{\kappa}
}\right),\quad
R=0.
$$

The hamiltonian flow (\ref{t-1-1}) after changing of the sign of the time
variable $t=-t^{1,1}$ and redenoting $v^1=v$, $v^2=\rho$ coincides with the
equations of motion of one-dimensional polytropic gas with the equation of
state $p={\kappa\over\kappa+1}\rho^{\kappa+1}$:
\eqa\label{fluid}
&&
v_t+\left( {v^2\over 2} + \rho^\kappa\right)_x=0
\nn\\
&&
\rho_t +(\rho\,v)_x =0.
\nn
\eeqa
The bihamiltonian structure (\ref{normalfm}) reads
\eqa\label{I2-pb}
&&\{v(x),v(y)\}^{[0]}_\lm=2\rho^{\kappa-1}(x)\,\delta'(x-y)
+\left( \rho^{\kappa-1}\right)_x\,\delta(x-y),\nn\\
&&\{v(x),\rho(y)\}^{[0]}_\lm=(v(x)-\lm)\,\delta'(x-y)+\frac1{\kappa}\,v'(x)\,
\delta(x-y),\nn\\
&&\{\rho(x),\rho(y)\}^{[0]}_\lm=\frac1{\kappa}\left(2\,\rho(x)\,\delta'(x-y)+
\rho'(x)\,\delta(x-y)\right).
\eeqa
This bihamiltonian structure for the polytropic gas equations has been found by
P.Olver \cite{olver2}. The above formulae remain valid also for the exceptional
values $\kappa=\pm 1$ where the expression for the potential of the Frobenius
manifold is to be modified. For the particular value $\kappa=3$ the Frobenius
manifold (\ref{I2-k}) corresponds to the $A_2$ topological minimal model
\cite{DVV1}.
\end{exam}

\begin{exam}\label{e-toda-2-7} The exceptional 2-dimensional Frobenius manifold 
with the charge $d=1$
\beq\label{potential-toda}
F=\frac12\,(v^1)^2\,v^2+e^{v^2}
\eeq
corresponds to the quantum cohomology of ${\bf CP}^1$ \cite{verdier} 
(it will also be
called ${\bf CP}^1$ sigma-model).
The deformed flat coordinates can also be expressed via modified Bessel
functions.
We have
\beq
\hat\mu=\left(\matrix{ -\frac12&0 \cr 0&\frac12}\right),\quad
R=\left(\matrix{0&0 \cr 2&0}\right).
\eeq
So
the normalized system (\ref{vtilde}) reads
\beq
({\tilde v}_1(v;z),{\tilde v}_2(v;z))=(\theta_1(v;z),\theta_2(v;z))\,
\left(\matrix{z^{-{1\over 2}}&0 \cr 0&z^{1\over 2}}\right)\left(\matrix{1&0 \cr
2\log{z}&1}\right)
\eeq
where
\eqa
&&\theta_1(v;z)=-2\,e^{z v^1}\left(K_0(2 z e^{\frac12\,v^2})
+(\log{z}+\gamma) I_0(2 z e^{\frac12\,v^2})\right)\nn\\
&&=-2 e^{z\, v^1} \sum_{m\geq 0} (\gamma-{1\over 2} v^2 + \psi(m+1)) e^{m\, v^2}
{z^{2\, m}\over (m!)^2},\label{theta-toda-1}\\
&&\theta_2(v;z)=z^{-1}\,e^{z v^1}\,I_0(2 z e^{\frac12\,v^2})-z^{-1}\nn\\
&&=z^{-1} \left(\sum_{m\geq 0}e^{m\, v^2 + z \, v^1} {z^{2\, m}\over
(m!)^2}-1\right).\label{theta-toda-2}
\eeqa
Here $\gamma$ denotes Euler's constant, $\psi(z)$ stands for the digamma
function.

The $-t^{1,1}$-flow has still the meaning of one-dimensional isentropic fluid
(polytropic gas)
with the equation of state $p=(\rho^2-2\rho+2) e^\rho$. It is more instructive to look
at the $t=t^{2,0}$-flow. Using (\ref{t-alpha-0}) we obtain the so-called long
wave limit of the Toda lattice equation
\beq\label{toda}
\rho_{tt} - (e^\rho)_{xx}=0.
\eeq
The bihamiltonian structure reads
\eqa\label{cp1-pb}
&&\{v(x),v(y)\}^{[0]}_\lm=2\,e^{\rho(x)}\,\delta'(x-y)+
\left(e^{\rho(x)}\right)_x\,\delta(x-y),\nn\\
&&\{v(x),\rho(y)\}^{[0]}_\lm=(v(x)-\lm)\,\delta'(x-y),\nn\\
&&\{\rho(x),\rho(y)\}^{[0]}_\lm=2\,\delta'(x-y).\nn
\eeqa
\end{exam}

\begin{remark} In the setting of two-dimensional topological field theory
\cite{DW, Dij5, Dij3, Witten2}
the Frobenius manifold is called {\it small phase space}. The basis of the Hamiltonian
densities $\theta_{\alpha,p}(v)$ coincides with particular two point
 tree level correlators 
$$
\theta_{\alpha,p}(v) =<\tau_p(\phi_\alpha) \tau_0(\phi_1)>_0
$$
as functions on the small phase space. Here $\phi_1$ corresponds to the identity
operator. Other two point tree level correlators are
$$
<\tau_p(\phi_\alpha)\tau_q(\phi_\beta)>_0 = \theta_{\alpha, p} *
\theta_{\beta,q}.
$$
To extend these formulae on the {\it big phase space} one is to evaluate these
functions on the topological solution $v=v({\bf t})$ 
(see (\ref{fp0}), (\ref{fp1}) below) of the Principal Hierarchy.  
\end{remark}

We are to clarify an important point about the bihamiltonian nature of the
hierarchy (\ref{F-hierarchy}). Although all these equations, as it has been
proved in Theorem \ref{t-2-7-1}, are bihamiltonian flows w.r.t. the Poisson pencil 
(\ref{normalfm}),
their hamiltonian densities in $A(M)$ not always belong to $A_0(M)$. For
example, the spatial translation flow on ${\cal L}(M)$, $M$ being the Frobenius
manifold of the Example \ref{e-toda-2-7} above
$$
{v_1}_t = {v_1}_x, ~~~{v_2}_t = {v_2}_x
$$
w.r.t. the first Poisson bracket has the hamiltonian $\bar{v_1 v_2}$ with
the density $v_1 v_2\in A_0(M)$. But the hamiltonian density of this flow
w.r.t. the second Poisson bracket 
$$
g={1\over 4} v_2^2 - \left[ {1\over 2} v_2 - \log{ v_1+\sqrt{ v_1^2 - 4 \,
e^{v_2}}}\right]^2
$$
does not belong to $A_0(M)$. We will show below that the hamiltonians of 
all but a finite number of 
the flows of the hierarchy (\ref{F-hierarchy}) w.r.t. the second Poisson
bracket belong to $A_0(M)$. 

\subsubsection{Periods of the Frobenius manifold and the bihamiltonian recursion
for the Principal Hierarchy}\label{sec-3-6-3}\par

The first step will be to show that the hamiltonians of the hierarchy
(\ref{F-hierarchy}) are obtained by a triangular transformation from those 
given by the bihamiltonian recursion
procedure presented in Section \ref{sec-3-1-2}
above. To do this we are to find the reducing
transformation for the Poisson pencil $\{~,~\}^{[0]}_\lambda=\{~,~\}^{[0]}_2
-\lambda\{~,~\}^{[0]}_1$ 
and to express the coefficients of the expansion of this reducing transformation
in terms of the hamiltonian densities $\theta_{\alpha,p}(v)$.

{\bf Definition.} The functions $p=p(v;\lambda)$ satisfying
\beq\label{2-7-55}
(\nabla^*-\lambda\, \nabla) dp=0
\eeq
are called {\it periods} of the Frobenius manifold. The system
(\ref{2-7-55}) is called {\it the Gauss-Manin system} on the Frobenius 
manifold \cite{D3}. 

Here $\nabla^*$ is the Levi-Civita connection for the metric $(~,~)$. The
connection is well-defined outside the discriminant $\Sigma\subset M$ (see
details in \cite{D3}).

Choosing a system of $n$ independent periods we obtain a system of flat
coordinates $p^1(v;\lambda)$, \dots, $p^n(v;\lambda)$ for the flat pencils of
metrics $(~,~)_\lambda:= (~,~)-\lambda\, <\ ,\ >$
\beq
\left(dp^i(v;\lambda), dp^j(v;\lambda)\right)_\lambda = G^{ij}
\eeq
for some constant nondegenerate matrix $G^{ij}$. According to the results
of Section \ref{sec-2-4-5}, 
choosing $p^i(v;\lambda)$, $i=1, \dots, n$ as a new 
system of depending variables we obtain a reduction of the Poisson pencil
to the canonical form
$$
\{ p^i(v(x);\lambda), p^j(v(y);\lambda)\}_\lambda = G^{ij} \delta'(x-y).
$$
 
Using results of \cite{duality}
we will produce a particular system of independent periods specified 
according to their behaviour for large $\lambda$. For technical reasons
it will be convenient to label these periods by {\it lower} indices; the roles
of upper and lower indices therefore will be interchanged in the subsequent
formulae.

\begin{theorem}\label{t-2-7-7} Let $M$ be a semisimple Frobenius manifold such that the
spectrum of ${\cal V}$ does not contain half-integers. Then the transformation
$p_\alpha=p_\alpha(v;\lambda)$ reducing the Poisson pencil to the constant form
$$
\{ p_\alpha(v(x);\lambda),p_\beta(v(y);\lambda)\}_\lambda =
G_{\alpha\beta}\delta'(x-y)
$$
with a constant matrix $G=(G_{\alpha\beta})$ 
\beq\label{gram}
G=-2\pi\,\eta\left[ e^{\pi\,i\,R} e^{\pi\,i\,\hat\mu} + 
e^{-\pi\,i\,R} e^{-\pi\,i\,\hat\mu}\right]^{-1}
\eeq
is given by the formula
\beq\label{period1}
{\bf p}(v;\lambda) = (p_1(v;\lambda), \dots, p_n(v;\lambda))=
(\pi_1(v;\lambda), \dots, \pi_n(v;\lambda)) \lambda^{-{1\over 2}-\hat\mu}
\lambda^{-R}
\eeq
where the
vector function ${  \pi}(v;\lambda)=(\pi_1(v;\lambda), \dots,
\pi_n(v;\lambda))$ is analytic for sufficiently large $|\lambda|$. It has the
following Taylor expansion at $\lambda=\infty$
\eqa\label{period2}
&&{  \pi}(v;\lambda)={  \pi}^{(0)}
(\lambda)+\sum_{m=0}^\infty \lambda^{-m} \sum_{p+q=m}
{\bf \theta}_p(v) 
\Gamma_q(R,\hat\mu +m+{1\over
2}),\\
&&
{ \pi}^{(0)}=\lambda \sum_{q\geq 0} \omega_1\, 
\Gamma_q(R,\hat\mu +q-{1\over 2}), ~~~\omega_1:=(\eta_{11}, \eta_{12}, \dots,
\eta_{1n}). \label{period21}
\eeqa
In this formula ${ \theta}_p(v) = (\theta_{1,p}(v), \dots, \theta_{n,p}(v))$.
The decomposition of the matrix polynomial $e^R = [e^R]_0 + [e^R]_1 +\dots$
was defined in (\ref{R8}), (\ref{R8-b}). 
The matrices $\Gamma_q(R, \hat\mu+a)$ for an arbitrary
complex number $a$ such that the spectrum of the matrix $\hat\mu + a\, {\rm id}$
does not contain negative integers are defined by
\beq
\Gamma_q(R,\hat\mu + a) = \left[ e^{R\, \pal_\nu}\right] _q \Gamma(\hat\mu+a
+\nu)_{\nu=0}.
\eeq
\end{theorem}

\pf As we already know, the reducing transformation for the Poisson pencil
$\{~,~\}_2-\lambda\{~,~\}_1$ is given by a system of the flat coordinates for the flat
pencil of metrics $(~,~)-\lambda <\ ,\ >$. The  latter can be obtained \cite{D3} by applying
the Laplace-type integrals to the deformed flat coordinates
$$
p_\alpha(v;\lambda) =\oint e^{-\lambda\, z} \tilde v_\alpha(v;z) \,{dz\over \sqrt{z}}, ~~\alpha=1,
\dots, n.
$$
Here the symbol of loop integral means just the possibility of integration by
parts dropping the boundary terms. On a semisimple Frobenius manifold the above  
loop integral can be regularized as follows. Consider the integral
\beq\label{laplace} 
p^{(\nu)}_\alpha(v;\lambda) =\int_0^{\infty e^{i\varphi}} 
e^{-\lambda\, z} \tilde v_\alpha(v;z) \,{dz\over z^{{1\over 2}-\nu}}
\eeq
along the ray $\arg z = \varphi$ on the complex $z$-plane. Here $\nu$ is a complex
parameter. The integral converges at $z=0$ for ${\rm Re}\, \nu >> 0$. It also
converges at $z= e^{i\varphi}\infty$ for sufficiently large $|\lambda|>r$ for
some $r = r(v)$
due to the exponential behaviour of the deformed flat coordinates at $z=\infty$
(here the semisimplicity of the Frobenius manifold plays the crucial role!).
Rotating the argument $\varphi$ we obtain an analytic continuation of the
integral onto the universal covering of the disc $r<|\lambda|<\infty$. It
remains to analytically continue the integral into the point of interest $\nu=0$
to obtain the needed functions
$$
p_\alpha(v;\lambda) := p^{(0)}_\alpha (v;\lambda).
$$
We will do the needed
analytic continuation just integrating the terms of the expansion (\ref{vtilde}),
(\ref{theta-alpha-p})
and then setting $x$ to zero. This can be done using the following calculations
$$
\int_0^{\infty e^{i\varphi}}e^{-\lambda\, z} z^{\hat \mu + p + \nu -{1\over 2}}
z^R dz = \int_0^\infty e^{-t} t^{\hat\mu + p + \nu-{1\over 2}}
t^{R_0 +{R_1\over \lambda} +{R_2\over \lambda^2}+\dots} dt \, \lambda^{-(\hat\mu
+p + \nu +{1\over 2})} \lambda^{-R}
$$
$$
=\sum_{q\geq 0}\sum_{k\geq 0} \int_0^\infty e^{-t} t^{\hat\mu + p + \nu-{1\over 2}}
\lambda^{-q}{[R^k]_q \log^k t\over k! } dt \, \lambda^{-(\hat\mu
+p + \nu +{1\over 2})} \lambda^{-R}
$$
$$
=\sum_{q\geq 0}\sum_{k\geq 0} {1\over k!} \pal_\nu^k \Gamma(\hat\mu+p+\nu+{1\over
2}) {[R^k]_q\over \lambda^q}\lambda^{-(\hat\mu
+p + \nu +{1\over 2})} \lambda^{-R}
$$
$$= \sum_q [e^{R\pal_\nu}]_q
\Gamma(\hat\mu+p+q+\nu+{1\over 2}) \, \lambda^{-(\hat\mu
+p+q + \nu +{1\over 2})} \lambda^{-R} 
$$
where we used the following obvious integral
$$
\int_0^\infty e^{-t} t^{s-1} \log^k t\, dt = \partial_s^k \Gamma(s)
$$
and also the commutation relation
$$
f(\hat\mu) [P(R_0, R_1, \dots)]_q = [P(R_0, R_1, \dots)]_q f(\hat\mu+q)
$$
valid for an arbitrary polynomial of the matrices $R_i$ and for an arbitrary
analytic function $f$. 

We see that the coefficients of the above series are meromorphic functions on
the complex $\nu$-plane with the poles at
\beq\label{poles}
\nu\in \cup_{\mu\in Spec \, \hat \mu} \cup_{k>0} (-\mu -k-{1\over 2}).
\eeq
For every pole in (\ref{poles}) only finite number of the coefficients
of the series become infinite. Therefore the sum of the series is a meromorphic
function in $x$ with poles at (\ref{poles}). Due to the assumption about the spectrum
of ${\cal V}$ = spectrum of $\hat\mu$ the value $\nu=0$ is not a pole of this 
series. Setting $\nu=0$ we obtain the proof of (\ref{period1}), (\ref{period2}).
Note that the term $\pi^{(0)}_\alpha(\lambda)$ does not depend on $v$. It drops
from the Gauss - Manin system (\ref{2-7-55}). We choose this term in such a way
to have the identity
\beq
{\pal p(v;\lambda)\over \pal \lambda} = - {\pal p(v;\lambda)\over \pal v^1}
\eeq
valid.
\epf

\begin{remark}
The Poisson brackets $\{ \pi_\alpha(v(x);\lambda),
\pi_\beta(v(y); \lambda)\}_\lambda$ are also constant but they depend on
$\lambda$:
$$
\{ \pi_\alpha(v(x);\lambda),
\pi_\beta(v(y); \lambda)\}_\lambda =\tilde g_{\alpha\beta}(\lambda) \delta'(x-y)
$$
$$
\left(\tilde G_{\alpha\beta}\right) =- 2\pi\, \lambda\,\eta 
\left[ e^{\pi\, i \, ( R_0 +{R_1\over \lambda} +\dots) }e^{\pi\,i\,\hat \mu}
+e^{-\pi\,i\,( R_0 +{R_1\over \lambda} +\dots)} e^{-\pi\,i\,\hat\mu}\right]^{-1}.
$$
\end{remark}

We are now able to write explicitly down the bihamiltonian recursion relation
for the Principal Hierarchy.

\begin{theorem} Under the assumptions of the theorem \ref{t-2-7-7} the following
recursion relation holds true
\beq\label{recur-2-7}
{\cal R} {\pal\over\pal t^{p-1}} = {\pal\over\pal t^p} \left( p+\hat\mu+{1\over
2}\right) +\sum_{k=0}^p {\pal\over\pal t^{p-k}} R_k.
\eeq
Here
$$
{\pal\over\pal t^p}:= \left( {\pal\over\pal t^{1,p}}, \dots,
{\pal\over\pal t^{n,p}}\right).
$$
\end{theorem}

\pf Applying the recursion operator to the hamiltonian flow
$$
v_t = \{ v(x), \bar f\}_1 =\pal_x \nabla f(v), ~~ \bar f := \int f(v)\, dx
$$
with an arbitrary Hamiltonian density $f(v)$
one obtains
$$
{\cal R}^\alpha_\beta \, v_t^\beta = {\nabla^*}^\alpha \pal_\gamma f(v)\,
v^\gamma_x.
$$
For the generating function of the flows of Principal Hierarchy we take
$$
f= \tilde v(v;z)
$$
and use the identity \cite{D3}
\beq
\nabla^* d\tilde v = (\pal _z -{1\over 2\, z}) \nabla d\tilde v
\eeq
to arrive at the needed recursion relation. \epf 

If the spectrum of $\hat\mu$ contains half-integers then the Gauss - Manin
system has solutions polynomial in $\lambda$. The recursion operator
(\ref{recur-2-7}) becomes degenerate. In other words, although the flows of the
Principal Hierarchy remain bihamiltonian, their Hamiltonians are not described
by the bihamiltonian recursion procedure.

More specifically, if ${1\over 2}\in \, {\rm Spec}\, \hat\mu$ then
the two Poisson brackets have common Casimirs. E.g., in the particular example
of the ${\bf CP}^1$ model the variable $v_2$ is the density of a Casimir
for both Poisson brackets.

\subsubsection{Solutions to the Principal Hierarchy and their tau-functions}
\label{sec-3-6-4}\par

We will now describe a natural class of solutions to the hierarchy
(\ref{F-hierarchy}) and compute explicitly their tau-functions.

We will consider {\it analytic solutions} of the hierarchy, i.e.,  power 
series in the variables ${\bf t}:= (t^{\alpha,p})$ with the coefficients
in ${\mathbb C}[[\epsilon]]$
\beq\label{formalsolution}
v=v(x, {\bf t},\epsilon)=a_0(\epsilon)
+ \sum_{k>0}a_{\alpha_1, p_1; \dots;\alpha_k, p_k }(\epsilon)
t^{\alpha_1, p_1}\dots t^{\alpha_k, p_k}|_{t^{1,0}\mapsto t^{1,0}+x}.
\eeq
Analyticity of the function in infinite number of variables is to be 
understood as follows.
Setting $t^{\alpha,p}=0$ for all $\alpha=1, \dots, n$, $p\geq N$ for an 
arbitrary positive integer $N$ we must obtain a power series in the 
finite number of variables $x$ and $t^{\beta,q}$, 
$\beta=1, \dots, n$, $0\leq q < N$ with the coefficients in  ${\mathbb C}[[\epsilon]]$. Every term of the formal power series in $\epsilon$ must be a
{\it convergent} series in $x$ and $t^{\beta,q}$, 
$\beta=1, \dots, n$, $0\leq q < N$
in a ball near the origin of the space ${\mathbb C}^{n \,N+1}$ 
(the size of the ball may depend on $N$).

To make it possible the substitution of the solution to the equations of the 
hierarchy the vector $a_0(0)$ must be a point in the Frobenius manifold $M$. 
The vector $v_x$ for $\epsilon=0$, ${\bf t}=0$, $x=0$ can be considered as an element of 
the tangent space at this point:
\beq\label{velocity}
v_x(x=0, {\bf t}=0, \epsilon=0)\in T_{v=a_0(0)} M.
\eeq

{\bf Definition}. The solution (\ref{formalsolution}) is called 
{\it monotone} at the origin if the vector (\ref{velocity}) is an invertible element of the algebra $T_{v=a_0(0)}M$.

We will now construct a dense subset in the space of analytic monotone 
solutions. 

Let us fix a point $v_0\in M$ and a collection of formal power series 
$c^{\alpha,p}(\epsilon)\in {\mathbb C}[[\epsilon]]$ with constant coefficients,
$\alpha=1, \dots, n$, $p\geq 1$ with only finite number of them being nonzero
such that the multiplication operator
\beq\label{operator-m}
m_0:=\left( \cdot\, \nabla \sum_{p\ge 1} 
c^{\alpha, p}(0) \theta_{\alpha, p-1}(v)\right) _{v=v_0}: T_{v_0}M\to T_{v_0}M
\eeq
is invertible. We are to also fix $n$ series 
$c^{\alpha,0}(\epsilon)\in {\mathbb C}[[\epsilon]]$ such that
$$
c^{\al,0}(0) :=- \nabla^\al \sum_{p\ge 1} 
c^{\beta,p}(0)\theta_{\beta,p}(v_0).
$$
The solution to the hierarchy will be defined by the following system of 
equations
\beq\label{tsarev}
x\, e +\sum_{p\geq 0} \tilde t^{\al, p} \nabla \theta_{\al, p}(v)=0
\eeq
where
\beq\label{shift}
\tilde t^{\al, p}: = t^{\al, p}-c^{\al, p}(\epsilon).
\eeq

\begin{theorem} 1) There exists a unique solution to the system (\ref{tsarev})
in the form (\ref{formalsolution}) with $a_0(0) = v_0$. It satisfies the 
equations of the hierarchy (\ref{F-hierarchy}).
2) The solutions of the form (\ref{formalsolution}) are dense in the space of analytic
monotone solutions to the hierarchy.
\end{theorem}

The first part of the theorem is an analogue of the Tsarev's generalized 
hodograph transform in the theory of integrable systems of hydrodynamic type 
\cite{tsarev} adapted for the case of hierarchies of these systems. 
A construction of a dense set of solutions for certain particular classes of 
systems of hydrodynamic type has been obtained in \cite{pisa}, \cite{tsarev1}.

\pf Differentiating (\ref{tsarev}) w.r.t. $v^1$, \dots, $v^n$ and setting
$x=0$, ${\bf t}=0$, $\epsilon=0$ we obtain the Jacobi matrix of the system
coinciding with the nondegenerate matrix $-m_0$ of the operator 
(\ref{operator-m}). Therefore
existence and uniqueness of an analytic solution of (\ref{tsarev}) 
with $t^{\al,p}=0$ for $p\geq N$
$$
v=v_0+\sum_{k>0}v_0^{[k]}\epsilon^k+ \sum_{k>0}\sum_{0\leq p_1, \dots, p_k < N}
a_{\alpha_1, p_1; \dots;\alpha_k, p_k }(\epsilon)
t^{\alpha_1, p_1}\dots t^{\alpha_k, p_k}|_{t^{1,0}\mapsto t^{1,0}+x}
$$
for every positive $N$ is an immediate consequence of the implicit function theorem. Differentiating (\ref{tsarev}) w.r.t. $x$ and $t^{\al,p}$ 
we obtain 
$$
w\cdot v_x =-e
$$
$$
w\cdot \pal_{t^{\alpha,p}} v =- \nabla \theta_{\al,p}(v)
$$
where the operator of multiplication by the vector
$$
w:= \sum \tilde t^{\beta,q} \nabla \theta_{\beta,q-1}(v)
$$
is invertible for small $t^{\beta,q}$ and for $v$ close to $v_0$ due to our 
choice of  $v_0$ and of the constants $c^{\beta,q}(\epsilon)$. 
The equations of the hierarchy
$$
 \pal_{t^{\alpha,p}} v =\nabla \theta_{\al,p}(v)\cdot \pal_x v
$$
readily follow by dividing over $w$.

Let us prove density of the constructed solutions. Let 
$v=v(x, {\bf t},\epsilon)$ 
be a monotone analytic solution to the hierarchy s.t. $v(0,0,0)=v_0\in M$. 
From 
the monotonicity condition and from the level zero equations of the hierarchy
\beq\label{level0}
\pal_{t^{\al,0}}v^\beta = c_{\al\gamma}^\beta(v) \pal_x v^\gamma
\eeq
it follows that
$$
\det\left( \pal v^\alpha(0,0,0)\over \pal t^{\beta,0}\right) \neq 0.
$$  
Restricting this solution onto the subspace $t^{\al,p}=0$ for $p>0$ and 
using the nondegeneracy of the Jacobian we can rewrite this restriction in the 
following implicit form
\eqa
x+&t^{1,0}&=f^1(v,\epsilon)\ ,\nn\\
&t^{2,0}&=f^2(v,\epsilon)\ ,\nn\\
& &\dots\ ,\label{array}\\
&t^{n,0}&=f^n(v,\epsilon)\ .\nn
\eeqa
Here $f^1(v,\epsilon), \dots, f^n(v,\epsilon)$ are some formal power series in $\epsilon$ with the coefficients analytic in a ball near
$v=v_0$.

\begin{lemma} The functions $f^\al(v,\epsilon)$ have the form
\beq\label{grad}
f^\al(v,\epsilon) =\nabla^\al f(v,\epsilon)
\eeq
for some function $f(v,\epsilon)\in A(M)$ 
Conversely, every function $f(v,\epsilon)$ satisfying (\ref{star})
defines a solution to the Principal Hierarchy (\ref{F-hierarchy}) 
in the implicit
form
\beq\label{gen-sol}
x\,e+\sum t^{\al,p}\,\nabla\theta_{\al,p}(v)=\nabla f(v,\epsilon).
\eeq 
\end{lemma}

\pf Differentiating (\ref{array}) w.r.t. $x$ and $t^{\alpha,0}$ and using again
the level zero part (\ref{level0}) of the equations of the hierarchy
we obtain
$$
\delta_\beta^\alpha = {\pal f^\alpha\over \pal v^\epsilon} c_{\beta\gamma}^\epsilon {\pal v^\gamma\over \pal x}.
$$
The last equation can be recast into the form
$$
e_\alpha = \nabla f_\alpha \cdot v_x
$$
where we denote
$$
f_\al (v,\epsilon):= \eta_{\al\beta}f^\beta(v,\epsilon).
$$
In particular, 
$$
v_x = (\nabla f_1)^{-1}.
$$
Hence
$$
\nabla f_\al = e_\alpha \cdot \nabla f_1,
$$
i.e.,
$$
\pal_\beta f_\alpha = c_{\al\beta}^\gamma \pal_\gamma f_1.
$$
The symmetry in $\al$ and $\beta$ proves closedness of the one-form
$$
f_\alpha(v,\epsilon) dv^\al = d\, f(v,\epsilon)
$$
and also implies the equation (\ref{star}).
\epf

To finish the proof of the Theorem we just approximate near $v=v_0$ the coefficients of the $\epsilon$-expansion of
the solution $f(v,\epsilon)$
to the system (\ref{star}) by the coefficients of the expansion of 
another solution $\tilde f(v,\epsilon)$ to (\ref{star})
polynomial in $v^1$. Applying Lemma \ref{l-2-7-4}
to the function $\tilde f(v,\epsilon)$ (with possibly the adding to it 
a linear in $v^\al$ term) we obtain a finite linear combination
$$
\tilde f(v,\epsilon) =\sum_{p\ge 1} 
c^{\alpha, p}(\epsilon) \theta_{\alpha, p}
$$
for some constants $c^{\al,p}(\epsilon)$. 
We now use these constants together with 
$\nabla^\al \tilde f(v_0,0) =-c^{\al,0}(\epsilon)$ 
to produce a 
solution (\ref{tsarev})
to the hierarchy. Let us denote this solution  
$\tilde v(x, {\bf t},\epsilon)$. It will approximate the given number of the coefficients of the expansion of the given solution to 
$v(x, {\bf t},\epsilon)$ to (\ref{gen-sol}) in a power series in $\epsilon$ 
when restricted onto a finite-dimensional subspace
$t^{\alpha, p}=0$ for $p\geq N$ for sufficiently small $|x|$ and 
$|t^{\alpha,p}|$, $0\leq p < N$. The Theorem is proved.
\epf

The equation (\ref{tsarev}) can be rewritten as the following stationary 
point equation
\beq\label{tsarev1}
\nabla \Phi_{x,{\bf t},c(\epsilon)}(v)=0
\eeq
where the function $ \Phi_{x,{\bf t},c(\epsilon)}(v)$ on $M$ depending on 
the parameters
$x$, ${\bf t}= (t^{\al,p})$ and $c(\epsilon)=(c^{\al,p}(\epsilon))$ 
has the form
\beq\label{Fi}
\Phi_{x,{\bf t},c(\epsilon)}(v) =
\sum \tilde t^{\alpha,p} \theta_{\alpha,p}(v)|
_{\tilde t^{1,0}\mapsto \tilde t^{1,0}+x}.
\eeq
(The solution depends also on the choice of the point $v_0\in M$ such that 
$$
\nabla^\al \Phi_{0,{\bf 0}, c(0)}(v_0) =0, ~~\al=1, \dots, n.
$$
However, locally $v_0$ is uniquely determined by $c(0)$ due to invertibility 
of the vector $\nabla \sum c^{\al, p}(0) \theta_{\al, p-1}(v_0)\in T_{v_0}M$.)
As we have just proved the dependence of the stationary point that we denote
$v(x, {\bf t}, c(\epsilon))$ on the 
parameters $x$ and 
${\bf t}$ satisfies the equations of the hierarchy (\ref{F-hierarchy}). 
The representation
(\ref{tsarev1}) will be useful in all calculations with the solutions of the 
hierarchy and with their tau-functions. Observe that $\Phi_{x,{\bf t},c(\epsilon)}(v)$
can be considered as a vector of the space 
${\cal K}\otimes {\mathbb C}[[\epsilon]]$ of the densities of the conservation 
laws of the hierarchy depending explicitly on $x$, ${\bf t}$ and $\epsilon$.

\begin{theorem}
The tau-function of the solution $v(x,{\bf t},c(\epsilon))$ defined by 
(\ref{tsarev1}) 
has the form
\eqa\label{F0}
&&
{\cal F}_0 (x, {\bf t},c(\epsilon))=\epsilon^2\log \, \tau = {1\over 2} 
\Phi_{x,{\bf t},c(\epsilon)}(v)* \Phi_{x,{\bf t},c(\epsilon)}(v) 
|_{v=v(x,{\bf t},c(\epsilon))}
\nn\\
&&
={1\over 2} \sum \tilde t^{\alpha,p} \tilde t^{\beta,q}
\Omega_{\alpha,p;\beta,q}(v(x,{\bf t},c(\epsilon))). 
\nn\\
\eeqa
The first derivatives of the tau-function w.r.t. the times of the hierarchy are given by the formula
\beq\label{derF0}
\epsilon^2\pal_{t^{\alpha,p}} \log \, \tau = 
\theta_{\alpha,p}(v) *  \Phi_{x,{\bf t},c(\epsilon)}(v) 
|_{v=v(x,{\bf t},c(\epsilon))}=
\sum\tilde t^{\beta,q}\Omega_{\alpha,p;\beta,q}(v(x,{\bf t},c(\epsilon))) .
\eeq
\end{theorem}

Recall that the product map $*$ was defined in (\ref{productmap}), (\ref{productmap1}),
(\ref{productmap2}).

\pf Differentiating (\ref{F0}) w.r.t. $t^{\gamma, r}$ and using
$$
\nabla (\theta_{\alpha,p} * \theta_{\beta,q}) =\nabla \theta_{\alpha,p} \cdot \nabla \theta_{\beta,q}
$$
we obtain
$$
\epsilon^2\pal_{t^{\gamma,r}} \log \tau 
= \sum \tilde t^{\beta,q} \theta_{\gamma,r} * \theta_{ \beta,q}|_{\tilde t^{1,0} 
\mapsto \tilde t^{1,0} + x}
+ {1\over 2} < \nabla \Phi_{x,{\bf t}, c(\epsilon)}(v) \cdot \nabla 
\Phi_{x, {\bf t}, c(\epsilon)}(v) \cdot \nabla \theta_{\gamma,r}(v), v_x>.
$$
The second part of the formula vanishes for $v= v(x, {\bf t}, c)$ due to (\ref{tsarev1}). This proves (\ref{derF0}). 
Repeating the trick we obtain
$$
\epsilon^2{\pal^2 \log \tau\over \pal t^{\alpha, p} \pal t^{\beta, q}} =\Omega_{\alpha,p; \beta,q}.
$$
The theorem is proved. 
\epf

\begin{exam} The particular solution to the Principal Hierarchy
specified by the constants
\beq
c^{\alpha,p}=\delta^\alpha_1 \,\delta^p_1
\eeq
will be called {\rm topological solution}. It is specified by the following
fixed point equation
\beq\label{fp0}
v =\nabla\Phi_{x,{\bf t}}(v), ~~\Phi_{x,{\bf t}}(v)=\sum \bar
t^{\alpha,p}\theta_{\alpha,p}(v).
\eeq
The expansion of the topological solution has the form
\beq\label{fp1}
v^\al(t)=t^{\al,0}+\sum_{k\ge 1,\, p_i\ge 1} 
A^\al_{\beta_1, q_1;\dots;\beta_k, q_k}(t^{1,0},\dots, t^{n,0})\,
t^{\beta_1,q_1}\dots t^{\beta_k,q_k},
\eeq
the coefficients are determined recursively by (\ref{fp0}). For example,
we have
$$
A_{\beta,q}^\al=\left. 
\frac{\pal\theta_{\beta,q}}{\pal v_\al}\right|_{v_\gamma=t^{\gamma,0}},
\quad 
A^\al_{\beta_1,q_1;\beta_2,q_2}=\left.
\frac12\,\frac{\pal^2\theta_{\beta_1,q_1}}
{\pal v_\al \pal v_\gamma}\,\frac{\pal\theta_{\beta_2,q_2}}
{\pal v_\gamma}\right|_{v_\xi=t^{\xi,0}}.
$$
As it was shown in \cite{D92},
the logarithmic derivatives of the tau-function of the topological solution
satisfy the genus zero topological recursion relations \cite{DW}. 
In order to formulate these recursion relations we introduce the symbols
(``the genus zero correlation functions'')
\beq\label{correl0}
\dbl\tau_{p_1}(\phi_{\alpha_1}) \tau_{p_2}(\phi_{\alpha_2})\dots 
\tau_{p_k}(\phi_{\alpha_k})\dbr_0 := 
\epsilon^k {\pal^k \log\tau\over \pal
t^{\alpha_1, p_1} \pal t^{\alpha_2, p_2} \dots \pal t^{\alpha_k, p_k}}.
\eeq
They are functions of all the times $t^{\alpha,p}$. The following identities
hold true for these functions
\eqa
&&\dbl\tau_p(\phi_\alpha) \tau_q(\phi_\beta)\tau_r(\phi_\gamma)\dbr_0\nn\\
&&= \dbl\tau_{p-1}(\phi_\alpha) \tau_0(\phi_\nu)\dbr_0 \eta^{\nu\mu}\dbl
\tau_0(\phi_\mu)
\tau_q(\phi_\beta) \tau_r(\phi_\gamma)\dbr_0.\label{rec-rel-0}
\eeqa
In topological sigma-models expanding the function $\log \tau$ at the point of
classical limit one obtains from (\ref{rec-rel-0}) the corresponding identities
for the intersection numbers of the Gromov - Witten Mumford - Morita - Miller 
classes on $\bar {\cal M}_{0,k}$ \cite{Witten2, ruan-tian}.
 
On the {\rm small phase space} $t^{\alpha,p}=0$ for $p>0$ the logarithm of the 
tau-function
coincides \cite{D92} with the potential of the Frobenius manifold
\beq\label{prepo1}
\left. \log\tau\right|_{t^{\alpha,0}= v^\alpha, ~t^{\alpha,p}=0, ~p>0} 
= {1\over \epsilon^2} F(v).
\eeq
The formula
\beq\label{prepo2}
F(v)={1\over 2} \Omega_{1,1;1,1}(v) - v^\alpha \Omega_{\alpha,0;1,1}
+{1\over 2} v^\alpha v^\beta \Omega_{\alpha,0; \beta,0}(v)
\eeq
was used in the derivation of (\ref{prepo1}). This formula, together with
(\ref{productmap3})gives an expression of $F(v)$ via gradients of the functions
$\theta_{\alpha,p}(v)$, $0\leq p\leq 3$.

Note that the quasihomogeneity axiom was not used in the proofs of these statements.
So, the above relations remain valid also for degenerate Frobenius manifolds. 
\end{exam}

\begin{remark} The Principal Hierarchy appears also in the so-called
symplectic field theory of Ya. Eliashberg, A. Givental and H. Hofer \cite{egh}
at the genus zero approximation. For example, the long wave limit of the Toda
lattice essentially appears in their calculation of the genus zero Gromov -
Witten invariants of the projective plane.
We are going to consider the new problems of the theory of integrable
systems inspired by \cite{egh} in a subsequent publication.
\end{remark}

\subsubsection{Complete integrability of the Principal Hierarchy corresponding to a semisimple Frobenius
manifold}\label{sec-3-6-5}\par

We are now to prove {\it completeness} of the system 
$$
H_{\al,p}=\int\theta_{\al,p+1} dx
$$ 
of conservation laws
of the hierarchy (\ref{F-hierarchy}). 

Let us assume that 
$$
I=\int h(v; v_x,\dots,v^{(m)}) dx
$$ 
is a conservation law of the hierarchy,
$$
\{ I, \bar \theta_{\al,p}\}_1=0, ~~\al=1, \dots, n, ~~p=0, 1, 2, \dots .
$$

\begin{lemma}\label{completeness} $I$ is a conservation law of the hierarchy (\ref{F-hierarchy}) {\rm iff}
\beq\label{c.law}
{\pal\over \pal v^{\alpha,k}}{\delta I\over \delta v(x)}=0, ~~\alpha =1, \dots, n, ~~k>0.
\eeq
\end{lemma}

\pf Denote 
$$
W^\alpha=\eta^{\alpha\beta}{\delta I\over \delta v^\beta(x)}.
$$
Then $I$ gives a conservation law for the dispersionless hierarchy 
iff
$$
{\delta \over \delta v^\alpha(x)}\left(W^\gamma\,\partial_x 
{\partial\theta_\beta(v;z)\over\partial v^\gamma}\right)=0,\quad
\alpha, \beta=1,\dots,n
$$
identically in $z$. Using
$$
{\partial^2 \theta_\beta(v;z)\over \partial v^\lambda \partial v^\mu}=
z\, c_{\lambda\mu}^\nu(v){\partial \theta_\beta(v;z)\over \partial v^\nu}
$$
we obtain, after
division by $z$,
\eqa
&&{\partial \left(W^\gamma\,c^{\rho}_{\gamma\sigma}\right)
\over \partial v^\alpha}\,
{\partial\theta_\beta\over\partial v^{\rho}}\,v^\sigma_x
+z\,W^\gamma\,c^{\rho}_{\gamma\sigma}\,c^{\nu}_{\al\rho}
{\partial\theta_\beta\over\partial v^{\nu}}v^\sigma_x
-\partial_x\left(
{\partial W^\gamma\over \partial v^\al_x}
\,c^{\rho}_{\gamma\sigma}\,\frac{\pal\theta_\beta}{\pal v^{\rho}}\,v^\sigma_x
+W^\gamma\,c^{\rho}_{\gamma\al}\,\frac{\pal\theta_\beta}{\pal v^{\rho}}
\right)\nn\\
&& 
+\sum_{k=2}^{2m}(-1)^k\partial_x^k
\left({\partial W^\gamma\over \partial v^{\al,k}}\,
\,c^{\rho}_{\gamma\sigma}\,\frac{\pal\theta_\beta}{\pal v^{\rho}}\,v^\sigma_x
\right)
\nn\\
&&={\partial \left(W^\gamma\,c^{\rho}_{\gamma\sigma}\right)
\over \partial v^\alpha}\,
{\partial\theta_\beta\over\partial v^{\rho}}\,v^\sigma_x
-\partial_x \left(W^\gamma\,c^{\rho}_{\gamma\al}\right)\,
{\partial\theta_\beta\over\partial v^{\rho}}
+\sum_{k=1}^{2m}(-1)^k\partial_x^k
\left({\partial W^\gamma\over \partial v^{\al,k}}\,
\,c^{\rho}_{\gamma\sigma}\,\frac{\pal\theta_\beta}{\pal v^{\rho}}\,
v^\sigma_x\right)\nn
\eeqa
Multiplying by the inverse matrix of $\left(\frac{\pal \theta_\beta(v;z)}{\pal v^\xi}
\right)$
we arrive at a polynomial of degree $2m$ in $z$ with  the coefficient of $z^{2m}$
given by
\eqa
&&\frac{\pal W^\gamma}{\pal v^{\al,2m}}
c^{\rho_1}_{\gamma\sigma}\,c^{\rho_2}_{\sigma_1\rho_1}\,
c^{\rho_3}_{\sigma_2\rho_2}\dots c^{\rho_{2m}}_{\sigma_{2m-1}\rho_{2m-1}}\,
c^{\xi}_{\sigma_{2m}\rho_{2m}}\,v^\sigma_x\,v^{\sigma_1}_x\dots 
v^{\sigma_{2m}}_x\nn\\
&&=\frac{\pal W^\gamma}
{\pal v^{\al,2m}}\sum_i\psi_{i\gamma}\,\psi^\xi_i\,u_{i,x}^{2m+1}\nn
\eeqa
where $u_i$ are the canonical coodinates. The vanishing of the above 
expression yields
$$
\frac{\pal W^\gamma}
{\pal v^{\al,2m}}=0.
$$
In a similar way, we prove inductively that
$$
\frac{\pal W^\gamma}
{\pal v^{\al,k}}=0,\quad k=1,\dots 2m-1.
$$
Therefore $W^\gamma$ does not depend on the $x$-derivatives of 
$v^\al$. 
\epf

\begin{theorem}\label{completeness-theorem}
 Let $I=\int h(v; v_x,\dots,v^{(m)}) dx$ be a conservation law
of the hierarchy (\ref{F-hierarchy}) polynomial in $v^1$. Then
\beq
h(v; v_x,\dots,v^{(m)})=\sum c^{\al,p} \theta_{\alpha,p}(v) + {\rm total ~ derivative}
\eeq
where only finite number of the constant coefficients $c^{\al,p}$ is not equal to zero.
\end{theorem}

\pf This follows from the above lemma and from the lemma \ref{l-2-7-4}.
\epf

\setcounter{equation}{0}
\setcounter{theorem}{0}

\subsection{Quasitrivial bihamiltonian structures}\label{sec-3-7}\par

After having settled the problem of normal forms of the leading term of the
expansion of the Poisson pencil (\ref{mainclass1})-(\ref{mainclass2}) we now
address the problem of construction and classification of the higher order
terms. Recall that we want to classify the Poisson pencils up to the action
of the Miura-type transformations independent of the parameter $\lambda$ of the
pencil. More precisely,

{\bf Definition.} Two Poisson pencils 
\beq\label{2-8-1}
\{ u^i(x), u^j(y)\}_\lambda =\sum_{k=0}^\infty
\epsilon^k \left[ \{u^i(x), u^j(y)\}_2^{[k]}-\lambda \{u^i(x),u^j(y)\}_1^{[k]}
\right]
\eeq
and 
\beq\label{2-8-2}
\{ v^i(x), v^j(y)\}_\lambda =\sum_{k=0}^\infty
\epsilon^k \left[ \{v^i(x), v^j(y)\}_2^{[k]}-\lambda \{v^i(x),v^j(y)\}_1^{[k]}
\right]
\eeq
are called {\it equivalent} if there exists a Miura-type transformation
$$
v^i=\sum_{k=0}^\infty \epsilon^k F^i_k(u; u_x, \dots, u^{(k)}), ~~i=1, \dots, n
$$
transforming (\ref{2-8-1}) to (\ref{2-8-2}) for every $\lambda$. 
The Poisson pencil (\ref{2-8-1}) is called
{\it trivial} if it is equivalent to (\ref{2-8-2}) with 
$\{~,~\}_{1,2}^{[k]}=0$ for $k>0$.

As in Section \ref{sec-2-4} above, 
the infinitesimal description of the space of classes of
equivalence of Poisson pencils with  a given leading term can be
done in terms of  certain cohomology. More precisely, the two Poisson
brackets $\{~,~\}_{1,2}^{[0]}$ induce two anticommuting differentials $\pal_1$
and $\pal_2$ on multivectors,
$$
\pal_i=\left[ ~.~, \{~,~\}_i^{[0]}\right], ~~i=1,2,
$$
$$
\pal_1^2 =\pal_2^2 =\pal_1\pal_2+\pal_2\pal_1=0.
$$
As we already know both the differentials have trivial cohomology. 

\begin{lemma} Let us denote
\eqa
&&H^k({\cal L}(M); \pal_1,\pal_2):= 
{{\rm Ker}}\, \pal_1\pal_2|_{\Lambda^{k-1}} /\left(
{{\rm Im}}\pal_1+{{\rm Im}}\pal_2\right), ~~k>1,
\nn\\
&&
H^1({\cal L}(M); \pal_1,\pal_2):= {{\rm Ker}}\, \pal_1\pal_2 |_{\Lambda^0}
\label{2-8-3}
\\&&
H^0({\cal L}(M); \pal_1,\pal_2):= {{\rm Ker}}\, \pal_1|_{\Lambda^0}\cap 
{{\rm Ker}}\, \pal_2|_{\Lambda^0}.\nn
\eeqa
The zero cohomology coincides with the algebra of common Casimirs for the two
Poisson brackets.
The first cohomology coincides with the space of bihamiltonian vector
fields for the Poisson pencil (\ref{pbht}). The second cohomology coincides with
classes of equivalence of
the infinitesimal deformations of the Poisson pencil modulo infinitesimal
Miura-type transformations.
\end{lemma}

\pf The interpretation of the zero cohomology is straightforward by the
definition.
For a local functional $\bar h$, $h=h(u; u_x, \dots; \epsilon)\in {\cal A}$
the condition $\pal _1\pal_2 \bar h=0$ means that $\pal_2 \bar h\in {{\rm
Ker}}\, \pal_1$. Due to the triviality of the first cohomology of $\pal_1$ the
last condition implies existence of a local functional $\bar f$ such that 
$$
\pal_2 \bar h = \pal_1 \bar f.
$$
That is, the vector field $\pal_2\bar h$ is a bihamiltonian one. Obviously, the converse
statement is also true.

Let us now look at the infinitesimal deformations of the Poisson pencil.
Without loss of generality we may assume that the perturbation of
$\{~,~\}_1^{[0]}$ is trivial, due to triviality of the second cohomology of
$\pal_1$. The infinitesimal deformation of $\{~,~\}_2^{[0]}$ must be annihilated
by $\pal_2$ and also by $\pal_1$, due to the compatibility condition
of the Poisson brackets. So the deformation of the Poisson pencil
must be of the form
\beq\label{2-8-4}
\{~,~\}_1^{[0]}\mapsto \{~,~\}_1^{[0]}+{\cal {O}}(\epsilon^2), ~~
\{~,~\}_2^{[0]}\mapsto \{~,~\}_2^{[0]}+ \epsilon \,\pal_1 X +
{\cal {O}}(\epsilon^2),
~~\pal_2 \pal_1 X=0.
\eeq
This transformation is trivial if it can be generated by another vector field 
$Y$. This means that
$$
\pal_1 Y=0, ~~\pal_2 Y=\pal_1X.
$$
The first of the two equations implies $Y=-\pal_1 \bar a$ for some local functional
$\bar a$. The second one proves existence of another local functional $\bar b$
such that $X = \pal_2 \bar a +\pal_1 \bar b$. \epf

Similar arguments prove the following simple statement.

\begin{theorem} The classes of equivalence of bihamiltonian structures
on the loop space with the given $\{~,~\}_{1,2}^{[0]}$ are in one-to-one
correspondence with classes of equivalence of vector fields
$$
X, ~~X|_{\epsilon=0}=0
$$
satisfying
\beq\label{maur-cart}
\pal_1\left( -\pal_2 X +{1\over 2}[X,\pal_1X]\right)=0
\eeq
modulo shifts along the $\{~,~\}_1^{[0]}$-hamiltonian vector fields
$$
X\mapsto \exp[{\rm ad}_{\pal_1 h}]X
$$
with $\epsilon$ dependent Hamiltonian $h$.
\end{theorem}

We leave the proof of this statement to the reader.

We will call the groups (\ref{2-8-3}) 
{\it the bihamiltonian cohomology} of the pencil
(\ref{pbht}). The calculation of  the bihamiltonian cohomology seem to be a
nontrivial problem. Another problem to be fixed is the one of obstructions
to extension of a given infinitesimal deformation (\ref{2-8-4}) 
to a global one. It is
easy to see from (\ref{maur-cart}) that the first obstruction is the class of equivalence of the
cocycle
\beq\label{obstr}
[X, \pal_1 X]\in  H^3({\cal L}(M); \pal_1,\pal_2).
\eeq
The analysis of this and higher obstructions seems to be an interesting problem
of infinite dimensional Poisson geometry.

\begin{exam} For $n=1$ all deformations upto the fourth order of the Poisson
pencil
\beq\label{pencil-kdv}
\{ u(x), u(y)\}_\lambda =(u-\lambda)\delta'(x-y) +{1\over 2}u_x \delta(x-y)
\eeq
have been classified by P.Lorenzoni \cite{lorenzoni}. They are
parametrized by one arbitrary function $f=f(u)$ of one variable as follows
\eqa
&&
\{ u(x), u(y)\}_\lambda =(u-\lambda)\delta'(x-y) +{1\over 2}u_x \delta(x-y)
\nn\\
&&
+\epsilon^2\left[ -2f\delta'''(x-y)-3\pal_x f\delta''(x-y)-
\pal_x^2 f\delta'(x-y)\right]
\nn\\
&&
+\epsilon^4\left[ 4 g\delta^V(x-y) + 10 \pal_x g\delta^{IV}(x-y)
+(20\, \pal_x^2 g-8\,hu_{xx})\delta'''(x-y)\right.\nn\\
&&\left.
+(20 \pal_x g -12 \pal_x(hu_{xx}))\delta''(x-y)
+(6 \pal_x^4 g-4\, \pal_x^2(hu_{xx}))\delta'(x-y)\right]\nn\\
&& +{\cal{O}}(\epsilon^5).\label{paolo}
\eeqa
In the r.h.s. of this fomula
$$
u=u(x), ~u_x=u_x(x), ~ u_{xx}=u_{xx}(x), ~~g=f f', ~~
h=f f'' + {f'}^2, ~~ f=f(u(x)).
$$
In particular, the obstruction (\ref{obstr}) is trivial for an arbitrary
infinitesimal deformation of the order $\epsilon^2$.
All the above Poisson pencils are inequivalent for different $f(u)$. In
particular, for $f(u)=c$ one obtains the KdV Poisson pencil.
$$
\{ u(x), u(y)\}_\lambda = [u(x)-\lambda]\, \delta'(x-y) +{1\over 2} u'
\delta(x-y) + c \,\epsilon^2 \delta'''(x-y).
$$
\end{exam}

We will now impose the main restriction onto the class of Poisson pencils
that will allow us to get rid of the above unpleasant cohomological problems.
Let us extend the class of Miura-type transformations.

{\bf Definition.} The transformations of the form
\beq\label{qtt-2-8}
u^i\mapsto v^i = u^i + \sum_{k=1}^\infty \epsilon^k F^i_k(u; u_x, \dots,
u^{(n_k)}), ~~i=1, \dots, n
\eeq
where the coefficients $F^i_k$ are quasihomogeneous of the degree $k$ {\it
rational} functions in the derivatives $u_x$, \dots, $u^{(n_k)}$ will be called
{\it quasi-Miura transformations}.  The Poisson pencil (\ref{mainclass1})--
(\ref{mainclass3}) 
is called {\it
quasitrivial} if there exists a quasi-Miura transformation reducing the pencil
to its leading term (\ref{pbht}).

We emphasize that the coefficients of the Poisson pencil are still to be
polynomials in the derivatives. All the denominators must disappear after
the transformation.

\begin{remark}
As it was proved in \cite{lorenzoni}, the deformation (\ref{paolo}) 
is quasitrivial
up to the order 4. This suggests that {\bf all} Poisson pencils of the form
(\ref{mainclass2}) corresponding to a semisimple Frobenius structure in
$\{~,~\}_\lambda^{[0]}$ could be quasitrivial. To our opinion, the problem of
quasitriviality 
deserves further investigation.
\end{remark}

In the next section we will show quasitriviality of the KdV hierarchy. Even in
this simplest example the quasitriviality is something unobserved before.
 
\setcounter{equation}{0}
\setcounter{theorem}{0}

\subsection{Quasitriviality of KdV and genus expansion in topological gravity}
\label{sec-3-8}\par

To construct a transformation
that maps any solution $v$ of the Riemann hierarchy (\ref{burgers}) 
to a solution
$u$ of the KdV hierarchy we will proceed following \cite{dm}. 
Every solution $v$ of the Riemann
hierarchy can be represented in the standard implicit form (\ref{tsarev}), i.e.,
\beq\label{eqno1}
x+\tilde t_0 + \tilde t_1 v + \tilde t_2 {v^2\over 2}
+ \tilde t_3 {v^3\over 6} +\dots =0.
\eeq
Here
$$
\tilde t_p = t_p - c_p
$$
where the constants $c_p$ correspond to the choice of the solution (in this section we will systematically 
suppress the explicit dependence of the coefficients $c_p$ on $\epsilon$). 
Representing the equation (\ref{eqno1}) in the variational form (\ref{tsarev1})
we obtain
\beq\label{eq1}
\Phi_{x,{\bf t},c}'(v)=0, ~~
\Phi_{x,{\bf t},c}'(v)=(x+\tilde t_0) v + \tilde t_1 {v^2\over 2!} +\dots.
\eeq
Let $\bar h_p =\int h_p (u; u_x, \dots)dx$ be the Hamiltonians of the KdV
hierarchy normalized as in (\ref{h-kdv}), i.e. 
$$
h_p = {u^{p+2}\over (p+2)!} +\epsilon^2 ({\rm terms ~ with ~ derivatives}), ~~p\geq -1.
$$
Let us construct a functional depending on the same parameters $x$, ${\bf t}$, $c$ 
\beq\label{functional-kdv}
I_{x,{\bf t},c}[u] = \int \left( (x+\tilde t_0) u +\sum_{p>0}
 \tilde t_p h_{p-1}(u;
 u_x, \dots)\right).
\eeq
At the moment we do not care about the precise meaning of the integral. We will be interested only in the Euler - Lagrange equation 
\eqa
&&{\delta\over \delta u(x)}I_{x,{\bf t},c}[u] =
\sum_k (-1)^k \partial_x^k {\pal\over \pal u^{(k)}} 
\left( (x+\tilde t_0) u +\sum_{p>0}
 \tilde t_p h_{p-1}(u;
 u_x, \dots)\right)\nn\\
&&
=x+\tilde t_0 +\sum_{p>0} \tilde t_p {\delta \bar h_{p-1}\over \delta u(x)}
=0.\label{el}
\eeqa 
The first few terms of the Euler - Lagrange equation (\ref{el})
read
$$
x+\tilde t_0 + \tilde t_1 u + \tilde t_2\left( {u^2\over 2} -\epsilon^2{u''\over 6}
\right) 
+ \tilde t_3 \left( {u^3\over 6} - 
{\epsilon^2\over 12}({u'}^2 + 2 u\, u'') + \epsilon^4 {u^{IV}\over 60}\right)+\dots=0.
$$
Truncating $t_p=0$ for $p\geq N$ (and assuming that only finite number of the constants $c_p$ is distinct from zero) we obtain an ODE for the function $u=u(x)$ depending on the times $t_0$, \dots, $t_{N-1}$ and on the constants $c=(c_0, c_1, \dots)$ as on the parameters.

\begin{lemma}\label{l-2-9-1}
(cf. \cite{moore}) The space of solutions to the 
Euler - Lagrange equation (\ref{el}) is invariant w.r.t. the flows of the KdV hierarchy.
\end{lemma}

\pf Let $u_0(x)$ be a solution to the differential equation (\ref{el})
with $t_p=t_p^0$, $p=0$, $1$, \dots. We are to prove that the solution
to the Cauchy problem for the KdV hierarchy with the initial data
$$
u(x, {\bf t})|_{t_p = t_p^0, ~p=0, 1, \dots} =u_0(x)
$$
will satisfy the same ODE (\ref{el}).
\epf

Let $v=v(x, {\bf t}, c)$ be the solution to the Riemann hierarchy
determined by (\ref{eq1}) such that $v'(0,{\bf 0},c)\neq 0$. (The solution 
depends on the choice of a simple root $v_0$ of the polynomial
$\sum c_p {v^p_0\over p!} =0$. It will be understood that such a choice has already been done.)

\begin{lemma} There exists a unique 
solution to (\ref{el}) in the form of power series in $\epsilon^2$
\beq\label{expansion}
u= v+\epsilon^2 u^{[1]} + \epsilon^4 u^{[2]} +\dots.
\eeq
\end{lemma}

\pf 
We plug (\ref{expansion}) into the equation (\ref{el})) and compute recursively the terms of the
expansion. For example, for the first correction we obtain
$$
u^{[1]} =-{1\over 24}{2\,v'' \left( \tilde t_2 + \tilde t_3 v +\tilde t_4 {v^2\over 2}
+\dots\right) +{v'}^2\left( \tilde t_3 +\tilde t_4 v +
\dots\right)\over
\tilde t_1 + \tilde t_2 v +\tilde t_3 {v^2\over 2} + \dots}.
$$
\epf

\begin{cor} The solution (\ref{expansion}) to (\ref{el}) satisfies equations of the KdV hierarchy.
\end{cor}

Thus we obtain a map
\beq\label{correspond}
{\rm the ~stationary ~point ~(\ref{eq1})}\mapsto {\rm 
the ~stationary ~function~
(\ref{expansion})~ of ~ (\ref{el})}
\eeq
transforming solutions of the Riemann hierarchy to the solutions to the KdV hierarchy. We will now show that this is a quasitriviality transformation.

First we will prove

\begin{lemma}\label{l-2-9-4} 
There exist universal (i.e., independent on the choice of the 
solution $v$ to the Riemann hierarchy) polynomials 
$P^{[2k]}(v', v'', \dots v^{(3k)})$
quasihomogeneous of the degree $6\,k-2$ such that the transformation (\ref{correspond}) is given by
\beq \label{correspond1}
v\mapsto u = v + 
\sum_{k\geq 1} \epsilon^{2\,k} (v')^{2-4\,k} P^{[2k]} (v', v'', \dots, v^{(3\, k)}). 
\eeq
\end{lemma}

\pf
It is technically convenient to return to the original normalization 
of Example \ref{e-2-2-4} for the KdV
hierarchy 
$$
{\pal u\over \pal t_k} =\pal _x {\delta I_k\over \delta u(x)}
$$
where
the
generating function of the densities of the KdV integrals is to be determined from the Riccati
equation
$$
i\epsilon\,\chi'-\chi^2=u-\lambda,
$$
$$
\chi=k+\sum_{m=1}^\infty {\chi_m\over k^m}, ~~k=\sqrt{\lambda},
$$
$$
I_k = -4\, \int \chi_{2k+3}dx.
$$
We can rewrite the Euler - Lagrange equation (\ref{el})
in the following form. Introduce the series
$$
{\bf t}(\lambda) := \tilde t_0 + \sum_{k=1}^\infty {2^k\over (2k-1)!!} \tilde
t_k \lambda^k.
$$
Let us also introduce the linear operator $Res$ acting on the series in
inverse powers of $k=\sqrt{\lambda}$ by
$$
Res \, f : = res_{k=\infty} {\bf t}(\lambda) f \,dk.
$$
Then the Euler - Lagrange equation reads
\beq\label{el-res}
Res {\delta \int \chi\,dx\over \delta u(x)} =0.
\eeq
In a similar way, the variational equation (\ref{eq1}) can be written as
$$
Res {1\over \sqrt{\lambda-v}}=0.
$$
We will now expand the variational derivatives in powers of $\epsilon^2$.
Using the formula
$$
{\delta \int\chi\, dx\over \delta u(x)} = -{1\over 2 \chi_R},
$$
where $\chi_R$ is the real part of $\chi$, $\chi=\chi_R + i\chi_I$, $\chi_I=
{1\over 2}
{\chi_R'\over \chi_R}$, we rewrite (\ref{el-res}) as
\beq\label{el-res1}
Res\, {1\over \chi_R} =0.
\eeq
Using differential equation for $f:=1/\chi_R$,
$$
f^2 + {\epsilon^2\over \lambda-u} \left[ {1\over 2} f'' f -{1\over 4}
{f'}^2\right] = {1\over \lambda - u}
$$
(a consequence of Riccati)
we can expand $1/\chi_R$ in the series of the form
$$
{1\over \chi_R}=
\sigma + {1\over 8}\epsilon^2\,\left( \sigma\,{{\sigma'}}^2 - 2\,\sigma^2\,{\sigma''} \
\right) 
$$
$$
+ {1\over 128}\epsilon^4\,\left( 3\,\sigma\,{{\sigma'}}^4 - 
       12\,\sigma^2\,{{\sigma'}}^2\,{\sigma''} + 
       12\,\sigma^3\,{{\sigma''}}^2 + 
       16\,\sigma^3\,{\sigma'}\,{\sigma'''} + 8\,\sigma^4\,{\sigma^{IV}} \
\right) +O(\epsilon^6)
$$
where
$$
\sigma ={1\over \sqrt{\lambda-u}}.
$$ 
Let us now compute the first correction $u^{[1]}$ in the
expansion (\ref{expansion}). Within the order $\epsilon^2$ the equation
(\ref{el-res1}) reads
\beq\label{eqno4}
Res\left[  \sigma 
+ {1\over 8}\epsilon^2\,\left( \sigma\,{{\sigma'}}^2 
- 2\,\sigma^2\,{\sigma''} \right)
\right]=O(\epsilon^4).
\eeq
Denote
$$
\sigma_0:={1\over \sqrt{\lambda-v}}.
$$
We must expand the above equation and retain the linear in $\epsilon^2$ terms.
Substituting
$$
u=v+\epsilon^2 u^{[1]} +\dots
$$
in $\sigma$ we obtain
$$
\sigma = \sigma_0 +{\epsilon^2\over 2} u^{[1]} \sigma_0^3+\dots
$$
In the second term in (\ref{eqno4}) 
we may replace $\sigma\to\sigma_0$. Next, we are to
observe the following simple rules for differentiating $\sigma_0$:
$$
\sigma_0'={1\over 2} \sigma_0^3,~~\sigma_0''={1\over 2} v''\sigma_0^3 +{3\over
4}{v'}^2 \sigma_0^5, \dots.
$$
So, the equation (\ref{eqno4}) can be rewritten as
$$
Res \left[ {1\over 2} u^{[1]} \sigma_0^3-{1\over 8} v''\sigma_0^5+{5\over 32}
{v'}^2 \sigma_0^7\right]=0.
$$
It remains to calculate the residues of the form
$$
Res \,\sigma_0^{2k+1}={2^k\over (2k-1)!!}Q_k
$$
where the rational functions $Q_k$ in the derivatives are defined recursively
$$
Q_{k+1} ={1\over v'} Q'_k, ~~ Q_1 =-{1\over v'}.
$$
To prove the last formula it suffices to observe that
$$
Res {d^k\sigma_0\over dv^k} = \tilde t_k +\tilde t_{k+1} v+\dots = Q_k,
$$
and to compute
$$
{d^k \sigma_0\over dv^k} = {(2k-1)!!\over 2^k} \sigma_0^{2k+1}.
$$
Finally we obtain the needed formula in the form
$$  
u^{[1]} ={1\over Q_1} \left[ {1\over 6} v'' Q_2 -{1\over 12} {v'}^2 Q_3\right]=-{1\over 12} (\log v')''.
$$
It is clear how to proceed with higher corrections. We want to emphasize that
the expressions
$$
u^{[k]} ={P^{[2k]} (v', v'', \dots, v^{(3k-2)})\over {v'}^{4k-2}}, ~~k\geq 1
$$
do not depend on $v$ explicitly.
\epf

\begin{cor} The correspondence (\ref{correspond}) 
$$
\left\{\matrix{{\rm solutions ~ to ~(\ref{eq1})}\cr
v(x,{\bf t},\epsilon)=v_0(x, {\bf t})+ \epsilon v_1(x,{\bf t})+\dots\cr}\right\}
 ~
\mapsto ~\left\{\matrix{{\rm solution ~ to ~(\ref{el}) }\cr
u(x,{\bf t},\epsilon)=u_0(x, {\bf t})+ \epsilon u_1(x,{\bf t})+\dots\cr}\right\}
$$
is a quasi-Miura transformation
\eqa\label{kdv-miura}
&&u=F(v; v_x, v_{xx}, \dots;\epsilon)\nn\\
&&= v -{\epsilon^2\over 12} (\log v')'' + 
\epsilon^4 \left[ {v^{IV}\over 288 {v'}^2} -{7v'' v'''\over 480 {v'}^3}
+{{v''}^3\over 90 {v'}^4}\right]''+O(\epsilon^6).
\eeqa
\end{cor}

We are now to prove that the quasi-Miura transformation is one-to-one.
Loosely speaking we want to prove that an arbitrary monotone function $u(x,\epsilon)$
satisfies an ODE (\ref{el}) possibly of infinite order with ${\bf t} =0$ and appropriate coefficients $c_p$ that may depend on $\epsilon$. More precisely,

\begin{lemma}\label{l-2-9-6} 
Let $u=u(x,\epsilon)\in {\mathbb C}[[x,\epsilon]]$ be an arbitrary 
power series satisfying $u_x(0,0)\neq 0$. Then there exist unique power series
$$
c_p(\epsilon) = c_p^{(0)} + \epsilon\, c_p^{(1)}+\epsilon^2  c_p^{(2)}+\dots, ~~p=0, 1, 2, \dots
$$
such that the following identity in the ring $ {\mathbb C}[[x,\epsilon]]$
holds true
\beq\label{pinfinity}
x= c_0(\epsilon) + 
\sum_{p> 0}c_p(\epsilon) h_{p-2}(u;  u_x,\dots,\epsilon^{2p-2} u^{(2 p -2)}).
\eeq
\end{lemma}

\pf The leading coefficients $c^{(0)}_p$ must satisfy
$$
\sum_p c^{(0)}_p {u_0^p(x)\over p!}=x
$$
where $u_0(x) =u(x,0)$. Therefore they are equal to the derivatives of the inverse function
$$
c^{(0)}_p = \left.{d^p x\over d u_0^p}\right|_{u_0=0}, ~~p=0, 1, 2, \dots .
$$
Proceeding by induction we obtain
$$
\sum_p c^{(k)}_p {u_0^p(x)\over p!}=f_k(x)
$$
where $f_k(x)$ is a polynomial in 
$$
c^{(i)}_q, ~~i=0, \dots, k-1
$$
and in  
$$
u_j(x)=\left.
{d^j u(x,\epsilon)\over d\epsilon^j}\right|_{\epsilon=0}, ~~j=0, \dots, k
$$
and their derivatives in $x$. Therefore
$$
c^{(k)}_p =\left. {d^p f_k(x)\over d u_0^p}
\right|_{u_0=0}, ~~p=0, 1, 2, \dots .
$$
\epf

\begin{cor} The transformation (\ref{correspond}) establishes a one-to-one
correspondence between solutions $v(x,{\bf t}, \epsilon)$ to the Riemann hierarchy satisfying $v_x(0,0,0)\neq 0$ and solutions $u(x,{\bf t}, \epsilon)$ to the
KdV hierarchy satisfying $u_x(0,0,0)\neq 0$.
\end{cor}

\pf Let $c_p(\epsilon)$ be the coefficients determined according to 
Lemma \ref{l-2-9-6} by $u(x,0,\epsilon)$. According to 
Lemma \ref{l-2-9-1} the solution 
$u(x,{\bf t},\epsilon)$ to the KdV hierarchy satisfies
$$
x+t_0-c_0(\epsilon)+\sum_{p>0} (t_p-c_p(\epsilon)) {\delta \bar h_{p-1}\over \delta u(x)}=0
$$
identically in ${\bf t}$. Let $v=v(x, {\bf t}, c(\epsilon))$ be the solution (\ref{eq1}) to the Riemann hierarchy determined by these coefficients. By the construction the quasi-Miura transformation maps $v$ to $u(x,{\bf t}, \epsilon)$.
\epf

We will now prove

\begin{lemma}\label{l-2-9-8} 
The quasi-Miura (\ref{correspond1}) 
$v\mapsto u=F(v; v_x, \dots;\epsilon)$ transforms the vector fields of the Riemann hierarchy to those of the KdV hierarchy.
\end{lemma}

\pf Let
$$
\hat K_j(u; u_x, \dots; \epsilon)=\sum_m \epsilon^{2m} u_x^{-p_m} 
\hat K_j^{[2m]}(u; u_x, \dots, u^{(q_m)})
=\left( \sum_s {\pal F\over \pal v^{(s)}} \pal_x^{s+1}\right) {v^{j+1}\over (j+1)!}
$$
be the result of application of the quasi-Miura transform to the flows of the 
Riemann hierarchy. Here 
$K_j^{[2m]}(u; u_x, \dots, u^{(n_m)})$ are some 
polynomials in the derivatives.
The precise values of the positive numbers $p_m$ and $q_m$ (that also depend on $j$) is not important. 
Denote 
$$
K_j(u; u_x, \dots, u^{2j+1};\epsilon)=\sum_{m=0}^{2j} 
\epsilon^{2m} K_j^{[2m]}(u; u_x, \dots, u^{(2m+1)})
$$ 
the r.h.s.
of the $j$-th equation of the KdV hierarchy.
According to the Lemma \ref{l-2-9-6} the identity
$$
\sum_{m=0}^{2j} 
\epsilon^{2m} K_j^{[2m]}(u; u_x, \dots, u^{(2m+1)})=
\sum_m \epsilon^{2m} u_x^{-p_m} 
\hat K_j^{[2m]}(u; u_x, \dots, u^{(q_m)})
$$
holds true for an arbitrary monotone solution $u(x, {\bf t}, \epsilon) $
to the KdV hierarchy. From this it easily follows that 
$$
u_x^{-p_m} 
\hat K_j^{[2m]}(u; u_x, \dots, u^{(q_m)})=
 K_j^{[2m]}(u; u_x, \dots, u^{(2m+1)}).
$$
\epf

Denote
$$
\hat h_k(v; v_x, \dots; \epsilon) = h_k(u; \epsilon\, u_x, \dots,\epsilon^{2k+2} 
u^{(2k+2)}), ~~k=-1, 0, 1, \dots, 
$$
the functions in the derivatives obtained from the Hamiltonian densities of KdV
by the inverse to the quasi-Miura transformation (\ref{correspond1}).

\begin{lemma}\label{l-2-9-9} 
$$
\hat h_k= {v^{k+2}\over (k+2)!} + {\rm total ~ derivative}.
$$
\end{lemma}

\pf Applying the inverse to the quasi-Miura transformation to the infinitesimal form of the conservation law
$$
{\partial h_k(u, \epsilon\, u_x, \dots)\over \partial t_j} = 
{\partial \Omega_{k+1,j}(u, \epsilon\, u_x, \dots)\over \pal x}
$$
(here $\Omega_{k+1,j}(u, \dots)$ is the density of the flux of the conserved 
quantity along the $j$-th time defined in (\ref{stress})) 
we obtain, according to Lemma \ref{l-2-9-4}
$$
{\partial\hat h_k(v, v_x, \dots; \epsilon)\over \partial t_j} = 
{\partial \hat \Omega_{k+1,j}(v, v_x, \dots;\epsilon)\over \pal x}.
$$
In the last equation the time derivative is taken w.r.t. the Riemann hierarchy;
the functions $\hat \Omega_{k+1,j}(v, v_x, \dots;\epsilon)$ are obtained from 
$\Omega_{k+1,j}(u, \epsilon\, u_x, \dots)$ by the same inverse quasi-Miura. Therefore
$\hat h_k(v, v_x, \dots; \epsilon)$ is the density of a conservation law 
for the Riemann hierarchy. Due to the completeness theorem 
\ref{completeness-theorem}
it must coincide with the standard density $v^{k+2}/(k+2)!$ up to a total derivative in $x$.
\epf

We are now ready to prove the main result of this section.

\begin{theorem} \label{t-2-9-10}
The (inverse to) the quasi-Miura transformation
(\ref{correspond1})
transforms the Magri Poisson pencil (\ref{kdv-ham2}),
(\ref{kdv-ham2}) to the Poisson pencil (\ref{pencil-kdv})
for the Riemann hierarchy.
\end{theorem}

\pf Applying the inverse quasi-Miura to the first Poisson bracket of the KdV
we obtain a Poisson bracket
$$
\{ v(x), v(y) \}^{\hat{}}_1 =
\sum \epsilon^{2m}A_{2m, s}(v; v_x, \dots, ) \delta^{(2m-s+1)}(x-y).
$$
From Lemma \ref{l-2-9-9} it follows that 
$$
{\delta \int \hat h_k\over \delta v(x)} = {v^{k+1}\over (k+1)!}, ~~k=0, 1, \dots .
$$
From Lemma \ref{l-2-9-8} it follows that
$$
\sum \epsilon^{2m}A_{2m, s}(v; v_x, \dots, ) \pal_x^{2m-s+1} 
{v^{k+1}\over (k+1)!} = \pal_x {v^{k+1}\over (k+1)!}
$$
for all non-negative $k$. Multiplying the last equation by $z^{k+1}$, where $z$ is an indeterminate, and summing in $k$ we obtain
$$
\sum \epsilon^{2m}A_{2m, s}(v; v_x, \dots, ) \pal_x^{2m-s+1} e^{z\, v(x,\epsilon)} =\pal_x e^{z\, v(x, \epsilon)}
$$
for all $z$ and for an arbitrary function $v(x,\epsilon)$. From this it easily
follows that $\{~, ~\}^{\hat{}}_1 = \{~, ~\}_1$.

Applying the inverse quasi-Miura to the second Poisson bracket for KdV
we obtain a Poisson bracket
$$
\{ v(x), v(y)\}^{\hat{}}_2 =
\sum \epsilon^{2m}B_{2m, s}(v; v_x, \dots, ) \delta^{(2m-s+1)}(x-y).
$$
that must satisfy the recursion relation (\ref{recur-2-7}) 
for the Riemann hierarchy
$$
\sum \epsilon^{2m}B_{2m, s}(v; v_x, \dots, ) \pal_x^{2m-s+1}{v^{k}\over k!} 
=(k+{1\over 2}) \pal_x {v^{k+1}\over (k+1)!}.
$$
Multiplying by $z^{k+1}$ and summing in $k$ we obtain the identity
$$
\sum \epsilon^{2m}B_{2m, s}(v; v_x, \dots, ) \pal_x^{2m-s+1} 
e^{z\, v} =\left( z\, v +{1\over 2}\right) v_x e^{z\, v}
$$
valid for any $z$ and for an arbitrary function $v=v(x,\epsilon)$.
This proves that $\{~, ~\}^{\hat{}}_2 = \{~,~\}_2$. The Theorem is proved.
\epf

Our approach can easily be extended to prove quasitriviality of the Gelfand -
Dickey hierarchy (also called nKdV). We will study in a separate publication
the problem of quasitriviality of other hierarchies of integrable 1+1 PDEs, in
particular 
of the Drinfeld - Sokolov hierarchies of $D$ and
$E$ type and of Toda lattice.

We will now prove that, in addition to Lemma \ref{l-2-9-4}, the following
statement.

\begin{lemma} There exists a function
\beq\label{deltaf}
\Delta f = \Delta f(v', v'', \dots; \epsilon^2) = \sum_{k=1}^\infty
\epsilon^{2k-2}\Delta f^{[k]}(v', \dots, v^{(3k-2)}) 
\eeq
where
$$
\Delta f^{[1]}(v') =-{1\over 12} \log v'
$$
and $\Delta f^{[k]}(v', \dots, v^{(3k-2)}) $ is a quasi-homogeneous function in
the derivatives of the degree $2k-2$ such that the correspondence
(\ref{kdv-miura}) is represented as
\beq\label{kdv-miura2}
v\mapsto u=v+\epsilon^2 \pal_x^2 \Delta f(v', v'', \dots, ;\epsilon^2).
\eeq
\end{lemma}
\pf (cf. the proof of Theorem \ref{t-2-10-1} below). We already know from
Lemma \ref{l-2-9-9} that
$$
h_{p-1}(u; u', \dots, u^{(p-1)}; \epsilon) ={v^{p+1}\over (p+1)!} +\epsilon
\pal_x g_{p-1}(v, v', \dots;\epsilon)
$$
for some functions $g_k(v, v', \dots; \epsilon)$. Using the tau-symmetry
$$
{\pal h_{p-1}\over \pal t_q} = {\pal h_{q-1}\over \pal t_p}
$$
$$
{\pal\over \pal t_q} {v^{p+1}\over (p+1)!} = 
{\pal\over \pal t_p} {v^{q+1}\over (q+1)!}
$$
of the KdV hierarchy and of the Riemann hierarchy we obtain
$$
{\pal\over \pal x} \left[ {\pal g_{p-1}\over \pal t_q} - 
{\pal g_{q-1}\over \pal t_p}\right]=0.
$$
This implies existence of a function
$$
\Delta f=\sum_{k=1}^\infty
\epsilon^{2k-2}\Delta f^{[k]}(v, v', \dots, v^{(3k-2)}) 
$$
such that
$$
g_{p-1}(v, v', \dots; \epsilon) = \epsilon {\pal \Delta f\over \pal t_p}
=\epsilon \,\sum {\pal \Delta f\over \pal v^{(i)}} \left( {v^{p+1}\over
(p+1)!}\right)^{(i+1)}.
$$
In particular, the quasi-Miura transformation itself reads
$$
u=v + \epsilon \pal_x g_0 = v + \epsilon^2 \pal_x^2 \Delta f
$$
where we may assume, due to quasihomogeneity of the terms $u^{[k]}$ in the
derivatives that
$$
\Delta f =\Delta f(v', v'', \dots ;\epsilon^2)
$$
does not depend explicitly on $v$. Lemma is proved. \epf

\begin{exam} The topological solution to the Riemann hierarchy is determined by
the equation
\beq\label{g0}
v= t_0 + t_1 v + t_2 {v^2\over 2} + t_3 {v^3\over 6}+\dots
\eeq
(we omit $x$ identifying $x$ with $t_0$). The tau-function of this solution
\eqa
&&
\log\tau_0 = {1\over \epsilon^2} \left( \frac{t_0^3}{6} + \frac{t_0^3\,t_1}{6} + 
  \frac{t_0^3\,t_1^2}{6} + \frac{t_0^3\,t_1^3}{6} + 
  \frac{t_0^3\,t_1^4}{6} + \frac{t_0^4\,t_2}{24}+ 
  \frac{t_0^4\,t_1\,t_2}{8}\right.\nn\\
&&\quad\left. + 
  \frac{t_0^4\,t_1^2\,t_2}{4} + 
  \frac{t_0^5\,t_2^2}{40} + \frac{t_0^5\,t_3}{120} + 
  \frac{t_0^5\,t_1\,t_3}{30} + \frac{t_0^6\,t_4}{720}+\dots\right).
\nn\\
\eeqa
Applying the quasi-Miura transformation (\ref{kdv-miura}) we obtain, after
changing the normalization
$$
\epsilon^2 \mapsto -{\epsilon^2\over 2}
$$
a solution to the KdV hierarchy with the tau-function (\ref{tau-kw}).
We will show below in Section \ref{sec-3-10-4}
this series coincides with the Witten -
Kontsevich
generating function of the Mumford - Morita
- Miller intersection numbers on the moduli spaces $\bar{\cal M}_{g,n}$ of all
genera
\beq\label{g-exp}
{\cal F} = \sum_{g=0}^\infty \epsilon^{2g-2} {\cal F}_g
\eeq
where 
\beq
{\cal F}_g =\sum_{n=1}^\infty {1\over n!} t_{p_1} \dots t_{p_n} <\phi_{p_1}\dots
\phi_{p_n}>_g
\eeq
\beq
<\phi_{p_1} \dots \phi_{p_2}>_g =\int_{\bar {\cal M}_{g,n}} c_1^{p_1} ({\cal
L}_1)\wedge \dots \wedge c_1^{p_n} ({\cal L}_n)
\eeq
where ${\cal L}_i$ is the tautological line bundle over the moduli space 
$\bar {\cal M}_{g,n}$ of stable algebraic curves $C$ of genus $g$ with $n$ punctures
$x_1$, \dots, $x_n$, i.e., the fiber of ${\cal L}_i$ coincides with the
cotangent line $T^*_{x_i} C$. 
\end{exam}

The idea to express the terms ${\cal F}_1$, ${\cal F}_2$, \dots, of the genus
expansion (\ref{g-exp}) as functions of $v'$, $v''$, \dots where $v=v(t)$
is the solution (\ref{g0}) is due to Dijkgraaf and Witten \cite{DW}.
This idea proved to be fruitful also in other models of 2D topological field
theory \cite{eguchi1, eguchi2, eguchi3, zograf, manin-zograf, givental}.

\begin{exam} Applying (\ref{kdv-miura}) to the monotone at $x=0$
function $v=v(x)$
\beq\label{zog}
x = \sqrt{v} J_1(2\sqrt{v}) =\sum_{m=0}^\infty (-1)^m {v^{m+1}\over m! (m+1)!}
\eeq
one obtains 
\beq\label{zog1}
u(x, \epsilon)=\pi^6 \sum_{g=0}^\infty \left( {\epsilon\over \pi^3} \right)^{2\,
g} \sum_n Vol({\cal M}_{g,n}) \left({x\over \pi^2}\right)^n
\eeq
where $Vol({\cal M}_{g,n})$ is the Weil - Petersson volume of the moduli space
of punctured algebraic curves. This is a reformulation of the result of P.
Zograf \cite{zograf} (see also \cite{manin-zograf}).

\end{exam}

\setcounter{equation}{0}
\setcounter{theorem}{0}

\subsection{Properties of quasitrivial Poisson pencils}\par

Let 
$$
\{ u^\alpha(x), u^\beta(y)\}_\lambda = \sum_{k\geq 0}\epsilon^k\{ u^\alpha(x), u^\beta(y)\}_\lambda^{[k]}
$$
be a quasitrivial Poisson pencil written in the normal coordinates with the leading term
$$
\{u^\alpha(x),u^\beta(y)\}_\lambda^{[0]}
=\left( g^{\alpha\beta}(u(x))-\lambda\, \eta^{\alpha\beta}\right)\delta'(x-y)
+\Gamma_\gamma^{\alpha\beta}(u)u^\gamma_x \delta(x-y)
$$
determined by a $n$-dimensional rigid semisimple Frobenius manifold $M$ (see the formula (\ref{normalfm}) above). Let
\beq\label{qt}
u_\alpha =v_\alpha +\sum_{k>0} \epsilon^k F_\alpha^{[k]}(v; v_x, \dots, 
v^{(n_k+2)})
\eeq
be the quasitriviality transformation for the pencil:
$$
\sum_{k\geq 0}\epsilon^k\{ u^\alpha(x), u^\beta(y)\}_\lambda^{[k]}
=\left( g^{\alpha\beta}(v(x))-\lambda\, \eta^{\alpha\beta}\right)\delta'(x-y)
+\Gamma_\gamma^{\alpha\beta}(v)\,v^\gamma_x \,\delta(x-y).
$$
Here we lower the indices as usual by means of the constant matrix $\eta_{\alpha\beta}$,\newline $ F_\alpha^{[k]}(v; v_x, \dots, 
v^{(n_k)})$ are some functions rational in the derivatives.   

Our first result is

\begin{theorem} \label{t-2-10-1} Let the quasi-Miura transformation (\ref{qt}) depend polynomially on $v^1$. Then there exists a function
\beq\label{fe}
{\cal F}(v; v_x, \dots, ;\epsilon) =\sum_{k>0} \epsilon^k {\cal F}^{[k]}(v; v_x, \dots, v^{(n_k)})
\eeq
such that the quasitriviality has the form
\beq\label{qt1}
u_\alpha = v_\alpha  + \pal_x \pal_{t^{\al,0}} {\cal F}(v; v_x, \dots, ;\epsilon).
\eeq
Moreover, the tau-structure for the pencil $\{~,~\}_\lambda$  written in the coordinates $v$ has the form, up to an equivalence (\ref{2-6-26}),
(\ref{2-6-27}),
\beq\label{qt2}
h_{\alpha,p}(v; v_x, \dots ;\epsilon) = \theta_{\al,p}(v) +
\pal_x \pal_{t^{\al,p}} {\cal F}(v; v_x, \dots, ;\epsilon).
\eeq
\end{theorem}

We recall (see Section \ref{sec-3-6-2}
above) that polynomiality in $v^1$ means that
every coefficient $ F_\alpha^{[k]}(v; v_x, \dots, 
v^{(n_k)})$ is a polynomial in $v^1$ of the degree that may depend on $k$.

\pf By the definition of the normal coordinates $\bar u_\al$ is a Casimir
of $\{~,~\}_1$. Since $\{~,~\}_1$ is obtained from  $\{~,~\}_1^{[0]}$ 
by the change of coordinates (\ref{qt}), and $\bar v_\al$ is a Casimir of
$\{~,~\}_1^{[0]}$, it follows that
$\bar v_\al$
is also a  Casimir of $\{~,~\}_1$. Hence the difference
$u_\alpha-v_\alpha=O(\epsilon)$ is a conserved quantity for the Principal
Hierarchy. Due to Lemma \ref{completeness}
$u_\alpha-v_\alpha$ must be a total derivative (polynomiality in the derivatives
assumption can be eliminated by considering arbitrary functions in the
derivatives). Hence the quasitriviality transformation must have the form
$$
u_\alpha =v_\alpha +\pal_x f_{\alpha,0}(v; v_x, \dots;\epsilon)
$$
for some function $f_{\alpha,0}(v; v_x, \dots;\epsilon)$.

Let 
$h_{\al,p}$ be the densities of the commuting Hamiltonians corresponding to a choice of a tau-structure for the Poisson pencil
$\{~,~\}_\lambda$ satisfying a recursion relation
$$
\{~.~, \bar h_{\alpha,p}\}_2 = \sum_{q=0}^p A_{\alpha,p}^{\beta,q}
\{~.~,\bar h_{\beta,q+1}\}_1
$$
with some constant coefficients $A_{\alpha,p}^{\beta,q}$. We have
$h_{\al,0}=u_\al$ since $u_\alpha$ are the normal coordinates for the chosen tau-structure. Rewriting the densities in the $v$-coordinates
$$
h_{\alpha,p}=h_{\al,p}(v; v_x, \dots;\epsilon)=\sum_{k\ge 0} \epsilon^k 
h_{\al,p}^{[k]} (v; v_x, \dots)
$$
we obtain the same recursion relation
$$
\{ ~.~, \int h_{\al,p}(v; v_x, \dots;\epsilon)\, dx\}_2^{[0]} =
 \sum_{q=0}^p A_{\alpha,p}^{\beta,q}
\{~.~,\int h_{\beta,q+1}(v; v_x, \dots;\epsilon)\, dx \}_1^{[0]}
$$ 
with the initial data 
$$
\int h_{\al,0}(v; v_x, \dots;\epsilon)\, dx = \int v_\al\,dx.
$$
Therefore the Hamiltonians $\int h_{\al,p}(v; v_x, \dots;\epsilon)\, dx$
are linear combinations of the standard Hamiltonians 
$\int \theta_{\beta,0}(v)\, dx$, \dots, $\int \theta_{\beta,p}(v)\, dx$
 of the hierarchy (\ref{F-hierarchy}). It follows that the leading terms $h_{\al,p}(v)^{[0]}$ give a tau-structure of the standard hierarchy. It must be related to
the standard tau-structure by a transformation of the form (\ref{2-6-26}),
(\ref{2-6-27}). Modifying
the Hamiltonians $h_{\al,p}$ by the inverse transformation we obtain an equivalent tau-structure for the Poisson pencil $\{~,~\}_\lambda$ satisfying
$$
h_{\al,p}(v; v_x, \dots;\epsilon) = \theta_{\al,p}(v)+ \epsilon\,\pal_x 
f_{\al,p}(v; v_x, \dots; \epsilon)
$$
where 
$$
f_{\al,p}(v; v_x, \dots;\epsilon) = \sum_{k>0} \epsilon^{k-1} f_{\al,p}^{[k]}
(v; v_x, \dots)
$$
for some functions $  f_{\al,p}^{[k]}
(v; v_x, \dots)$. By the definition of the tau-structure we have
$$
{\pal f_{\al,p}(v; v_x, \dots;\epsilon) \over \pal t^{\beta,q}}
= {\pal f_{\beta,q}(v; v_x, \dots;\epsilon) \over \pal t^{\al,p}}.
$$
In particular, 
\beq\label{qt10}
{\pal f_{1,0}(v; v_x, \dots;\epsilon) \over \pal t^{\beta,q}}
= {\pal f_{\beta,q}(v; v_x, \dots;\epsilon) \over \pal x}.
\eeq
So $ f_{1,0}(v; v_x, \dots;\epsilon)$ is an integral of the hierarchy 
(\ref{F-hierarchy}). Due to the polynomiality assumption it must be a total derivative of some function that we denote ${\cal F}$
$$
 f_{1,0}(v; v_x, \dots;\epsilon)= \pal_x {\cal F}  (v; v_x, \dots;\epsilon).
$$
From (\ref{qt10}) we get
$$
\pal_x \pal_{t^{\beta,q}} {\cal F} = \pal_x f_{\beta,q}.
$$
This proves the theorem.
\epf 

We will next obtain upper estimates for the order of the highest derivative
in ${\cal F}^{[k]}$, and we will also describe the explicit form of the first
three terms of this expansion.

Let us first consider the infinitesimal deformation of the Poisson pencil
\newline
$\{v^\al (x),v^\beta (y)\}_\lambda^{[0]}$
caused by
a quasi-Miura transformation
\beq\label{Miura-1}
v_\al\mapsto
w_\al=v_\al+\epsilon\,\frac{\pal^2 {{\cal F}}(v; v_x, \dots, v^{(l)})}{\pal x\pal t^{\al,0}}+O(\epsilon^2).
\eeq
We have changed the notations for the dependent functions of the hierarchy
since the variables $u^i=u^i(v)$, $i=1, \dots, n$, will be reserved for
denoting the canonical coordinates on the Frobenius manifold. 
As in the Section \ref{sec-2}, we denote $v^{\al,s}$ and $w^{\al,s}$ the jet coordinates,
$$
v^{\al, s} ={\pal v^\al\over \pal x^s}, ~~
w^{\al, s} ={\pal w^\al\over \pal x^s},
$$
$$
v^{\al,0}=v^\al, ~~w^{\al,0}=w^\al, ~~v^{\al,1}=v^\al_x, ~~w^{\al,1}=w^\al_x, ~
\dots .
$$

\begin{lemma}\label{l-2-10-2}
The deformed Poisson pencil 
has the form
\eqa
&&\{w^\al(x),w^\beta(y)\}_1=\eta^{\al\beta}\,\delta'(x-y)\nn\\
&&\quad+\epsilon\left(W^{\al\beta}(w,w_x,\dots)\,\delta^{(K_l)}(x-y)
+R^{\al\beta}(w,w_x,\dots)\,\delta^{(K_l-1)}(x-y)+\dots\right)\nn\\
&&\quad+
{\cal O}(\epsilon^2),\label{deform-pb1}\\
&&\{w^\al(x),w^\beta(y)\}_2=g^{\al\beta}(w(x))\,\delta'(x-y)+
\Gamma^{\al\beta}_\gamma(w(x))\,w^\gamma_x\,\delta(x-y)\nn\\
&&\quad+\epsilon\left(S^{\al\beta}(w,w_x,\dots)\,\delta^{(K_l)}(x-y)
+Q^{\al\beta}(w,w_x,\dots)\,\delta^{(K_l-1)}(x-y)+\dots\right)\nn\\
&&\quad+
{\cal O}(\epsilon^2),\label{deform-pb}
\eeqa
where the integer $K_l$ is equal to $l+3$ when $l=2\,m$ and it is equal
to $l+2$ when $l=2m-1$. 
\end{lemma}

The proof can be obtained by a simple straightforward computation.

\begin{lemma}\label{l-2-10-3}
Let $l=2\,m$, then in the deformed Poisson bracket (\ref{deform-pb})
the term \newline
$S^{\al\beta}(w,w_x,\dots)\,\delta^{(l+3)}(x-y)$
does not appear {\rm iff} ${{\cal F}}$
does not depend on $v^{\al,2\,m}, \al=1,\dots,n$.
\end{lemma}

\pf From the form of the quasi-Miura transformation we see that the 
functions $S^{\al\beta}(v,v_x,\dots)$
are given by the formulae
$$
\left(g^{\al\nu}\,c^{\gamma\beta}_\nu+g^{\beta\nu}\,c^{\gamma\al}_\nu\right)
\frac{\pal{{\cal F}}}{\pal v^{\gamma,2m}}=2\,g^{\al\nu}\,c^{\gamma\beta}_\nu
\frac{\pal{{\cal F}}}{\pal v^{\gamma,2m}}
$$
So outside the discriminant $\det (g^{\al\beta})=0$ of the Frobenius manifold 
the above expression vanishes
{\it iff} $\frac{\pal{{\cal F}}}{\pal v^{\gamma,2m}}=0$.
The lemma is proved.\epf

\begin{lemma}\label{l-2-10-4}
Let $l$ be an odd positive integer, $l=2\,m-1$.
Denote $h_\al=\frac{\pal {\cal F}}{\pal v^{\al,2m-1}}$,
then
\eqa
&&W^{\al\beta}=-2m\,\eta^{\al\gamma}\,c^{\beta\xi}_\gamma\,\pal_x h_\xi+
\eta^{\al\gamma}\,\frac{\pal h_\gamma}{\pal t_{\beta,0}}+\eta^{\beta\gamma}\,
\frac{\pal h_\gamma}{\pal t_{\al,0}}+2\,\eta^{\al\gamma}\,c^{\beta\xi}_\gamma\,
\frac{\pal{{\cal F}}}{\pal v^{\xi,2m-2}}\nn\\
&&\quad+2m\,\pal_x(\eta^{\al\gamma}\,c^{\beta\xi}_\gamma)h_\xi.\label{Wab}
\eeqa
\eqa
&&S^{\al\beta}=-2m\,g^{\al\gamma}\,c^{\beta\xi}_\gamma\,\pal_x h_\xi+
g^{\al\gamma}\,\frac{\pal h_\gamma}{\pal t_{\beta,0}}+g^{\beta\gamma}\,
\frac{\pal h_\gamma}{\pal t_{\al,0}}+2\,g^{\al\gamma}\,c^{\beta\xi}_\gamma\,
\frac{\pal{{\cal F}}}{\pal v^{\xi,2m-2}}\nn\\
&&\quad+\left(2m\,\pal_x(g^{\al\gamma}\,c^{\beta\xi}_\gamma)+
\pal_\gamma g^{\al\beta}\,c^{\gamma\xi}_\nu\,v^\nu_x\right) h_\xi,\label{Sab}
\eeqa
and
\eqa\label{combination}
c^\lambda_{\al\beta}\left(S^{\al\beta}-{\cal U}^\al_\gamma\,W^{\gamma\beta}\right)
=(2m+1)\,c^\lambda_{\al\beta}\,c^{\al\gamma}_\nu\,c^{\beta\xi}_\gamma\,
v^\nu_x\,\frac{\pal {\cal F}}{\pal v^{\xi,2m-1}},\quad 1\le\lambda\le n,
\eeqa
where ${\cal U}^\al_\beta=E^\gamma\,c^\al_{\gamma\beta}$ is the matrix of the operator of multiplication by the Euler vector field.
In particular, the left hand sides of the last equalities vanish {\rm iff}
${{\cal F}}$ does not depend on $v^{\al,2\,m-1}, \al=1,\dots,n$.
\end{lemma}

Here we use notations
$$
{\pal\over \pal t_{\alpha,p}} := \eta^{\al\beta}{\pal\over \pal t^{\beta,p}}
$$
for the linear combinations of the flows of the hierarchy 
(\ref{F-hierarchy}).

\pf The identities (\ref{Wab}) and (\ref{Sab}) are deduced from their
definitions. 
The equalities in (\ref{combination}) follow directly from 
(\ref{Wab}) and (\ref{Sab}).
To prove the last statement of the lemma, let us rewrite the right hand 
side of the identity (\ref{combination}) in the canonical coordinates
$u^i=u^i(v^1, \dots, u^n)$, $1\le i \le n$ (see Section \ref{sec-3-10-4a} 
for the 
definition of the canonical coordinates and of the functions $\psi_{i,\al}$
to be used below).
It becomes equal to
$$
(2m+1)\eta^{\lambda\nu}\,\psi_{i\nu}\,\frac{u^i_x}{\psi_{i1}^3}\,\frac{\pal {\cal F}}
{\pal u^{i,2m-1}}.
$$
Since $\det(\eta^{\lambda\nu}\psi_{i\nu})\ne 0$, we deduce that the left hand sides 
of the above equalities vanish {\it iff}
${{\cal F}}$ does not depend on $u^{i,2m-1}, \ 1\le i\le n$ (equivalently, ${{\cal F}}$ does 
not depend on $v^{\al,2m-1}, \ 1\le \al\le n$). Lemma is proved.\epf

\begin{theorem}\label{t-2-10-5}
Let  
\eqa
&&\{w^\al(x),w^\beta(y)\}_1=\eta^{\al\beta}\delta'(x-y)\nn\\
&&\quad +\sum_{i\ge 1}\epsilon^i\,\{w^\al(x),w^\beta(y)\}_1^{[i]}\label{genus-exp-pb1}
\eeqa
\eqa
&&\{w^\al(x),w^\beta(y)\}_2=g^{\al\beta}(w(x))\,\delta'(x-y)+
\Gamma^{\al\beta}_\gamma(w(x))\,w^\gamma_x\,\delta(x-y)\nn\\
&&\quad +\sum_{i\ge 1}\epsilon^i\,\{w^\al(x),w^\beta(y)\}_2^{[i]}
\label{genus-exp-pb}
\eeqa
be a quasitrivial Poisson pencil. Here
$\{w^\al(x),w^\beta(y)\}_1^{[i]} $, $\{w^\al(x),w^\beta(y)\}_2^{[i]}$ 
have the form
\eqa
&&\{w^\al(x),w^\beta(y)\}_1^{[i]}=\sum_{l=0}^{i+1} 
H^{\al\beta}_{i,l}(w;w_x,w_{xx},\dots, w^{(l)})\,
\delta^{(i+1-l)}(x-y),\nn\\
&&\{w^\al(x),w^\beta(y)\}_2^{[i]} =\sum_{l=0}^{i+1} 
K^{\al\beta}_{i,l}(w;w_x,w_{xx},\dots, w^{(l)})\,
\delta^{(i+1-l)}(x-y),
\eeqa
and $H^{\al\beta}_{i,l}, K^{\al\beta}_{i,l}$ are quasihomogeneous 
polynomials in the derivatives of the degree $l$.
Then the quasi-triviality transformation
\beq
v_\al\mapsto
w_\al=v_\al+\sum_{k\ge 1}\epsilon^k\,\frac{\pal^2 {\cal F}^{[k]}
(v^1,\dots,v^n;v^{1,1},\dots, v^{n,1};\dots,v^{1,m_k},\dots,v^{n,m_k})}
{\pal x\pal t^{\al,0}},\label{Miura-2}
\eeq
must have the property
\beq\label{3m-2}
m_{2g},\quad m_{2g+1}\le 3g-2.
\eeq
Moreover, modulo additive constants, we 
\eqa
&&{\cal F}^{[1]}=0,\label{F12}\\
&&{\cal F}^{[2]}=\sum_{i=1}^n a_i\,\log(u^i_x)+f(v),\label{F1}\\
&&{\cal F}^{[3]}=\sum_{i=1}^n f_i(v)\,u^i_x.
\label{F32}
\eeqa
where $a_i$ are constants, $f(v)$ and $f_i(v)$ are some functions, 
and $u^1(v)$, \dots, $u^n(v)$ are the canonical coordinates on $M$.
\end{theorem}

\pf From Lemma \ref{l-2-10-3} and Lemma \ref{l-2-10-4} 
we see that ${\cal F}^{[1]}$ must be a 
constant, and ${\cal F}^{[2]}$, ${\cal F}^{[3]}$ only depend on 
$v^1,\dots,v^n;
v^1_x,\dots, v^n_x$. 
We are to prove that ${\cal F}^{[2g]}, {\cal F}^{[2g+1]}$  depend at most
on $v^1,\dots,v^n; v^{1,1},\dots,v^{n,1};\dots;v^{1,3g-2},\dots,v^{n,3g-2}$ 
for $g\ge 2$. 
This can be done
by induction. Assume 
that the statement is true for $g\le N-1$, we need to prove the validity of 
the above
property for $g=N$. Express the left and right hand sides
of (\ref{genus-exp-pb1}) and 
(\ref{genus-exp-pb}) in the $v$-coordinates by using the quasi-Miura
transformation (\ref{Miura-2}) and compare the $\epsilon^{2N}$ terms
of both sides. We denote, as in Lemma \ref{l-2-10-2}, 
by $W^{\al\beta}$
and 
$S^{\al\beta}$
the coefficients of the highest order derivatives of the delta function 
in the $\epsilon^{2N}$ term
that are contributed
by ${\cal F}^{[2]},\dots,{\cal F}^{[2N-1]}$ in 
the left hand sides of (\ref{genus-exp-pb1}) and 
(\ref{genus-exp-pb}) respectively. When $N$ is even this
is the coefficients of $\delta^{(3N+1)}(x-y)$ and when $N$ is odd this is the 
coefficient of $\delta^{(3N)}(x-y)$. While on the right hand side of 
(\ref{genus-exp-pb1}) and
(\ref{genus-exp-pb}) the highest order derivatives of the delta function 
in the $\epsilon^{2N}$ term is $\delta^{(2N+1)}(x-y)$ due to the form 
of $\{~,~\}_2^{[2N]}$. So, in order that the equality 
(\ref{genus-exp-pb}) holds true, the ${\cal F}^{[2N]}$ term in the quasi-Miura
transformation (\ref{Miura-2}) must be responsible for the killing of the 
$S^{\al\beta}\,\delta^{(3N+1)}(x-y)$ term when $N$ is even and the 
$S^{\al\beta}\,\delta^{(3N)}(x-y)$ term when $N$ is odd in the 
left hand sides
of (\ref{genus-exp-pb1}) and 
(\ref{genus-exp-pb}). We first consider the
case when $N$ is odd. In this case, by using 
Lemma \ref{l-2-10-3}  and Lemma \ref{l-2-10-4} we
immediately deduce that ${\cal F}^{[2g]}$ depends at most on 
$v^1,\dots,v^n; v^{1,1},\dots,v^{n,1},\dots,v^{1,3N-2},\dots,v^{n,3N-2}$. 
Now let us consider
the case when $N$ is even. From Lemma \ref{l-2-10-3} 
we see that ${\cal F}^{[2N]}$
depends at most on $v^\al, v^{\al,1},\dots,v^{\al,3N-1}$, $\al=1, \dots, n$.
We are to prove that actually it does not depend on
$v^{1,3N-1},\dots,v^{n,3N-1}$. To prove this let us 
compute $W^{\al\beta}$ and $S^{\al\beta}$ to obtain
\eqa
&&S^{\al\beta}=\sum_{k=1}^{N-1}\frac{\pal{\cal F}^{[2k]}}
{\pal v^{\gamma,3k-2}}\,\frac{\pal{\cal F}^{[2N-2k]}}
{\pal v^{\nu,3(N-k)-2}}\,c^{\al\gamma}_{\sigma}\,c^{\beta\nu}_{\rho}\,
g^{\sigma\rho},\nn\\
&&W^{\al\beta}=\sum_{k=1}^{N-1}\frac{\pal{\cal F}^{[2k]}}
{\pal v^{\gamma,3k-2}}\,\frac{\pal{\cal F}^{[2N-2k]}}
{\pal v^{\nu,3(N-k)-2}}\,c^{\al\gamma}_{\sigma}\,c^{\beta\nu}_{\rho}\,
\eta^{\sigma\rho},\nn
\eeqa
and 
$$
c^\lambda_{\al\beta}\,\left(S^{\al\beta}-{\cal U}^\al_\gamma\,W^{\gamma\beta}\right)
=0.
$$
So from Lemma \ref{l-2-10-4} 
we see that ${\cal F}^{[2N]}$ indeed does not depend
on $v^{1,3N-1},\dots,v^{n,3N-1}$. In a similar way,
we can prove that ${\cal F}^{[2N+1]}$ depends at most on
$v^1,\dots,v^n$; $v^{1,1},\dots,v^{n,1}$;\newline
$\dots;v^{1,3M-2},\dots,v^{n,3M-2}$. 
We have thus finished the procedure of induction and proved (\ref{3m-2}). 

We next prove that ${\cal F}^{[2]}$ must have the form
(\ref{F1}). Our assumption on $\{~,~\}_1^{[2]}$,  
$\{~,~\}_2^{[2]}$ implies that 
$H^{\al\beta}_{2,0}$ and $K^{\al\beta}_{2,0}$ are functions of 
 $u^1,\dots,u^n$. Rewrite both sides of (\ref{genus-exp-pb1}) and
(\ref{genus-exp-pb}) in the $v$-coordinates by using the quasi-Miura
transformation (\ref{Miura-2}) and compare the coefficients of 
$\epsilon^{2}\delta'''(x-y)$, we obtain by using (\ref{combination})
the following identities:
\beq\label{system-F}
3\,c^\lambda_{\al\beta}\,c^{\al\gamma}_\nu\,c^{\beta\xi}_\gamma\,
v^\nu_x\,\frac{\pal {\cal F}^{[2]}}{\pal v^{\xi,1}}
=c^\lambda_{\al\beta}\left(K^{\al\beta}_{2,0}-{\cal U}^\al_\gamma\,
H^{\gamma\beta}_{2,0}\right),
\quad 1\le\lambda\le n.
\eeq
In the canonical coordinates the left hand sides of the  above identities 
have the expressions
$$
\eta^{\lambda\nu} \sum_{i=1}^n \frac{\psi_{i\nu}}{\psi_{i1}^3}\,u^i_x\,
\frac{\pal {\cal F}^{[2]}}
{\pal u^i_x}.
$$
So from (\ref{system-F}) we see that there exist functions $a_1(u),
\dots,a_n(u)$ of $u^1,\dots,u^n$ such that 
\beq\label{expa-F1}
u^i_x\,\frac{\pal {\cal F}^{[2]}}{\pal u^i_x}=a_i(u),\quad 1\le i\le n,
\eeq
which yields
\beq\label{for-F}
{\cal F}^{[2]}=\sum_{i=1}^n a_i(u)\,\log(u^i_x)+f(u)
\eeq
for certain function $f$ of $u^1,\dots, u^n$. Now let's prove that 
$a_i(u)$ are constants. Indeed, the coefficients of 
$\epsilon^2\,\delta'''(x-y)$ in the left hand side of 
(\ref{genus-exp-pb}) written in the $v$-coordinates have the expressions
\eqa
&&-2\,g^{\al\gamma}\,c^{\beta\xi}_\gamma\,\pal_x\left(
\frac{\pal {\cal F}^{[2]}}
{\pal v^{\xi}_x}\right)+
g^{\al\gamma}\,\frac{\pal}{\pal t_{\beta,0}}\left(\frac{\pal {\cal F}^{[2]}}
{\pal t^{\gamma,1}}\right)+g^{\beta\gamma}\,
\frac{\pal}{\pal t_{\al,0}}\left(\frac{\pal {\cal F}^{[2]}}
{\pal t^{\gamma,1}}\right)\nn\\
&&\quad +\left(2m\,\pal_x(g^{\al\gamma}\,c^{\beta\xi}_\gamma)+
\pal_\gamma g^{\al\beta}\,c^{\gamma\xi}_\nu\,v^\nu_x\right) 
\frac{\pal {\cal F}^{[2]}}{\pal v^{\xi}_x}
+2\,g^{\al\gamma}\,c^{\beta\xi}_\gamma\,
\frac{\pal{\cal F}^{[2]}}{\pal v^{\xi}}\label{exp-F}
\eeqa
Substituting the formula (\ref{for-F}) into the above expressions,
we see that the first four summands are rational polynomials in the 
$x$-derivatives of
$u^1,\dots,u^n$,
and the last summand is a linear combination of $\log(u^1_x),
\dots,\log(u^n_x)$. Since (\ref{exp-F}) should be functions of $u^1,
\dots, u^n$ only, we deduce
that 
$$
g^{\al\gamma}\,c^{\beta\xi}_\gamma\,
\frac{\pal{a_i}}{\pal v^{\xi}}=0,\quad 1\le i\le n.
$$
Putting $\beta=n$ in the above equations we obtain 
$$
\frac{\pal{a_i}}{\pal v^{\xi}}=0,\quad 1\le i, \xi\le n
$$
for generic point $v$ when $\det(g^{\al\beta}(v))\ne 0$.
So $a_1(u), \dots, a_n(u)$ are constants.

We will now prove that ${\cal F}^{[3]}$ has the form (\ref{F32}). Since
${\cal F}^{[3]}$ only depends on $v^1,\dots,v^n; v^1_x,\dots, v^n_x$, 
the highest 
order of the derivatives of the delta-function in $\{~,~\}_2^{[3]}$
is $3$, and the coefficients of $\delta'''(x-y)$ depend linearly on 
$v^1_x,\dots,v^n_x$ since $\deg K^{\al\beta}_{3,1}=1$. An approach
similar to the one given in the derivation of (\ref{expa-F1})
yields
\beq
u^i_x\,\frac{\pal {\cal F}^{[3]}}{\pal u^i_x}=\sum_{k=1}^n b^i_k(u)\,u^k_x,
\quad 1\le i\le n.
\eeq
By using the compatibility condition 
$\frac{\pal}{\pal u^i_x}\left(\frac{\pal}{\pal u^j_x}  {\cal F}^{[3]}\right)=
\frac{\pal}{\pal u^j_x}\left(\frac{\pal}{\pal u^i_x} {\cal F}^{[3]}\right)
$
we have
$$
b^i_k(u)=0,\quad i\ne k
$$ 
which yields
$$
 {\cal F}^{[3]}=\sum_{i=1}^n b^i_i(u)\,u^i_x+h(u).
$$
Here $h(u)$ is certain function of $u^1,\dots,u^n$. To prove that $h(u)$
is a constant we use the explicit expression for $K^{\al\beta}_{3,1}$
which is given by the right hand side of (\ref{Sab}) with $m=1$ and
${\cal F}= {\cal F}^{[3]}$. From this expression we get 
$$
g^{\al\gamma}\,c^{\beta\xi}_\gamma\,
\frac{\pal{h(u)}}{\pal v^{\xi}}=0,\quad 1\le \al, \beta\le n
$$
which implies that $h(u)$ is a constant.
The theorem is proved.
\epf

Since ${\cal F}^{[1]}$ is a constant and $ {\cal F}^{[3]}$ is a polynomial 
in the $x$-derivatives of $v^\al$,
the quasitrivial bihamiltonian structure (\ref{genus-exp-pb1}) and
(\ref{genus-exp-pb}) is equivalent to a  quasitrivial bihamiltonian structure
whose quasitriviality transformation does not contain the $\epsilon$ and
$\epsilon^3$ terms.
The equivalence is established by the Miura transformation
$$
u_\al\mapsto u_\al-\epsilon^3\frac{\pal^2  {\cal F}^{[3]}(v,v_x)}
{\pal x\pal t^{\al,0}}.
$$
We were not be able to prove that, for an arbitrary quasitriviality
transformation all the terms with odd powers of $\epsilon$ can be gauged
out by a Miura transformation. In the next section we will prove that this is the case under 
an additional assumption about the structure of the symmetry algebra of the hierarchy.

\setcounter{equation}{0}
\setcounter{theorem}{0}

\subsection{Virasoro symmetries}\par

In this section we will show that the Principal Hierarchy 
(\ref{F-hierarchy}) on ${\cal L}(M^n)$ of a $n$-dimensional 
Frobenius manifold
$M^n$ always admits a rich algebra of symmetries isomorphic to the half of the Virasoro algebra
with the central charge $n$. The Virasoro algebra is constructed in terms of the spectrum of the Frobenius
manifold. The operators of the Virasoro algebra act by nonlinear first order differential operators
on the tau-cover of the hierarchy.
We will characterize general solution of the Principal Hierarchy in terms of the action of the Virasoro
algebra. Due to quasitriviality the action of the Virasoro algebra can be extended to the full
hierarchy (\ref{pde00}). 
Our last condition requires linearity of this action of Virasoro onto the tau-function
of the full hierarchy. We prove that, for a semisimple Frobenius manifold $M^n$ this condition
uniquely determines the quasitriviality transformation (\ref{qtt-2-8}). 
For the semisimple Frobenius manifolds
coinciding with quantum cohomology of a smooth projective variety $X$ we identify our condition 
of linearization
of the Virasoro symmetries with the Virasoro constraints conjectured by T.Eguchi {\it et al.} \cite{eguchi1}
in the theory of the higher genus Gromov - Witten invaraints of $X$.

\subsubsection{From Galilean invariance to Virasoro symmetries of the Principal Hierarchy}\label{sec-3-10-1}\par

This subsection is based on the paper \cite{selecta}. We also find a nice generating formula
for the Virasoro symmetries that will be useful in subsequent calculations.
 
Let us begin with the following 

{\bf Definition.} The PDE 
\beq
{\pal v^\alpha\over \pal s} = B^\alpha(x, {\bf t}, v; v_x, \dots, ; \epsilon)
\eeq
is called (infinitesimal) {\it symmetry} of the hierarchy 
(\ref{F-hierarchy}) if it commutes with all the flows of the hierarchy
$$
{\pal \over \pal s} {\pal v\over \pal t^{\alpha,p}} = {\pal\over \pal t^{\alpha, p}} {\pal v\over \pal s}.
$$
\medskip
According to our definition the flows of the hierarchy themselves are symmetries. In this case the r.h.s.
does not depend on $x$, ${\bf t}$. Less trivial example is given by

\begin{lemma} The flow
\beq\label{galilean}
{\pal v\over \pal s} = e + \sum_{p=1}^\infty t^{\alpha, p} {\pal v\over \pal t^{\alpha, p-1}}
\eeq
is a symmetry of the Principal Hierarchy (\ref{F-hierarchy}).
\end{lemma}

Here $e$ is the unity vector field on the Frobenius manifold. 

\pf (cf. \cite{selecta}) Let us consider first the flow
$$
{\pal v\over \pal s} = e.
$$
Using 
$$
\pal_1 \theta_{\alpha, p+1} = \theta_{\alpha,p}, ~~p\geq 0, ~~\pal_1 \theta_{\alpha,0}= \eta_{\alpha \,1}
$$
we obtain 
$$
{\pal \over \pal s} {\pal v\over \pal t^{\alpha,p}} - {\pal\over \pal t^{\alpha, p}} {\pal v\over \pal s}
= {\pal v\over \pal t^{\alpha, p-1}}, ~p>0, ~~ 
{\pal \over \pal s} {\pal v\over \pal t^{\alpha,0}} = {\pal\over \pal t^{\alpha, 0}} {\pal v\over \pal s}.
$$
The term $ \sum_{p=1}^\infty t^{\alpha, p} {\pal v\over \pal t^{\alpha, p-1}}$ in (\ref{galilean}) compensates the
noncommutativity of the above flow with the equations of the hierarchy (cf. \cite{fuchssteiner}). \epf

\begin{exam} For the Riemann equation
$$
v_t = v\, v_x
$$
the above symmetry coincides with the infinitesimal form 
$$
v_s = 1 + t\, v_x
$$
of the Galilean transformation 
$$
x\mapsto x+ c\, t, ~~t\mapsto t, ~~ v\mapsto v+c.
$$
Here $c$ is an arbitrary parameter.
\end{exam}

Also in the general case we will call (\ref{galilean}) 
the Galilean symmetry of the hierarchy (\ref{F-hierarchy}).
It is natural to produce an infinite chain of other symmetries by applying the recursion operator
\beq\label{recur-pc}
{\cal R}{\pal\over \pal  s_{m-1}} = {\pal \over \pal s_m}, ~~m\geq 0
\eeq
where we redenote $ \pal / \pal s \mapsto \pal / \pal s_{-1}$ the Galilean symmetry (\ref{galilean}). Such symmetries
were discovered in \cite{chen-lee-lin} for the case of KdV (also the idea  appeared already in \cite{shabat}).
It was shown in \cite{zubelli-magri} that, for the symplectic bihamiltonian structures the above chain
satisfies the Virasoro commutation  relations
\beq\label{vc-sym}
\left[ {\pal\over \pal s_i}, {\pal\over \pal s_j}\right] =(j-i) {\pal\over \pal s_{i+j}}, ~~i, j\geq -1
\eeq
if 
$$
Lie_{\pal / \pal s_0} {\cal R} = {\cal R}.
$$

We cannot apply directly this result to the Principal Hierarchy. Indeed, the Poisson pencil (\ref{normalfm}) is not
symplectic. From practical point of view the recursion procedure 
(\ref{recur-pc}) produces nonlocal flows, starting
from $m=1$. For the simplest example of Riemann hierarchy
$$
{\pal v\over \pal s_0} = v + {x\over 2} v_x + \sum_{k\geq 0} \left ( k+{1\over 2}\right)  
t_k {\pal v\over \pal t_k}
$$
$$
{\pal v \over \pal s_1} = v^2 + {3\over 4} x^2 {\pal v\over \pal t_1} + {1\over 4} \pal_x^{-1} v +\sum_{k\geq
0}  \left( k+{1\over 2}\right) \left( k+{3\over 2}\right) t_k {\pal v\over \pal t_{k+1}}
$$
etc. 

The general problem of dealing with the nonlocalities in the bihamiltonian recursion procedure has been 
analyzed, in
the style of formal variational calculus, in \cite{khorkova, sergeyev}. For the case of Principal
Hierarchy it was shown in \cite{selecta} that the bihamiltonian recursion procedure can be used to produce
Virasoro symmetries of the tau-cover (\ref{2-4-16}) of the Principal Hierarchy. We present here this result in a
slightly modified form, using generating functions.

Let us consider first the nonresonant case, i.e., assuming that the spectrum
of $\hat\mu$ contains no half-integers.
Let $p_\alpha(v; \lambda)$, $\alpha=1, \dots, n$ be a basis of independent periods (e.g., the one described
in Theorem \ref{t-2-7-7}). Introduce the constant Gram matrix
$$
G^{\alpha\beta} := \left( {\pal\over\pal p_\alpha}, {\pal\over\pal p_\beta}\right)_\lambda=
{\pal v^\sigma\over \pal p_\alpha} g_{\sigma\epsilon}(v; \lambda) 
{\pal v^\epsilon\over  \pal p_\beta}, ~~ 
(g_{\sigma\epsilon}(v;\lambda)) = \left(
g^{\sigma\epsilon}(v)-\lambda \eta^{\sigma\epsilon}\right)^{-1}.
$$
Recall that, if the basis of independent periods is chosen as in 
Theorem \ref{t-2-7-7}, then
\beq
\left( G^{\alpha\beta}\right) = -{1\over 2\, \pi} (e^{\pi\, i\, R} e^{\pi\, i \, \hat \mu} + 
e^{-\pi\, i\, R} e^{-\pi\, i \, \hat \mu} )\eta^{-1}.
\eeq
Introduce {\it action functions} 
$$
s_\alpha=s_\alpha(x,f, \pal f/ \pal {\bf t}; \lambda)
$$ 
\eqa
&&
s_\alpha =\int^x p_\alpha(v;\lambda)dx 
= x\, \omega_1 \sum_{q\geq 0} {\Gamma_q (R, \hat\mu+q-{1\over
2})\over \lambda^{q-1}} \lambda^{-{1\over 2} -\hat\mu} \lambda ^{-R}\nn\\
&&\qquad
 +\sum_{m\geq 0} \sum_{p+q=m} {\partial
f\over \pal t^p} {\Gamma_q(R, \hat\mu + m +{1\over 2})\over \lambda^m} 
\lambda^{-{1\over 2} -\hat\mu} \lambda ^{-R}.\nn\\
\eeqa
Put 
\eqa
&&
S_\alpha=S_\alpha(x, {\bf t}, f, \pal f/ \pal {\bf t}; \lambda)\nn\\
&&=s_\alpha
+ \sum_m \lambda^{m+1} \sum_{p-q=m} (-1)^p t_p 
\Gamma_q (R, \hat\mu-m-{1\over 2}) 
\lambda^{-{1\over 2} -\hat\mu} \lambda ^{-R}\nn\\
&&
= \int_0^\infty {dz\over \sqrt{z}} e^{-\lambda\, z} \left[ \sum_{p\geq 0} {\pal f \over \pal t^p} z^p 
+\sum_{q\geq 0} (-1)^q \bar t_q z^{-q-1}\right] z^{\hat\mu} z^R.
\label{laplace-2-11}
\eeqa
Here we use short notations for the following row vectors
$$
{\pal\over \pal t^p} := \left( {\pal\over \pal t^{1,p}}, \dots, {\pal \over \pal t^{n,p}}\right),
$$
$$
t_q : = (t_{1,q}, \dots, t_{n,q}), ~~ t_{\alpha,q}:= \eta_{\alpha\beta} t^{\beta,q},
$$
and we denote
$$
\bar t_q = t_q, ~q>0, ~~ \bar t_0 = t_0 + x\, \omega_1.
$$
The Laplace integrals in (\ref{laplace-2-11}) 
are defined as in Theorem \ref{t-2-7-7}.

It will be always assumed that the periods $p=p(v; \lambda)$ are chosen in such a way that
\beq
\pal_\lambda p = - \pal_1 p.
\eeq
The functions $S_\alpha$, $\alpha = 1, \dots, n$ satisfy the following simple identities
\beq
\pal_x {\pal\over \pal \lambda} S_\alpha = {\pal p_\alpha\over \pal \lambda} = - \pal_1 p_\alpha
\eeq
\beq\label{2-11-9}
{\pal\over \pal t^{\gamma,0}} {\pal \over \pal \lambda} S_\alpha = - \pal_\gamma p_\alpha.
\eeq
In the proof of the second formula one is to use the equation
$$
\pal_\alpha \pal_\beta p = c_{\alpha\beta}^\gamma(v) \pal_\gamma\pal_1 p.
$$

\begin{theorem} \label{t-2-11-3}
The flows $\pal / \pal s_m$, $m\geq -1$, are  defined on the tau-covering of the Principal
Hierarchy by the following generating formula
\eqa
&&{\pal \over \pal s} = \sum_{m\geq -1} {1\over \lambda^{m+2}} 
{\pal \over \pal s_{m}}
\\
&&
{\pal f\over \pal s} = \left[ -{1\over 2} {\pal S_\alpha\over \pal \lambda} G^{\alpha\beta}
{\pal S_\beta\over \pal \lambda}\right]_+
\label{2-11-10a}\\
&&
{\pal f_{\gamma,p}\over \pal s } = {\pal\over \pal t^{\gamma,p}} {\pal f\over
\pal s}
\label{2-11-10b}\\
&&
{\pal v_\gamma\over \pal s} = \left[ \,\pal_x \pal_\gamma p_\alpha \,  G^{\alpha\beta} {\pal
S_\beta\over \pal \lambda}\right]_+ - \eta _{\gamma\epsilon} \left((E-\lambda\,
e)^{-1}\right)^\epsilon
\label{2-11-10c}
\eeqa
\end{theorem}

In these formulae $[ ~~]_+$ means the regular part of the expansion of the 
function vanishing at $\lambda = \infty$.

\pf We first derive the part of the symmetries not containing explicitly the times.

\begin{lemma} Let $B_{-1}=(B_{-1,\alpha})$, $B_0=(B_{0,\alpha})$, $B_1=(B_{1,\alpha})$, $B_2=(B_{2,\alpha})$, \dots, be the r.h.s. of the flows 
$$
{\pal v_\alpha\over \pal s_m} = B_{m,\alpha}
$$
defined recursively by
\beq
B_m = {\cal R} B_{m-1}, ~m\geq 0, ~~ B_{-1, \alpha} 
= \eta_{\alpha,1}.
\eeq 
Here
\beq
{\cal R}^\alpha_\beta = {\cal U}^\alpha_\beta + \eta_{\beta\gamma}\Gamma^{\alpha\gamma}_\epsilon v^\epsilon_x
\pal_x^{-1}
\eeq
is the recursion operator,
$$
{\cal R} = \{~,~\}_2 \{~, ~\}_1^{-1}.
$$
Then the generating function
$$
B_\lambda := {B_{-1}\over \lambda} + {B_0\over \lambda^2} + 
{B_1\over \lambda^3} +\dots
$$
reads
\beq
(B_\lambda)_\alpha =-{1\over 2} \pal_{t^{\alpha,0}} \pal_x \left[ \pal_\lambda \int^x  p^a dx \,G_{ab}
\pal_\lambda \int^x  p^b\right].
\eeq
Here $p^a=p^a(v; \lambda)$, $a=1, \dots, n$ is a system of independent periods of the Frobenius manifold, the
constant Gram matrix $G_{ab}$ is defined by 
$$
G_{ab} := {\pal v^\al\over\pal p^a} {\pal v^\beta\over\pal p^b} 
\left( g^{\alpha\beta}(v)-\lambda \,
\eta^{\alpha\beta}\right)^{-1}.
$$
\end{lemma}

\pf Denote $\varpi_1$ and $\varpi_2$ the tensors of the first and the second Poisson brackets respectively. Then
$$
B_\lambda = {1\over \lambda}\left( 1-{{\cal R}\over \lambda}\right)^{-1}B_{-1}
 =-\varpi_1 (\varpi_2 -\lambda\, \varpi_1)^{-1} B_{-1}.
$$
Let us first compute 
$$
b:= (\varpi_2 -\lambda\, \varpi_1)^{-1} B_{-1}.
$$
To this end we are to solve the following linear inhomogeneous equation
\beq
(g^{\alpha\beta}-\lambda\, \eta^{\alpha\beta}) \pal_x b_\beta + \Gamma^{\alpha\beta}_\epsilon v^\epsilon_x
b_\beta = \delta_1^\alpha.
\eeq
Denote
$$
\phi_\alpha^a:=\pal_\al p^a, ~~ a=1, \dots, n
$$
the basis of solutions of the corresponding linear homogeneous equation. Applying variation of constants we
obtain
$$
b_\alpha = \phi_\alpha^a G_{ab}\int^x \pal_1 p^b dx.
$$
Therefore
\beq
(B_\lambda)_\alpha =-\pal_x \left( \pal_\alpha p^a G_{ab} \int^x \pal_1 p^b dx\right).
\eeq
Using
$$
\pal_{t^{\alpha,0}} \int^x p\, dx =\pal_\alpha p
$$
valid for an arbitrary period $p=p(v;\lambda)$
(cf. (\ref{2-11-9}) above) we complete the proof of the lemma. \epf

To complete the proof of the theorem we are to apply the operator 
$$
{1\over \lambda}\left( 1-{{\cal R}\over \lambda}\right)^{-1}
 =-\varpi_1 (\varpi_2 -\lambda\, \varpi_1)^{-1}
$$
at the sum $\sum t^{\epsilon, r} {\pal v\over \pal t^{\epsilon, r-1}}$.
We leave this part of calculations as an exercise to the reader.

The final step in the proof of the theorem is to check that, indeed, the flows (\ref{2-11-10a})--(\ref{2-11-10c}) are symmetries
of the tau-covering of the Principal Hierarchy satisfying the Virasoro commutation relations. This can be done
by identifying the coefficients of the expansion of the flows 
(\ref{2-11-10c}) with the symmetries obtained in
\cite{selecta}). The theorem is proved. \epf

In the resonant case we can regularize the formula 
(\ref{2-11-10a})--(\ref{2-11-10c}) as follows:
Introduce
\beq\label{2-11-16}
S_\alpha^{(\nu)} =\int_0^\infty {dz\over z^{{1\over 2}-\nu}} \, 
e^{-\lambda\, z}
\left[ \sum_{p\geq 0} {\pal f\over \pal t^p} z^p 
+ \sum_{q\geq 0} (-1)^q \bar t_q z^{-q-1}\right] z^{\hat\mu} z^R.
\eeq
Define the deformed Gram matrix by
\beq \label{2-11-17}
G^{\alpha\beta}(\nu) := -{1\over 2\, \pi} \left[ \left( e^{\pi\, i\, R}
 e^{\pi\, i\, (\hat\mu+\nu)} +
e^{-\pi\, i\, R} e^{-\pi\, i\, (\hat\mu+\nu)}\right) 
\,\eta^{-1}\right]^{\alpha\beta}.
\eeq
Then we define the regularized symmetry by
\beq
{\pal f\over \pal s} =\lim_{\nu\to 0}\left[ -{1\over 2} 
{\pal S_\alpha^{(\nu)}\over \pal \lambda} G^{\alpha\beta}(\nu) 
{\pal S_\beta^{(-\nu)}\over \pal \lambda}\right]_+.
\eeq
Existence of the limit and the Virasoro commutation relations
can be proved representing the r.h.s. in the form similar to (\ref{2-11-28}).

\subsubsection{Free field realization of the Virasoro algebra}
\label{sec-3-10-2}\par

Let $({\cal L}, <\ ,\ >, \hat\mu, R)$ be the spectrum of a $n$-dimensional Frobenius manifold. We will construct here a
representation of the Virasoro algebra in the ring of functions of infinite number of variables
${\bf t} = (t^{\alpha,p})$. The construction is isomorphic to that 
of \cite{selecta} but the formulae
are much simpler.

Let us choose a basis $e_1$, \dots, $e_n$ in ${\cal L}$, denote $\eta_{\alpha\beta}=<e_\alpha, e_\beta>$.
We introduce the Heisenberg algebra with the generators
$$
a_{\alpha,p}, ~~\alpha=1, \dots, n, ~~p\in {\bf Z}+{1\over 2}
$$
and the commutation relations
\beq\label{heis}
[a_{\alpha,p},a_{\beta,q}]=(-1)^{p-{1\over 2}} \eta_{\alpha\beta} \delta_{p+q, 0}.
\eeq
Introduce the row vectors
$$
{\bf a}_p = (a_{1,p}, \dots, a_{n,p})
$$
and put
\beq
\phi_\alpha (\lambda) =\left( \int_0^\infty {dz\over z} e^{-\lambda\, z} \sum_{p\in {\bf Z}+{1\over 2}} {\bf
a}_p z^{p+\hat\mu} z^R\right)_\alpha, ~~\alpha=1, \dots, n.
\eeq
We consider first the nonresonant case where the spectrum of $\hat\mu$ contains
no half-integers. In this case
the Virasoro operators are given by the following generating function
\beq\label{2-11-21}
T(\lambda) = \sum_{m\in {\bf Z}} {L_m\over \lambda^{m+2}} =-{1\over 2}
: \pal_\lambda \phi_\alpha \, G^{\alpha\beta} \pal_\lambda \phi_\beta \, : 
+ {1\over 4\, \lambda^2} {\rm tr}\,
\left( {1\over 4} - \hat\mu^2\right).
\eeq
Here
\beq
G^{\alpha\beta}= -{1\over 2\, \pi} \left[ \left( e^{\pi\, i\, R} e^{\pi\, i\, \hat\mu} +
e^{-\pi\, i\, R} e^{-\pi\, i\, \hat\mu}\right) \,\eta^{-1}\right]^{\alpha\beta},
\eeq
the normal ordering is defined by
\eqa
&&:a_{\alpha,p} a_{\beta,q}: = a_{\beta,q} a_{\alpha,p} \quad
{\rm if} ~q<0, ~p>0,
\nn\\
&&
:a_{\alpha,p} a_{\beta,q}:= a_{\alpha,p} a_{\beta,q}\quad
{\rm otherwise}.\nn
\eeqa

\begin{lemma} The operators $L_m$ satisfy Virasoro commutation relations
\beq
[L_i, L_j] = (i-j) L_{i+j} + n\, {i\,(i^2-1)\over 12} \delta_{i+j,0}.
\eeq
\end{lemma}

In more general resonant case we regularize the formula (\ref{2-11-21}) as follows. Introduce
the operator-valued functions
\beq
\phi_\alpha^{(\nu)}(\lambda)=\left( \int_0^\infty {dz\over z^{1-\nu}} \, 
e^{-\lambda\, z} \sum_{p\in {\bf Z}+{1\over 2}} {\bf
a}_p z^{p+\hat\mu} z^R\right)_\alpha, ~~\alpha=1, \dots, n.
\eeq
and the Gram matrix $G^{\alpha\beta}(\nu)$ as in (\ref{2-11-17}).
Here $\nu$ is an arbitrary complex parameter. Put
\beq\label{regular-T}
T^{(\nu)}(\lambda) =\sum_{m\in {\bf Z}} {L_m^{(\nu)}\over \lambda^{m+2}}
=-{1\over 2}
: \pal_\lambda \phi^{(\nu)}_\alpha \, G^{\alpha\beta}(\nu) \pal_\lambda 
\phi_\beta^{(-\nu)} \, : + {1\over 4\, \lambda^2} {\rm tr}\,
\left( {1\over 4} - \hat\mu^2\right).
\eeq
\begin{lemma} Let $k$ be the minimal positive integer such that $R^k=0$. Then
there exist the limits
\eqa
&&
L_m:= \lim_{\nu\to 0} L_m^{(\nu)}, \quad m\geq -1
\nn\\
&&
L_m:= \lim_{\nu\to 0} \nu^k L_m^{(\nu)}, \quad m<-1.
\nn\\
\eeqa
These operators satisfy the following commutation relations
\eqa
&&
[L_i, L_j] =0, ~~i, j <-1, ~{\rm or}~ i+j\geq -1, ~{\rm but} ~(i+1)(j+1)<0,
\nn\\
&&
[L_i, L_j]=(i-j) L_{i+j}, ~~i+j<-1, ~(i+1)\, (j+1)<0, ~~{\rm or} ~i, j\geq -1.
\nn
\eeqa
\end{lemma}

The proof easily follows form the explicit formula
\eqa
&&
T^{(\nu)}(\lambda) 
= {1\over 2\pi}\times \nn\\
&&\times
\sum_{p, q\in {\bf Z}+{1\over 2}} 
\sum_{r\geq 0} :\,{\bf a}_p \left[ e^{R\, \pal_\nu}\right]_r {\Gamma(\hat\mu + \nu
+ p +r +1) \cos\pi(\hat\mu+\nu) \Gamma(-\hat\mu-\nu +q+1)\over
\lambda^{p+q+r+2}} {\bf a}^q\,: 
\nn\\
&&
+{1\over 4\, \lambda^2} \tr \left({1\over 4}
-\hat\mu^2\right).\label{2-11-28}
\eeqa
In this formula 
$$
{\bf a}^q = \left( a^{1,q}, \dots, a^{n,q}\right)^T
$$
is a column vector with the entries
$$
a^{\alpha,q}:= \eta^{\alpha\beta} a_{\beta,q}.
$$
\epf

To prove the commutation relations of the Virasoro operators it suffices to
compute the commutator
$$
[T(\lambda_1),T(\lambda_2)] ={1\over 2} \sum_{p,q} \sum_s 
(-1)^{s-{1\over 2}} :{\bf a}_p \left( M^p_{-s}(\lambda_2) M^s_q(\lambda_1) 
- M^p_{-s}(\lambda_1) M^s_q(\lambda_2)\right) \, {\bf a}^q:
$$
$$
+{1\over 2} \sum_{p,q>0} (-1)^{p+q+1} (M^p_q(\lambda_1) M^{-q}_{-p}(\lambda_2)
- M^p_q(\lambda_2) M^{-q}_{-p}(\lambda_1).
$$
Here
\beq\label{def-M}
M^p_q(\lambda) =\sum_{r\ge 0} N^p_q(r)\,\lm^{-p-q-r-2}
\eeq
with 
\eqa
\hskip -1.9 truecm
&&N^p_q(r) =\nn\\
\hskip -1.5truecm 
&&\left.{1\over \pi}\, \left[ e^{R\, \pal_\nu}\right]_r
\left( \Gamma(\hat\mu+\nu +p+r+1) \cos \pi(\hat\mu+\nu)
\Gamma(-\hat\mu-\nu+q+1)\right)\right|_{\nu=0}.
\eeqa
This can be easily done using the identity
\eqa
&&\sum_{r_1+r_2=m+n-p-q}
(-1)^{p+r_2-n-\frac12}\,N^p_{n-p-r_2}(r_2)\,N^{-n+p+r_2}_q(r_1)\nn\\
&&+\sum_{r_1+r_2=m+n-p-q}
(-1)^{m-p-r_2-\frac12}\,N^p_{-m+p+r_2}(r_2)\,N^{m-p-r_2}_q(r_1)\nn\\
&&
=(m-n)\,N^p_q(m+n-p-q).\nn
\eeqa

A natural representation of the Heisenberg algebra (\ref{heis}) is obtained as follows
\eqa 
&&
a_{\alpha,p} =\epsilon{\pal\over \pal t^{\alpha, p-{1\over 2}}}, ~~p>0,\nn\\
&&
a_{\alpha,p} = \epsilon^{-1}(-1)^{p+{1\over 2}} \eta_{\alpha\beta} t^{\beta, -p-{1\over 2}}, ~~p<0.\nn\\
\eeqa
In this representation $T(\lambda)$ and $L_m$ become linear second order 
differential operators
\beq
T(\lambda) = T(\epsilon^{-1}{\bf t}, \epsilon{\pal/\pal{\bf t}};\lambda)
\eeq
\eqa 
&&\hskip -1.2truecm
L_m = L_m(\epsilon^{-1}{\bf t}, \epsilon{\pal/\pal{\bf t}})\nn\\
&&\hskip -1.2truecm
=\epsilon^2 \sum a_m^{\alpha,p; \beta,q}
{\pal^2\over \pal t^{\alpha,p} \pal t^{\beta,q}} 
+ \sum {b_m}^{\alpha,p}_{\beta,q}\,  t^{\beta,q} {\pal\over \pal
t^{\alpha,p}} 
+ \epsilon^{-2} c^m_{\alpha,p; \beta,q} \, t^{\alpha,p}\,  t^{\beta,q}
\label{2-11-32}
\eeqa
for some constant coefficients $a_m^{\alpha,p; \beta,q}$, 
${b_m}^{\alpha,p}_{\beta,q}$,
$c^m_{\alpha,p; \beta,q}$ depending on $m$ and on 
the spectrum of the Frobenius manifold.

\begin{exam} For $n=1$ (the KdV theory) the Virasoro operators 
(\ref{2-11-32})
coincide, up to normalization of the independent variables, with the well
known
in the KdV theory \cite{DVV, FKN, KacSchwarz} realization of the Virasoro
algebra
\eqa\label{kdv-vir}
&&
L_m={\epsilon^2\over 2} \sum_{k+l=m-1} {(2k+1)!!\,(2l+1)!!\over 2^{m+1}}
{\pal^2\over \pal t_k \pal t_l}\nn\\
&&\qquad + \sum_{k\geq 0} {(2k+2m+1)!!\over 2^{m+1}
(2k-1)!!}\, t_k\,{\pal\over\pal t_{k+m}}+{1\over 16} \delta_{m,0}, 
\quad m\geq 0,
\nn\\
&&
L_{-1}=\sum_{k\geq 1} t^k{\pal\over \pal t_{k-1}} 
+{1\over 2\epsilon^2} t_0^2,
\nn\\
&&
L_{-m} ={1\over 2\,\epsilon^2}
\sum_{k+l=m-1} {2^{m-1}\over (2k-1)!!(2l-1)!!} t_k\, t_l\nn\\
&&\qquad
+\sum_{k\geq 0} {2^{m-1} (2k+1)!!\over (2m+2k-1)!!}\, t_{k+m} 
{\pal\over \pal
t_k},\quad m>1.
\eeqa
\end{exam}

\begin{exam} For the Frobenius manifold (\ref{potential-toda}) 
(the ${\bf CP}^1$-model)
the regularized operators (\ref{regular-T}) read
\eqa
&&
 L_m={\epsilon^2\over 2} \sum_{k=1}^{m-1} k!\, (m-k)! 
{\pal^2\over \pal t^{2,k-1}
\pal t^{2,k-m-1}}\nn\\
&&+\sum_{k\geq 1} {(m+k)!\over (k-1)!}\left(
t^{1,k} {\pal\over \pal t^{1, m+k}}
+ t^{2,k-1}{\pal\over \pal t^{2,k-1}}\right)
+2\,\sum_{k\geq 0} \alpha_m(k) t^{1,k} {\pal \over \pal t^{2, m+k-1}}, 
\quad m>0
\nn\\
&&
L_0 =\sum_{k\geq 1} k\, \left( t^{1,k} {\pal\over \pal t^{1,k}}+t^{2,k-1}
{\pal\over \pal t^{2,k-1}}\right)
\nn\\
&&
L_{-1}=\sum_{k\geq 1} t^{\alpha,k} {\pal\over \pal t^{\alpha, k-1}}+{1\over
\epsilon^2 } t^{1,0} t^{2,0}
\nn\\
&&
L_{-m} = {1\over \epsilon^2} \sum_{k=1}^{m-1} {t^{1,k}\, t^{1, m-k}\over
(k-1)!(m-k-1)!},\quad m>1.
\nn\\
\eeqa
Here the integer coefficients $\alpha_m(k)$ are defined by
$$
\alpha_m(0)= m!, ~~\alpha_m(k) = {(m+k)!\over (k-1)!} \sum_{j=k}^{m+k} {1\over
j}, ~k>0.
$$
\end{exam}

\begin{lemma} The generating symmetry flow (\ref{2-11-10a}) 
can be represented as follows
\beq
{1\over \epsilon^2}{\partial f\over \pal s} = 
\left[ T(\epsilon^{-1}\bar{\bf t}, \epsilon{\pal/\pal{\bf t}};\lambda) 
\exp({f\over\epsilon^{2} })\right]_+ + O(1),
\eeq
$$
\bar {\bf t}= (\bar t^{\alpha,p}).
$$
\end{lemma}

The proof is straightforward.

\begin{remark} After this paper was finished a very interesting work of A.
Givental appeared \cite{givental1}. In particular, an elegant realization
of our Virasoro operators $L_m$ with $m\geq -1$ was obtained in
\cite{givental1}.
\end{remark}

\subsubsection{Virasoro symmetries and solutions of the Principal Hierarchy}

Let $v=v(x,{\bf t}, c(\epsilon))$ be a solution to the Principal Hierarchy specified by
the   series  $c(\epsilon) =(c^{\alpha,p}(\epsilon))$ in $\epsilon$ with constant coefficients as in 
(\ref{formalsolution}).
Put
\beq
f(x,{\bf t},c(\epsilon)) ={\cal F}_0(x,{\bf t})=\log\, \tau, ~~ f_{\alpha,p}(x,{\bf t},c(\epsilon))=
\pal_{t^{\alpha,p}} \log\,\tau.
\eeq  
Here the tau-function and its first derivatives for the solution $v=v(x,{\bf t}, c(\epsilon))$
are defined by the formulae (\ref{F0}), (\ref{derF0}).

\begin{theorem}\label{t-2-11-10} 
The functions $v=v(x,{\bf t}, c(\epsilon))$, $f=f(x,{\bf t},c(\epsilon))$,
\newline
$f_{\alpha,p}=f_{\alpha,p}(x,{\bf t},c(\epsilon))$ are the stationary points of the following symmetries
of the Principal Hierarchy
\eqa
&&
{\pal f\over \pal \tilde s} = \left[ -{1\over 2} {\pal S_\alpha(\lambda, {\bf \tilde t})\over \pal \lambda} 
G^{\alpha\beta}
{\pal S_\beta(\lambda, {\bf\tilde t})\over \pal \lambda}\right]_+ =0
\label{2-11-37a}\\
&&
{\pal f_{\gamma,p}\over \pal \tilde s } = {\pal\over \pal t^{\gamma,p}} 
\left[ -{1\over 2} {\pal S_\alpha(\lambda, {\bf \tilde t})\over \pal \lambda} 
G^{\alpha\beta}
{\pal S_\beta(\lambda, {\bf\tilde t})\over \pal \lambda}\right]_+ =0
\label{2-11-37b}\\
&&
{\pal v_\gamma\over \pal \tilde s} = \left[ \pal_x \pal_\gamma p_\alpha \,  G^{\alpha\beta} {\pal
S_\beta(\lambda, {\bf\tilde t})\over \pal \lambda}\right]_+ - \eta _{\gamma\epsilon} 
\left((E-\lambda\,
e)^{-1}\right)^\epsilon =0
\label{2-11-37c}
\eeqa
where
${\bf \tilde t} =(\tilde t^{\alpha,p})$,
\beq
\tilde t^{\alpha,p} = t^{\alpha,p}-c^{\alpha,p}(\epsilon).
\eeq
\end{theorem}

Observe that the difference between the ``shifted symmetries'' 
(\ref{2-11-37c}) and the original ones (\ref{2-11-10c})
is a linear combination of the flows of the Principal Hierarchy. In other words, the solution $v=v(x, {\bf t},
c(\epsilon))$ satisfies the following infinite family of constraints
\beq\label{2-11-39}
{\pal v\over \pal s} = {b_m}^{\alpha,p}_{\beta,q}\, c^{\beta,q}\,
{\pal v\over \pal t^{\alpha,p}}, ~~m\geq -1.
\eeq

The theorem was proved in \cite{selecta} for the particular case $c^{\alpha,p}=\delta^\alpha_1 \,\delta^p_1$.
The proof can be repeated also in the general case without major changes.
The crucial point in the proof is the identity
\beq\label{omega-sym}
\pal_{E^{m+1}}\Omega_{\al,p;\beta,q}
=2\sum a^{\lm,k;\epsilon,l}_m\,\Omega_{\al,p;\lm,k}\,
\Omega_{\epsilon,l;\beta,q}+b^{\lm,k}_{m;\al,p}\,\Omega_{\lm,k;\beta,q}+
b^{\lm,k}_{m;\beta,q}\,\Omega_{\lm,k;\al,p}+2\,c^m_{\al,p;\beta,q}.
\eeq
valid for any $m\ge -1$. 

The expanded form of the stationary equations (\ref{2-11-37a}) reads 
\beq 
\sum a_m^{\alpha,p;\beta,q} {\pal f\over \pal t^{\alpha,p}} {\pal f\over \pal t^{\beta,q}}
+ \sum {b_m}^{\alpha,p}_{\beta,q} \tilde t^{\beta,q} {\pal f\over \pal t^{\alpha,p}}
+\sum c^m_{\alpha,p;\beta,q} \tilde t^{\alpha,p} \tilde t^{\beta,q} =0, ~~ m\geq -1.
\eeq
The constant coefficients are the same as in (\ref{2-11-32}).
The stationary equations for the first derivatives of the tau-function and for $v_\alpha$ are obtained from
the above by differentiation. 

\subsubsection{Linearization of Virasoro symmetries}
\label{sec-3-10-4}\par

Let us consider a quasitrivial bihamiltonian tau-symmetric hierarchy 
\beq\label{2-11-41}
{\pal w_\alpha\over \pal t^{\alpha,p}} = \{ w_\alpha(x), H_{\alpha,p}\}_1
\eeq
obtained from the Principal 
Hierarchy by a 
quasi-Miura transformation  
\beq\label{2-11-42}
v_\alpha\mapsto w_\alpha=v_\alpha + \pal_x \pal_{t^{\alpha,0}} 
\sum_{k\geq 1} \epsilon^k
{\cal F}^{[k]}(v; v_x, \dots, v^{(n_k)}).
\eeq
(For technical reasons we have changed the notations of the dependent variables
of the hierarchy from $u^1$, \dots, $u^n$ to $w^1$, \dots, $w^n$.) 

\begin{theorem} The flows ${\pal\over \pal\hat s_m}$ defined by the generating function
\eqa
&&
{\pal\over \pal\hat s} =\sum_{m\geq -1}  {1\over \lambda^{m+2}} {\pal\over \pal\hat s_m}
\\
&&
{\pal\log \tau\over \pal \hat s} = {\pal f\over \pal s} + \sum_{k\geq 1} \epsilon^k
\sum_{r=0}^{n_k} {\pal {\cal F}^{[k]}\over \pal v^{\alpha,r}}\,
\pal_x^r {\pal v^\alpha\over \pal s}
\label{2-11-43a}\\
&&
{\pal\over\pal\hat s} {\pal \log \tau\over \pal t^{\alpha,p}} = {\pal \over \pal t^{\alpha,p}}
{\pal\log \tau\over \pal \hat s} 
\label{2-11-43b}\\
&&
{\pal w_\alpha\over\pal\hat s} =\epsilon^2 \pal_x \pal_{t^{\alpha,0}} {\pal\log \tau\over \pal \hat s} 
\label{2-11-43c}\\
\eeqa
are symmetries of the tau-cover of the hierarchy (\ref{2-11-41}). They satisfy the Virasoro commutation relations
(\ref{vc-sym}). Every solution of the hierarchy, its tau-function and the first derivatives of it is a stationary 
point of the symmetries with the shifted times as in 
(\ref{2-11-37a})--(\ref{2-11-37c}).
\end{theorem}

\pf This is obtained by applying the ``change of coordinates'' 
(\ref{2-11-42}) to Theorems \ref{t-2-11-3} and \ref{t-2-11-10}.
\epf

{\bf Definition.} We say that the quasitriviality (\ref{2-11-42}) {\it linearizes the Virasoro symmetries}
if there exists linear differential operators 
\beq\label{hat-L}
\hat L_m=\hat L_m (\epsilon^{-1} {\bf t}, \epsilon {\pal\over\pal {\bf t}}), ~~m\geq -1
\eeq
with coefficients polynomial in $\{\epsilon^{-1} t^{\alpha,p}\}$
such that
\beq
{\pal \tau\over\pal\hat s_m} = \hat L_m \tau, ~~m\geq -1.
\eeq

\begin{exam} The quasi-triviality (\ref{kdv-miura}) 
\eqa
&&v\mapsto u=v+\epsilon^2 \pal_x^2 \Delta f(v', v'', \dots; \epsilon^2)\nn\\
&&=v +{\epsilon^2\over 24} (\log v')'' + \epsilon^4 \left(
{v^{IV}\over 1152\, {v'}^2} - {7\, v'' v'''\over 1920\, {v'}^3}
+{{v''}^3\over 360\, {v'}^4}\right)''+O(\epsilon^6).\label{kdv-miura0}
\eeqa
for the KdV hierarchy 
satisfies the condition of linearization of Virasoro constraints.
\end{exam}

Observe that we have changed the normalization
$$
\epsilon^2 \mapsto -{\epsilon^2\over 2}.
$$

\pf From independence of $\Delta f$ of $v$ and from Lemma \ref{l-2-9-8}
it follows that (\ref{kdv-miura0}) transforms the Galilean symmetry
\beq\label{gal-riem}
{\pal v\over \pal s_{-1}} = 1+\sum_{p\geq 1} t_p {\pal v\over \pal t_{p-1}}
\eeq
of the Riemann hierarchy to the Galilean symmetry
\beq\label{gal-kdv}
{\pal u\over \pal s_{-1}} = 1+\sum_{p\geq 1} t_p {\pal u\over \pal t_{p-1}}
\eeq
of the KdV hierarchy. For higher symmetries we use, due to Theorem
\ref{t-2-9-10}, that the quasi-Miura (\ref{kdv-miura0}) transforms the recursion
$$
{\pal v\over \pal s_m} = \left( v+{1\over 2} v' \pal_x^{-1} \right) 
{\pal v\over
\pal s_{m-1}}, ~~m\geq 0
$$
for the symmetries of the Riemann hierarchy to the recursion
\beq\label{dvv}
{\pal u\over \pal s_m} = \left({1\over 8} \epsilon^2 \pal_x^2+ u
+{1\over 2} u' \pal_x^{-1} \right) {\pal u\over
\pal s_{m-1}}, ~~m\geq 0
\eeq
for the symmetries of the KdV hierarchy. As it was shown in \cite{DVV}, the
symmetries generated by the recursion procedure (\ref{dvv}), (\ref{gal-kdv})
have the form
\beq\label{dvv1}
{\pal u\over \pal s_m} = \epsilon^2 \pal_x^2 {L_m \tau\over \tau}, ~~m\geq -1
\eeq
where the Virasoro operators $L_m$ are defined in (\ref{kdv-vir}). Observe that
our normalization of the flows of the KdV hierarchy differs from that of
\cite{DVV}. Our flows satisfy the recursion relation
$$
\left(p+{1\over 2}\right) {\pal u\over \pal t_p} =\left[ {\epsilon^2\over 8}
\pal_x^2 + u +{1\over 2} u' \pal_x^{-1}\right] {\pal u\over \pal t_{p-1}},
~~p\geq 1.
$$
\epf

Let us assume that the quasitriviality transformation (\ref{2-11-42}) is such that $F_1=0$. This can always be done
according to Theorem \ref{t-2-10-5}. 

\begin{lemma} The operators $\hat L_m$ in (\ref{hat-L}) must have the form
\beq
\hat L_m = L_m + \kappa_0 \, \delta_{m,0}, ~~m\geq -1
\eeq
where the Virasoro operators $L_m$ were defined in Section \ref{sec-3-10-2}
and $\kappa_0$ is a constant.
\end{lemma}

\pf By definition the principal symbol of the operator $\hat L_m$ coincides with
$$
\sum a_m^{\alpha,p;\beta,q} {P_{\alpha,p}} {P_{\beta,q}}
+ \sum {b_m}^{\alpha,p}_{\beta,q} \tilde t^{\beta,q} {P_{\alpha,p}}
+\sum c^m_{\alpha,p;\beta,q} \tilde t^{\alpha,p} \tilde t^{\beta,q}, ~m\geq -1.
$$
Here the variables $P_{\alpha,p}$ correspond to ${\pal\over \pal t^{\alpha,p}}$.
Due to vanishing of the terms of the order $\epsilon^{-1}$ in the expansion of the tau-function
we conclude that
$$
\hat L_m = \epsilon^2 \sum a_m^{\alpha,p; \beta,q}
{\pal^2\over \pal t^{\alpha,p} \pal t^{\beta,q}} + \sum {b_m}^{\alpha,p}_{\beta,q} \bar t^{\beta,q} 
{\pal\over \pal
t^{\alpha,p}} + \epsilon^{-2} c^m_{\alpha,p; \beta,q} \bar t^{\alpha,p} \bar t^{\beta,q}
+\kappa_m=L_m +\kappa_m
$$
for some constants $\kappa_m$. Using the Virasoro commutation relations
$$
\left[ {\pal\over\hat \pal s_i}, {\pal\over\pal\hat s_j}\right] = (j-i) {\pal\over\pal\hat s_{i+j}}, 
~~i,j\geq
-1
$$
and
$$
[L_i, L_j] =(j-i) L_{i+j}, ~~i,j\geq -1
$$
we derive that $\kappa_m=0$ for $m\neq 0$. \epf

Let $v=v(x,{\bf t},c(\epsilon))$ be a solution to the Principal Hierarchy
described in Section \ref{sec-3-6-4}. We call it {\it admissible} w.r.t. the
quasitriviality transformation if all the denominators of the rational functions
$\pal_x \pal_{t^{\alpha,0}} 
{\cal F}^{[k]}(v; v_x, \dots, v^{(n_k)})$ do not vanish at the point
$x=0$, ${\bf t}=0$. We will see below that any monotone solution
to the Principal Hierarchy is admissible for our class of quasitriviality
transformations.
 
Denote $w=w(x,{\bf t}, c(\epsilon))$  the solution to the full hierarchy 
(\ref{2-11-41}) obtained from an admissible
solution $v=v(x,{\bf t},c(\epsilon))$  by the quasitriviality transformation
(\ref{2-11-42}). All these will be called {\it admissible solutions to the full
hierarchy}. Applying to an admissible solution $w(x,{\bf t}, c(\epsilon))$
the quasi-Miura inverse to (\ref{2-11-42}) one obtains a solution $v(x,{\bf t},
c(\epsilon))$ of the form (\ref{formalsolution}).
The tau-function of the solution $w(x,{\bf t}, c(\epsilon))$ has the form
\beq\label{tau-w}
\tau(x,{\bf t}, c(\epsilon)) = \exp \, [ \epsilon^{-2} {\cal F}_0(x,{\bf t}, c(\epsilon))
+ \sum_{k\geq 1} \epsilon^{k-2}
{\cal F}^{[k]}(v; v_x, \dots, v^{(n_k)})|_{v=v(x,{\bf t},c(\epsilon))}]
\eeq
where the function ${\cal F}_0(x,{\bf t}, c(\epsilon))$ has been defined 
in (\ref{F0}). 

\begin{lemma} The tau-function (\ref{tau-w}) of any admissible  solution 
\newline
$w=w(x,{\bf t}, c(\epsilon))$ to the full hierarchy satisfies the following system of 
{\rm Virasoro constraints}
\beq\label{vc-2-11}
\hat L_m(\epsilon^{-1} {\bf \tilde t}, \epsilon {\pal\over \pal {\bf t}})\,
  \tau(x,{\bf t}, c(\epsilon)) =0,
~~m\geq -1
\eeq
where 
$$
\tilde t^{\alpha,p} = t^{\alpha,p} -c^{\alpha,p}(\epsilon) + x\, 
\delta^\alpha_1 \delta^p_0.
$$
\end{lemma}

\pf Substituting the expansion (\ref{tau-w}) 
into the definition of the linearization
$$
\hat L_m \tau = {\pal \tau\over \pal \hat s_m}
$$
and collecting the coefficients of $\epsilon^{k-2}$ we obtain, for every $k\geq 0$
\eqa
&&
\sum_{l=0}^k a_m^{\alpha,p;\beta,q}\, 
{\pal {\cal F}^{[l]}\over \pal t^{\alpha,p}}\,  
{\pal {\cal F}^{[k-l]}\over \pal t^{\beta,q}} + a_m^{\alpha,p;\beta,q}\, 
{\pal^2 {\cal F}^{[k-2]}\over \pal t^{\alpha,p}\pal t^{\beta,q}} 
+{b_m}_{\beta,q}^{\alpha,p}\, t^{\beta,q}\,
{\pal {\cal F}^{[k]}\over \pal t^{\alpha,p}} 
\nn\\
&&
+\delta_{k,0}\,\sum c^m_{\alpha,p; \beta,q}\, \bar t^{\alpha,p}\,
 \bar t^{\beta,q} 
+\kappa_0\, \delta_{m,0}\, 
\delta_{k,2}={\pal {\cal F}^{[k]}\over \pal s_m}.
\eeqa
In this equation we identify ${\cal F}_0 = :{\cal F}^{[0]}$.
The equality with $k=0$ holds true due to the definition of the 
symmetry (\ref{2-11-43a})--(\ref{2-11-43c}). For $k>0$
we can use (\ref{2-11-39}) in order to rewrite the r.h.s. in the form
$$
{\pal {\cal F}^{[k]}\over \pal s_m}
=\sum {\pal {\cal F}^{[k]}\over \pal v^{\gamma,r}}\, \pal_x^r 
{\pal v^\gamma\over \pal s_m} = \sum {\pal {\cal F}^{[k]}\over \pal v^{\gamma,r}} \pal_x^r 
{b_m}^{\alpha,p}_{\beta,q}\, c^{\beta,q}\,{\pal v^\gamma\over \pal t^{\alpha,p}}=
{b_m}^{\alpha,p}_{\beta,q}\, c^{\beta,q}\, {\pal {\cal F}^{[k]}\over \pal t^{\alpha,p}}.
$$
This proves the lemma.\epf

From the above statements it follows

\begin{theorem} The logarithm of the tau-function 
$$
{\cal F} : = \log \tau
$$
of any admissible solution of the full hierarchy satisfies the following
system of equations
\eqa
&&
\sum_{p,q>0}\left[ {\epsilon^2\over 2} 
{\pal\over\pal t^{p-{1/2}}}M^p_q(\lambda) {\pal\over\pal t_{q-{1/2}}} {\cal F}
+{\pal{\cal F}\over\pal t^{p-{1/2}}}M^p_q(\lambda)
{\pal{\cal F}\over\pal t_{q-{1/2}}}\right.
\nn\\
&&
\left. (-1)^{p-{1\over 2}} \tilde t^{p-{1\over 2}} M^q_{-p} (\lambda)
{\pal{\cal F}\over\pal t_{q-{1/2}}}
+{1\over 2\epsilon^2} (-1)^{p+q+1}\tilde t^{p-{1\over 2}} M^{-q}_{-p}(\lambda)
\tilde t_{q-{1\over 2}}+{\kappa_0\over \lambda^2}\right]_+=0
\nn\\ \label{loop-eq}
\eeqa
for any $\lambda$.
\end{theorem}

\begin{exam} For $n=1$ the equation (\ref{loop-eq}) coincides with the loop
equation in topological gravity \cite{DVV, FKN} (for the particular solution
with $c_j = \delta_{j,1}$) or in the double scaling limit \cite{BK, DS, GM}
of the Hermitean matrix model at the $k$-th multicritical point (for the
particular solution with $c_j=\delta_{j, k+1}$, $k\geq 1$) \cite{David90}
(see also \cite{AmbjornMakeenko, Makeenko}).
\end{exam}

Motivating by this example, we introduce

{\bf Definition.} The equation (\ref{loop-eq}) is called {\it generalized loop
equation} in the theory of integrable hierarchies.
\smallskip
 
 In Section \ref{sec-3-10-5} we will develop a technique to obtain a universal
``perturbative
solution'' to the loop equation. To this end we need first to recall
some technical tricks of using canonical coordinates on semisimple 
Frobenius manifolds.
 
\subsubsection{Working with Frobenius manifolds in the canonical
coordinates}\label{sec-3-10-4a}\par

In this Section we summarize, following \cite{D92, D3} 
a very efficient technique of working with
semisimple Frobenius manifolds based on introduction of {\it canonical
coordinates}.
 
Let $M$ be a semisimple Frobenius manifold. Denote $M_{s\, s}\subset M$
the open dense subset in $M$ consisting of all points $v\in M$ s.t.
the operator of multiplication by the Euler vector field
$$
E(v)\cdot \, : T_v M\to T_v M
$$
has simple spectrum. Denote $u_1(v)$, \dots, $u_n(v)$ the eigenvalues
of this operator, $v\in M_{s\, s}$. The mapping
$$
M_{s\, s} \to \left({\mathbb C}^n\setminus \cup_{i<j}(u_i=u_j) \right)
/S_n, ~~v\mapsto (u_1(v), \dots, u_n(v))
$$
is an unramified covering. Therefore one can use the eigenvalues as local
coordinates on $M_{s\, s}$. The vectors $\pal/\pal u_i$, $i=1, \dots, n$
are basic idempotents of the algebra $T_vM$ for any $v\in M_{s\, s}$
$$
{\pal \over \pal u_i} \cdot {\pal \over \pal u_j} = 
\delta_{ij} {\pal \over \pal u_i}.
$$
Orthogonality
$$
\left< {\pal \over \pal u_i} , {\pal \over \pal u_j}\right> =0, ~i\neq j
$$
readily follows from the above multiplication table.
We call $u_1$, \dots, $u_n$ canonical coordinates on $M_{s\, s}$.
We will never use summation over repeated indices when working in the canonical
coordinates.

Choosing locally branches of the square roots
\beq\label{psi-1}
\psi_{i1}(u):= \sqrt{<\pal/\pal u_i, \pal/\pal u_i>}, ~i=1, \dots, n
\eeq
we obtain a transition matrix $\Psi =(\psi_{i\alpha}(u))$
from the basis $\pal /\pal v^\alpha$ to the
orthonormal basis
\beq\label{basis}
\psi_{11}^{-1}(u) {\pal \over \pal u_1} , \psi_{21}^{-1}(u) {\pal \over \pal u_2} ,
\dots , \psi_{n1}^{-1}(u) {\pal \over \pal u_n} 
\eeq
of the normalized idempotents 
\beq\label{mat-psi}
{\pal\over \pal v^\alpha} =\sum_{i=1}^n {\psi_{i\alpha}(u)\over
\psi_{i1}(u)}{\pal\over \pal u_i}.
\eeq
Equivalently, the Jacobi matrix has the form
\beq
{\pal u_i\over \pal v^\alpha} ={\psi_{i\alpha}\over\psi_{i1}}.
\eeq
The matrix $\Psi(u)$ satisfies orthogonality condition
$$
\Psi^T(u) \Psi(u) \equiv \eta, ~~\eta = (\eta_{\alpha\beta}),
~~\eta_{\alpha\beta}:= \left< {\pal \over \pal v^\alpha}, {\pal \over \pal
v^\beta}\right>.
$$
The lengths (\ref{psi-1}) coincide with the first column of this matrix. So,
the metric $<\ ,\ >$ in the canonical coordinates reads
\beq
<\ ,\ >=\sum_{i=1}^n \psi_{i1}^2(u) du_i^2.
\eeq

Denote $V(u)=\left( V_{ij}(u)\right)$ the matrix of the antisymmetric operator 
${\cal V}$ 
(\ref{2-7-4b}) w.r.t. the orthonormal frame
\beq\label{mat-v}
V(u) := \Psi(u)\, {\cal V}\, \Psi^{-1}(u).
\eeq
It is a solution to the following system of commuting time-dependent
Hamiltonian flows on the Lie algebra $so(n)$ equipped with the standard Lie -
Poisson bracket (\ref{lie-poiss})
\beq\label{eq-per-v}
{\pal V\over \pal u_i} =  \{ V, H_i(V; u)\}, ~~i=1, \dots, n
\eeq
with quadratic Hamiltonians
\beq\label{iso-ham}
H_i(V; u) ={1\over 2} \sum_{j\neq i} {V_{ij}^2\over u_i -u_j}.
\eeq
The matrix $\Psi(u)$ satisfies
\beq\label{eq-psi}
{\pal\Psi\over \pal u_i} = V_i(u) \Psi, ~~ V_i(u):= {\rm ad}_{E_i} {\rm ad}_U
^{-1} (V(u)).
\eeq
Here $U={\rm diag}\, (u_1, \dots, u_n)$, the matrix unity $E_i$ has the entries
$$
(E_i)_{ab}= \delta_{ai} \delta_{ib}.
$$

The {\it isomonodromic
tau-function} of the semisimple Frobenius manifold is defined by
\beq\label{tau-i}
d\log\tau_I(u) =  \sum_{i=1}^n H_i(V(u); u) du_i.
\eeq
It is an analytic function on a suitable unramified covering of $M_{s\, s}$.

The system (\ref{eq-per-v}) coincides with the equations of isomonodromy
deformations of the following linear differential operator with rational
coefficients
\beq\label{isomono}
{dY\over  dz} = \left( U+{V\over z}\right) \, Y.
\eeq
The latter is nothing but the last
component of the deformed flat connection (\ref{def-con}) written in the
orthonormal frame (\ref{basis}). The integration of (\ref{eq-per-v}),
(\ref{eq-psi}) and, more generally, the
reconstruction of the Frobenius structure can be reduced to a solution
of certain Riemann - Hilbert problem \cite{D4}.

In the canonical coordinates
the intersection form $(~,~)_\lambda$, having in the flat coordinates the Gram matrix
$$
\left( g^{\alpha\beta}(v)-\lambda\,\eta^{\alpha\beta}\right)^{-1},
$$
becomes equal to
\beq\label{2-11-56}
(~,~)_\lambda = \sum_{i=1}^n {\psi_{i1}^2(u)\over u_i-\lambda} du_i^2.
\eeq
The differential equations for the periods $p=p(v;\lambda)$ (i.e., the Gauss -
Manin system in the terminology of \cite{D3}) can be recast into the form
\eqa
&&
{\pal \phi_i\over\pal u_j} = -{V_{ij}\over u_i-u_j} \phi_j, ~~j\neq i
\nn\\
&&
{\pal \phi_i\over\pal u_i} ={1\over \lambda-u_i} \left[{1\over 2} \phi_i
+\sum_s V_{is} \phi_s\right] + \sum_s {V_{is}\over u_i-u_s} \phi_s
\label{2-11-57}
\eeqa
where the functions $\phi_i=\phi_i(v;\lambda)$ are defined by
\beq
{\pal p(v(u);\lambda)\over \pal u_i} =\psi_{i1}(u) \phi_i(v(u);\lambda).
\eeq
These functions also satisfy the following equations which will
be used later:
\beq\label{2-11-58b}
\frac{\pal \phi_i}{\pal \lm}=\frac{\phi_i}{2(u_i-\lm)}+\sum_k \frac{V_{ik}\,\phi_k}
{u_i-\lm}.
\eeq
For the basis of periods $p_\alpha=p_\alpha(v;\lambda)$ such that
$$
G^{\alpha\beta} =\left( {\pal\over\pal p_\alpha}, {\pal\over\pal
p_\beta}\right)_\lambda
$$
we put
\beq
{\pal p_\alpha(v(u);\lambda)\over \pal u_i} =\psi_{i1}(u) \phi_{i\alpha}
(u;\lambda).
\eeq
From (\ref{2-11-56}) and from the tensor law for the metric $(~,~)_\lambda$ the 
orthogonality condition follows
\beq\label{2-11-60}
\phi_{i\alpha}(v(u);\lambda) G^{\alpha\beta}\phi_{j\beta}(v(u);\lambda)
={\delta_{ij} \over u_i-\lambda}.
\eeq

\subsubsection{Loop equation on the jet space of a semisimple Frobenius manifold}\label{sec-3-10-5}\par

In this section we develop the main tool for computing the ``perturbative
solution'' of the loop equation. In particular we
will prove that the loop equation (\ref{loop-eq}) 
uniqely determine
the quasitriviality transformation (\ref{2-11-42}). 

The coefficients ${\cal F}^{[k]}$ for $k>0$ are to be determined from the recursion relations
\eqa
\hskip -1.5truecm
&&
2\,a_m^{\alpha,p;\beta,q} {\pal {\cal F}_0\over \pal t^{\alpha,p}}  
{\pal {\cal F}^{[k]}\over \pal t^{\beta,q}}
+{b_m}_{\beta,q}^{\alpha,p} \tilde t^{\beta,q}{\pal {\cal F}^{[k]}\over \pal t^{\alpha,p}} \nn\\
\hskip -1.5truecm
&&
= -\sum_{l=1}^{k-1} a_m^{\alpha,p;\beta,q} {\pal {\cal F}^{[l]}\over \pal t^{\alpha,p}}  
{\pal {\cal F}^{[k-l]}\over \pal t^{\beta,q}} - a_m^{\alpha,p;\beta,q} 
{\pal^2 {\cal F}^{[k-2]}\over \pal t^{\alpha,p}\pal t^{\beta,q}} 
-\kappa_0 \delta_{m,0} \delta_{k,2}, ~~m\geq -1\label{vc-Fk}
\eeqa
that must hold true for an arbitrary solution $v=v(x,{\bf t},c(\epsilon))$. The l.h.s. of this equation
is a linear differential operator acting on ${\cal F}^{[k]}$. 
We will call it {\it linearized Virasoro constraint}.
The coefficients of the linear 
differential operator
 depend on the choice of the solution
$v=v(x, {\bf t}, c(\epsilon))$.  Let us compute the
generating function (\ref{loop-eq}) of the linearized Virasoro constraints.

\begin{lemma} Let ${\cal F}(v; v_x, \dots, v^{(N)})$ be an arbitrary function 
on the 
jet space, the functions
${\cal F}_0(x,{\bf t}, c(\epsilon))$ be as in (\ref{F0}). Then
\eqa
&&
\sum_{m=-1}^\infty {1\over \lambda^{m+2}} \left[
2\,a_m^{\alpha,p;\beta,q} {\pal {\cal F}_0\over \pal t^{\alpha,p}}  
{\pal {\cal F}\over \pal t^{\beta,q}}
+{b_m}_{\beta,q}^{\alpha,p} \tilde t^{\beta,q}{\pal {\cal F}\over 
\pal t^{\alpha,p}}\right] 
\nn\\
&&
=\sum_{r=0}^N {\pal {\cal F}\over \pal v^{\gamma,r}}
\pal_x^r \left( {1\over E-\lambda}\right)^\gamma
+\sum_{r=1}^N {\pal {\cal F}\over \pal v^{\gamma,r}} \sum_{k=1}^r \left(\matrix{r \cr k\cr}\right)
\pal_x^{k-1} \pal_1 p_\alpha G^{\alpha\beta} \pal_x^{r-k+1} \pal^\gamma p_\beta.
\nn\\
\eeqa
\end{lemma}

\pf According to (\ref{2-11-37a})--(\ref{2-11-37c}), 
the l.h.s. of the last equation can be written in the form
\eqa
&&
\left[  {\pal S_\alpha({\bf \tilde t}, \lambda)\over \pal \lambda} \, G^{\alpha\beta}\,
\int_0^\infty {dz\over \sqrt{z}} \, e^{-\lambda\, z} \left( {\pal {\cal F}\over \pal t^p} z^{p+1} \,
z^{\hat\mu} z^R \right)_\beta\right]_+
\nn\\
&&
=\left[  {\pal S_\alpha({\bf \tilde t}, \lambda)\over \pal \lambda} \, G^{\alpha\beta}\,
\int_0^\infty {dz\over \sqrt{z}} \, e^{-\lambda\, z} \left( \sum_{r=0}^N{\pal {\cal F}\over \pal v^{\gamma,r}}
\pal_x^r {\pal v^\gamma\over \pal t^p}\, z^{p+1} \,
z^{\hat\mu} z^R \right)_\beta\right]_+
\nn\\
&&
=\left[  {\pal S_\alpha({\bf \tilde t}, \lambda)\over \pal \lambda} \, G^{\alpha\beta}\,
\sum_{r=0}^N
{\pal {\cal F}\over \pal v^{\gamma,r}}
\pal_x^r \pal_x \pal^\gamma p_\beta\right]_+
\nn\\
&&
=\left[- \sum_{r=0}^N
{\pal {\cal F}\over \pal v^{\gamma,r}}
\pal_x^r \left( {\pal S_\alpha({\bf \tilde t},\lambda)\over \pal \lambda} \, G^{\alpha\beta} \pal_x\pal_{t_{\gamma,0}}
{\pal S_\beta ({\bf \tilde t},\lambda)\over \pal \lambda} \right)\right.\nn\\
&&
\quad \left. -  \sum_{r=1}^N
{\pal {\cal F}\over \pal v^{\gamma,r}}\sum_{k=1}^r \left(\matrix{r\cr k\cr}\right) \pal_x^k {\pal S_\alpha
({\bf
\tilde t},\lambda)\over \pal\lambda} G^{\alpha\beta}\, \pal_x^{r-k+1}\pal^\gamma p_\beta\right]_+
\nn\\
&&
=\sum_{r=0}^N
{\pal {\cal F}\over \pal v^{\gamma,r}}\, \pal_x^r \left({1\over E-\lambda}\right)^\gamma
+\sum_{r=1}^N
{\pal {\cal F}\over \pal v^{\gamma,r}}\sum_{k=1}^r \left(\matrix{r\cr k\cr}\right) 
\pal_x^{k-1} \pal_1
p_\alpha\,G^{\alpha\beta}\, \pal_x^{r-k+1}\pal^\gamma p_\beta.
\nn\\
\eeqa
The lemma is proved.\epf

We will now rewrite the linearized Virasoro constraint operator in the canonical
coordinates $u_1, \dots, u_n$ on the Frobenius manifold. 

For any integer $r>0$ denote
\beq\label{2-11-61}
K^\gamma_r(u;u_x, \dots, u^{(r)};\lambda) =
\pal_x^r \left({1\over{E-\lm}}\right)^\gamma+
  \sum_{k=1}^r 
\left(\matrix{r\cr k\cr}\right)
\pal_x^{k-1} \pal_1
p_\alpha\,G^{\alpha\beta}\, \pal_x^{r-k+1}\pal^\gamma p_\beta
\eeq
where $p_\alpha=p_\alpha(v(u);\lambda)$, $p_\beta=p(v(u);\beta)$.

\begin{lemma} The differential polynomial $K_r(u;u_x, \dots, u^{(r)};\lambda)$
is a rational function in $\lambda$ with poles of the order $r+1$
at $\lambda=u_1$, \dots,
$\lambda=u_n$. It is regular at $\lambda =\infty$. The coefficient of the
highest order pole $(\lambda-u_i)^{-r-1}$ is equal to
\beq
- \psi_{i1}(u){\psi_i}^\gamma(u)\,(u_i')^r \left[ r! + 2^{-r} \sum_{k=1}^r \left(\matrix{
r\cr k}\right) (2k-3)!! (2r-2k+1)!! \right].
\eeq
\end{lemma}

\pf For $v\in M_{s\, s}$ the vector functions $\phi_i(v;\lambda)=(\phi_1(v;\lambda), \dots,
\phi_n(v;\lambda))^T$ are solutions of a Fuchsian system in $\lambda$
with regular singularities at $\lambda=u_1$,  \dots, $\lambda=u_n$,
$\lambda=\infty$. The monodromy group of this Fuchsian system preserves
the invariant bilinear form (\ref{2-11-60}). 
The monodromy does not depend on the point
of the Frobenius manifold. Therefore every term 
$$
\pal_x^{k-1} \pal_1
p_\alpha\,G^{\alpha\beta}\, \pal_x^{r-k+1}\pal^\gamma p_\beta
$$
in the expression (\ref{2-11-61}) 
is a single valued function in $\lambda$ with
regular singularities at the same points. Hence this is a rational function
with poles at $u_1, \dots, u_n, \infty$. From the expansion (\ref{period1}),
(\ref{period2}) it easily
follows that the pole at infinity disappears. To compute the pole of the highest
order it suffices to derive from (\ref{2-11-57}) that
\beq
\pal_x^m \phi_{i} ={(2m-1)!!\over 2^m} {(u_i')^m\over (\lambda-u_i)^m}\phi_i
+\dots
\eeq
where dots denote the terms with poles of the lower order. Using the
orthogonality we obtain the highest order pole at $\lambda=u_i$ in
\eqa
&&
\pal_x^{k-1}(\sum_i\psi_{i1}(u) \phi_{i\alpha}) \, G^{\alpha\beta}\pal_x^{r-k+1}
(\sum_j{\psi_j}^\gamma(u) \phi_{j\beta})
\nn\\
&&=-
{(2k-3)!! (2r-2k+1)!!\over 2^r } {(u_i')^r\over (\lambda-u_i)^{r+1}}
\psi_{i1}(u) {\psi_i} ^\gamma(u)+\dots
\nn\\
\eeqa
where dots denote poles of lower order. To finish the proof of the lemma
it remains to observe that
\beq
\left({1\over E-\lambda}\right)^\gamma
=\sum_{i=1}^n {\psi_{i1}(u) {\psi_i}^\gamma(u)\over u_i-\lambda}.
\eeq
So, the highest order pole at $\lambda=u_i$ in the first term in 
(\ref{2-11-61})
equals
$$
- {\psi_{i1}(u)\,\psi_i^\gamma(u)\,r!\,(u_i')^r\over (\lambda-u_i)^{r+1}}.
$$
The lemma is proved.\epf

\begin{theorem} The solution of the system of Virasoro constraints 
(\ref{vc-2-11}) 
for a semisimple
Frobenius manifold is unique.
\end{theorem}

\pf It suffices to prove that, for any $N$ the linear homogeneous equation 
\beq
\sum_{r=0}^N
{\pal {\cal F}\over \pal v^{\gamma,r}}\, \pal_x^r 
\left({1\over E-\lambda}\right)^\gamma
+ \sum_{r=1}^N
{\pal {\cal F}\over \pal v^{\gamma,r}}\sum_{k=1}^r \left(\matrix{r\cr k\cr}\right) 
\pal_x^{k-1} \pal_1
p_\alpha\,G^{\alpha\beta}\, \pal_x^{r-k+1}\pal^\gamma p_\beta=0
\eeq
has only trivial solution ${\cal F} = {\cal F}(v; v_x, \dots, v^{(N)})$.
According to the previous Lemma the l.h.s. of the equation is a rational
function in $\lambda$ with poles of the order $N+1$ at the points $\lambda=u_1$,
\dots, $\lambda=u_n$. The highest order pole at the point $\lambda=u_i$ reads
$$
-\sum_\gamma {\pal {\cal F}\over \pal v^{\gamma,N}} 
\frac{\psi_{i1}(u)\,\psi_i^\gamma(u)\,(u_i')^N}{(\lambda-u_i)^{N+1}}
\left[ N! + 2^{-N} \sum_{k=1}^N \left(\matrix{
N\cr k}\right) (2k-3)!! (2N-2k+1)!! \right].
$$
Due to nondegeneracy of the matrix $({\psi_i}^\gamma(u))$ we derive that
$$
{\pal {\cal F}\over \pal v^{\gamma,N}}=0, ~~\gamma=1, \dots, n.
$$
The theorem is proved.\epf

\begin{cor} If the quasitriviality transformation satisfies the assumption
of linearization of the Virasoro symmetries, then it is equivalent, up to
an action of the Miura group, to the transformation of the form
\beq\label{2-11-67}
v_\alpha\mapsto w_\alpha=v_\alpha + \pal_x\pal_{t^{\alpha,0}}
\sum_{g=1}^\infty \epsilon^{2\,g} {\cal F}_g(v; v_x, \dots, v^{(3g-2)})
\eeq
where the functions  ${\cal F}_g$ are uniquely determined from the system of
Virasoro constraints (\ref{vc-2-11}).
\end{cor}

\pf Using Theorem \ref{t-2-10-5}, 
we can kill the first term ${\cal F}^{[1]}$ in the
quasitriviality transformation (\ref{2-11-42}). 
All subsequent trems ${\cal F}^{[k]}$
for odd $k$ must vanish due to the uniqueness  theorem. Theorem \ref{t-2-10-5}
implies that  the highest order of jets in ${\cal F}_g$ is equal to $3\,g-2$.
The corollary is proved.\epf

In the framework of the theory of Gromov - Witten invariants, 
the generating function of the genus $g$ Gromov-Witten invariants
is conjectured to have the form ${\cal F}_g(v, v_x,\dots,v^{3g-2})$
where $v^\al$ are some special two point correlation functions,
it makes part of the so-called {\it Virasoro conjecture}
of T.Eguchi {\it et. al.} \cite{eguchi1}. For the case where Frobenius manifold $M$
coincides with the quantum cohomology of a smooth projective variety $X$
with $H^{\rm odd}(X)=0$ the Virasoro conjecture suggests a system of 
differential equations for the generating function of Gromov - Witten
invarants of $X$ and their descendents that coincides, in this particular case,
with our Virasoro constraints (see details in \cite{selecta}).

\begin{cor} The partial derivatives $\pal {\cal F}_g/\pal u^{i,k}$
with respect to the jet coordinates $u^{i,k}:= \pal_x^k u^i$, 
$k=1, \dots, 3g-2$, are rational functions
of the jets with the denominator  containing only powers of $u^1_x$, \dots,
$u^n_x$.
\end{cor}

\begin{cor} Any monotone solution $w=w(x,{\bf t}, \epsilon)$ satisfying
$w(0,{\bf 0},0)\in M_{s\, s}$ is admissible.
The quasitriviality transformation (\ref{2-11-67}) establishes a
one-to-one correspondence between monotone solutions $v(x,{\bf t}, \epsilon)$
to the Principal Hierarchy satisfying 
and admissible solutions to the full hierarchy.
\end{cor}

Let us now derive the generating function of the remaing part of the system
of Virasoro constraints (\ref{loop-eq}). 
The following statement will be convenient
for such a derivation.

\begin{lemma} Let $b'_{\alpha,p}$ and $b''_{\beta,q}$ be two sets of elements
of a commutative algebra, $\alpha, \beta = 1, \dots, n$, $p, q =
0, 1, \dots$. Denote
$$
b'_p=(b'_{1,p}, \dots, b'_{n,p}), ~~ b''_p=(b''_{1,p}, \dots, b''_{n,p})
$$
and put
\eqa
&&\sigma'_\alpha(\lambda):= \left( \int_0^\infty {dz\over \sqrt{z}} e^{-\lambda\, z} \left(
\sum_{p=0}^\infty b'_p z^p\right) z^{\hat\mu} z^R\right)_\alpha,\nn\\
&&
\sigma''_\beta(\lambda):= \left( \int_0^\infty {dz\over \sqrt{z}} e^{-\lambda\, z} \left(
\sum_{p=0}^\infty b''_p z^p\right) z^{\hat\mu} z^R\right)_\beta.
\nn
\eeqa
Then
\beq\label{2-11-68}
-{1\over 2} {\pal \sigma'_\alpha(\lambda)\over\pal \lambda}
G^{\alpha\beta} {\pal \sigma''_\beta(\lambda)\over\pal \lambda}
= \sum_{m=1}^\infty \sum_{p+q=m-1} {a_m^{\alpha,p; \beta,q}\over \lambda^{m+2}}
b'_{\alpha,p} b''_{\beta,q}.
\eeq
\end{lemma}

Proof follows immediately from the free field realization of the Virasoro
operators.

Let us denote 
$$
\Delta {\cal F}(v; v_x, \dots; \epsilon) := \sum_{k=1}^\infty \epsilon^{k-2} {\cal F}^{[k]}.
$$

\begin{lemma} The loop equation (\ref{loop-eq}) is equivalent to the
following differential  equations for the function $\Delta {\cal F}$ on the jet
space
\eqa
&&
{\pal \Delta{\cal F}\over \pal v^{\gamma,r}}\, \pal_x^r \left({1\over E-\lambda}\right)^\gamma
+\sum_{r\geq 1}
{\pal \Delta {\cal F}\over \pal v^{\gamma,r}}\sum_{k=1}^r \left(\matrix{r\cr k\cr}\right) 
\pal_x^{k-1} \pal_1
p_\alpha\,G^{\alpha\beta}\, \pal_x^{r-k+1}\pal^\gamma p_\beta
\nn\\
&&
= {1\over 2} {\pal p_\alpha(v;\lambda)\over \pal \lambda} *
{\pal p_\beta(v;\lambda)\over \pal \lambda} \, G^{\alpha\beta}
-{\kappa_0\over\lambda^2}
\nn\\
&&
+{\epsilon^2\over 2}\sum \left( {\pal^2 \Delta{\cal F}\over \pal v^{\gamma,k} 
\pal
v^{\rho,l}}+{\pal\Delta {\cal F}\over \pal v^{\gamma,k}} 
{\pal\Delta {\cal F}\over \pal v^{\rho,l}} 
\right)\pal_x^{k+1} \pal^\gamma p_\alpha G^{\alpha\beta} \pal_x^{l+1} \pal^\rho p_\beta
\nn\\
&&
+{\epsilon^2\over 2} \sum  {\pal\Delta {\cal F}\over \pal v^{\gamma,k}} 
\pal_x^{k+1} \left[\nabla{\pal p_\alpha(v;\lambda)\over \pal \lambda} \cdot 
\nabla{\pal p_\beta(v;\lambda)\over \pal \lambda} \cdot v_x\right]^\gamma
G^{\alpha\beta}.
\label{2-11-69}
\eeqa
\end{lemma}

\pf Direct substitution of $\exp ( \epsilon^{-2} {\cal F}_0 + \Delta {\cal F})$
into the system of Virasoro constraints gives
\eqa
&&
{\pal \Delta{\cal F}\over \pal v^{\gamma,r}}\, \pal_x^r \left({1\over E-\lambda}\right)^\gamma
+ \sum_{r\geq 1}
{\pal \Delta {\cal F}\over \pal v^{\gamma,r}}\sum_{k=1}^r \left(\matrix{r\cr k\cr}\right) 
\pal_x^{k-1} \pal_1
p_\alpha\,G^{\alpha\beta}\, \pal_x^{r-k+1}\pal^\gamma p_\beta
\nn\\
&&
+\epsilon^2 \sum{a^{\alpha,p;\beta,q}_m\over \lambda^{m+2}}
\left( {\pal^2 \Delta {\cal F}\over \pal t^{\alpha,p} \pal t^{\beta,q}}
+{\pal \Delta {\cal F}\over \pal t^{\alpha,p}} 
{\pal \Delta {\cal F}\over \pal t^{\beta,q}}\right) 
\nn\\
&&
+\sum {a^{\alpha,p;\beta,q}_m\over \lambda^{m+2}}\Omega_{\alpha,p;\beta,q}(v)
+{\kappa_0\over\lambda^2}=0.
\nn\\
\eeqa
Applying the formula (\ref{2-11-68}) to
$$
b'_{\alpha,p}=b''_{\alpha,p}= {\pal \Delta {\cal F}\over \pal t^{\alpha,p}}
$$
we obtain
$$
\sigma'_\alpha(\lambda)=\sigma''_\alpha(\lambda)
=\{ \Delta {\cal F}(v(x); \dots; \epsilon), \left[\int_0^\infty 
{dz\over z^{3/2}}e^{-\lambda\, z} 
\left( \sum_{p=0}^\infty \bar \theta_p z^p\right) 
z^{\hat\mu} z^R\right]_\alpha\}_1.
$$
So
$$
{\pal\sigma_\alpha(\lambda)\over \pal \lambda} 
= - \{ \Delta{\cal F}(v(x); \dots; \epsilon), \int_0^\infty {dz\over \sqrt{z}}
e^{-\lambda\, z} \left( \int \tilde v_\alpha(v(x); z) dx\right)\}_1
$$
$$
=- \{ \Delta{\cal F}(v(x); \dots; \epsilon), \bar p_\alpha(\lambda)\}_1.
$$
Therefore
$$
\sum {a_m^{\alpha,p; \beta,q}\over \lambda^{m+2}}
{\pal \Delta{\cal F}\over \pal t^{\alpha,p}} 
{\pal \Delta{\cal F}\over \pal t^{\beta,q}}=-{1\over 2} 
\{\Delta{\cal F}(v(x); \dots; \epsilon), \bar p_\alpha(\lambda)\}_1
G^{\alpha\beta}
\{\Delta{\cal F}(v(x); \dots; \epsilon), \bar p_\beta(\lambda)\}_1.
$$
This gives the second term in the third line in the equation (\ref{2-11-69}). 
Similarly, to calculate
the two terms of the expression
$$
\sum {a_m^{\alpha,p; \beta,q}\over \lambda^{m+2} }
{\pal^2 \Delta{\cal F}\over \pal t^{\alpha,p}\pal t^{\beta,q}}
$$
$$
= \sum {a_m^{\alpha,p; \beta,q}\over \lambda^{m+2}}
{\pal^2 \Delta{\cal F}\over \pal v^{\gamma,k}\pal v^{\rho,l}}
\pal_x^k {\pal v^\gamma\over \pal t^{\alpha,p}} \pal_x^l {\pal v^\rho\over \pal
t^{\beta,q}}
+ \sum {a_m^{\alpha,p; \beta,q}\over \lambda^{m+2}}
{\pal\Delta{\cal F}\over\pal v^{\gamma,k}} \pal_x^{k+1}
{\pal \Omega_{\alpha,p;\beta,q}(v)\over \pal t_{\gamma,0}}
$$
we are to use the same trick and also the formula
$$
{\pal \Omega_{\alpha,p;\beta,q}(v)\over \pal t_{\gamma,0}}=
\left[ \nabla \theta_{\alpha,p}\cdot \nabla\theta_{\beta,q}\cdot
v_x\right]^\gamma.
$$
This will give the first term of the third line and also the last line of the 
equation (\ref{2-11-69}).
Finally, applying the same trick to the calculation of
$$
\sum {a_m^{\alpha,p; \beta,q}\over \lambda^{m+2}}{\pal^2 {\cal F}_0\over \pal
t^{\alpha,p} \pal t^{\beta,q}} = 
\sum {a_m^{\alpha,p; \beta,q}\over \lambda^{m+2}}
\theta_{\alpha,p}* \theta_{\beta,q}
$$
we obtain the first term of the second line of the equation. The lemma is proved.
\epf

\begin{exam}
For $n=1$ we have
$$
p=\sqrt{v-\lambda}, ~~ G=4,
$$
$$
\pal_\lambda p \, * \, \pal_\lambda p = {1\over 32\, \lambda^2} - {1\over 32(v-\lambda)^2}.
$$
So the spelling of the loop equation (\ref{loop-eq}) reads
\eqa
&&
\sum_r {\pal \Delta{\cal F}\over \pal v^{(r)}} \pal_x^r {1\over v-\lambda} 
+\sum_{r\geq 1} {\pal \Delta{\cal F}\over \pal v^{(r)}}\sum_{k=1}^r \left(
\matrix{ r \cr k\cr}\right)
\pal_x^{k-1} {1\over \sqrt{v-\lambda}}
\pal_x^{r-k+1} {1\over \sqrt{v-\lambda}}
\nn\\
&&
={1\over 16 \,\lambda^2} -{1\over 16(v-\lambda)^2} -{\kappa_0\over \lambda^2}
\nn\\
&&
\quad
+{\epsilon^2\over 2} \sum\left[ {\pal^2 \Delta{\cal F}\over \pal v^{(k)} \pal v^{(l)}}
+ {\pal \Delta{\cal F}\over \pal v^{(k)}}{\pal \Delta{\cal F}\over \pal v^{(l)}}\right]
\pal_x^{k+1} {1\over \sqrt{v-\lambda}} \pal_x^{l+1} {1\over \sqrt{v-\lambda}}
\nn\\
&&\quad 
-{\epsilon^2\over 16} \sum {\pal \Delta{\cal F}\over \pal v^{(k)}}\pal_x^{k+2} {1\over (v-\lambda)^2}.
\label{2-11-71}
\eeqa
Now we can determine recursively each term of the expansion
$$
\Delta {\cal F} = {\cal F}_1 (v; v_x) + \epsilon^2 
{\cal F}_2 (v; v_x, v_{xx}, v_{xxx}, v_{xxxx})+\dots
$$
substituting it into equation (\ref{2-11-71}). 
For ${\cal F}_1$ we obtain 
$$
{1\over v-\lambda} {\pal {\cal F}\over \pal v} -{3\over 2} {v'\over (v-\lambda)^2} {\pal {\cal F}\over \pal
v'} = {1\over 16 \, \lambda^2} -{1\over 16 (v-\lambda)^2} -{\kappa_0\over \lambda^2}.
$$
This implies that
\beq
\kappa_0={1\over 16}, ~~{\cal F}_1 ={1\over 24} \log v'.
\eeq
For the next term
$$
{\cal F} :={\cal F}_2 (v; v', v'', v''', v^{IV})
$$
we obtain from (\ref{2-11-71})
\eqa
&&
{1\over (v-\lambda)^5} \left( {105 \over 2048}{v'}^2 
-{945\over 16} {v'}^4 {\partial {\cal F}\over \partial v^{IV}}\right)
\nn\\
&&
+{1\over (v-\lambda)^4} \left( -{49\over 1536} v'' +{735\over 8} {v'}^2 v''
{\partial {\cal F}\over \partial v^{IV}}+{105\over 8} {v'}^3 
{\partial {\cal F}\over \partial v'''}\right)
\nn\\
&&
+{1\over (v-\lambda)^3} \left[{1\over 192} {v'''\over v'} -{23\over 4608} {{v''}^2\over {v'}^2}
-\left(16 {v''}^2 +{87\over 4} v' v'''\right)
{\partial {\cal F}\over \partial v^{IV}}-{55\over 4} v' v'' 
{\partial {\cal F}\over \partial v'''} \right.\nn\\
&&\left.-{15\over 4} {v'}^2 
{\partial {\cal F}\over \partial v''}\right]
+{1\over (v-\lambda)^2}\left( 3 v^{IV} {\partial {\cal F}\over \partial v^{IV}}
+{5\over 2} v''' {\partial {\cal F}\over \partial v'''}
+2 v'' {\partial {\cal F}\over \partial v''}
+{3\over 2} v' {\partial {\cal F}\over \partial v'}\right)\nn\\
&&
-{1\over v-\lambda} {\partial {\cal F}\over \partial v}=0.
\eeqa
Solving this system we easily obtain
\beq
{\cal F}_2 = {v^{IV}\over 1152\, {v'}^2} - {7\, v'' v'''\over 1920\, {v'}^3}
+{{v''}^3\over 360\, {v'}^4}.
\eeq
\end{exam}

\medskip
\begin{exam}
For the two-dimensional Frobenius manifold (\ref{potential-toda})
(i.e., for the ${\bf CP}^1$ sigma-model) one can choose the following system
of independent periods
\eqa
&&
p_1= v_2 - 2 \,\log \left(v_1-\lambda +  
\sqrt{ (v_1-\lambda)^2 - 4\, \exp{v_2}}
\right)
\nn\\
&&
p_2= v_2.
\nn\\
\eeqa
The Gram matrix is equal to
$$
\left( G^{\alpha\beta}\right) ={1\over 2} \left( \matrix{ -1 & 0\cr 0 & 1
\cr}\right).
$$
A simple calculation gives
$$
\pal_\lambda p_\alpha * \pal_\lambda p_\beta G^{\alpha\beta}
= {2 e^{v_2} (4\, e^{v_2} + (v_1-\lambda)^2)\over 
(4\, e^{v_2} - (v_1-\lambda)^2)^3}.
$$
So the loop equation (\ref{loop-eq}) reads
\eqa
&&
\sum_{r\geq 0}\left({\pal \Delta{\cal F}\over \pal v_1^{(r)}} \pal_x^r
{v_1-\lambda\over D} - 2 {\pal \Delta{\cal F}\over \pal v_2^{(r)}} 
\pal_x^r {1\over D}\right)
\nn\\
&&
+\sum_{r\geq 1} \sum_{k=1}^r \left(\matrix{r\cr k}\right) \pal_x^{k-1}
{1\over \sqrt{D}}\left( {\pal \Delta{\cal F}\over \pal v_1^{(r)}}
\pal_x^{r-k+1} {v_1-\lambda\over \sqrt{D}} -2 
{\pal \Delta{\cal F}\over \pal v_2^{(r)}}\pal_x^{r-k+1} {1\over \sqrt{D}}\right)
\nn\\
&&
=D^{-3} e^{v_2} \left( 4\, e^{v_2}+(v_1-\lambda)^2\right) -{\kappa_0^2\over
\lambda^2}
\nn\\
&&
+\sum_{k,l}{\epsilon^2\over 4} \left[ -\left( 
{\pal^2 \Delta{\cal F}\over \pal v_1^{(k)}\pal v_1^{(l)}}+
{\pal \Delta{\cal F}\over \pal v_1^{(k)}}
{\pal \Delta{\cal F}\over \pal v_1^{(l)}}\right)
\, \pal_x^{k+1} {v_1-\lambda\over \sqrt{D}} \pal_x^{l+1} {v_1-\lambda\over
\sqrt{D}}\right.
\nn\\
&&
+4\, \left( 
{\pal^2 \Delta{\cal F}\over \pal v_1^{(k)}\pal v_2^{(l)}}+
{\pal \Delta{\cal F}\over \pal v_1^{(k)}}
{\pal \Delta{\cal F}\over \pal v_2^{(l)}}\right)\,\pal_x^{k+1} {v_1-\lambda\over
\sqrt{D}} \pal_x^{l+1} {1\over \sqrt{D}}
\nn\\
&&
\left. -4\, \left( 
{\pal^2 \Delta{\cal F}\over \pal v_2^{(k)}\pal v_2^{(l)}}+
{\pal \Delta{\cal F}\over \pal v_2^{(k)}}
{\pal \Delta{\cal F}\over \pal v_2^{(l)}}\right)\,
\pal_x^{k+1} {1\over \sqrt{D}} \pal_x^{l+1} {1\over \sqrt{D}}\right]
\nn\\
&&
-{\epsilon^2\over 2}\sum_k \left\{ {\pal \Delta{\cal F}\over \pal v_1^{(k)}}
\pal_x^{k+1} e^{v_2}{ -[(v_1-\lambda)^2 + 4\, e^{v_2}]\, v_1' +4 \, e^{v_2}
(v_1-\lambda) v_2'\over D^3}
\right.
\nn\\
&&
\left.  +{\pal \Delta{\cal F}\over \pal v_2^{(k)}} \pal_x^{k+1}
e^{v_2} {4 \, (v_1-\lambda) \, v_1' - [ (v_1-\lambda)^2 + 4\, e^{v_2}]\, v_2'
\over D^3}\right\}
\nn\\
\eeqa
where
$$
D= (v_1-\lambda)^2 -4\, e^{v_2}.
$$
\end{exam}

In the next sections we will solve the system of Virasoro constraints for low
genera $g\leq 2$ and compare this solution with the topological one.

\subsubsection{Genus one case and the final form of the loop equation}
\par

\begin{lemma} For an arbitrary semisimple Frobenius manifold and an arbitrary
system of independent periods $p_\alpha(v; \lambda)$, $\alpha=1, \dots, n$,
$G^{\alpha\beta}=(\pal/\pal p_\alpha, \pal/\pal p_\beta)_\lambda$, the following
identity holds true
\eqa
&&\pal_\lambda p_\alpha (v(u);\lambda)*\pal_\lambda p_\beta (v(u); \lambda)
G^{\alpha\beta}\nn\\
&& = -{1\over 8} \sum_{i=1}^n {1\over (\lambda-u_i)^2}
+\sum_{i<j} {V_{ij}^2\over (\lambda-u_i)(\lambda-u_j)} +{1\over 2\, \lambda^2}
\tr\left( {1\over 4} -\hat\mu^2\right).
\eeqa
\end{lemma}

\pf Introducing, as above, the functions $\phi_{i\,\alpha}(v; \lambda)$ by
$$
\pal_i p_\alpha(v(u); \lambda) =\psi_{i\,1}(u) \phi_{i\,\alpha}(v(u);\lambda)
$$
we obtain
$$
\pal_i (\pal_\lambda p_\alpha * \pal_\lambda p_\beta G^{\alpha\beta})
= \pal_\lambda \phi_{i\, \alpha} \pal_\lambda \phi_{i\,\beta} G^{\alpha\beta}.
$$
Using differential equations (\ref{2-11-58b}) we rewrite the r.h.s. as
$$
= {1\over (u_i-\lambda)^2} \left[ {1\over 2} \phi_{i\, \alpha} +
\sum_k V_{ik} \phi_{k\, \alpha} \right] G^{\alpha\beta}
\left[ {1\over 2} \phi_{i\, \beta} +
\sum_l V_{il} \phi_{l\, \beta} \right]
$$
$$
= {1\over 4 (u_i-\lambda)^3} +\sum_j {V_{ij}^2\over (u_i-\lambda)^2 
(u_j-\lambda)}.
$$
To derive the last formula we have used the orthogonality condition 
(\ref{2-11-60}).
Using the differential equations for the matrix $V_{ij}$ 
\cite{D3, D4} we integrate
the last formula to obtain
\beq\label{2-11-78}
\pal_\lambda p_\alpha * \pal_\lambda p_\beta G^{\alpha\beta}
=-{1\over 8} \sum_i {1\over (u_i-\lambda)^2} + \sum_{i<j} {V_{ij}^2\over
(u_i-\lambda)(u_j-\lambda)} + c(\lambda)
\eeq
where the rational function $c(\lambda)$ is an integration constant. It can have
poles only at $\lambda=\infty$. To determine this integration constant
we will use the basis of the regularized periods (see Section 
\ref{sec-3-6-3} above) 
as follows 
$$
\lim_{\nu\to 0}
\pal_\lambda p_\alpha^{(\nu)} G^{\alpha\beta}(\nu) * \pal_\lambda
p_\beta^{(-\nu)} = -{1\over \pi}\times
$$
$$
\times \sum_{p,q\geq 1} \sum_{r\geq 0}
{\theta_{p-1} \left[ e^{R\,\pal_\nu}\right]_r \left( \Gamma(\hat\mu +\nu
+p-r+{1\over 2}) \cos {\pi(\hat\mu +\nu)} \Gamma(-\hat\mu -\nu+q+{1\over
2})\right) * \theta ^{q-1}\over \lambda^{p+q+r+1}}.
$$
Here
$$
\theta_k =(\theta_{1,k}, \dots, \theta_{n,k}), ~~\theta^k = (\theta^{1,k},
\dots, \theta^{n,k})^t, ~~\theta^{\alpha,k}:=
\eta^{\alpha\beta}\theta_{\beta,k}.
$$
Therefore the integration constant $c(\lambda)$ must be chosen in such a way
that the r.h.s. of (\ref{2-11-78}) = $O(1/\lambda^3)$. Hence
$$
c(\lambda) = {n\over 8\, \lambda^2} -\sum_{i<j}{V_{ij}^2\over \lambda^2}
= {1\over 2\, \lambda^2} \tr \left( {1\over 4} - V^2\right).
$$
To complete the proof it remains to observe that
$$
\tr\, V^2 = \tr \,(\hat\mu +R_0)^2 = \tr\, \hat\mu^2
$$
due to nilpotency of $R_0$ and commutativity $[R_0,\hat\mu]=0$.
\epf

\begin{theorem} For an arbitrary semisimple Frobenius manifold 
the system\newline (\ref{2-11-69}) implies 
\beq\label{F1-fm}
{\cal F}_1= \log {\tau_I(u)\over J^{1/24}(u)} +{1\over 24}\sum_{i=1}^n
\log u_i',
\eeq  
\beq
\kappa_0 ={1\over 4} \tr \left( {1\over 4} -\hat\mu^2\right).
\eeq
Here the isomonodromic tau-function $\tau_I(u)$ of the Frobenius manifold
is defined by (\ref{tau-i}),
$J(u)$ is the Jacobian of the transformation from canonical to the flat
coordinates
\beq
J(u) = \det\left( {\pal v^\alpha\over \pal u_i}\right) = \pm \prod_{i=1}^n
\psi_{i\,1}(u).
\eeq
\end{theorem}

We proved the formula (\ref{F1-fm}) 
in \cite{cmp} using results \cite{getzler1} on topology
of the moduli space $\bar{\cal M}_{1,4}$. Remarkably, the very same formula
follows from our axioms of integrable PDEs!

\pf
We are to solve the following equation for the function ${\cal F}_1={\cal
F}_1(v; v_x)$
\eqa
&&
\sum_{i=1}^n {\pal {\cal F}_1\over \pal u_i} {1\over u_i-\lambda}
-\sum_{i=1}^n {\pal {\cal F}_1\over \pal u_i'} {u_i'\over (u_i-\lambda)^2}
+\sum {\pal {\cal F}_1\over \pal v^\gamma_x} \pal_1 p_\alpha G^{\alpha\beta}
\pal_x \pal^\gamma p_\beta
\nn\\
&&
= -{1\over 16} \sum_{i=1}^n {1\over (\lambda-u_i)^2}
+{1\over 2}\sum_{i<j} {V_{ij}^2\over (\lambda-u_i)(\lambda-u_j)} 
+{1\over 4\, \lambda^2}
\tr\left( {1\over 4} -\hat\mu^2\right) -{\kappa_0\over \lambda^2}.
\nn\\ \label{2-11-83}
\eeqa
Using the formulae
$$
{\pal {\cal F}_1\over \pal v^\gamma_x} =\sum_i {\psi_{i\, \gamma}(u)\over
\psi_{i\, 1} (u)} {\pal{\cal F}_1\over \pal u_i'}
$$
$$
\pal_\sigma p_\alpha =\sum \psi_{k\, \sigma}(u) \phi_{k\, \alpha}(v(u); \lambda)
$$
 we obtain
$$
\sum {\pal {\cal F}_1\over \pal v_x^\gamma} \pal_1 p_\alpha G^{\alpha\beta}
\pal_x \pal^\gamma p_\beta =
\sum{\pal {\cal F}_1\over \pal u_i'} \left( {\psi_{j\, 1}\over \psi_{i\,1}}
{\psi_{i\,\gamma} \pal_x {\psi_j}^\gamma\over u_j-\lambda}
+{\psi_{k\, 1}\over \psi_{i\, 1}} \phi_{k\,\alpha} G^{\alpha\beta} \pal_x
\phi_{i\,\beta}\right).
$$
With the help of the differential equations for the functions $\psi_{i\, \sigma}$ and
$\phi_{k\, \alpha}$ and the orthogonality conditions (\ref{2-11-60}) and 
$$
\psi_{i\, \gamma} {\psi_j}^\gamma = \delta_{ij}
$$
we rewrite the equation (\ref{2-11-83}) in the form
$$
\sum_i {\pal{\cal F}_1\over \pal u_i} {1\over u_i-\lambda}
-{3\over 2} \sum {\pal{\cal F}_1\over \pal u_i'} {u_i'\over (u_i-\lambda)^2}
-\sum {\pal{\cal F}_1\over \pal u_i'}{\psi_{j\, 1}\over \psi_{i\, 1}}
 {V_{ij}\over (u_i-\lambda)(u_j-\lambda)}
$$
$$
= -{1\over 16} \sum_{i=1}^n {1\over (\lambda-u_i)^2}
+{1\over 2}\sum_{i<j} {V_{ij}^2\over (\lambda-u_i)(\lambda-u_j)} 
+{1\over 4\, \lambda^2}
\tr\left( {1\over 4} -\hat\mu^2\right) -{\kappa_0\over \lambda^2}.
$$
The formulae (\ref{F1-fm}) and (\ref{tau-i}) easily follow from the last equation.
\epf

As in (\ref{correl0}) we define ``genus one correlators''
\beq\label{correl1}
\dbl\tau_{p_1}(\phi_{\alpha_1}) \tau_{p_2}(\phi_{\alpha_2})\dots 
\tau_{p_k}(\phi_{\alpha_k})\dbr_1 := 
\epsilon^k {\pal^k {\cal F}_1(v; v_x)\over \pal
t^{\alpha_1, p_1} \pal t^{\alpha_2, p_2} \dots \pal t^{\alpha_k, p_k}}
\eeq 
where instead of $v$, $v_x$ one is to substitute the topological solution
(\ref{fp0}), (\ref{fp1}) and its $x$-derivative.

\begin{theorem} For an arbitrary semisimple Frobenius manifold the genus 1
solution
(\ref{F1-fm}) of the loop equation evaluated on the topological solution
(\ref{fp0}), (\ref{fp1}) satisfies the following identities.

1). The genus one topological recursion relations 
\eqa
&&\dbl\tau_p(\phi_\alpha)\dbr_1 =\dbl
\tau_{p-1}(\phi_\alpha) \tau_0(\phi_\nu)\dbr_0
\eta^{\nu\mu} \dbl\tau_0(\phi_\mu)\dbr_1\nn\\
&&\qquad +{1\over 24} \eta^{\nu\mu}
\dbl\tau_{p-1}(\phi_\alpha) \tau_0(\phi_\nu)\tau_0(\phi_\mu)\dbr_0.
\label{rec-rel-1}
\eeqa

2). The restriction
$$
G(v): = {\cal F}_1(v; v_x)
$$
onto the small phase space $t^{\alpha,p}=0$ for $p>0$ evaluated onto the
topological solution satisfies 
\eqa
&&\sum_{1\le \al_1,\al_2,\al_3,\al_4\le n} z_{\al_1} z_{\al_2}
z_{\al_3} z_{\al_4} \left(3\,c^{\mu}_{\al_1\al_2}\,c^{\nu}_{\al_3\al_4}
\,\frac{\pal^2 G}{\pal v^\mu\pal v^{\nu}}-4\,
c^{\mu}_{\al_1\al_2}\,c^{\nu}_{\al_3\mu}
\,\frac{\pal^2 G}{\pal v^{\al_4}\pal v^{\nu}}\right.
\nn\\
&&
-c^{\mu}_{\al_1\al_2}\,c^{\nu}_{\al_3\al_4\mu}
\,\frac{\pal G}{\pal v^{\nu}}+2\,
c^{\mu}_{\al_1\al_2\al_3}\,c^{\nu}_{\al_4\mu}
\,\frac{\pal G}{\pal v^{\nu}}+\frac16 
c^{\mu}_{\al_1\al_2\al_3}\,c^{\nu}_{\al_4\mu\nu}
\nn\\
&&\left.
+\frac1{24} c^{\mu}_{\al_1\al_2\al_3\al_4}\,c^{\nu}_{\mu\nu}-
\frac14 c^{\mu}_{\al_1\al_2\nu}\,c^{\nu}_{\al_3\al_4\mu}\right)=0.
\label{Getz}
\eeqa   
Here the coefficients $c^\al_{\beta\de}=c^\al_{\beta\de}(v)$ are structure
functions of the Frobenius multiplication on $T_vM$,
$c^{\al\beta\ga\de}=c^{\al\beta\ga\de}(v)$ and $c^\mu_{\alpha\beta\gamma\delta}=
c^\mu_{\alpha\beta\gamma\delta}(v)$
are linear
combinations of the 4th and 5th order derivatives of the potential $F(v)$, 
\beq\label{Def-ofc}
c^{\al\beta}_{\ga\mu}=\eta^{\al{\al'}}\,\eta^{\beta{\beta'}}\,
\frac{\pal^4 F(v)}{\pal v^{\al'}\,\pal v^{\beta'}\,\pal v^\ga\,
\pal v^\mu}, ~~ c^\mu_{\alpha\beta\gamma\delta}=\eta^{\mu\nu} {\pal^5 F(v)\over
\pal v^\nu \pal v^\alpha\pal v^\beta \pal v^\gamma \pal v^\delta}. 
\eeq
It is
understood that all the coefficients of the degree 4 polynomial in $z_1$, \dots,
$z_n$ in (\ref{Getz}) are equal to zero.
\end{theorem}

Proof of the Theorem can be obtained just inverting the arguments of \cite{cmp}
where, vice versa, the formula (\ref{F1-fm}) has been derived starting from
the ``topological'' equations (\ref{rec-rel-1}), (\ref{Getz}).

Similarly to the genus zero relations (\ref{rec-rel-0}) the above equations
(\ref{rec-rel-1}), \newline(\ref{Getz}) can be spelled out as an infinite system
of identities for the genus one and genus zero Gromov - Witten and Mumford -
Morita - Miller classes in the case when the Frobenius manifold comes from
quantum cohomology. For this case
the genus one topological recursion relations (\ref{rec-rel-1}) were derived by
E. Witten in \cite{Witten2},
the equations (\ref{Getz}) were derived by E. Getzler \cite{getzler1}. They are
known to be the defining relations for the genus one classes (see details
in \cite{getzler1}).

We arrive at the main 
\begin{theorem} The tau-symmetric quasitrivial (0,n) Poisson pencil with the 
leading term (\ref{normalfm})
corresponding to a semisimple Frobenius manifold and satisfying the axiom of
linearization is obtained from (\ref{normalfm}) by
the transformation
$$
v_\alpha\mapsto w_\alpha v_\alpha + \epsilon^2 \pal_x \pal_{t^{\alpha,0}} \Delta {\cal
F}(v; v_x, v_{xx}, \dots; \epsilon^2), ~~\alpha=1, \dots, n
$$
where the function 
$$
\Delta {\cal F} = \sum_{g\geq 1} \epsilon^{2g-2} {\cal F}_g (v; v_x, \dots,
v^{(3g-2)})
$$
is uniquely determined from the following universal loop equation
\eqa\label{glavnoe}
&&
{\pal \Delta{\cal F}\over \pal v^{\gamma,r}}\, \pal_x^r \left({1\over E-\lambda}\right)^\gamma
+ \sum_{r\geq 1}
{\pal \Delta {\cal F}\over \pal v^{\gamma,r}}\sum_{k=1}^r \left(\matrix{r\cr k\cr}\right) 
\pal_x^{k-1} \pal_1
p_\alpha\,G^{\alpha\beta}\, \pal_x^{r-k+1}\pal^\gamma p_\beta
\nn\\
&&
= -{1\over 16} \tr \left( {\cal U}-\lambda\right)^{-2} -{1\over 4} \tr \left[
\left( {\cal U}-\lambda\right)^{-1} {\cal V}\right]^2 
\nn\\
&&
+{\epsilon^2\over 2}\sum \left( {\pal^2 \Delta{\cal F}\over \pal v^{\gamma,k} 
\pal
v^{\rho,l}}+{\pal\Delta {\cal F}\over \pal v^{\gamma,k}} 
{\pal\Delta {\cal F}\over \pal v^{\rho,l}} 
\right)\pal_x^{k+1} \pal^\gamma p_\alpha G^{\alpha\beta} \pal_x^{l+1} \pal^\rho p_\beta
\nn\\
&&
+{\epsilon^2\over 2} \sum  {\pal\Delta {\cal F}\over \pal v^{\gamma,k}} 
\pal_x^{k+1} \left[\nabla{\pal p_\alpha(v;\lambda)\over \pal \lambda} \cdot 
\nabla{\pal p_\beta(v;\lambda)\over \pal \lambda} \cdot v_x\right]^\gamma
G^{\alpha\beta}.
\eeqa
\end{theorem}

It remains an open problem to prove existence of solution to the loop equation
(\ref{glavnoe}) -- this sounds plausible since the number of equations is equal
to the number of unknowns -- and also to prove polynomiality of the resulting
Poisson pencil in every order $g$. For $g=1$ this follows from the results of
\cite{cmp}. In the next section we will consider the $g=2$ terms.

\subsubsection{Genus two}

Let us now proceed to the genus two case. After long but straightforward
calculations we obtain the following formula for ${\cal F}_2(u; u_x, \dots,
u_{xxxx})$. Denote 
$$
u_{ij}=u_i-u_j, ~i\neq j,
$$
$$
u_i' :=u_{i,x}, ~~u_i'' := u_{i,xx} ~~{\rm etc.}
$$
$$
h_i=h_i(u) := \psi_{i\,1}(u), ~i=1, \dots, n.
$$

\begin{theorem} For an arbitrary semisimple Frobenius manifold the following
formula holds true 
\eqa
&&
{\cal F}_2=\frac{1}{1152}\frac{u_i^{IV}}{{u_i'}^2\,h_i^2}
-\frac{7}{1920}\,\frac{{u_i''}\,u_i'''}{{u_i'}^3\,h_i^2}+
\frac{1}{360}\frac{{u_i''}^3}
{{u_i'}^4\,h_i^2}
+\frac1{40}\frac{V_{ij}^2\,u_i'''}{u_{ij}\,{u_i'}\,h_i^2}\nn\\
&&+\frac1{640}\frac{V_{ij}\,h_j\,{u_j'}\,u_i'''}
{u_{ij}\,{u_i'}^2\,h_i^3}-
\frac{19}{2880}\frac{V_{ij}\,u_i'''\,h_j}{u_{ij}\,{u_i'}\,
h_i^3}+
\frac1{1152}\frac{V_{ij}\,u_i'''\,h_i}{u_{ij}\,{u_j'}\,
h_j^3}\nn\\
&&+\frac7{40}\frac{V_{ij}^2\,V_{ik}^2\,{u_i''}}{u_{ij}\,u_{ik}\,h_i^2}
-\frac1{240}\frac{V_{ij}^2\,V_{ik}\,{u_i''}\,h_k
\left(32\,{u_i'}-7\,{u_k'}\right)}
{u_{ij}\,u_{ik}\,{u_i'}\,h_i^3}\nn\\
&&+\frac{1}{40}\frac{V_{ij}\,V_{jk}^2\,{u_i''}\,h_i}
{u_{ij}\,u_{jk}\,h_j^3}-
\frac{1}{48}\frac{V_{ij}\,V_{jk}^2\,{u_j'}\,{u_i''}}
{u_{ij}\,u_{jk}\,{u_i'}\,h_i\,h_j}-
\frac{3}{64}\frac{V_{ij}^2\,{u_i''}}{u_{ij}^2\,h_i^2}\nn\\
&&-\frac{11}{480}\frac{V_{ij}^2\,{u_i''}^2}{u_{ij}\,{u_i'}^2\,
h_i^2}+
\frac{29}{5760}\frac{V_{ij}\,V_{jk}\,{u_i''}\,h_i\,h_k
\left({u_k'}-2\,{u_j'}\right)}
{u_{ij}\,u_{jk}\,{u_j'}\,h_j^4}\nn\\
&&+
\frac{1}{1920}\frac{V_{ij}\,V_{ik}\,{u_i''}\,h_j\,h_k
\left(54\,{u_i'}^2-25\,{u_i'}\,{u_j'}-{u_j'}\,{u_k'}\right)}
{u_{ij}\,u_{ik}\,{u_i'}^2\,h_i^4}\nn\\
&&+
\frac{1}{384}\frac{V_{ij}\,V_{ik}\,{u_i''}\,h_k
\left({u_i'}-{u_k'}\right)}
{u_{ij}\,u_{ik}\,{u_j'}\,h_j^3}-
\frac{1}{384}\frac{V_{ik}\,V_{jk}\,{u_k'}\,{u_i''}\,h_i}
{u_{ik}\,u_{jk}\,{u_j'}\,h_j^3}\nn\\
&&+
\frac{1}{576}\frac{V_{ij}\,V_{jk}\,{u_i''}\,h_k
\left(2\,{u_j'}-{u_k'}\right)}
{u_{ij}\,u_{jk}\,{u_i'}\,h_i\,h_j^2}\nn\\
&&-
\frac{1}{5760}\frac{V_{ij}\,V_{jk}\,{u_k'}\,{u_i''}\,h_k
\left(27\,{u_i'}+{u_k'}\right)}
{u_{jk}\,u_{ik}\,{u_i'}^2\,h_i^3}-
\frac{19}{1920}\frac{V_{ij}\,V_{jk}\,{u_i''}\,h_k}
{u_{ij}\,u_{ik}\,h_i^3}\nn\\
&&+
\frac{1}{5760}\frac{V_{ij}\,V_{jk}\,h_k
\left(27\,{u_i'}\,{u_k'}-{u_j'}^2+2\,{u_j'}\,{u_k'}\right)\,{u_i''}}
{u_{ij}\,u_{jk}\,{u_i'}^2\,h_i^3}\nn\\
&&+
\frac{1}{288}\frac{V_{ij}\,V_{jk}\,{u_i''}\,h_i}
{u_{jk}\,u_{ik}\,h_k^3}+
\frac{1}{384}\frac{V_{ij}\,V_{jk}\,{u_i'}\,{u_i''}\,h_i}
{u_{ij}\,u_{ik}\,{u_k'}\,h_k^3}-
\frac{1}{576}\frac{V_{ij}\,V_{jk}\,{u_k'}\,{u_i''}}
{u_{jk}\,u_{ik}\,{u_i'}\,h_i\,h_k}\nn\\
&&+
\frac{1}{1920}\frac{V_{ij}\,{u_i''}^2\,h_j
\left(11\,{u_i'}-5\,{u_j'}\right)}
{u_{ij}\,{u_i'}^3\,h_i^3}-
\frac{1}{5760}\frac{V_{ij}\,{u_i''}\,{u_j''}\,h_j}
{u_{ij}\,{u_i'}^2\,h_i^3}\nn\\
&&+
\frac{1}{5760}\frac{V_{ij}\,{u_i''}\,h_j
\left(57\,{u_i'}^2-27\,{u_i'}\,{u_j'}-{u_j'}^2\right)}
{{u_i'}^2\,h_i^3}\nn\\
&&+
\frac{1}{1152}\frac{V_{ij}\,{u_i''}\,h_i
\left(4\,{u_j'}-3\,{u_i'}\right)}
{{u_j'}\,h_j^3}-
\frac{1}{576}\frac{V_{ij}\,{u_j'}\,{u_i''}}
{{u_i'}\,h_i\,h_j}\nn\\
&&-
\frac{1}{1152}\frac{V_{ij}\,{u_i''}\,{u_j''}}
{{u_i'}\,{u_j'}\,h_i\,h_j}+
\frac1{10}\frac{V^2_{ij}\,V^2_{ik}\,V^2_{il}\,{u_i'}^2}
{u_{ij}\,u_{ik}\,u_{il}\,h_i^2}\nn\\
&&
-\frac7{20}\frac{V^2_{ij}\, V^2_{ik}\, V_{il}\,h_l\,{u_i'}^2}
{u_{ij} u_{ik} u_{il}\,h_i^3}
+\frac7{40}\frac{V^2_{ij}\, V^2_{ik}\, V_{il}\, 
h_l\,{u_i'}\,{u_l'}}
{u_{ij}\, u_{ik}\, u_{il}\, h_i^3}
-\frac1{8} \frac{V^2_{ij}\,V_{ik}\,V^2_{kl}\,{u_i'}\,{u_k'}}
{u_{ij}\,u_{ik}\,u_{kl}\,h_i\,h_k}\nn\\
&&+
\frac1{40} \frac{V^2_{ij}\,V_{ik}\,V_{kl}\,h_l
\left({u_k'}^2-3\,{u_i'}^2-2\,{u_k'}\,{u_l'}\right)}
{u_{ij}\,u_{ik}\,u_{kl}\,h_i^3}
+\frac3{40}\frac{V^2_{ij}\,V_{ik}\,V_{kl}\,{u_i'}\,{u_l'}\,h_l}
{u_{ij}\,u_{ik}\,u_{il}\,h_i^3}\nn\\
&&
+\frac1{40}\frac{V^2_{ij}\,V_{ik}\,V_{kl}\,h_l\left(
3\,{u_i'}^2+{u_l'}^2\right)}
{u_{ij}\,u_{kl}\,u_{il}\,h_i^3}
+\frac1{48}\frac{V^2_{ij}\,V_{ik}\,V_{kl}\,h_l\,{u_i'}
\left(2\,{u_k'}-{u_l'}\right)}{u_{ij}\,u_{ik}\,u_{kl}\,h_i\,
h_k^2}\nn\\
&&+\frac5{96} \frac{V^2_{ij}\,V_{ik}\,V_{il}\,h_k\,h_l\,
\left(4\,{u_i'}^2-4\,{u_i'}\,{u_k'}+{u_k'} {u_l'}\right)}
{u_{ij} u_{ik} u_{il}\,h_i^4}
-
\frac{83}{480}\frac{V^2_{ij}\,V^2_{ik}\,{u_i'}^2}
{u_{ij}\,u_{ik}^2\,h_i^2}\nn\\
&&+
\frac1{144}\frac{V_{ij}\,V_{ik}\,V_{jl}\,V_{kl}\,{u_i'}^2}
{u_{ik}\,u_{jl}\,u_{il}\,h_i^2}
-\frac1{144}\frac{V_{ij}\,V_{ik}\,V_{jl}\,V_{kl}\,{u_i'}^2}
{u_{ij}\,u_{ik}\,u_{kl}\,h_i^2}-
\frac1{48}\frac{V_{ij}^2\,V_{ik}\,V_{kl}\,{u_i'}\,{u_l'}}
{u_{ij}\,u_{kl}\,u_{il}\,h_i\,h_l}
\nn\\
&&+\frac{29}{1920} \frac{V_{ij}\,V_{ik}\,V_{jl}\,h_k\,h_l
\left({u_k'}\,{u_l'}-{u_i'}\,{u_k'}+2\,{u_i'}^2-{u_i'}\,{u_l'}\right)}
{u_{ij}\,u_{ik}\,u_{il}\,h_i^4}\nn\\
&&-\frac{29}{5760} \frac{V_{ij}\,V_{ik}\,V_{jl}\,h_k\,h_l\,
{u_l'}^2\left(2\,{u_i'}-{u_k'}\right)}
{u_{ik}\,u_{jl}\,u_{il}\,h_i^4\,{u_i'}}\nn\\
&&-
\frac{29}{5760} \frac{V_{ij}\,V_{ik}\,V_{jl}\,h_k\,h_l\,
{u_j'}\left(2\,{u_k'}\,{u_l'}+2\,{u_i'}\,{u_j'}-{u_j'}\,{u_k'}
-4\,{u_i'}\,{u_l'}\right)}
{u_{ij}\,u_{ik}\,u_{jl}\,h_i^4\,{u_i'}}\nn\\
&&-
\frac{1}{1152} \frac{V_{ij}\,V_{ik}\,V_{jl}\,h_k\,h_l\,
\left(4\,{u_i'}\,{u_j'}-4\,{u_i'}\,{u_l'}+{u_k'}\,{u_l'}\right)}
{u_{ij}\,u_{ik}\,u_{jl}\,h_i^2\,h_j^2}\nn\\
&&-\frac1{384}\frac{V_{ij}\,V_{ik}\,V_{jl}\,h_l
\left({u_i'}\,{u_j'}^2-2\,{u_j'}\,{u_i'}\,{u_l'}\right)}
{u_{ij}\,u_{ik}\,u_{jl}\,{u_k'}\,h_k^3}\nn\\
&&+
\frac1{1152}\frac{V_{ij}\,V_{ik}\,V_{jl}\,h_l\,{u_i'}^2
\left({u_i'}-3\,{u_l'}\right)}
{u_{ij}\,u_{ik}\,u_{il}\,{u_k'}\,h_k^3}
-
\frac1{384}\frac{V_{ij}\,V_{ik}\,V_{jl}\,h_l\,{u_i'}\,{u_l'}^2}
{u_{ik}\,u_{jl}\,u_{il}\,{u_k'}\,h_k^3}\nn\\
&&-
\frac1{1152}\frac{V_{ij}\,V_{ik}\,V_{jl}\,h_l\,{u_j'}^2
\left(3\,{u_l'}-2\,{u_j'}\right)}
{u_{ij}\,u_{jl}\,u_{jk}\,{u_k'}\,h_k^3}-
\frac1{288}\frac{V_{ij}\,V_{ik}\,V_{jl}\,h_l\,{u_j'}
\left({u_j'}-2\,{u_l'}\right)}
{u_{ik}\,u_{jl}\,u_{jk}\,h_k^3}\nn\\
&&+
\frac1{576}\frac{V_{ij}\,V_{ik}\,V_{jl}\,h_l\,{u_k'}
\left(2\,{u_k'}-3\,{u_l'}\right)}
{u_{ik}\,u_{jk}\,u_{kl}\,h_k^3}-
\frac1{1152}\frac{V_{ij}\,V_{ik}\,V_{jl}\,h_l\,{u_l'}^3}
{u_{jl}\,u_{kl}\,u_{il}\,{u_k'}\,h_k^3}\nn\\
&&+
\frac1{288}\frac{V_{ij}\,V_{ik}\,V_{jl}\,h_l\,{u_l'}^2}
{u_{ik}\,u_{jl}\,u_{kl}\,h_k^3}-
\frac1{576}\frac{V_{ij}\,V_{ik}\,V_{jl}\,h_k\,{u_l'}
\left({u_k'}-2\,{u_i'}\right)}
{u_{ik}\,u_{jl}\,u_{il}\,h_i^2\,h_l}\nn\\
&&-
\frac1{1152}\frac{V_{ij}\,V_{ik}\,V_{jl}\,{u_k'}\,{u_l'}}
{u_{ik}\,u_{jl}\,u_{kl}\,h_k\,h_l}\nn\\
&&-
\frac7{1440}\frac{V_{ij}\,V_{ik}\,V_{il}\,h_j\,h_k\,
h_l\left(8\,{u_i'}^3-12\,{u_i'}^2\,{u_j'}-
{u_j'}\,{u_k'}\,{u_l'}+6\,{u_i'}\,{u_j'}\,{u_k'}\right)}
{u_{ij}\,u_{ik}\,u_{il}\,h_i^5\,{u_i'}}\nn\\
&&-
\frac{29}{1152}\frac{V_{ij}\,V_{ik}\,V_{jk}\,{u_i'}^2}
{u_{ij}\,u_{ik}^2\,h_i^2}-
\frac{1}{320}\frac{V_{ij}^2\,V_{ik}\,h_k
\left(3\,{u_i'}^2-8\,{u_k'}^2\right)}
{u_{ij}\,u_{ik}^2\,h_i^3}-
\frac{53}{1920}\frac{V_{ij}^2\,V_{ik}\,h_k\,{u_i'}\,{u_k'}}
{u_{ij}\,u_{ik}\,u_{jk}\,h_i^3}\nn\\
&&-
\frac{V_{ij}^2\,V_{ik}\,{u_i'}\,h_k}{u_{ij}^2\,u_{jk}\,h_i^3}
\left(\frac{27}{640}\,{u_k'}-\frac{233}{2880}\,{u_i'}\right)\nn\\
&&-
\frac{V_{ij}^2\,V_{ik}\,{u_i'}\,h_k}
{u_{ik}^2\,u_{jk}\,h_i^3}
\left(\frac{233}{2880}\,{u_i'}-\frac{67}{960}\,{u_k'}\right)
+
\frac{1}{1152}\frac{V_{ij}^2\,V_{ik}\,h_i\,{u_i'}^3}
{u_{ij}\,u_{ik}^2\,{u_k'}\,h_k^3}\nn\\
&&-
\frac{1}{576}\frac{V_{ij}^2\,V_{ik}\,h_i\,{u_i'}^3}
{u_{ij}^2\,u_{ik}\,{u_k'}\,h_k^3}-
\frac{1}{48}\frac{V_{ij}^2\,V_{ik}\,{u_i'}\,{u_k'}}
{u_{ij}\,u_{ik}^2\,h_i\,h_k}\nn\\
&&+
\frac{233}{1440}\frac{V_{ij}^3\,h_j\,{u_i'}^2}
{u_{ij}^3\,h_i^3}-
\frac{43}{384}\frac{V_{ij}^3\,h_j\,{u_i'}\,{u_j'}}
{u_{ij}^3\,h_i^3}-
\frac{1}{12}\frac{V_{ij}^3\,{u_i'}\,{u_j'}}
{u_{ij}^3\,h_i\,h_j}\nn\\
&&+
\frac{29}{5760}\frac{V_{ij}\,V_{ik}\,{u_j'}\,{u_k'}\,
h_j\,h_k
\left({u_k'}-6\,{u_i'}\right)}
{u_{ij}\,u_{ik}^2\,{u_i'}\,h_i^4}\nn\\
&&+
\frac{29}{5760}\frac{V_{ij}\,V_{ik}\,
h_j\,h_k
\left(3\,{u_i'}\,{u_k'}+3\,{u_j'}\,{u_k'}+6\,{u_i'}\,{u_j'}
-6\,{u_i'}^2-2\,{u_j'}^2\right)}
{u_{ij}^2\,u_{ik}\,h_i^4}\nn\\
&&+
\frac{1}{576}\frac{V_{ij}\,V_{ik}\,{u_j'}\,
h_k\left(2\,{u_i'}-{u_k'}\right)}
{u_{ij}^2\,u_{ik}\,h_i^2\,h_j}\nn\\
&&+
\frac{1}{1152}\frac{V_{ij}\,V_{ik}\,u_{ij}\,h_k
\left(3\,{u_i'}^2\,{u_k'}-3\,{u_i'}\,\,{u_k'}^2+
{u_k'}^3-{u_i'}^3\right)}
{u_{ik}^2\,u_{jk}^2\,{u_j'}\,h_j^3}\nn\\
&&+
\frac{1}{576}\frac{V_{ij}\,V_{ik}\,u_{ik}\,h_k
\left(-{u_i'}^3+3\,{u_j'}^2\,{u_k'}-4\,{u_i'}\,{u_j'}\,{u_k'}+
2\,{u_i'}^2\,{u_j'}-2\,{u_j'}^3\right)}
{u_{ij}^2\,u_{jk}^2\,{u_j'}\,h_j^3}\nn\\
&&+
\frac{1}{384}\frac{V_{ij}\,V_{ik}\,h_k
\left(-{u_i'}\,{u_k'}^2+{u_i'}^3-6\,{u_j'}^2\,{u_k'}\right)}
{u_{ij}\,u_{jk}^2\,{u_j'}\,h_j^3}\nn\\
&&+
\frac{1}{288}\frac{V_{ij}\,V_{ik}\,h_k
\left(4\,{u_i'}\,{u_j'}\,{u_k'}+
{u_j'}\,{u_k'}^2-2\,{u_i'}^2\,{u_j'}+3\,\,{u_j'}^3\right)}
{u_{ij}\,u_{jk}^2\,{u_j'}\,h_j^3}\nn\\
&&+
\frac{1}{384}\frac{V_{ij}\,V_{ik}\,h_k
\left(2\,{u_i'}\,{u_k'}^2-
{u_i'}^2\,{u_k'}-{u_k'}^3\right)}
{u_{ik}\,u_{jk}^2\,{u_j'}\,h_j^3}\nn\\
&&+
\frac{1}{288}\frac{V_{ij}\,V_{ik}\,h_k
\left({u_j'}\,{u_k'}^2-2\,{u_i'}\,{u_j'}\,{u_k'}
+{u_i'}^2\,{u_j'}\right)}
{u_{ik}\,u_{jk}^2\,{u_j'}\,h_j^3}\nn\\
&&+
\frac{1}{384}\frac{V_{ij}\,V_{ik}\,h_k
\,{u_i'}^2\,{u_k'}}
{u_{ij}^2\,u_{jk}\,{u_j'}\,h_j^3}-
\frac{1}{576}\frac{V_{ij}\,V_{ik}
\,{u_j'}\,{u_k'}}
{u_{ik}\,u_{jk}^2\,h_j\,h_k}
\nn\\
&&+
\frac{1}{1152}\frac{V_{ij}^2\,{u_i'}
\left(37\,{u_i'}\,{u_j'}\,h_j^2+
10\,{u_i'}\,{u_j'}\,h_i^2-3\,{u_i'}^2\,h_i^2
+11\,{u_j'}^2\,h_j^2\right)}
{u_{ij}\,{u_j'}\,h_i^2\,h_j^2}\nn\\
&&-
\frac{1}{576}\frac{V_{ij}\,h_j
\left(4\,{u_i'}^3+
4\,{u_i'}\,{u_j'}^2-6\,{u_i'}^2\,{u_j'}-{u_j'}^3\right)}
{u_{ij}^3\,{u_i'}\,h_i^3}+
\frac{1}{576}\frac{V_{ij}\,{u_i'}\,{u_j'}}
{u_{ij}^3\,h_i\,h_j}.
\nn\\
\eeqa
A summation over repeated indices is assumed in each term
of the formula provided the denominators do not vanish. 
\end{theorem}

\begin{exam} 
For the ${\bf CP}^1$ model defined by the potential (\ref{potential-toda})
the canonical coordinates are
\beq\label{canonical-cp1}
u_1=v_1+2\,\exp(\frac{v_2}2),\quad 
u_2=v_1-2\,\exp(\frac{v_2}2).
\eeq
The functions $h_1, h_2$ and the matrix $V$ have the form
\eqa
&&h_1=\frac{\sqrt{2}}{\sqrt{u_1-u_2}},\quad
h_2=-\frac{\sqrt{2}\, i}{\sqrt{u_1-u_2}},\nn\\
&&V=\frac{i}2\left(\matrix{0 & -1\cr 1 & 0\cr}\right).\nn
\eeqa
By substituting the above formulae for the functions $h_i, V_{ij}$ into 
the formula of genus two free energy 
for general semisimple Frobenius manifolds, we get the following expression
for the genus two free energy of the ${\bf CP}^1$ model:
\eqa
&&24^2\,{\cal F}_2=
\frac{{4\,u_1''}^3\,u_{12}}{5\,{u_1'}^4} - 
  \frac{{4\,u_2''}^3\,u_{12}}{5\,{u_2'}^4}-
\frac{u_1''\,u_2''}{4\,u_1'\,u_2'}
\nn\\
&&\quad+
\frac{3\,u_1''}{4\,{u_1'}^3}\left(\frac12\,u_1''\,u_2'-\frac75\,u_1'''\,u_{12}
\right)
+
\frac{3\,u_2''}{4\,{u_2'}^3}\left(\frac12\,u_2''\,u_1'+\frac75\,u_2'''\,u_{12}
\right)\nn\\
&&\quad+
\frac1{4\,{u_1'}^2}\left(\frac{33}{10}\,{u_1''}^2-\frac9{10}\,u_1'''\,u_2'+
\frac1{10}\,u_1''\,u_2''+u_1^{IV}\,u_{12}\right)\nn\\
&&\quad+
\frac1{4\,{u_2'}^2}\left(\frac{33}{10}\,{u_2''}^2-\frac9{10}\,u_2'''\,u_1'+
\frac1{10}\,u_1''\,u_2''-u_2^{IV}\,u_{12}\right)\nn\\
&&\quad-
\frac1{4\,u_1'}\left(\frac{17}{5}\,u_1'''+\frac1{2}\,u_2'''\right)-
\frac1{4\,u_2'}\left(\frac{17}{5}\,u_2'''+\frac1{2}\,u_1'''\right)\nn\\
&&\quad-
\frac1{10\, u_{12}^2}\left(\frac{{u_1'}^3}{u_2'}+\frac{{u_2'}^3}{u_1'}
\right)
-\frac{1}{u_{12}^2}\left({u_1'}^2-\frac{11}{5}\,u_1'\,u_2'+{u_2'}^2\right)
\nn\\
&&\quad+
\frac{u_1''-u_2''}{u_{12}}\left(\frac{u_2'}{5\,u_1'}+\frac{u_1'}{5\,u_2'}+
1\right).
\eeqa
Here we denote $u_{12}=u_1-u_2$ as above. 

Let us now explain how to use the above formula for computation of the genus 2
Gromov - Witten invariants and their descendents in the ${\bf CP}^1$ model.
We must substitute into the
formula for ${\cal F}_2$ the genus zero two point correlation
functions
\beq\label{correlation-zero}
v_\al=\frac{\pal^2}{\pal t^{1,0}\pal t^{\alpha,0}}
 \langle\exp(\sum \tau_p(\phi_\beta) t^{\beta,p})\rangle_0,\quad \al=1,2,
\eeq
where $\tau_p(\phi_\al), p\ge 0, \al=1,2$ are the gravitational descendents 
of the primary fields $\phi_1=1, \phi_2$ of the ${\bf CP}^1$ model.
The expansion of these two point correlation functions in power serieses
of $t^{\al,p}$ can be obtained by solving the following equations 
(see Section \ref{sec-3-6-4})
\beq\label{fp}
v_\al(t)=\sum_{p\ge 0} \sum_{\beta=1}^2 t^{\beta,q} 
\frac{\pal\theta_{\beta,q}}{\pal v_\gamma},\quad \al=1,2.
\eeq
The generating function for $\theta_{\beta,q}$ is given in 
(\ref{theta-toda-1}), (\ref{theta-toda-2}), 
we list few of them 
\eqa
&&\theta_{1,0}=v_2,\quad\theta_{2,0}=v_1,\nn\\
&&\theta_{1,1}=v_1\,v_2,\quad
\theta_{2,1}=e^{v_2}+\frac12\,v_1^2,\nn\\
&&\theta_{1,2}=\frac12\,v_1^2\,v_2+v_2\,e^{v_2}
-2\,e^{v_2},\nn\\
&&\theta_{2,2}=\frac16\,v_1^3+v_1\,e^{v_2},\nn\\
&&\theta_{1,3}=\frac16\,v_1^3\,v_2+v_1\,v_2\,e^{v_2}
-2\,v_1\,e^{v_2},\nn\\
&&\theta_{2,3}=\frac1{24}\,v_1^4+\frac12\,v_1^2\,e^{v_2}+
\frac14\,e^{2\,v_2}.
\eeqa
The expansion of the correlation functions of (\ref{correlation-zero})
is the unique solution of the form 
$$
v_\al(t)=t^{\al,0}+\sum_{k\ge 1,\, p_i\ge 1} 
A^\al_{\beta_1, q_1;\dots;\beta_k, q_k}(t^{1,0},t^{2,0})\,
t^{\beta_1,q_1}\dots t^{\beta_k,q_k}
$$
of the equations given in (\ref{fp}),
the coefficients are determined recursively by (\ref{fp}). For example,
we have
$$
A_{\beta,q}=\left. 
\frac{\pal\theta_{\beta,q}}{\pal v_\al}\right|_{v_\gamma=t^{\gamma,0}},
\quad 
A_{\beta_1,q_1;\beta_2,q_2}=\left.
\frac12\,\frac{\pal^2\theta_{\beta_1,q_1}}
{\pal v_\al \pal v_\gamma}\,\frac{\pal\theta_{\beta_2,q_2}}
{\pal v_\gamma}\right|_{v_\xi=t^{\xi,0}}.
$$
Taking the expansions of $v_\al$ up to the 6-th order of $t^{\al,p}$,
we obtain the  following expansions for $u_1, u_2$ defined in 
(\ref{canonical-cp1}):  
\eqa
&&u_1'=1 + t^{1,1} + (t^{1,1})^2 + t^{1,0}\,t^{1,2} - (t^{1,2})^2 - t^{1,3} - 
  2\,t^{1,1}\,t^{1,3} - 2\,t^{1,2}\,t^{1,3} + (t^{1,3})^2\nn\\
&&\quad + t^{1,2}\,t^{2,0} + t^{2,1}
+ 2\,t^{1,1}\,t^{2,1} + t^{1,2}\,t^{2,1} - 3\,t^{1,3}\,t^{2,1} + 
  \frac{t^{2,0}\,t^{2,1}}{2} + (t^{2,1})^2 + t^{2,2} \nn\\
&&\quad+ t^{1,0}\,t^{2,2} + 
  2\,t^{1,1}\,t^{2,2} - 4\,t^{1,3}\,t^{2,2} + t^{2,0}\,t^{2,2} + 
  \frac{7\,t^{2,1}\,t^{2,2}}{2} + 2\,(t^{2,2})^2 
+ t^{2,3} + t^{1,0}\,t^{2,3}\nn\\
&&\quad + 
  2\,t^{1,1}\,t^{2,3} - \frac{t^{1,2}\,t^{2,3}}{2} - 5\,t^{1,3}\,t^{2,3} + 
  \frac{3\,t^{2,0}\,t^{2,3}}{2} + 3\,t^{2,1}\,t^{2,3} + 4\,t^{2,2}\,t^{2,3}
\nn\\ 
&&\quad+ 
  \frac{3\,(t^{2,3})^2}{2}+{\cal O}(t^3),\nn\\
&&u_1''=t^{1,2} + 3\,t^{1,1}\,t^{1,2} + t^{1,0}\,t^{1,3} - 4\,t^{1,2}\,t^{1,3} - 
  2\,(t^{1,3})^2 + t^{1,3}\,t^{2,0} + 3\,t^{1,2}\,t^{2,1} \nn\\
&&\quad+ t^{1,3}\,t^{2,1}+ 
  \frac{(t^{2,1})^2}{2} + t^{2,2} + 3\,t^{1,1}\,t^{2,2} + 3\,t^{1,2}\,t^{2,2} - 
  3\,t^{1,3}\,t^{2,2} + \frac{t^{2,0}\,t^{2,2}}{2} + 3\,t^{2,1}\,t^{2,2} \nn\\
&&\quad+ 
  \frac{7\,(t^{2,2})^2}{2} + t^{2,3} + t^{1,0}\,t^{2,3} + 3\,t^{1,1}\,t^{2,3} + 
  2\,t^{1,2}\,t^{2,3} - \frac{9\,t^{1,3}\,t^{2,3}}{2} + t^{2,0}\,t^{2,3} + 
  5\,t^{2,1}\,t^{2,3}\nn\\
&&\quad + 7\,t^{2,2}\,t^{2,3} + 4\,(t^{2,3})^2+{\cal O}(t^3),\nn\\
&&u_1'''=3\,(t^{1,2})^2 + t^{1,3} + 4\,t^{1,1}\,t^{1,3} - 4\,(t^{1,3})^2 + 
  4\,t^{1,3}\,t^{2,1} + 6\,t^{1,2}\,t^{2,2} + 4\,t^{1,3}\,t^{2,2}\nn\\
&&\quad + \frac{3\,t^{2,1}\,t^{2,2}}{2}
+ 3\,(t^{2,2})^2 + t^{2,3} + 4\,t^{1,1}\,t^{2,3} + 
  6\,t^{1,2}\,t^{2,3} - t^{1,3}\,t^{2,3} + \frac{t^{2,0}\,t^{2,3}}{2} + 
  4\,t^{2,1}\,t^{2,3} \nn\\
&&\quad + 12\,t^{2,2}\,t^{2,3} + 7\,(t^{2,3})^2+{\cal O}(t^3) ,\nn\\
&&u_1^{IV}=10\,t^{1,2}\,t^{1,3} + 10\,t^{1,3}\,t^{2,2} + \frac{3\,(t^{2,2})^2}{2} + 
  10\,t^{1,2}\,t^{2,3} + 10\,t^{1,3}\,t^{2,3} + 2\,t^{2,1}\,t^{2,3} \nn\\
&&\quad+ 
  10\,t^{2,2}\,t^{2,3} + 12\,(t^{2,3})^2+{\cal O}(t^3),\nn\\
&&u_2'=1 + t^{1,1} + (t^{1,1})^2 + t^{1,0}\,t^{1,2} - (t^{1,2})^2 - t^{1,3} - 
  2\,t^{1,1}\,t^{1,3} + 2\,t^{1,2}\,t^{1,3} + (t^{1,3})^2 \nn\\
&&\quad- t^{1,2}\,t^{2,0}- 
  t^{2,1} - 2\,t^{1,1}\,t^{2,1} + t^{1,2}\,t^{2,1} + 3\,t^{1,3}\,t^{2,1} - 
  \frac{t^{2,0}\,t^{2,1}}{2} + (t^{2,1})^2 + t^{2,2}\nn\\
&&\quad - t^{1,0}\,t^{2,2} + 
  2\,t^{1,1}\,t^{2,2} - 4\,t^{1,3}\,t^{2,2} + t^{2,0}\,t^{2,2} - 
  \frac{7\,t^{2,1}\,t^{2,2}}{2} + 2\,(t^{2,2})^2 - t^{2,3} + t^{1,0}\,t^{2,3}
\nn\\
&&\quad - 
  2\,t^{1,1}\,t^{2,3} - \frac{t^{1,2}\,t^{2,3}}{2} + 5\,t^{1,3}\,t^{2,3} - 
  \frac{3\,t^{2,0}\,t^{2,3}}{2} + 3\,t^{2,1}\,t^{2,3} \nn\\
&&\quad- 4\,t^{2,2}\,t^{2,3} + 
  \frac{3\,(t^{2,3})^2}{2}+{\cal O}(t^3) ,\nn\\
&&u_2''=t^{1,2} + 3\,t^{1,1}\,t^{1,2} + t^{1,0}\,t^{1,3} - 4\,t^{1,2}\,t^{1,3} + 
  2\,(t^{1,3})^2 - t^{1,3}\,t^{2,0} - 3\,t^{1,2}\,t^{2,1}\nn\\
&&\quad + t^{1,3}\,t^{2,1} - 
  \frac{(t^{2,1})^2}{2} - t^{2,2} - 3\,t^{1,1}\,t^{2,2} + 3\,t^{1,2}\,t^{2,2} + 
  3\,t^{1,3}\,t^{2,2} - \frac{t^{2,0}\,t^{2,2}}{2} + 3\,t^{2,1}\,t^{2,2}\nn\\
&&\quad - 
  \frac{7\,(t^{2,2})^2}{2} + t^{2,3} - t^{1,0}\,t^{2,3} + 3\,t^{1,1}\,t^{2,3} - 
  2\,t^{1,2}\,t^{2,3} - \frac{9\,t^{1,3}\,t^{2,3}}{2} + t^{2,0}\,t^{2,3} - 
  5\,t^{2,1}\,t^{2,3}\nn\\
&&\quad + 7\,t^{2,2}\,t^{2,3} - 4\,(t^{2,3})^2+{\cal O}(t^3) ,\nn\\
&&u_2'''=3\,(t^{1,2})^2 + t^{1,3} + 4\,t^{1,1}\,t^{1,3} - 4\,(t^{1,3})^2 - 
  4\,t^{1,3}\,t^{2,1} - 6\,t^{1,2}\,t^{2,2} + 4\,t^{1,3}\,t^{2,2}\nn\\
&&\quad - 
  \frac{3\,t^{2,1}\,t^{2,2}}{2} + 
3\,(t^{2,2})^2 - t^{2,3} - 4\,t^{1,1}\,t^{2,3} + 
  6\,t^{1,2}\,t^{2,3} + t^{1,3}\,t^{2,3} - \frac{t^{2,0}\,t^{2,3}}{2} + 
  4\,t^{2,1}\,t^{2,3}\nn\\
&&\quad - 12\,t^{2,2}\,t^{2,3} + 7\,(t^{2,3})^2+{\cal O}(t^3) ,\nn\\
&&u_2^{IV}=10\,t^{1,2}\,t^{1,3} - 10\,t^{1,3}\,t^{2,2} - \frac{3\,(t^{2,2})^2}{2} - 
  10\,t^{1,2}\,t^{2,3} + 10\,t^{1,3}\,t^{2,3}\nn\\
&&\quad - 2\,t^{2,1}\,t^{2,3} + 
  10\,t^{2,2}\,t^{2,3} - 12\,(t^{2,3})^2+{\cal O}(t^3) ,\nn\\
&&u_{12}=t^{1,2} + 3\,t^{1,1}\,t^{1,2} + t^{1,0}\,t^{1,3} - 4\,t^{1,2}\,t^{1,3} - 
  2\,(t^{1,3})^2 + t^{1,3}\,t^{2,0} + 3\,t^{1,2}\,t^{2,1} + t^{1,3}\,t^{2,1} 
\nn\\
&&\quad+ 
  \frac{(t^{2,1})^2}{2} + t^{2,2} + 3\,t^{1,1}\,t^{2,2} + 3\,t^{1,2}\,t^{2,2} - 
  3\,t^{1,3}\,t^{2,2} + \frac{t^{2,0}\,t^{2,2}}{2} + 3\,t^{2,1}\,t^{2,2} + 
  \frac{7\,(t^{2,2})^2}{2}\nn\\
&&\quad + t^{2,3}
 + t^{1,0}\,t^{2,3} + 3\,t^{1,1}\,t^{2,3} + 
  2\,t^{1,2}\,t^{2,3} - \frac{9\,t^{1,3}\,t^{2,3}}{2} + t^{2,0}\,t^{2,3} + 
  5\,t^{2,1}\,t^{2,3}\nn\\
&&\quad + 7\,t^{2,2}\,t^{2,3} + 4\,(t^{2,3})^2+{\cal O}(t^3) .\nn
\eeqa
In the r.h.s. of the above formulae we only keep the terms up to the 
quadratic ones in $t^{\al,p}$ with $p\le 3$. The function ${\cal F}_2$ thus
has the following expansion:
\eqa
&&{\cal F}_2=-\frac{t^{1,3}}{240} + \frac{7\,t^{2,2}}{5760}
-\frac{5\,(t^{1,2})^2}{576} - \frac{t^{1,1}\,t^{1,3}}{80} + 
  \frac{29\,(t^{1,3})^2}{2880} + \frac{7\,t^{1,3}\,t^{2,0}}{5760} + 
  \frac{7\,t^{1,2}\,t^{2,1}}{1920}\nn\\
&&\qquad + \frac{7\,t^{1,1}\,t^{2,2}}{1920} + 
  \frac{t^{1,3}\,t^{2,2}}{192} + \frac{(t^{2,2})^2}{1152} + 
  \frac{7\,t^{1,0}\,t^{2,3}}{5760} + \frac{t^{1,2}\,t^{2,3}}{192} + 
  \frac{25\,(t^{2,3})^2}{2304}+{\cal O}(t^3).\nn
\eeqa
Again the above  formula for $F_2$ is at the approximation 
up to quadratic terms in $t^{\al,p}$ with $p\le 3$. The above expansion
of the function ${\cal F}_2$ coincides, to the extend of the above
mentioned approximation, to the generating function
of the genus two Gromov-Witten invariants for the ${\bf CP}^1$ model 
defined by
\eqa
&&\tilde {\cal F}_2=
\sum_{n\ge 0} \frac1{k!} \langle (\sum_{p\ge 0}\sum_{\al=1}^2 
\tau_p(\phi_\al) t^{\al,p})^k\rangle_2
\nn\\
&&\quad=1+\sum_{p\ge 0}\sum_{\al=1}^2
\langle \tau_p(\phi_\al)\rangle_2 t^{\al,p}\nn\\
&&\qquad+
\frac12\,\sum_{p, q\ge 0}\sum_{\al,\beta=1}^2
\langle \tau_p(\phi_\al) \tau_q(\phi_\beta)\rangle_2 t^{\al,p}\,t^{\beta,q}
+{\cal{O}}(t^3).
\eeqa
A list of some of the invariants can be found in \cite{song}, 
the numbers $\langle \tau_{k,i}\rangle_2$ and $\langle \tau_{k,i}\,\tau_{m,j}\rangle_2$ there 
correspond to 
$\langle \tau_k(\phi_{i+1})\rangle_2$ and $\langle \tau_k(\phi_{i+1})\tau_m(\phi_{j+1})\rangle_2$ 
respectively. For example,
\eqa
&&\langle \tau_3(\phi_1)\rangle_2=-\frac1{240},\quad \langle \tau_2(\phi_2)\rangle_2=\frac7{5760},\nn\\
&&\langle \tau_2(\phi_1)\tau_2(\phi_1)\rangle_2=-\frac5{288},
\quad \langle \tau_3(\phi_1)\tau_3(\phi_1)\rangle_2=\frac{29}{1440},\nn\\
&&\langle \tau_2(\phi_1)\tau_1(\phi_2)\rangle_2=\frac7{1920},\quad
\langle \tau_3(\phi_1)\tau_2(\phi_2)\rangle_2=\langle \tau_2(\phi_1)\tau_3(\phi_2)\rangle_2
=\frac1{192},\nn\\
&&\langle \tau_2(\phi_2)\tau_2(\phi_2)\rangle_2=\frac1{576},\quad
\langle \tau_3(\phi_2)\tau_3(\phi_2)\rangle_2=\frac{25}{1152}.
\eeqa
\end{exam}

\begin{exam} 
The expansion of the function ${\cal F}_2$ for the ${\bf CP}^2$ model
can be obtained in a similar way as we did for the ${\bf CP}^1$ model.
The canonical coordinates $u_i$ are the roots of the characteristic polynomial
of the matrix correspondent to the operation of multiplication by the
Euler vector field, express $u_i$ in the form
$$
u_i=v_1+v_3^{-1}\,z_i(y),\quad i=1,2,3,
$$
where $y=v_2+3\,\log v_3$.
Then the functions $z_i(y)$ are roots of the following cubic polynomial:
\eqa
&&z^3 - f''\,z^2+\left(6\,f - 15\,f' - 9\,f''\right)\,z
-54\,f + 243\,f' 
\nn\\
&&+ 4\,{f'}^2 - 243\,f'' - 6\,f\,f'' + 
  3\,f'\,f'' + 9\,{f''}^2=0.
\eeqa
Here $f=f(y)$ is the generating function of the genus zero Gromov - Witten invariants
of ${\bf CP}^2$, i.e.
$$
f(y)=\sum_{k\ge 1}\frac{N_k^{(0)}}{(3 k-1)!} e^{k y}
$$
with $N^{(0)}_k$ being the number of rational curves of degree $k$ on
${\bf CP}^2$ passing through $3k-1$ generic points.
The matrix $\Psi=(\psi_{i,\al})$ is given by
\eqa
&&\psi_{i,1}=a_i\,v_3\,\sqrt{f''+9-z_i},\nn\\
&&\psi_{i,2}
=\frac{a_i\left(3 z_i+2 f'-6 f''\right)}
{\sqrt{f''+9-z_i}},
\nn\\
&&\psi_{i,3}=\frac{a_i\left(z_i^2-f''\,z_i+6 f'-18 f''\right)}
{v_3\,\sqrt{f''+9-z_i}},
\eeqa
where the coefficents $a_i$ are defined by
$$
a_i=\frac{1}{\sqrt{\prod_{k\ne i}(z_i-z_k)}}
$$
The matrix $V=(V_{ij})$ has the following expression
$$
V=\Psi\,\mu\,\Psi^{-1},\quad \mu=\mbox{diag}(-1,0,1).
$$
The genus zero two point functions $v_1, v_2, v_3$ can be obtained by solving
(\ref{tsarev}). 
Thus we have all the data to compute the expansion of the function
${\cal F}_2$ for the ${\bf CP}^2$ model. 
\end{exam}

Let us now consider the deformations of the Principal Hierarchy induced by the
quasitriviality transformation (\ref{2-11-42}) satisfying the loop equation
(\ref{loop-eq}).
It remains an open problem to prove that, at every order in the
$\epsilon^2$-expansion the coefficients of the bihamiltonian structure
are {\it polynomials in the derivatives}. For the first term in the
$\epsilon^2$-expansion this was proved in \cite{cmp}. At the next approximation
the calculations are more complicated. We consider here 
bihamiltonian structures corresponding to  two-dimensional Frobenius
manifolds (\ref{I2-k}).

As we already know from Example \ref{e-I2-k-2-7} of Section 
\ref{sec-3-6-2} for generic 
two-dimensional Frobenius manifold  
$$
F(v_1, v_2)= {1\over 2} v_1^2 v_2 + {v_2^{\kappa+1}\over \kappa^2-1}
$$
the $t=-t^{1,1}$-flow 
$$
v_t + v\cdot v_x=0
$$
of the Principal Hierarchy
coincides with the equations of motion of 1D polytropic gas
with the equation of state
of the form
$$
p  ={\kappa\,\rho^{k+1}\over \kappa+1}
$$
(we redenote as above $v_1 = u$, $v_2 =\rho$).

Our construction produces the following
deformation of the polytropic gas equations
\eqa\label{poly-def}
&&{\partial u\over \partial t}+\partial_x\left\{
{u^2\over 2}+{\rho^\kappa}\right.
\nn\\
&&+\epsilon^2\left[
\frac{\kappa-2}{8}\rho^{\kappa-3}\,\rho_x^2
+\frac{\kappa}{12}\,\rho^{\kappa-2}\rho_{xx}\right]
\nn\\
&&+\epsilon^4 (\kappa-2)(\kappa-3)\left[
a_{1}\, \rho^{-4} \,u_x^2\,\rho_x^2+
        a_{2}\, \rho^{\kappa-6} \,\rho_x^4\right.
\nn\\
&&\quad+
        a_{3}\, \rho^{-3} \,u_{xx}\,u_x\,\rho_x+
        a_{4}\, \rho^{-2} \,u_{xx}^2+
        a_{5}\, \rho^{-3} \,u_x^2\,\rho_{xx}
\nn\\
&&\quad+
        a_{6}\, \rho^{\kappa-5} \,\rho_x^2\,\rho_{xx}+
        a_{7}\, \rho^{\kappa-4} \,\rho_{xx}^2+
        a_{8}\, \rho^{-2} \,u_x\,u_{xxx}
\nn\\
&&\quad\left. \left.+
        a_{9}\, \rho^{\kappa-4} \,\rho_x\,\rho_{xxx}\right]+
        \epsilon^4 \frac{\kappa\,(\kappa^2-1)\,(\kappa^2-4)}{360} \rho^{\kappa-3}
        \,\rho_{xxxx}\right\} =O(\epsilon^6),
\nn
\eeqa
\eqa
&&{\partial \rho\over \partial t}+\partial_x\left\{
\rho\, u+\epsilon^2\,\left(\frac{(2-\kappa)(\kappa-3)}{12\,\kappa\,\rho}\,u_x\,\rho_x+
\frac16\,u_{xx}\right)\right.
\nn\\
&&\nn\\&&
+\epsilon^4 (\kappa-2)(\kappa-3)\left[
b_{1}\, \rho^{-4} \,u_x\,\rho_x^3+
        b_{2}\, \rho^{-3} \,\rho_x^2\,u_{xx}+
        b_{3}\, \rho^{-3} \,u_x\,\rho_x\,\rho_{xx}\right.
\nn\\
&&\left.\left.+
        b_{4}\, \rho^{-2} \,u_{xx}\,\rho_{xx}+
        b_{5}\, \rho^{-2} \,u_{xxx}\,\rho_x+
        b_{6}\, \rho^{-2} \,u_x\,\rho_{xxx}+
        b_{7}\, \rho^{-1} \,u_{xxxx}\right]\right\}
	\nn\\
	&&=O(\epsilon^6).
\nn
\eeqa
integrable up to corrections of the order $O(\epsilon^6)$. The coefficients
are given below
\eqa
&&a_{1}={\frac{
     36 + 144\,\kappa - 59\,{\kappa^2} + 19\,{\kappa^3}  }{5760\,{\kappa^3}}},
~~a_{2}={\frac{
     60 + 176\,\kappa + 433\,{\kappa^2} - 182\,{\kappa^3} + 17\,{\kappa^4}  }{5760\,
     {\kappa^3}}}
\nn\\&&\nn\\
&&a_{3}={\frac{
     6 - 19\,\kappa - 11\,{\kappa^2} - 4\,{\kappa^3}  }{1440\,{\kappa^3}}},~~
a_{4}={\frac{
      -6 - 5\,\kappa + 13\,{\kappa^2}  }{1440\,{\kappa^3}}},~~
a_{5}={\frac{
      -42 + 13\,\kappa - 7\,{\kappa^2}  }{2880\,{\kappa^2}}}
\nn\\&&\nn\\
&&a_{6}={\frac{
      -36 - 72\,\kappa - 245\,{\kappa^2} - 61\,{\kappa^3} + 30\,{\kappa^4}  }{2880\,
     {\kappa^2}}},~~
a_{7}={\frac{
      6 + 5\,\kappa + 15\,{\kappa^2} + 5\,{\kappa^3} + 5\,{\kappa^4} }{1440\,{\kappa^2}}}
\nn\\&&\nn\\
&&a_{8}={\frac{1 }{120\,\kappa}}, ~~
a_{9}={\frac{ 2 + 5\,\kappa  }
     {240}}
     \nn
     \eeqa
    
\eqa
&&b_{1}={\frac{
      108 + 192\,\kappa - 97\,{\kappa^2} + 17\,{\kappa^3} 
      }{2880\,{\kappa^3}}},~~
b_{2}={\frac{
     -18 - 75\,\kappa + 47\,{\kappa^2} - 10\,{\kappa^3}  }{1440\,{\kappa^3}}}
\nn\\&&\nn\\
&&b_{3}=-{\frac{
      6 + 17\,\kappa - 5\,{\kappa^2} + 2\,{\kappa^3}  }{288\,{\kappa^3}}},~
b_{4}={\frac{
      6 - 4\,\kappa + {\kappa^2}  }{180\,{\kappa^2}}},~
b_{5}={\frac{
      6 + \kappa + {\kappa^2} }{720\,{\kappa^2}}}\nn\\
&&\nn\\
&&b_{6}={\frac{
      6 + \kappa + {\kappa^2}  }{720\,{\kappa^2}}},~
b_{7}=-{\frac{1 }{360\,\kappa}}
\eeqa
The polytropic gas equations can be considered as the dispersionless limit
of the above integrable system.

The bihamiltonian structure of the above system will be written in the
coordinates where the first Poisson bracket has the canonical form
(\ref{poly-def})
up to $O(\epsilon^6)$. These coordinates $\bar v_1$, $\bar v_2$ are given by the
following Miura-type transformation (recall $v_1=u$, $v_2=\rho$; we do not
distinguish between upper and lower indices in these formulae)
\beq
{\bar v}_\al=v_\al+\epsilon^2\,A_\al(v,v',v'')+\epsilon^4\,
B_\al(v,v',\dots,v^{IV}), ~\al=1,\, 2\nn
\eeq
where
$$
A_1=\frac{(\kappa-2)(\kappa-3)}{24\,\kappa}\,\frac{\pal}{\pal x}\left(\frac{
v_1'}{v_2}\right),\quad
A^2=\frac{(\kappa-2)(\kappa-3)}{24\,\kappa}\,\frac{\pal}{\pal x}\left(\frac{
v_2'}{v_2}\right)
$$
\eqa
&&B_1=\frac{(\kappa-2)(\kappa-3)}{240}\,\pal_x
\left[\frac {{\kappa}^{2}-13\,\kappa+6}{{\kappa}^{2}}\left(
\frac14\,
\frac{v_1'''}{v_2^2}
-\frac12\,
\frac{v_1''\,v_2'}{v_2^3}
-{\frac {7}{12}}\,
\frac{v_1'\,v_2''}{v_2^3}\right)\right.
\nn\\
&&\quad\left.+
\frac16\,{\frac {5\,{\kappa}^{3}-54\,{\kappa}^{2}+25\,\kappa-6}{{{\kappa}}^{3}}}\,
\frac{v_1'\,{v_2'}^2}{v_2^4}
-\frac1{18}\,{\frac {(\kappa+3) (\kappa+2) (\kappa+1)}{
\left (\kappa-1\right ){\kappa}^{4}}}\,
\frac{{v_1'}^3}{v_2^{\kappa+2}}\right].\nn\\
&&B_2=0\nn
\eeqa
In the new coordinates our deformed bihamiltonian 
structure reads
\eqa\label{1st-const}
&&\{{\bar v}^\al(x),{\bar v}^\beta(y)\}_1
=\{{\bar v}^\al(x),{\bar v}^\beta(y)\}_1^{[0]}+{\cal O}(\epsilon^6),\nn\\
&&\{{\bar v}^\al(x),{\bar v}^\beta(y)\}_2
=\{{\bar v}^\al(x),{\bar v}^\beta(y)\}_2^{[0]}+
\epsilon^2\,(d_1 X)^{\al\beta}+\epsilon^4\,(d_1 Y)^{\al\beta}+
{\cal O}(\epsilon^6).\nn
\eeqa
with the vector fields $X$ and $Y$ defined by
\eqa
&&X^1=\frac{1}{24}\left[(3\,\kappa^2-7\,\kappa+2)\, v_2^{\kappa-3}\,{v_2'}^2
+(4 \kappa-2) v_2^{\kappa-2} v_2''\right]\nn\\
&&X^2=\frac{\kappa+1}{12\,\kappa} v_1''\nn
\eeqa
and 
\eqa
&&Y^1=\frac{(\kappa-2)(\kappa-3)}{240}\left({\frac {g_{{1}}v_1'''v_1'}{{v_2}^{2}}}
+{\frac {g_{{2}}{v_1''}^{2}}{{v_2}^{2}}}
+{\frac {g_{{3}}{v_1'}^{2}\,v_2''}{{v_2}^{3}}}
+{\frac {g_{{4}}{v_1'}^{2}{v_2'}^{2}}{{v_2}^{4}}}
+{\frac {g_{{5}}{v_1'}^{4}}{{v_2}^{\kappa+2}}}\right.\nn\\
&&\quad\left.
+\frac {g_6\,{{v_2}^{\kappa}\,v_2^{IV}}}{{{v_2}^{3}}}
+\frac {g_7\,{{v_2}^{\kappa}{v_2''}^{2}}}{{{v_2}^{4}}}
+\frac {g_{8}\,{{v_2}^{\kappa}\,v_2''\,{v_2'}^{2}}}{{{v_2}^{5}}}
+\frac {g_{9}\,{{v_2}^{\kappa}\,{v_2'}^{4}}}{{{v_2}^{6}}}\right)\nn\\
&&
Y^2=\frac{(\kappa-2)(\kappa-3)}{240}\,\left({\frac {f_{{1}}\,v_1^{IV}}{v_2}}
+{\frac {f_{{2}}\,v_1'''\,v_2'}{{v_2}^{2}}}
+{\frac {f_{{3}}v_1''\,v_2''}{{v_2}^{2}}}
-2\,f_3\,\frac{v_1''\,{v_2'}^{2}}{{v_2}^{3}}
+{\frac {f_{{4}}v_1'\,v_2'''}{{v_2}^{2}}}\right.\nn\\
&&\quad\left.-6\,f_4\,\frac{v_1'\,v_2''\,v_2'}{{v_2}^{3}}
+6\,f_4\,\frac{v_1'\,{v_2'}^{3}}{{v_2}^{4}}\right),\nn
\eeqa
where
\eqa
&&g_1=-\frac1{24}\,{\frac {9\,{\kappa}^{3}-120\,{\kappa}^{2}+101\,\kappa+30}
{{\kappa}^{3}}}\nn\\
&&g_2=-\frac1{12}\,{\frac {6\,{\kappa}^{3}-81\,{\kappa}^{2}
+79\,\kappa+6}{{\kappa}^{3}}}\nn\\
&&g_3=-\frac1{12}\,{\frac {8\,{\kappa}^{3}-109\,{\kappa}^{2}
+117\,\kappa-6}{{\kappa}^{3}}}\nn\\
&&g_4=\frac1{24}\,{\frac {37\,{\kappa}^{4}-478\,{\kappa}^{3}
+461\,{\kappa}^{2}+28\,\kappa+12}{{\kappa}^{4}}}\nn\\
&&g_5=-{\frac {1}{72}}\,{\frac { (\kappa
+3)(\kappa+2)(\kappa^2-1)}{{\kappa}^{5}}}\nn\\
&&g_6=\frac14\,{\frac {9\,{\kappa}^{3}-12\,{\kappa}^{2}
+41\,\kappa-18}{{\kappa}^{2}\left (\kappa-3\right )}}\nn\\
&&g_7=-\frac1{24}\,{\frac {
43\,{\kappa}^{3}-30\,{\kappa}^{2}+197\,\kappa-90}
{{\kappa}^{2}}}\nn\\
&&g_8=-\frac{1}{12}\,{\frac{
63\,{\kappa}^{4}-287\,{\kappa}^{3}+492\,{\kappa}^{2}-1258
\,\kappa+540}{{\kappa}^{2}}}\nn\\
&&g_9=-{\frac {1}{72}}\,{\frac {78\,{\kappa}^{6}-732\,{\kappa}^{5}+2268\,{\kappa}
^{4}-4393\,{\kappa}^{3}+8340\,{\kappa}^{2}-3395\,\kappa-6}{{\kappa}^{3}}}\nn\\&&\nn\\
&&
f_1=-\frac{7\,\kappa^2-35\,\kappa+18}{12\,\kappa^3}\nn\\
&&f_2=-\frac{21\,\kappa^3-312\,\kappa^2+529\,\kappa-138}{24\,\kappa^3}\nn\\
&&f_3=-\frac{21\,\kappa^3-291\,\kappa^2+424\,\kappa-84}{12\,\kappa^3}\nn\\
&&f_4=-\frac{21\,\kappa^3-284\,\kappa^2+389\,\kappa-66}{24\,\kappa^3}\nn
\eeqa

For the particular value $\kappa=2$ the above formulae truncate at the order
$\epsilon^2$. We obtain a system of two uncoupled KdVs
$$
\dot u_\pm + u_\pm u_\pm' \pm {\epsilon^2\over 12 \sqrt{2}} u_\pm'''=0
$$
for
$$
u_\pm = u\pm \rho \sqrt{2}.
$$

For $\kappa=3$ the formulae (\ref{1st-const}) define a Poisson pencil. That is, the Jacobi
identity becomes exact but not just modulo $O(\epsilon^6)$.
One obtains the bihamiltonian structure for the Boussinesq hierarchy 
(see details in \cite{cmp}). The Boussinesq equation itself for the unknown
function $v_2=v_2(x,t)$ is the spelling of the $t^{2,0}$-flow of the full
hierarchy.

The last interesting particular value of the parameter $\kappa$ is $\kappa=1$.
The formula for the potential of the Frobenius manifold must be
modified
\beq\label{fr-nls}
F={1\over 2} v_1^2 v_2 +{1\over 2} v_2^2 \log v_2.
\eeq
However, the formula for the $-t^{1,1}$-flow of the hierarchy and the
bihamiltonian structure are well defined for this particular value of the
parameter. The PDE (\ref{poly-def}) coinicides, in the leading approximation in
$\epsilon$, with the nonlinear Schr\"odinger equation (NLS)
\beq\label{nls}
i\, \psi_t =-{1\over 2} \psi_{xx} + |\psi|^2 \psi
\eeq
when using
\eqa\label{nls-change}
&&
u={1\over 2i} \left(\log \psi/\bar\psi\right)_x
\nn\\
&&
\rho=|\psi|^2
\nn\\
\eeqa
as the new dependent variables
\eqa\label{nls-conserv}
&&
u_t + \pal_x \left[ {u^2\over 2} +\rho+\epsilon^2 \left( {{\rho'}^2\over 8\,
\rho^2} -{\rho''\over 4\, \rho}\right)\right]=0
\nn\\
&&
\rho_t + \pal_x (\rho\, u) =0.
\eeqa
The relationship of the bihamiltonian structure of NLS to topological field
theory was discussed in \cite{bonora-xiong} (see also recent paper
\cite{chang-tu1, chang-tu2}). 
To reduce the bihamiltonian structure of NLS to our normal form
one is to perform a Miura-type transformation that contains an infinite number
of terms in the $\epsilon$-expansion. We will consider this transformation
in a subsequent publication.

The above formulae work also for the exceptional value $\kappa=-1$ where
the potential of the Frobenius manifold reads
$$
F={1\over 2} v_1^2 v_2 -{1\over 2} \log v_2.
$$

It remains to consider the Frobenius manifold (\ref{potential-toda}) of the ${\bf CP}^1$
model. The transformation reducing the first Poisson structure to the canonical
form in this case becomes very simple
\beq
v_\al\mapsto \bar v_\al = v_\al +{\epsilon^2\over 24} v_\al''
+{\epsilon^4\over 1920} v_\al^{IV} +O(\epsilon^6), ~~\al=1, \, 2.
\eeq
The Poisson pencil has the form (\ref{1st-const}) with
\eqa
&&X^1=\left(\frac16\,v_2''+\frac18\,{v_2'}^2\right)e^{v_2},\quad
X^2=\frac1{12} v_1''\nn\\
&&Y^1=\left(\frac1{120}\,v_2^{IV}+\frac1{96} v_2'\,v_2'''+\frac1{288}\,
{v_2''}^2-
\frac1{1152}\,{v_2'}^4\right)e^{v_2}\nn\\
&&Y^2=-\frac1{720}\,v_1^{IV}.\nn
\eeqa
The resulting Poisson pencil is equivalent, within the order $\epsilon^4$
approximation, to the bihamiltonian structure of Toda lattice \cite{YZ}.

\def\theequation{\thesubsection.\arabic{equation}}
\setcounter{equation}{0}
\setcounter{theorem}{0}

\section{Conclusions}\par

In this paper we presented a new approach to the problem of classification
of hierarchies of 1+1-dimensional integrable  PDEs of certain class. We have 
formulated
a system of simple axioms that can be used as defining relations
in the theory of such hierarchies. We developed an efficient tools for a
perturbative reconstruction of the hierarchy starting from a given semisimple
Frobenius manifold. We computed the first few terms of the perturbative expansion
and showed that, for low genera  the identities for the tau-function of a
particular solution to the hierarchy correctly reproduce the universal
identities among Gromov - Witten classes and their descendents in the cohomology
of the moduli spaces of stable algebraic curves. 

Of course, our system of axioms can probably be improved (in particular,
independence of the quasitriviality axiom is questionable). However, as the 
reader may agree, it is something more than just putting Gromov - Witten 
invariants into the
game in the beginning in order to get them back at the very end.

Some important problems of our classification programme remain to be fixed.
We did not prove polynomiality in the derivatives in every term of the perturbative
expansion of the 
integrable hierarchy corresponding
to an arbitrary semisimple Frobenius manifold.  For $g\leq 1$ this follows
from the results of \cite{cmp}. We are to also check whether main examples
of integrable hierarchies (e.g., Drinfeld - Sokolov hierarchies
\cite{drinfeld-sokolov}
of the ADE type,
Toda lattice etc.) satisfy our four axioms. In the present paper we did it only
for the illustrative example of KdV.

A more challenging problem is to find a ``non-perturbative construction'' 
of the integrable hierarchies in question performing a summation over all
genera. We hope that the structure theory of semisimple Frobenius manifolds
based on a Riemann - Hilbert problem \cite{D3} will be important for such a
summation.
A related problem is to study relationships between our classification
approach and topological results of the very recent paper of A. Givental
\cite{givental1}. It might happen that a reasonable non-perturbative 
construction of the
hierarchy can be obtained only for certain particularly ``nice'' Frobenius 
manifolds.
Recall \cite{D3}, that generic $n$-dimensional Frobenius manifold depends on
${n(n-1)\over 2}$ complex parameters. Mathematically or physically interesting
integrable hierarchies might correspond to special points in the space of
parameters. For example, from results of C. Hertling \cite{hertling} it can 
probably
be derived that integrable hierarchies satisfying our axioms that are polynomial
also in $u$
must correspond
to the Frobenius structures on the orbit space of Weyl groups of the ADE type.
However, just the universality of our approach suggests to apply it to other
ingredients of the theory of integrable systems, first of all to the theory
of classical and quantum W-algebras \cite{zam, FatLuk, DIZ} and, more
generally, of vertex Poisson algebras and vertex algebras \cite{frenkel}
in order to unravel the eventual role of the topology of moduli spaces of stable
algebraic curves in the theory of these algebraic structures. 

From the point of view of applications to the theory of dispersive waves
the integrable systems we are
constructing may look ugly. Indeed, they typically contain infinite number
of terms of the dispersion in the r.h.s. However, one can obtain an {\it
approximately integrable system} just truncating the expansion at some order
in $\epsilon$. We expect that solutions to such an approximately integrable 
system
exhibit an integrable behaviour within a certain range of physical parameters.
For $n=1$ such approximately integrable systems were studied in \cite{lorenzoni}
(only the bihamiltonian property has been used). In numerical experiments
elastic scattering of solitons was observed in such systems. It could be of
interest to study numerically our approximately integrable deformation of the
one-dimensional gas dynamic equations constructed in Section 3.10.8.

Finally, the problem of classification of Poisson pencils of
$(p,q)$-brackets with \newline
$(p,q)\neq (0,n)$ and study of the properties of
associated bihamiltonian hierarchies is still completely open. The $(n,0)$ case was
recently studied by P.~Lorenzoni (unpublished) where preliminary results 
were obtained. Such a classification is needed to extend our general
classification strategy onto other classes of integrable systems not admitting
dispersionless limit (e.g., integrable systems of the Sine-Gordon type).

We plan to study these problems in subsequent publications.
\smallskip

{\bf Note added.}  After finishing this work we have received the paper
of S. Barannikov, Semi-infinite variations of Hodge structures and integrable
hierarchies of KdV-type, math.AG/0108148. The author also presents a
construction of a bihamiltonian  hierarchy associated to an arbitrary semisimple
Frobenius manifold. Virasoro symmetries of this hierarchy were not discussed.
Relations of the Barannikov's  hierarchy to Gromov - Witten invariants
is still to be investigated.

\medskip

{\bf Acknowledgments.} We wish to thank O. Babelon, Ya. Eliashberg, G. Falqui, G. Felder, E. Frenkel, 
E. Getzler, 
V. Kac, F. Magri, A. Maltsev, V. Sokolov,
P. Wiegmann and C.-S. Xiong  
for fruitful discussions. Part of this work was done during our
visit to MSRI. We wish to thank the organizers of the MSRI activity
``Random Matrix Models and Their Applications'' for the invitation and generous
support. One of the authors (B.D.) also acknowledges generous support
of LPTHE, of Isaac Newton Mathematical Institute, of ETH, of Tsinghua University where some parts of this work were done; some parts of Y.Z.'s work
were done during his visit to the International Center for Theoretical 
Physics (ICTP), he acknowledges the generous support of ICTP and 
of the Commission on Development and Exchanges of the International 
Mathematical Union.
The researches of B.D. were
partially supported by Italian Ministry of Education research grant ``Geometry
of Integrable Systems''. The researches of Y.Z. were partially supported
by the NSFC  
grants No.10025101, 10041002 and
the Special Funds of Chinese Major Basic Research Project 
``Nonlinear Sciences'', he also acknowledges the partial support by the 
startup grants from 
the Chinese 
Ministry of Education and Tsinghua University.

\def\thetheorem{A.\arabic{theorem}}
\def\theprop{A.\arabic{prop}}
\def\thelemma{A.\arabic{lemma}}
\def\thecor{A.\arabic{cor}}
\def\theexam{A.\arabic{exam}}
\def\theremark{A.\arabic{remark}}
\def\theequation{A.\arabic{equation}}

\newpage
\appendix
\makeatletter
\renewcommand{\@seccntformat}[1]{{Appendix:}\hspace{-2.3cm}}
\makeatother
\renewcommand{\thesection}{Appendix:}
\section{\quad\qquad \ \ Degenerate Frobenius manifolds, spectral curves,
and Poisson pencils}

Let $M$ be a degenerate Frobenius manifold, $F(v)$ the corresponding solution to
the equations of associativity, $e=\partial_1$ and $E=E^\alpha\partial_\alpha$
the unity and the Euler vector
fields on $M$ resp. satisfying
\beq
\partial_EF(v) = k\, F(v) +{\rm quadratic}, ~~[e,E]=0.
\label{(1)}
\eeq
We introduce the following linear operator on $TM$
\beq
{\cal V} := {k\over 2} -\nabla \,E
\eeq
(cf. (\ref{2-7-4b}) above).
This (1,1)-tensor is constant in the flat coordinates on $M$.

\begin{lemma}\label{l-a-2-12-1} The following equations hold true 
\beq\label{a-2-12-3}
Lie_E <\ ,\ > = k <\ ,\ >
\label{(2)}
\eeq
\beq\label{a-2-12-4}
Lie_E(a\cdot b) = Lie_Ea\cdot b+a\cdot Lie_Eb.
\label{(3)}
\eeq
\end{lemma}

\pf Differentiating (\ref{(1)}) along $\partial_1$, $\partial_\alpha$,
$\partial_\beta$ one obtains
\beq
<{\cal V}x, y> +<x,{\cal V}y>=0
\label{(4)}
\eeq
for any two vectors $x$, $y$. This implies (\ref{(2)}). Differentiating the same
equation along $\partial_\alpha$, $\partial_\beta$, $\partial_\gamma$ we obtain
$$
\partial_E c_{\alpha\beta\gamma}=-{k\over 2}\, c_{\alpha\beta\gamma} 
+{\cal V}_\alpha^\epsilon
c_{\epsilon\beta\gamma} +{\cal V}^\epsilon_\beta c_{\alpha\epsilon\gamma}
+{\cal V}^\epsilon_\gamma c_{\alpha\beta\epsilon}.
$$
Raising the index $\gamma$ and using (\ref{(4)}) we arrive at
$$
Lie_E\, c_{\alpha\beta}^\gamma =0.
$$
This coincides with (\ref{(3)}). \epf

Denote
\beq
L(v;z): = z\, {\cal U} + {\cal V}.
\eeq

\begin{lemma} The operator $L(v;z)$ is horizontal w.r.t. to the deformed flat
connection
\beq 
\tilde \nabla(z) \, L(v;z) =0.	
\eeq
\end{lemma}

The proof easily follows from Lemma \ref{l-a-2-12-1}.

\begin{cor} The spectral curve does not depend on the point $v\in M$ of the
degenerate Frobenius manifold. The differentials of the deformed flat connection $\xi = d\tilde v$,
$\tilde v= \tilde v(v;z)$ are meromorphic sections of the  vector bundle of
eigenvectors of $L(v;z)$ over the
spectral curve
\beq\label{a-2-12-8}
{\cal C}:~\det (L(v;z) - w) =0.
\eeq
\end{cor}

\pf For the covectors $\xi$ we obtain the compatible system
\eqa
&&
\pal_\alpha \xi = z\, C_\alpha(v)\,\xi
\\
&&
L(v;z)\, \xi = w\, \xi.
\label{line-bundle}
\eeqa
This proves the corollary.
\epf

Let us assume that the spectral curve has no singularities at $|z|<\infty$
and that (\ref{line-bundle}) is a line bundle on ${\cal C}$. We will see below
that this is the case under the assumption of semisimplicity. In this case
one can choose a meromorphic trivialization of the line bundle. 
Denote $\tilde v(v;P)$  the corresponding flat function for the deformed
connection $\tilde \nabla$ considered as a meromorphic function on
the spectral curve depending on the point $v\in M$. Let $z=0$ be not a
ramification point of the spectral curve. Denote
$$
P_\alpha := (z=0, w=\mu_\alpha), ~~\alpha = 1, \dots, n
$$
the points of the spectral curve over $z=0$. Expanding the branches of the
function $\tilde v(v; z)$ near $P\to P_\alpha$ we obtain, after a multiplication
by an appropriate normalization factor depending only on $z$, the basis
of the deformed flat coordinates $\tilde v_\alpha(v; z)$, $\alpha=1, \dots, n$.
Recall that the coefficients of expansion of these functions in $z$ are
Hamiltonian densities of the analogue of the Principal Hierarchy for the case
of degenerate Frobenius manifold.
This procedure can be modified in an obvious way in the case when $z=0$ is a
ramification point.

Due to symmetry of ${\cal U}$ and antisymmetry of ${\cal V}$ w.r.t. the bilinear
form $<\ ,\ >$ the spectral curve admits a holomorphic involution
\beq\label{a-2-12-10}
\sigma:\, (z,w)\mapsto (-z, -w).
\eeq
The orthogonality condition in this case
reads as follows: the inner product\newline
 $<\nabla \tilde v(v;P), \nabla
\tilde v(v;\sigma(P))>$ does not depend on $v\in M$. Considering the function
on the Cartesian square of the spectral curve
$$
{<\nabla \tilde v(v;P), \nabla
\tilde v(v;\sigma(P'))>-<\nabla \tilde v(v;P), \nabla
\tilde v(v;\sigma(P))>\over z(P)+z(P')}
$$
and expanding it into a power series in $z=z(P)$, $z'=z(P')$ near
$P\to P_\alpha$, $P'\to P_\beta$ we obtain, after an appropriate
renormalization, the matrix $\Omega(v;z, z')$.

The system for the differentials of periods $p(v;\lambda)$ 
(i.e., of the flat coordinates
of the flat pencils of metrics $(~,~)_\lambda$) reads 
$$
({\cal U} -\lambda)^{-1} \partial_\alpha \phi +C_\alpha \mu \phi =0, ~~
\phi =\nabla p(v;\lambda).
$$
From here it follows that
$$
\partial_{E-\lambda e} \phi =0
$$
and, more important, that the operator $({\cal U}-\lambda)^{-1} {\cal V}$
is horizontal w.r.t. the connection $\nabla^* -\lambda\, \nabla$. 
Therefore $\phi$ is meromorphic on the spectral curve
$$
\det [\nu({\cal U} -\lambda) +\mu]=0, ~~
({\cal U}-\lambda)^{-1} \mu \phi +\nu \phi =0
$$
The new spectral curve is birationally equivalent to the old one:
\beq\label{linea}
w=\lambda\nu, ~z=\nu.
\eeq
Therefore the system of independent periods $p_1(v;\lambda)$, \dots,
$p_n(v;\lambda)$ is obtained as follows
\beq
p_\alpha(v;\lambda) = \tilde v(v; P_\alpha(\lambda))
\eeq
where $P_\alpha(\lambda)$ are the intersection points of the line (\ref{linea})
with the spectral curve. For large $\lambda$ they can be ordered in such a form
that 
$$
P_\alpha(\infty) = P_\alpha, ~~\alpha=1, \dots, n.
$$

\begin{exam} An arbitrary two-dimensional degenerate Frobenius manifold has the
potential
\beq
F={1\over 2} v_1^2 v_2 + a^2 v_2 \log v_2, ~~E=v_2 \partial_2
\eeq
where $a$ is a parameter. Then
$$
{\cal V}={\rm diag} \,(1/2, -1/2),
$$
$$
L(v;z) =\left(\matrix{ 1/2 & -a^2 z/v_2\cr v_2 z & -1/2 \cr}\right).
$$
This gives the spectral curve
$$
w^2 + a^2 z^2 ={1\over 4}.
$$
Uniformizing it
$$
 z = {1\over a} {s\over 1+s^2}, ~~~
w = {1\over 2} {1-s^2\over 1+s^2}
$$
we obtain the eigenvector in the form
$$
\xi =\phi \left(\matrix{ a\cr s v_2\cr}\right)
$$
where $\phi$ is a normalizing factor. The dependence of it on $v_1$ and
$v_2$ is completely determined by the above differential equations.
This gives
$$
\tilde v = e^{s\, v_1\over a\, (1+s^2)} v_2^{1\over 1+s^2}.
$$
Therefore the basis of the deformed flat coordinates reads
$$
\tilde v_1 = z^{-1} \left( e^{z\, v_1} v_2^{1-\sqrt{1-4a^2 z^2}\over 2}
-1\right)
$$
$$
\tilde v_2 = e^{z\, v_1} v_2 ^{1+\sqrt{1-4 a^2 z^2}\over 2}.
$$
The coefficients of expansion of these functions in powers of  $z$
are hamiltonians of the hierarchy related to the degenerate polytropic gas
equations
\eqa
&&
\rho_t +(\rho\, u)_x =0
\nn\\
&&
u_t + u\, u_x -a^2 (\log \rho)_x =0.
\nn\\
\eeqa
Here $u=v_1$, $\rho = v_2$, $t=-t^{2,1}$.
  
The above involution becomes $s\mapsto -1/s$. The orthogonality gives
$$
<\nabla \tilde v(s), \nabla\tilde v(-{1\over s})> = -2zw.
$$
Using this one can easily compute the tau-fucntion of the hierarchy.
\end{exam}

Let us now add the assumption of semisimplicity. Recall that this means
that the eigenvalues of the operator of multiplication by $E$ are pairwise
distinct. One can locally introduce a system of canonical coordinates $u_1,
\dots, u_n$ on $M$ such that the vector fields
$$
\pi_i :={\partial\over \partial u_i}, ~~i=1, \dots, n
$$
are the basic idempotents of the Frobenius algebra on $T_uM$ at each point
$(u_1, \dots, u_n)=u\in M$,
$$
\pi_i \cdot \pi_j =\delta_{ij} \pi_i.
$$

\begin{lemma} The Euler vector field in the canonical coordinates equals
\beq\label{a-2-12-14b}
E=\sum_{i=1}^n c_i {\partial\over \partial u_i}
\eeq
where $c_1$, \dots, $c_n$ are pairwise distinct constants.
\end{lemma}

\pf Define the matrix $(b_{ij})$ by
$$
Lie_E\,\pi_i =\sum_{j=1}^n b_{ij}\pi_j.
$$
Using (\ref{a-2-12-4}) one obtains, for $i\neq j$
$$
0=Lie_E (\pi_i\cdot \pi_j) = b_{ij}\pi_j +b_{ji}\pi_i.
$$
Hence $b_{ij}=b_{ji}=0$. Using the same trick for the square of $\pi_i$ we have
$$
b_{ii}\pi_i =Lie_E(\pi_i\cdot\pi_i) = 2b_{ii}\pi_i.
$$
Hence
\beq\label{a-2-12-14c}
Lie_E\pi_i=0, ~~i=1, \dots, n.
\eeq
This means that $E$ has the form (\ref{a-2-12-14b}) with some constants $c_i$. These constants
clearly are the eigenvalues of the operator of multiplication by $E$,
$$
E\cdot \pi_i = c_i \pi_i.
$$
Therefore
$$
c_i\neq c_j ~~{\rm for }~i\neq j.
$$

Introduce, as usual, the basis of normalized idempotents
$$
f_i={\pi_i\over \sqrt{<\pi_i, \pi_i>}}
$$
and the transition matrix
$$
\Psi=(\psi_{i\alpha}(u)), ~~\partial_\alpha =\sum_{i=1}^n \psi_{i\alpha}f_i.
$$
The matrix $\Psi$ satisfies the linear system
\beq
\partial_i \Psi =V_i \Psi, ~~i=1, \dots, n
\label{(7)}
\eeq
where the antisymmetric matrix $V_i$ is defined by this equation. The matrices
$V_i$ satisfy
$$
[E_i, V_j] =[E_j,V_i], ~~
\partial_jV_i -\partial_jV_j +[V_i, V_j]=0, ~~i\neq j.
$$
The last system coincides with the equations of flateness of the deformed
connection $\tilde \nabla$
\beq
[\partial_i -(zE_i +V_i), \partial_j-(zE_j+V_j)]=0
\label{(8)}
\eeq
depending on the parameter $z$.
Here $E_i$ is the $i$-th
matrix unity, i.e.
$$
(E_i)_{pq}=0 ~~{\rm for}~ (p, q)\neq (i,i), ~~(E_i)_{ii}=1.
$$
The flat coordinates and the structure constants can be reconstructed from an
arbitrary solution to (\ref{(7)}), (\ref{(8)}) by the formulae of \cite{D92}
\beq
v_\alpha =\Omega_{\alpha,0; 1,0}, ~~dv_\alpha =\sum_{i=1}^n
\psi_{i1}\psi_{i\alpha} du_i
\label{(9)}
\eeq
\beq
\eta_{\alpha\beta} =  \sum_{i=1}^n \psi_{i\alpha}\psi_{i\beta}
\label{(10)}
\eeq
\beq
F={1\over 2} \Omega_{1,1; 1,1} - v^\alpha \Omega_{\alpha,0; 1,1} +{1\over 2}
v^\alpha v^\beta \Omega_{\alpha,0; \beta,0}
\label{(11)}
\eeq
\beq
e=\sum_{i=1}^n {\partial\over\partial u_i}.
\label{(12)}
\eeq
All these statements were proved in \cite{D92} for an arbitrary solution to the
associativity equations (the quasihomogeneity has not been used in their
derivation).

Let us now use the Euler vector field of the form (\ref{a-2-12-14b}).

\begin{lemma}\label{l-a-2-12-6} The matrices $\Psi$, $V_i$ satisfy
\beq
\partial_E \Psi =\Psi\, {\cal V}
\label{(13)}
\eeq
\beq
\partial_E V_i =0, ~~i=1, \dots, n.
\label{(14)}
\eeq
\end{lemma}

\pf Using (\ref{a-2-12-3}) and also
$$
<\ ,\ >=\sum_{i=1}^n \psi_{i1}^2\, du_i^2
$$
we obtain
$$\partial_E \psi_{i1}=-{k\over 2} \psi_{i1}.
$$
From this using commutativity (\ref{a-2-12-14c}) we obtain 
(\ref{(13)}). Therefore
$$
\partial_E V_i = \partial_E \left[ \partial_i \Psi \cdot \Psi^{-1}\right] =0.
$$
\epf

\begin{lemma} The operators $\partial_i-(zE_i+V_i)$ commute with the
operators of multiplication by the matrix
$$
L(z):= zC+V, ~~C={\rm diag} (c_1, \dots, c_n), ~~V=\sum_{i=1}^n c_i V_i =-
\Psi{\cal V}\Psi^{-1}.
\eqno(15)
$$
\end{lemma}

\pf This is just the spelling of the statement of Lemma \ref{l-a-2-12-6}.
\epf

We arrive at

\begin{theorem} For a semisimple degenerate Frobenius manifold
the antisymmetric matrix $V\in so(n)$ satisfies the Euler
equations of free rotations of $n$-dimensional rigid body
\beq
\partial_iV = [V_i, V], ~~V_i = {\rm ad}_{E_i} {\rm ad}_C^{-1}(V), ~~i=1, \dots,
n.
\label{(17)}
\eeq
Conversely, let $V=V(u_1, \dots, u_n)$ be any solution to (\ref{(17)}), and
$\Psi=\Psi(u_1, \dots, u_n)$ be a solution to the linear system (\ref{(7)}). 
Then the
formulae (\ref{(9)}) - (\ref{(12)}), (\ref{a-2-12-14b}) 
determine a structure of a 
semisimple degenerate
Frobenius manifold on an open subset in $M$ consisting of the points $u=(u_1,
\dots, u_n)$ where
$$
\prod_{i=1}^n \psi_{i1}(u)\neq 0.
$$
\end{theorem}

The nonlinear system (\ref{(17)}) can be integrated \cite{matrix} in terms of Prym
theta-functions of the spectral curve (\ref{a-2-12-8})
with the involution (\ref{a-2-12-10}).
The solution $\Psi$ to the
linear system (\ref{(7)}) and also all the ingredients of the degenerate Frobenius
structure can be given in terms of Baker - Akhiezer functions on the spectral
curve. In the generic situation the spectral curve is a smooth plane 
algebraic curve of the degree $n$. It has the genus equal 
$$
g={(n-1)(n-2)\over 2}.
$$
The fixed points of the involution $\sigma$ are
$$
\infty_k = (z, \, w \to \infty, {w\over z} \to c_k), ~~k=1, \dots, n
$$
and, for odd $n$, also the point 
$$
P_0 = (z=0, w=0).
$$
Baker-Akhiezer function is a
vector function $Y=Y(u, P)$, $u=(u_1, \dots, u_n)$, meromorphic in $P=(z,w)\in
{\cal C}\setminus (\infty_1, \dots, \infty_n)$ such that
\beq
Y=e^{z\, u_k} (e_k +O({1\over z}), ~~P\to \infty_k, ~~k=1, \dots, n
\eeq
with a nonspecial divisor of poles $D$, $\deg D= {n(n-1)\over 2}$ that must
belong to the odd part of the generalized Jacobian $J({\cal C}, \infty_1, \dots,
\infty_n)$ of the spectral curve with identified infinite points. Here 
$$
(e_k)_i =\delta_{ik}.
$$
We will give
elsewhere the explicit formulae in terms of Prym theta functions for the
corresponding degenerate Frobenius manifold and for the
analog
of the Principal Hierarchy for this manifold. We also postpone for a subsequent
publication the study of the perturbations of the hierarchy.

\begin{remark} The method of constructing solutions to the associativity
equations using theta-functions has been suggested in \cite{D92}. The explicit
solutions for an arbitrary $n$-sheeted Riemann surface with an involution were
obtained by I.Krichever \cite{krich}. The relationship of these solutions to
bihamiltonian structures of hydrodynamic type was not observed before. This
relationship takes place only for the case where the spectral curve is a plane
algebraic curve of the degree $n$ or its degeneration.
\end{remark}

\begin{exam}
$n=3$. According to the above theory, three-dimensional
semisimple degenerate Frobenius manifolds are expressed via solutions to the
classical Euler equations of free rotations of a rigid body. The potential
of the degenerate Frobenius manifold reads
\beq
F = {1\over 2} v_1^2 v_3 + {1\over 2} v_1 v_2^2 + f(v_2, v_3)
\eeq
where
\beq
f(x,y) =x^2 g(y\, x^{-2}) + a\, x^2 \log x + b\, y\, \log y.
\eeq
Here $a$, $b$ are arbitrary constants, $g=g(s)$ is given by an elliptic quadrature
\beq
s^3 g'' ={1\over 8} \left( 1-4\, a\, s + \sqrt{ 1-8\, a\, s+16 (a^2+b)\, s^2 + 
8(c-8\, a\, b)\, s^3}\right)
\eeq 
where $c$ is another constant. The spectral curve is a plane cubic
\beq
w^3 - 2a\, w^2 z +4b\, w\, z^2 - c \, z^3 -w + 2a\, z=0.
\eeq
\end{exam}

\medskip

E-mail addresses: dubrovin@sissa.it, yzhang@math.tsinghua.edu.cn
\end{document}